\documentclass[12pt]{amsart}
\usepackage{amssymb,amscd}
\usepackage[russian]{babel}
\usepackage{graphicx}
\usepackage{pdfpages}
\usepackage{placeins}
\newtheorem{theorem}{Теорема}[section]
\newtheorem{proposition}{Предложение}[section]
\newtheorem{corollary}{Следствие}[section]
\newtheorem{lemma}{Лемма}[section]

\def\Ad{\operatorname{Ad}}
\def\DOLL{\operatorname{DOLL}}

\theoremstyle{definition}
\newtheorem{definition}{Определение}[section]

\theoremstyle{remark}

\def \x#1{\begin{tabular}{c}#1\end{tabular}}

\ifx\pdfpagewidth\undefined \else
\fi
 \evensidemargin \oddsidemargin

\advance\hoffset by -1in
\advance\voffset by -1in

\usepackage[cp1251]{inputenc}
\newcount \lemmacount
\def \lemma #1.{\global \advance \lemmacount by 1
{\bf Лемма \No \thesubsection.\the \lemmacount \ #1.}%
 \begingroup \sl
 \def \par {\endgroup \par}}

\author{И.~А.~Иванов-Погодаев, А.~Я.~Канель-Белов}

\title{Конструкция бесконечной конечно определенной нильполугруппы}

\begin{document}

\let \mathbf=\texttt
\tabcolsep 2pt

\begin{abstract}

Работа посвящена решению проблемы Л.~Н.~Шеврина и М.~В.~Сапира. Строится конечно определенная бесконечная нильполугруппа, удовлетворяющая тождеству $x^9=0$.
Доказательство основано на изучении свойств апериодических мозаик, аналогов теоремы Гудмана-Штраусса на равномерно-эллиптическом пространстве. Пространство называется {\it равномерно-эллиптическим}, если любые две точки $A$ и $B$ на расстоянии $D$ соединяются системой геодезических,  образующих диск ширины $\lambda\cdot D$ для некоторой глобальной константы $\lambda>0$.

Исследование поддержано грантом РФФИ N 14-01-00548

\end{abstract}

\maketitle

\medskip
\tableofcontents

\medskip

\section{Введение} \label{nachalo}

Работа посвящена построению конечно определенных нильполугрупп. Доказана следующая

\medskip
{\bf Теорема.} {\it Существует конечно определенная бесконечная нильполугруппа, удовлетворяющая тождеству $x^9=0$.}
\medskip

Теории полугрупп посвящено ряд обзоров (см. например, \cite{Shevrin,Sapirbook}).
 Построение таких полугрупп затрагивает, с одной стороны, проблемы Бернсайдовского типа, а с другой -- проблемы построения конечно определенных объектов. Поэтому следует, хотя бы кратко, упомянуть историю вопроса в этой связи.

\subsection{Проблемы бернсайдовского типа} \label{vvedenie}

Проблемы Бернсайдовского типа оказали огромное значение в развитие современной алгебры. Эта проблематика охватила большой круг вопросов,
как в теории групп, так и в смежных областях, стимулировала алгебраические исследования. Проблемам бернсайдовского типа
 посвящена обзорная статья Е.~И.~Зельманова \cite{Zelmanov}.

Важнейшее значение имеет теория групп. Ясно, что группы, удовлетворяющие тождеству $x^2=1$ коммутативны. Для групп с тождеством $x^3=1$ доказать локальную конечность несколько труднее (и она была доказана самим Бернсайдом \cite{Burnside}). Вопрос о локальной конечности групп с тождеством $x^4=1$ стоял открытым чуть меньше 40 лет \cite{Sanov}, а с тождеством $x^6=1$~-- свыше 50 \cite{HallM}.

Получению контрпримеров предшествовало построение бесквадратных слов над трехбуквенным алфавитом, а также бескубных слов над алфавитом их двух букв, т.е. построению бесконечных ниль-полугрупп с тождествами $x^2=0$ и $x^3=0$.  Вопрос о локальной  конечности групп с тождеством $x^n = 1$ был решен отрицательно в знаменитых работах П.~С.~Новикова и С.~И.~Адяна  \cite{Novikov-Adyan}: было доказано существование для любого нечетного $n > 4381$  бесконечной группы с $m>1$ образующими, удовлетворяющей тождеству $X^n = 1$. Результаты Новикова-Адяна А.~И.~Мальцев рассматривал как основное событие алгебры 20 века.
 Эта оценка  была улучшена до $n > 665$ С.~И.~Адяном \cite{Adyan}. Совсем недавно С.~И.~Адяном был совершен неожиданный прорыв: эта оценка была улучшена до $n \ge 101$ (см.  \cite{Adyan1}).

 Работы П.~С.~Новикова и С.~И.~Адяна оказали огромное влияние на творчество И.~А.~Рипса, который в дальнейшем разработал метод канонической формы и также построил примеры бесконечных периодических групп.  Позднее А.~Ю.~Ольшанский предложил геометрически
наглядный вариант доказательства для нечетных $n>10^{10}$ \cite{Olshansky}. Для всех достаточно больших четных $n$ примеры бесконечных 2-порожденных групп периода $n$ были построены независимо И.~Лысенком и С.~Ивановым. Подробная библиография и история вопроса -- см. \cite{Adyan}.

 В $PI$-случае вопросы локальной конечности алгебраических алгебр решаются положительно. В ассоциативном случае
соответствующий результат был получен И.~Капланским и Д.~Левицким. Чисто комбинаторное доказательство для ассоциативного случая получается из теоремы Ширшова о высоте \cite{Shirshov1}, \cite{Shirshov2}.
 Теореме Ширшова о высоте посвящена обширная библиография (см., например, работы \cite{BelovRowenShirshov,Kem09,BelovHeightObzor,BelovKharitonov}).
  Для $PI$-алгебр Ли соответствующий результат был получен  А.~И.~Кострикиным (для нулевой характеристики, кроме того в характеристике $p$ им вместе с Е.~И.~Зельмановым показана локальная нильпотентность алгебр Ли с тождеством $y\circ \Ad(x)^n=0$) и в общем случае Е.~И.~Зельмановым. Это позволило решить так называемую ``ослабленную проблему Бернсайда'', то есть доказать наличие максимальной конечной группы среди $k$-порожденных групп, удовлетворяющих тождеству $x^n=1$.   Подробная библиография по этому вопросу изложена в монографии  \cite{Kostrikin}.
Первый контрпример к неограниченной проблеме был найден Е.~С.~Голодом в 1964 году на  основе универсальной конструкции Е.~С.~Голода -- И.~Р.~Шафаревича \cite{Golod,GolodShafarevich}. Эта конструкция позволила также построить не локально конечные периодические группы неограниченной экспоненты.

\subsection{Проблемы конечной определенности}

Все имеющиеся примеры бесконечных периодических групп бесконечно порождены. Чрезвычайно глубоким и вдохновляющим является следующий открытый вопрос (входящий список основных алгебраических проблем в теории групп):

\medskip

{\bf Вопрос.}\
{\it Существует ли конечно определенная бесконечная периодическая группа?
}
\medskip

Известна классическая теорема Хигмана о вложении рекурсивно-определенных групп в конечно определенные.

В работе А.~Ю.~Ольшанского и М.~В.~Сапира \cite{OlshanskiiSapir} была построена конечно определенная группа являющаяся расширением конечно порожденной бесконечной периодической группы с помощью циклической.

На проблематику, связанную с построением разного рода
 экзотических объектов с помощью конечного числа определяющих соотношений обратил
внимание В.~Н.~Латышев. Им же была поставлена проблема существования конечно
определенного ниль-кольца~\cite{Dnestrovsk}.

\medskip
{\bf Вопрос (В.~Н.~Латышев).} {\it Существует ли конечно определенное бесконечномерное нилькольцо?
}
\medskip

В качестве продвижения в решении этого вопроса можно
 рассматривать результаты Г.~П.~Кукина, В.~Я.~Беляева о вложениях
 рекурсивно определенных объектов в конечно определенные \cite{Kukin,Belyaev}. В.~А.~Уфнаровским был построен пример конечно определенной алгебры промежуточного роста \cite{Ufnar1}. В работе В.~В.~Щиголева \cite{Schigolev} была изучена связь между понятиями ниль и нильпотентности конечно определённых алгебр в зависимости от количества определяющих соотношений и порождающих. Также построен пример алгебр с малым количеством определяющих соотношений, у которых все слова длины два нильпотентны.

\medskip

Фундаментальную проблему существования конечно определенной нильполугруппы поставили Л.~Н.~Шеврин и М.~В.~Сапир в Свердловской Тетради (3.61б) \cite{Sverdlovsk},  а также вопрос 3.8 в \cite{Obzor}.
%Философия проблемы -- в возможности добиться локальными средствами глобального эффекта (в данном случае, свойства нильности).

\medskip
{\bf Вопрос (Л.~Н.~Шеврин, М.~В.~Сапир).} {\it Существует ли конечно определенная бесконечная нильполугруппа?
}
\medskip

Авторам представляется, что полугрупповой вопрос предшествует кольцевому, а кольцевой -- групповому. Этот вопрос привлекал внимание авторов в течение многих лет. Ряд результатов, таких как конструкция конечно-определенной полугруппы с нецелой размерностью Гель\-фан\-да--Ки\-рил\-ло\-ва, построение алгебр с конечным базисом Гребнера но неразрешимой проблемой делителей нуля возникли из работы над этой проблемой. В своей кандидатской диссертации \cite{ivadis} один из авторов построил пример конечно определенной полугруппы, содержащий ненильпотентый ниль-идеал.

\begin{theorem}[\cite{ivadis}]
Существует конечно-определенная полугруппа H, с множеством слов $G$ от образующих, удовлетворяющая
следующим свойствам:

\begin{itemize}
  \item Существует ненильпотентный идеал $I = LH$,
где $L$~-- буква в $H$;
  \item Если слово $A\in G$ представляется в виде $A
= XYYZ$, где $X,Y,Z\in G$, тогда $LA = 0$.
\end{itemize}

\end{theorem}

Основной результат данной работы заключается в следующем:

\begin{theorem}
  Существует конечно определенная бесконечная нильполугруппа, удовлетворяющая тождеству $x^9=0$.
\end{theorem}

{\bf Благодарности.} Авторы признательны руководителям семинара <<Теория колец>> на кафедре Высшей Алгебры механико-математического факультета МГУ В.~Н.~Латышеву и А.~В.~Михалеву за полезные обсуждения и постоянное, в течение ряда лет, внимание к работе. Мы также благодарны  И.~А.~Рипсу, Л.~Н.~Шеврину, А.~Х.~Шеню за полезные обсуждения, Ф.~Дюранду, Ц.~Селле, Л.~А.~Бокутю, Ю. ~Чэну, Т.~Ферни за поддержку в участии на конференциях,  А.~С.~Малистову за помощь в оформлении статьи.

\section{Конечная определенность, связь с информатикой, геометрические методы}      \label{ScFnDefInfGen}

\subsection{Построение конечно определенных объектов и системы конечных автоматов}

При построении разного рода объектов трудности вызывает контроль над следствиями из вводимых соотношений, в особенности тогда, когда
следует доказать, что некое соотношение не является следствием заданных.
Зачастую используются три метода контроля над соотношениями:

\begin{enumerate}

\item Базис Гребнера и бриллиантовая лемма;
\item Теория малых сокращений;
\item Реализация машины Тьюринга или машины Минского.

\end{enumerate}

\medskip

В конечно определенном случае  вопросы, связанные с построением объектов, обладающих заданными свойствами, сильно усложняются и наибольшее значение приобретает третий метод.
При этом буква интерпретируется как состояние конечного автомата, а слово -- как цепочка взаимодействующих конечных автоматов. Если число соотношений конечно, то это взаимодействие {\it локально} и мы получаем связь с задачами самоорганизации, типа задачи Майхилла о стрелк\'ах и т.д. На этом пути были решены задача
о построении конечно определенных полугрупп с рекурсивной размерностью Гельфанда-Кириллова. Отметим, что подобная техника довольно громоздка для построений, требующих малого роста. Например, не удалось пока построить конечно определенную полугруппу с размерностью Гельфанда-Кириллова равной $2.5$.

 На этом же пути был получен ответ на известный открытый вопрос, поставленный В.~Н.~Латышевым -- была  построена  алгебра с неразрешимой проблемой делителей нуля и конечным базисом Гр\"еб\-не\-ра \cite{ivadis}. Отметим, что для автоматных мономиальных алгебр (в частности, конечно определенных) а также {\it алгебр с ограниченной переработкой} аналогичный вопрос (также поставленный В.~Н.~Латышевым) решается положительно \cite{BBL,Belov,Iudu}.

\medskip

Задача о построении конечно определенной бесконечной нильполугруппы имеет интерпретацию в этих терминах. Рассмотрим цепочку локально взаимодействующих конечных автоматов. У них есть цвета корпусов. Если автомат объявляет себя нулем (совершает самоубийство), то вся цепочка погибает. Можно ли добиться того, чтобы преобразования были обратимы (если $u=v$, то $v=u$),  при этом существовали сколь угодно длинные живые цепочки, и чтобы любая цепочка, у которой несколько раз подряд повторились цвета корпусов, погибала. Мы задаем локальный закон взаимодействия, а наш враг -- задает внутренние состояния автоматов. Требуется достичь адекватного поведения цепочки автоматов.

Хотя решение проблемы построения бесконечной нильполугруппы было достигнуто геометрическими методами, данная интерпретация демонстрирует связь с самоорганизующимися системами и могут быть интересными с точки зрения получения результата в этой области.

\subsection{Схема доказательства. Мозаики}

Пусть $W$~-- бесквадратное слово над конечным алфавитом из трех букв. Если каждое его неподслово (т.е. антислово) объявить нулем, то естественно возникающая полугруппа слов обладает тождеством $x^2=0$. Однако естественная конструкция, связанная с заданием множества нулевых слов как подслов слов из некоторого семейства в конечно определенном случае работает плохо.

Невозможно задать конечно определенную ненильпотентную нильполугруппу только {\it мономиальными} соотношениями, т.е. соотношениями типа $v=0$. Ибо если есть конечный список запретов и бесконечное слово без запрещенных подслов, то есть и бесконечное {\em периодическое} слово также без запрещенных подслов.

Что если от одномерных расположений (слов) перейти к двумерным?
Известно, что  существуют конечные наборы многоугольников (плиток) которыми можно
замостить плоскость лишь непериодически. Впервые такой набор был построен Робертом Бергером \cite{Berger}. В
дальнейшем были построены более простые примеры, например, Рафаэлем Робинсоном \cite{Robinson}. Широко известна также мозаика Пенроуза.
Итак, имеются контактные правила, для которых: \begin{itemize}
  \item Существует замощение всей плоскости, удовлетворяющее запретам.
  \item Однако таких периодических замощений не существует.
\end{itemize}
Поэтому если бы можно было бы умножать слева-справа-сверху-снизу, то такого рода объекты можно было бы построить.

Но как придать всему этому смысл? Будем интерпретировать элементы полугруппы как пути на мозаике (дальнейший анализ показывает, что удобно иметь дело с кратчайшими путями -- иначе можно много раз проходить один и тот же цикл). Буквы кодируют плитки и переходы между ними. Если локальный не порядок (два символа плитки без символа перехода между ними и т.д.) -- то произведение ноль. Кроме того, если локальный участок не располагается на мозаике, или не располагается как участок кратчайшего пути -- то он также нулевой.

Если же любую пару узлов, соединяемую кратчайшим путем с кодом $s_1$ можно соединить
кратчайшим путем с кодом $s_2$ и наоборот, то $s_1=s_2$.

\medskip

Итак, пусть есть периодическое слово $U=W^n$. Начинаем добавлять клетки к слову $U$, локально перебрасывая пути. Оно окружается мозаикой. Поскольку $U$ периодично, то не может быть расположено на нашей мозаике. Поэтому в какой то момент вставлять клетки не получится и мы доберемся до локального участка, несовместимого с мозаикой. Тем самым устанавливается равенство слова $U=W^n$ нулю.

Таким образом, возникает мозаика со своей глобальной структурой, которая и обеспечивает апериодичность. Локальные правила задают эту самую структуру, и путь, ``перекидывание'' которого задает область на мозаике. Возникают три языка: геометрический язык, описывающий глобальное поведение комплекса, комбинаторно геометрический язык контактов (локальных правил) и полугрупповой язык соотношений -- переброски путей. Для решения задачи слеует научиться переводить с одного языка на другой и, главное, обеспечить саму эту возможность.

\medskip

Возможность перевода с языка контактов на глобальный язык связана с утверждениями типа теоремы Х.~Гудмана-Штраусса \cite{GoodmanStrauss}.  Теорема Х.~Гудмана-Штрауса утверждает, что любую плоскую мозаику, связанную с подстановочной системой (плоской $\DOLL$-системой) можно задать локальными правилами. Наш подход, использующий реберную структуру, ближе к подходу Тома Ферни и Николя Бедериде \cite{Fernique}.

\newpage

\subsection{Апериодические мозаики. ``Демо-версия'' доказательства}

\subsubsection{Апериодические замощения}
Наличие апериодических мозаик связано со следующим обстоятельством. Рассмотрим работу машины Тьюринга на клетчатой ленте.  Движение головки и обмен сигналами можно перевести на язык укладки плиток на фазовой диаграмме ``пространство-время''. Из этого выводится алгоритмическая неразрешимость проблемы дополнения заданного набора плиток до замощения плоскости.

Алгоритмическая неразрешимость проблемы замощения, когда изначально ничего не положено, доказывается технически сложнее. Здесь ``пространство'' и ``время'' перемешиваются. Строятся замощения, в некоторых местах которых вынуждается один шаг машины Тьюринга, в других -- два шага и т.д. Линии, отвечающие работе машины, образуют некую структуру. Дальнейшее совершенствование и упрощение конструкций приводит к апериодическим мозаикам, ставшими классическими (плитки Аммана, мозаики Пенроуза и др.).

В $1961$ году Хао Вангом \cite{Wang} были рассмотрены квадратные плитки с разноцветными сторонами. Разные плитки можно прикладывать
друг к другу сторонами одного цвета. Был поставлен вопрос, существуют ли конечные наборы таких плиток,
с помощью которых могут быть получены только непериодические замощения плоскости. Первым такой набор был построен
Робертом Бергером \cite{Berger}, основная идея состояла в том, что замощения моделировали работу машины Тьюринга. При этом
использовалось несколько тысяч плиток. Позднее были придуманы наборы, содержащие небольшое количество плиток. Например,
интересна конструкция Рафаэля Робинсона \cite{Robinson}. В этой связи можно также упомянуть мозаики Пенроуза и Амана, состоящие
из многоугольников, которыми, при заданных условиях,  можно замостить плоскость лишь непериодически.
Различные примеры апериодических мозаик есть также в работе \cite{FerniqueIvanovBlvMitrafnv}.

%Впервые такой набор был построен Робертом Бергером \cite{Berger}. В
%дальнейшем были построены более простые примеры, например, Рафаэлем Робинсоном \cite{Robinson}. Широко известна также мозаика Пенроуза.

\begin{figure}[hbtp]
\centering
\includegraphics[width=0.4\textwidth]{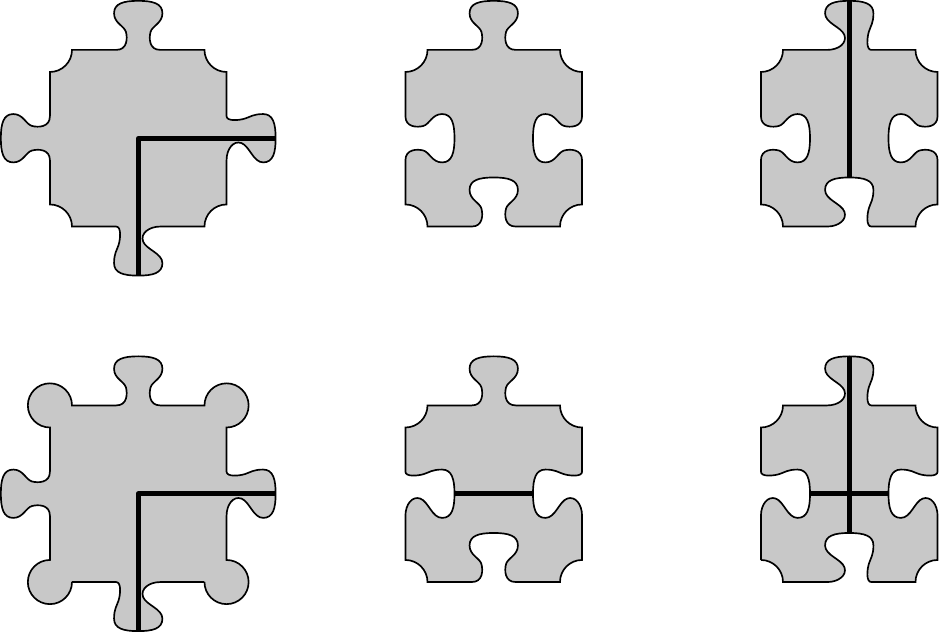}
\caption{Набор плиток Робинсона.}
\label{fig:robinson_bw}
\end{figure}

\subsubsection{Иерархия и апериодичность} \label{dol}

Рассмотрим другой способ получения непериодических замощений. Пусть имеется конечное число типов плиток и мы задаем
правила, по которым из нескольких маленьких плиток можно составлять большие макроплитки тех же типов.

\begin{quote}

{\bf Пример.} {\it Плитки могут быть квадратами $A$ и $B$, при этом, чтобы составить квадрат $A$ второго уровня, нужно взять
четыре квадрата $A$, $A$, $A$, $B$. А чтобы составить квадрат $B$, нужно взять четыре квадрата $B$ $A$ $B$ $A$.}

\end{quote}

Таким образом, получается иерархическая система. Каждую плитку можно разбить на требуемое число уровней иерархии.
Можно показать, что получаемое замощение будет непериодично. Аналогичный способ построения используется в подстановочных
системах, например, c помощью подстановок $1\rightarrow 10$, $0\rightarrow 01$ получается бескубное слово
$$1001011001101001 \dots$$

\medskip

Оказывается, язык граничных условий и язык иерархий схожи. А именно, любую иерархическую систему можно задать конечным числом граничных условий.

Пусть имеется конечное число типов плиток, причем заданы правила иерархии, по которым плитка уровня $n$ составляется из
нескольких плиток уровня $n-1$. Тогда для начального набора плиток первого уровня можно задать конечную систему
граничных условий так, чтобы задавалась мозаика, получаемая при иерархическом способе задания.

Иерархичность системы плиток гарантирует непериодичность замощения. Вследствие этого можно задавать с помощью граничных
условий мозаики, которые будут с гарантией непериодическими.

\subsubsection{Демо-версия доказательства}

\begin{figure}[hbtp]
\centering
\includegraphics[width=0.8\textwidth]{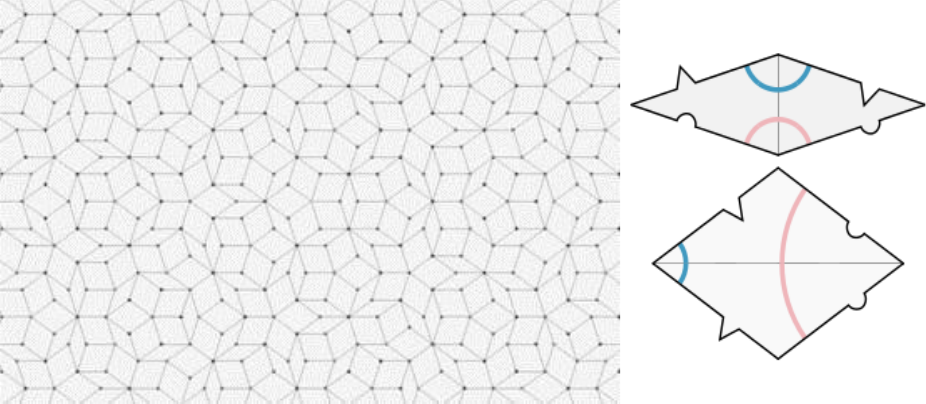}
\caption{Мозаика Пенроуза.}
\label{penrose}
\end{figure}

Рассмотрим одну из классических мозаик -- знаменитую {\it мозаику Пенроуза} (см. рисунок~\ref{penrose}).
\medskip

\paragraph{\bf Конструкция мозаики Пенроуза.}  Используются плитки двух видов -- толстый и тонкий ромбы. Есть граничные условия: стороны каждого ромба раскрашены в две пары цветов. Соприкасаться два ромба могут только сторонами двух цветов, образующих пару. На рисунке~\ref{penrose} цвета обозначены внешними и внутренними насечками двух типов.
\medskip

\medskip

\paragraph{\bf Конструкция полугруппы мозаики Пенроуза.}
Легко видеть, возможно конечное количество типов узлов-вершин где
сходятся несколько плиток. Кроме того, ребра в мозаике могут иметь десять возможных направлений. Можно выписать все возможные типы узлов и обозначить
их буквами алфавита. Теперь последовательность букв (слово) будет кодировать
последовательность узлов, которые мы проходим вдоль пути.

Для каждого узла введем также два параметра: по какому ребру мы в него вошли и по какому ребру мы вышли (включая информацию о цветах ребер). Теперь можно расширить алфавит, добавив буквы для всевозможных сочетаний параметров для разных типов узлов. Некоторые последовательности букв с гарантией не смогут представлять путь на мозаике (например, если ребро, по которому мы пришли в узел не соответствует по цвету ребру, по которому мы вышли из предыдущего узла). Такие последовательности мы будем заносить в список запрещенных.

Некоторые пути на мозаике можно объявить {\it эквивалентными}: например путь по двум соседним ребрам ромба эквивалентен пути по другой паре сторон. Можно выписать
все такие эквивалентности для разных вариантов получающихся узлов и получить
полный конечный список. Тогда, с помощью локальных замен, можно переводить одни пути в другие.

Таким образом, с мозаикой Пенроуза можно связать конечно определенную полугруппу, где словам соответствуют пути. Запрещенные пути -- это нулевые слова. Можно также запретить пути, которые не являются кратчайшими между двумя узлами. Для этого нужно сначала внести в список запрещенных короткие такие пути. После этого можно показать что с помощью локальных замен любой некратчайший путь приводится к виду содержащему запрещенный короткий подпуть.

В случае, если бы можно было бы показать, что любой путь, не вкладывающийся в мозаику, приводится к нулю с помощью указанных локальных правил, получающаяся полугруппа была бы нильполугруппой, так как на мозаике не лежит периодических путей.

\medskip

\paragraph{\bf Почему нужна другая мозаика.}
Проблема заключается в том, что на мозаике Пенроуза есть пути, которые недостаточно сильно меняются локальными заменами, то есть, недостаточно ``извиваются''. Это приводит к тому, что можно сконструировать путь, каждый локальный кусок которого может быть вложен в мозаику, но весь путь не может быть вложен. Локальные замены меняют его незначительно и преобразовать его в достаточной мере, чтобы диагностировать несоответствие мозаике, не получается.

В целом, каждый путь можно трактовать как массив информации об его окрестности. Когда мы производим локальные замены, происходит перенос информации вдоль пути. В случае, если пути можно шевелить незначительно, канал переноса информации будет ограничен, что не позволит выявить ситуацию, когда длинный кусок пути не может являться частью мозаики. (Например, когда путь -- это степень разрешенного слова.)

\subsection{Язык контактов vs язык соотношений. Подклейки.}
Итак, если осуществлять нашу программу на базе классических мозаик, то возникают трудности, связанные с тем, что в некоторых направлениях геодезические пути не поддаются изгибу и вокруг них ничего из переброски не наращивается. А именно достаточно протяженный в двух измерениях кусок мозаики обеспечивает вычислительный процесс.

Перевод с языка соотношений (т.е. движения геодезического пути) на язык контактов вешь более сложная и обеспечить его возможность не так просто. Чтобы локальное шевеление геодезической заполнило достаточную область, пространство должно быть {\it равномерно-эллиптическим}. Пространство называется {\it равномерно-эллиптическим}, если любые две точки $A$ и $B$ на расстоянии $D$ соединяются системой геодезических,  образующих диск ширины $\lambda\cdot D$ для некоторой глобальной константы $\lambda>0$. (Априори не очевидно, что такое пространство существует.) Такое пространство следует собрать из мозаики и установить аналог теоремы Гудмана-Штраусса для этой сборки.

Равномерно-эллиптическое пространство строится из подстановочной системы указанной на рис. \ref{fig:rule} при этом следует озаботится тем, чтобы степени узлов были бы ограничены, что усложняет конструкцию, вынуждая сделать композицию разных подстановок \ref{fig:subst}.

Далее, когда путь проходит через узел и по обе стороны от него есть пара плиток более низкой ступени иерархии, сходящиеся в этом узле, необходимы {\it подклейки}. Бесконечное число иерархических уровней приводит к тому, что к ребру может примыкать неограниченное число подклеек. При этом необходимо избежать следующей ситуации: мы вышли в подклейку и вернулись. Далее продеформировали участок пути в подклейку так, что его ``голова'' снова коснулась исходного ребра, а вход и выход в ``подклейку'' не изчезли. Далее возникший участок по ребру может уйти в иную подклейку там снова изогнуться и коснуться исходного ребра и т.д. Конструкция специально построена таким образом, чтобы этого избежать.
То есть {\em все организовано так, что путь может вернуться из подклейки ТОЛЬКО в исходную точку}. Далее размеры подклек к плиткам, подклеек к подклейкам и т.д. экспоненциально убывают  достаточно быстро. Тем самым, помимо всего прочего, обеспечивается гомотопическая тривиальность комплекса, образованного ``подклейками''.

\subsection{Последовательная канонизация. Завершение доказательства.}
Доказательство завершается так. Рассмотрим слово. Надо его либо привести к пути на комплексе (тем самым оно непериодично), либо к нулю. 
Оно последовательно приводится к {\it $k$-каноническому виду} и при этом $k$ растет, т.е. оно состоит из учасков границ плиток $k$-го уровня иерархии, кроме начала и конца, целиком содержащихся в плитке уровня $k-1$ (возможно в подклейке $(k-1)$-го уровня) и при этом являющихся $(k-1)$-каноническими. 

Далее с промощью соотношений (в том числе отвечающих выходу в подклейку) проверяется согласование соседних участков пути, располагаемых на границах плиток $k$-го уровня, их нахождение в $(k+1)$-ом уровне или возможность находится в двух соседних плитках $(k+1)$-го уровня. Отсутствие согласование означает локальную невозможность располагаться преобразованного пути на комплексе -- т.е. вывод равенства нулю. 
Вначале рассматривается плоский процесс. Если имеются команды входа в подклейку и  выхода из нее -- то проверяется ситуация между ними и производится сокращение.

Согласование означает возможность преобразования к $(k+1)$-каноническому виду, после чего процесс повторяется. Он заканчивается канонической формой слова (в нашем случае каноничность не означает однозначности, поскольку можно выбирать разные стороны плиток). Наверное, существует такой выбор -- но он более сложен и ситуация похожа на теорию канонической формы элементов в гиперболических (или более обще, в группах с неположительной кривизной) развитую И.~А.~Рипсом, когда вначале строится предканоническая форма (для данного уровня иерархии), потом осуществляется выбор.

\subsection{Специфика геометрических методов в полугрупповом случае}
Интерпретеция соотношений в полугруппах через ``rewriting diagram'' похожа на применение диаграмм Ван-Кампена в теории групп.
Рассмотрение мозаик также родственно работе с диаграммами Ван-Кампена в теории групп.
Использование путей по ребрам в плиточных замощениях и групп связанных с ними, встречается в работах Дж.~Конвея \cite{Conway}.
С помощью групповой техники он показал отстутствие различных разбиений. Вот модельный

\medskip
{\bf Пример.} {\it Рассмотрим шахматную доску, из которой вырезаны уголки. Тогда
ее нельзя разбить на доминошки $1\times 2$}.
\medskip

Предположим, что такое разбиение существует. Рассмотрим доминошки как клетки в диаграмме Ван-Кампена. Букве $a$ сопоставим вертикальную стрелку, букве $b$~-- горизонтальную, обратным буквам -- обратные стрелки. Вертикальная доминошка отвечает соотношению $a^2ba^{-2}b^{-1}=e$, а горизонтальная -- $b^2ab^{-2}a^{-1}=1$. Границе квадрата с вырезанными углами отвечает путь $a^7bab^7a^{-7}b^{-7}a^{-1}$. Если сопоставить элементу $a$ транспозицию элементов $1$ и $2$, а элементу $b$ -- транспозицию $2$ и $3$, то $a^2=b^2=e$ и соотношения выполнены, в то же время как $a^7bab^7a^{-7}b^{-7}a^{-1}=(ab)^4=ab\ne e$ есть трехчленный цикл.

Тем же методом решаются многие другие олимпиадные задачи на раскраску, но в то же время и такие, которые раскраской не делаются (подробности см. \cite{BelvIvanvMalstvMitrfnvKharitnv}).

Тем не менее мозаики и диаграммы Ван-Кампена существенно различаются. Например, свойства величин углов имеют далеко не полное отражение. Аналогия между мозаиками и группами довольно глубокая, но не вполне формализована. Имеется классический результат об алгоритмической разрешимости проблемы равенства слов в группе с одним соотношением. Имеется и аналогичный вопрос об алгоритмической разрешимости проблемы разбиения плоскости транслятами одной фигуры. Для связных фигур такая алгоритмическая разрешимость доказана, а в общем случае вопрос остается открытым.

Далее, при разбиении квадрата на домино возникают подквадратики $2\times 2$, разбитые на пары доминошек. Такую пару можно повернуть на $90^o$ и сделать {\it флип}. Такими флипами можно от одного разбиения перейти к любому другому. Аналогичный факт верен и разбиения на $k$-миношки и для многих других разбиениях. Данный факт родственен конечной порождаемости группы $\pi_2$ комплекса (с учетом действия фундаментальной группы $\pi_1$). Хорошо бы прояснить аналогию.

В полугруппах геометрические идеи работают не совсем так, как в группах, и в нашем случае эффекты ``углов'' лучше представлены.   По всей видимости, полугрупповая теория может пролить свет на мозаики. Более того, есть некая близость теории полугрупп к кольцам как к ``квантовой взвеси'', так что и кольцевая теория, если будет построена, даст дополнительное понимание.
Возможно, наш  подход окажется полезным и для других построений в полугруппах и кольцах.

\subsection{Апериодические замощения в нашей конструкции}
Как уже говорилось,  существуют конечные наборы многоугольников (плиток) которыми можно
замостить плоскость лишь непериодически. Методы построения таких мозаик, как правило, опираются на иерархические правила построения: задаются универсальные правила построения плиток уровня $n+1$ из плиток уровня $n$, для нескольких типов плиток $A_1,\dots,A_k$.
Х.~Гудман-Штраусс \cite{GoodmanStrauss}  показал, что иерархические мозаики можно получить, задав конечное число локальных правил.
Таким образом, локальные условия (конечность набора) могут приводить
к глобальному эффекту (непериодичности замощения).

Итак, основной задачей является конструирование мозаики, в которой указанные выше трудности были решены. Для этого она должна обладать несколькими свойствами:

\medskip
\begin{description}
  \item[1. Локальная конечность] Речь идет конечности возможных типов узлов, конечности изначально задаваемых пар эквивалентных путей, а также конечности списка запрещенных путей.
  \item [2. Возможность ``шевеления'' любого пути] Любой достаточно длинный путь, соединяющий узлы $A$ и $B$ может быть переведен локальными заменами в другой путь, отличающийся достаточно сильно от начального. Например, мерой отличия путей может являться максимальное расстояние между соответствующими точками путей.
  \item [3. Апериодичность] На мозаике не должно быть путей, отвечающих периодическим словам.
\end{description}

 \medskip
Как уже говорилось, мы используем геометрическую интерпретацию для алгебраических построений. Запрет для
двух (или более) плиток находиться рядом друг с другом схож с запретом для двух букв стоять рядом в разрешенном слове.
Возникает интерпретация слова как последовательности плиток на выложенной мозаике.
Апериодичность мозаики приводит к непериодичному характеру таких ``плиточных слов''.  В свою очередь,
непериодических замощений можно добиться, если применять иерархический способ построения.

В связи с этим используются  языки плиточных примыканий и иерархий. Эти языки во многом схожи, например, заданные
правила иерархии, когда плитки $A_1,\dots,A_k$ уровня $n+1$ составляются из плиток $A_1,\dots,A_k$ уровня $n$, можно выразить
с помощью конечного множества локальных правил для набора $A_1,\dots,A_k$. Эти локальные правила будут порождать те же непериодические мозаики чти и исходные правила иерархии.

\medskip

Дальнейшее развитие этих языков приводит к появлению более универсального языка путей на графе. То есть плиточная мозаика
рассматривается как граф, где вершины это узлы мозаики, а ребра -- границы плиток. В этом смысле понятие плитки
можно обобщить, рассматривая их уже как локальные подграфы из которых, с помощью локальных правил, можно составлять граф,
покрывающий плоскость.

Аналогом буквы будет тип вершины графа, аналогом слова -- путь,
проходящий через несколько вершин. Аналогом соотношения будет эквивалентность между путями с общими концами: например,
в простом $4$-цикле $ABCD$ выполнено соотношение $ABC=ADC$. Помимо таких, есть также мономиальные соотношения, выражающие
идею о невозможности существования какого-то пути на мозаике. Также, для обеспечения необходимого контроля над множеством
ненулевых слов, вводятся мономиальные соотношения, обнуляющие слова, соответствующие некратчайшим путям.
Оказывается, что при этом можно обойтись конечным числом соотношений. Немономиальные соотношения, при этом, не меняют длины слова.

\medskip

Язык путей на мозаике-графе позволяет выразить те же концепции и определить те же мозаики, что и языки иерархических плиток или граничных условий.
Таким образом возникает связанная с мозаикой конечно определенная полугруппа с набором свойств, индуцированных мозаикой.

\begin{figure}[hbtp]
\centering
\includegraphics[width=0.3\textwidth]{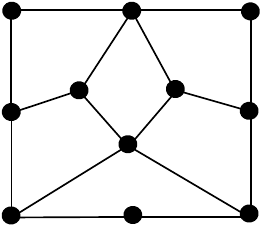}
\caption{Правило разбиения.}
\label{fig:rule}
\end{figure}

Для построения нильполугруппы используется мозаика сгенерированная с помощью следующего иерархического правила разбиения (рисунок~\ref{fig:rule}).

От мозаики нужно потребовать ряд дополнительных свойств. В частности, любой длинный путь должен допускать
возможность ``шевеления'' то есть, локального преобразований над ним, позволяющих в достаточной мере менять его. Для
достижения этого свойства к плоской мозаике производятся ``подклейки'', представляющие собой небольшие плоские подграфы, не лежащие
в исходной плоскости, и позволяющие обходить ``узкие места'' на исходном плоском графе. Структура подклеек также имеет иерархическую
природу и так же может быть задана на языке преобразований путей конечным образом.

В итоге мы получаем граф, обладающий набором важных для нас свойств, в котором любой путь, соединяющий произвольные точки $A$ и $B$,
является членом семейства геодезических эквивалентных друг другу путей, соединяющих эти точки. Причем эквивалентность двух путей
из этого семейства может быть получена путем цепочки локальных преобразований переводящих один путь в другой. ``Шевеление''
пути играет роль передачи информации. Фактически, определяющие соотношения задают правила передачи информации по пути.
Если задан произвольный длинный путь, мы можем начать работать над ним, совершая локальные преобразования. При этом возможны две альтернативы:

\begin{enumerate}
  \item В результате этой работы можно получить внутри некратчайший подпуть, либо подпуть указанный в числе запрещенных.
В этом случае наш путь представляет нулевое слово.
  \item В результате этой работы мы восстанавливаем некоторый кусок мозаики, внутри которого лежит пучок геодезических путей, эквивалентных нашему.
\end{enumerate}

\medskip

Мозаика не может содержать в себе путей, выражаемых периодическим словом. То есть, все периодические слова могут быть приведены к нулю
 локальными преобразованиями. При этом геодезические пути, лежащие на мозаике не приводятся к нулю, и могут иметь любую длину.
Таким образом полугруппа, соответствующая построенному графу-мозаике, будет конечно определенной нильполугруппой.

\subsection{Конечно определенные полугруппы} \label{semigroups}
%\medskip
Рассмотрим алфавит $\Omega$, состоящий из конечного числа букв. Слова, составленные из букв образуют полугруппу
относительно операции приписывания одного слова к другому. Кроме того, есть специальная буква $0$ (ноль), такая что
$W0=0W=0$ для любого слова $W$ из полугруппы. Определяющим соотношением в полугруппе считается равенство вида
$W_1=W_2$, где $W_1$ и $W_2$ -- некоторые слова, причем одно из них может быть нулем (или нулевым словом). В этом случае,
соотношение называется {\it мономиальным}. Для натурального $n$ степенью слова $W^n$ называется слово $WW\cdots W$, где $W$ выписано
$n$ раз подряд. Элемент полугруппы (слово) называется {\it ниль-элементом}, если существует натуральное $n$, такое что
$W^n=0$. Если все элементы полугруппы являются нильэлементами, то вся полугруппа называется нильполугруппой.
Разумеется, наша полугруппа содержит {\it ноль}.

Мы пользуемся геометрической интерпретацией: {\it буквам} отвечают вершины специального
графа, а {\it словам} -- пути в этом графе. И если слово не может быть представлено путем на графе, оно всегда будет
приводиться к нулю.

%\medskip

\subsubsection{Плитки и непериодические замощения} \label{plitki}

Начнем с чистой геометрии.
Мы задаем конечный набор граничных условий (как можно прикладывать плитки друг к другу),
и через задание локальных условий достигается глобальный эффект. Задание граничных условий схоже с заданием определяющих
соотношений в полугруппе, этим объясняется интерес к плиткам и замощениям.

\subsubsection{Узлы на мозаике и пути по границам плиток} \label{ways}

Помимо языков иерархических систем и граничных условий есть еще один подход. Можно рассматривать мозаики с
точки зрения путей по границам плиток. Это естественный шаг для перехода к полугруппе, так как пути по мозаике
логично отождествляются со словами в подходящем алфавите.

Назовем {\it узлами} вершины мозаики, где сходится несколько плиток. Ясно, что в плоской мозаике с конечным набором
типов плиток возможно конечное число видов узлов. Обозначим их буквами конечного алфавита. Последовательность
букв (слово) соответствует пути на мозаике. Некоторые слова могут вообще не встречаться на мозаике, а другие встречаться в разных местах.

Третий подход к заданию мозаик состоит в определении конечного списка невозможных путей, а также задание конечного
числа пар эквивалентных путей. С помощью такого языка можно задавать мозаики также, как и с использованием
иерархических систем плиток или граничных условий на плитки.

Мы будем строить такую мозаику в два приема. На первом этапе построим ее {\it плоскую часть}, которую будем считать базовой. Она будет строиться в виде иерархического бесконечного графа. На втором этапе мы покажем, как к построенному графу производятся {\it подклейки} -- плоские и конечные подграфы, лежащие в другой плоскости, но имеющие с базовым графом несколько общих ребер. В результате получится трехмерный комплекс.

\medskip

Подклейки нужны для того, чтобы появился маршрут, альтернативный проходу по ребру, к которому производится подклейка. В этом случае, становится возможным перенос информации вдоль пути и мы можем выявить ситуации с невозможным на мозаике путем.

\medskip

В дальнейших главах мы построим мозаику, где выполнены указанные выше свойства.

%\medskip

\subsection{Схема зависимости материала и состав глав}

Доказательство основного результата достаточно объемное. В этом параграфе мы кратко опишем состав глав, составляющих работу.

\medskip

Во {\bf введении} обсуждаются проблемы бернсайдовского типа и конечной определенности и связанные с этим вопросы.

В главе~\ref{ScFnDefInfGen} Обсуждается общая стратегия доказательства а также смежные вопросы.

В главе~\ref{geom} (``Геометрическая структура комплекса'') обсуждаются чисто геометрические свойства комплекса, лежащего в основе построения.

В главе~\ref{pathsemigroup} (``Строение полугруппы путей'') и главе~\ref{coding_section} (``Кодировка вершин и путей на комплексе'') вводится параметризация, благодаря которой структура построенного геометрического комплекса оказывается связанной со структурой полугруппы. В частности, вводится алфавит букв, кодирующих вершины и ребра на комплексе.

В главе~\ref{functions} (``Функции на узлах и структура цепей'') обсуждаются особенности комплекса, благодаря которым можно вычислять параметры вершин, по известным другим вершинам, в случае если вершины занимают друг относительно друга определенное положение. В этой главе мы готовим возможности вычисления близлежащих вершин комплекса.

 В главе~\ref{count} (``Оценка количества букв'') проводится оценка, сколько букв используется в алфавите. Если читателю этот вопрос не так интересен, эту главу можно пропустить без ущерба для дальнейшего изложения.

 Глава~\ref{flip_section} (``Разбор случаев расположения путей'') является самой объемной в работе. В ней доказывается, что по параметрам трех вершин, лежащих в вершинах некоторой плитки на комплексе, всегда можно вычислить параметры четвертой вершины. Это основное свойство комплекса, на котором основано построение полугруппы. Объем обусловлен большим количеством случаев расположения вершин. При первоначальном прочтении можно ознакомиться с несколькими из них. При рассмотрении локальных преобразований $7$, $8$, $9$, $10$ более подробно рассмотрен случай $\mathbb{C}1$-цепи, в частности, в явной форме выписана таблица характеризации локального преобразования для  $\mathbb{C}1$-цепи по виду кода.

 В главе~\ref{pasting_section} (``Локальные преобразования при выходе в подклееную макроплитку'') также разбираются случаи расположения путей, но в этот раз в ситуации, когда происходит выход в подклееную часть (предыдущая глава посвящена плоской ситуации). Аналогично, при первоначальном прочтении можно ознакомиться с несколькими случаями.

 Глава~\ref{canonic} (``Приведение к стандартной форме'') завершает доказательство. С помощью свойств комплекса, рассмотренных в предыдущих главах, показывается как можно преобразовать слово. В частности, показывается, что слово, содержащее подслово в девятой степени, приводится к нулю.

 Глава~\ref{Appendix} (``Приложение. Подсчет числа окружений'') носит технический характер. Выделяются возможные значения четверок-комбинаций типов ребер.

 В главе~\ref{future} (``Вопросы и перспективы'') обсуждаются смежные вопросы, представляющие интерес в связи с полугрупповой конструкцией.

\medskip

%В настоящей статье развивается новый подход к контролю над вводимыми соотношениями. Он основан на
%следующих соображениях.

%\subsection{Плитки, мозаики и полугруппы, уточнение плана доказательства}

\section{Геометрическая структура комплекса} \label{geom}

\subsection{Формальное задание полугруппы}

В этой главе мы обсудим формальную сторону записи полугруппы, то есть какие алфавиты используются для записи слов. Структура алфавита полугруппы сильно зависит от геометрии мозаики (главы~\ref{coding_section} и~\ref{functions}).

\medskip

Алфавит, отвечающий образующим полугруппы (в дальнейшем, просто {\it алфавит}) будет состоять из букв, кодирующих вершины и букв, кодирующих входящие и выходящие ребра в вершины. Слово из букв кодирует путь -- последовательность $X_1e_1e_2X_2e_3e_4X_3 \dots$, где $X_i$ -- буквы, кодирующие вершины, а $e_i$ -- буквы, кодирующие входящие и выходящие ребра.

Для вершины определяется несколько параметров: тип, уровень, окружение, информация, флаг подклейки. В целом эти параметры характеризуют геометрическое положение вершины внутри построенного комплекса. Для каждого параметра мы задаем конечное число значений. Буква, кодирующая вершину, представляет собой сочетание значений всех параметров. В некоторых случаях мы будем рассматривать множество букв, имеющих фиксированное значение некоторых параметров. В частности, определяющие соотношения вводятся для пар слов заданной формы, где мы фиксируем часть параметров для некоторых букв.

Для входящих и выходящих ребер используется свой алфавит, содержащий конечное число букв. Удобно ввести свой набор таких букв для вершины каждого типа, это будет сделано ниже.

\medskip

В следующих параграфах мы последовательно определим каждый из параметров вершин, вместе с объяснением геометрического смысла. Основная задача состоит в реализации одной из следующих возможностей:

{\bf 1.} По выданному слову установить, что данное слово не кодирует никакой путь (или кодирует запрещенный) и в этом случае привести слово к нулю, применяя определяющие соотношения;

{\bf 2.} По выданному слову установить, что данное слово кодирует некоторый путь на построенном геометрическом комплексе. При этом применение определяющих соотношений будет приводить к преобразованию пути в эквивалентный ему.

\medskip

Смысл задания параметров для вершин (тип, уровень, окружение, информация, флаг подклейки) состоит в том, чтобы знание набора значений этих параметров для простого пути $P_1$ из трех вершин и двух ребер между ними позволяло восстановить набор значений этих параметров для другого пути $P_2$ из трех вершин и двух ребер, имеющего начало и конец такие же как у $P_1$. Введение определяющих соотношений представляет собой перебор различных случаев расположения таких пар путей и приравнивание их кодов.

Другая часть определяющих соотношений состоит в обнулении слов, которые либо не представляют никакой путь, либо кодируют путь, который мы хотим сделать нулевым.

После введения таких определяющих соотношений появляется возможность преобразования пути через операцию с его кодом. То есть геометрическому преобразованию пути будет соответствовать преобразование слова в полугруппе. После этого можно будет показать, что слова, содержащие девятую степень любого слова, можно привести к нулю. Кроме того, можно будет убедиться, что на комплексе будут существовать сколь угодно длинные пути (и, соответственно, слова), не приводимые к нулю.

\medskip

Определим {\it алфавит типов вершин}, состоящий из символов

$$\{ \mathbb{A},\mathbb{B}, \mathbb{C}, \mathbb{U}, \mathbb{R}, \mathbb{D}, \mathbb{L}, \mathbb{UL}, \mathbb{LU}, \mathbb{UR}, \mathbb{RU}, \mathbb{DL}, \mathbb{LD}, \mathbb{DR}, \mathbb{RD} \}. $$

\medskip

Обозначим его как $\{ \mathbf{nodes} \}$.

Пусть $\{\mathbf{edges}\}$ -- алфавит из $20$ символов:

$$\{\mathbf{left}, \mathbf{top}, \mathbf{right}, \mathbf{bottom}, \mathbf{1A}, \mathbf{1B}, \mathbf{2A}, \mathbf{2B}, \dots, \mathbf{8A}, \mathbf{8B} \}.$$

Будем называть его {\it алфавит типов ребер}.

Рассмотрим множество упорядоченных четверок $( x_1, x_2, x_3, x_4)$, где $x_i\in \{ \mathbf{edges}\}$. Всего в нем $20^{4}$ элементов. Мы будем использовать только $210$ из них. Назовем эти $210$ четверок {\it алфавитом окружений макроплиток} и будем обозначать его как $\{ \mathbf{fours} \}$.

Список используемых четверок получается из комбинаторных свойств конструируемой мозаики. Подробно про это написано в Приложении, глава~\ref{Appendix}, предложения~\ref{combos_edges} и~\ref{combos_inner}.

Остальные параметры мы определим в связи с соответствующей геометрической структурой в последующих параграфах.

\medskip

\subsection{Иерархическое построение}

Построим {\it комплекс уровня $n$} с помощью итерационного процесса. На каждом шаге мы будем иметь дело со графом, состоящим из вершин и ребер, причем простые циклические пути из четырех ребер будем называть {\it плитками}.

\medskip

\begin{figure}[hbtp]
\centering
\includegraphics[width=0.4\textwidth]{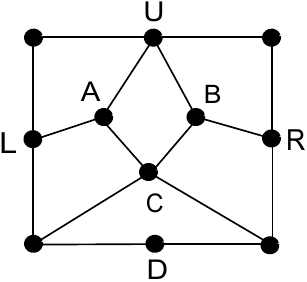}
\caption{Комплекс второго уровня.}
\label{fig:firststep}
\end{figure}

{\it Комплекс уровня $1$} представляет собой простой цикл из четырех вершин, {\it Комплекс уровня $2$} это граф на рисунке~\ref{fig:firststep}. Вершины $A$, $B$, $C$ будем называть {\it внутренними} или
{\it черными}, а вершины $U$, $L$, $R$, $D$ -- {\it боковыми}. То есть,  на втором этапе у нас имеется шесть плиток. Будем называть их, согласно их положению, левой верхней, левой нижней, правой нижней, правой верхней, средней, нижней.

Будем считать, что вершины $A$, $B$, $C$, $U$, $L$, $R$, $D$ имеют {\it нулевой} уровень.

\medskip

С комплексом будут производиться итерации двух типов: {\it разбиение} и {\it подклейка}.

\begin{figure}[hbtp]
\centering
\includegraphics[width=0.5\textwidth]{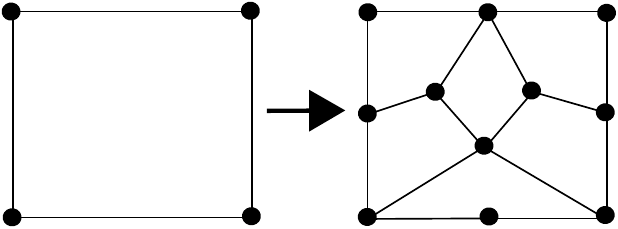}
\caption{Разбиение.}
\label{fig:subst}
\end{figure}

\begin{figure}[hbtp]
\centering
\includegraphics[width=0.5\textwidth]{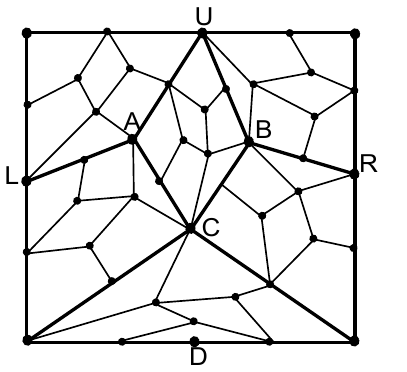}
\caption{Макроплитка третьего уровня.}
\label{fig:level2}
\end{figure}

\begin{definition}
При операции {\it Разбиения} каждая плитка комплекса разбивается на шесть плиток согласно правилу на рисунке~\ref{fig:subst}. При этом эта старая плитка перестает быть просто плиткой и становится макроплиткой второго уровня, создаются новые четыре боковые и три внутренние вершины. Одна из четырех возможных ориентаций определяется согласно положению (одному из шести) разбиваемой плитки в ее родительской макроплитке второго уровня (рисунок~\ref{fig:level2}).
\end{definition}

\medskip

\begin{definition}
{\it Макроплиткой} уровня $n$ назовем: для $n=1$ -- обычную плитку (простой $4$-цикл), для $n>1$ -- результат применения операции Разбиения к макроплитке уровня $n-1$.
\end{definition}

Все создаваемые при разбиении вершины получают уровень на $1$ больший, чем максимальный уровень до разбиения, в частности, если создается вершина в  середине ребра, она получает уровень на $1$ больший, чем максимальный уровень двух концов ребра.

\begin{definition}
При {\it подклейке} мы рассматриваем все пути $X_1X_2YZ_2Z_1$ из пяти  вершин (и четырех ребер), такие что:

1) Вершины $X_1$, $Y$, $Z_1$ {\bf не} являются тремя углами из четырех никакой макроплитки;

2) Вершины $X_1$, $Z_1$ являются боковыми вершинами  уровня $k-1$, где $k$ -- максимальный уровень вершины в комплексе;

3) Вершины $X_2$ является серединой $X_1Y$, а $Z_2$ -- серединой $Z_1Y$, то есть,  $X_2$ и $Z_2$  являются боковыми вершинами уровня $k$, созданными в комплексе самыми последними, при последней операции {\it Разбиения}.

4) Вершина $Y$ имеет уровень $k-2$.

\end{definition}

\medskip

\begin{figure}[hbtp]
\centering
\includegraphics[width=0.5\textwidth]{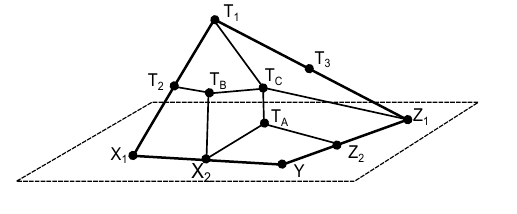}
\caption{Подклейка.}
\label{fig:pasting}
\end{figure}

Далее, для каждого такого пути создаются шесть новых вершин $T_1$, $T_2$, $T_3$, $T_A$,  $T_B$,  $T_C$, не лежащих в плоскости $X_1YZ_1$, а также проводятся новые ребра $X_1 T_2$, $X_2 T_A$, $X_2 T_B$, $T_2 T_B$, $T_C T_B$, $T_C T_A$, $T_2 T_1$, $T_3 T_1$, $T_3 Z_1$, $T_C Z_1$, $T_A Z_2$  и появляется новая {\it подклееная } макроплитка $X_1YZ_1T_1$ (рисунок~\ref{fig:pasting}).

\medskip

Созданные вершины будем называть {\it подклеенными}, причем $T_1$ присваивается уровень $k-1$, а вершинам $T_2$, $T_3$, $T_A$,  $T_B$,  $T_C$ -- уровень $k$.
Вершину $T_1$ будем считать угловой, вершины $T_2$, $T_3$ -- боковыми, а вершины $T_A$,  $T_B$,  $T_C$ -- внутренними.

\medskip

{\bf Примечание 1.} Подклееная макроплитка выглядит так же, как если бы к плитке $X_1YZ_1T_1$ применили Разбиение, считая сторону $X_1Y$ верхней. Также можно считать, что в этой макроплитке $T_2$ середина стороны $T_1X_1$, а $T_3$ -- середина стороны $T_1Z_1$.

\medskip

{\bf Примечание 2.} Вершину $Y$ будем называть {\it ядром подклейки}.

\medskip

\begin{definition}
 Комплексом $n+1$ уровня становится комплекс $n$ уровня, к которому применены {\it Разбиение} и затем {\it Подклейка}.
 \end{definition}

\medskip

{\bf Замечание 1.} После разбиения, примененного к комплексу первого уровня, образуется граф на рисунке~\ref{fig:level2}, на котором нет пяти вершин, где можно было бы выполнить подклейку. То есть, подклейки начинают проводиться после второго разбиения.

\medskip

{\bf Замечание 2.} Начиная с комплекса второго уровня, подклейки начинают применяться и к подклееным плиткам. Таким образом, возникают подклейки к частям подклееных когда-то макроплиток и т.д. Кроме того, начиная с третьего уровня комплекса путь $X_1YZ_1$ из
определения операции подклейки может частично лежать в базовой плоскости, а частично в подклееной части (если $Z_1$ появилась при предыдущей подклейке в качестве вершины $T_2$ или $T_3$).

\medskip

\begin{proposition}[Об ограниченности роста степени вершины] \label{growth_bound}

Для каждой вершины $X$ существует такое натуральное $N$, что
начиная с уровня комплекса $N$, степень (число входящих ребер) вершины $X$ не меняется, то есть она одинакова для комплексов уровня $N$ и $N+k$ для любого натурального $k$.

\end{proposition}

\medskip

\begin{figure}[hbtp]
\centering
\includegraphics[width=0.5\textwidth]{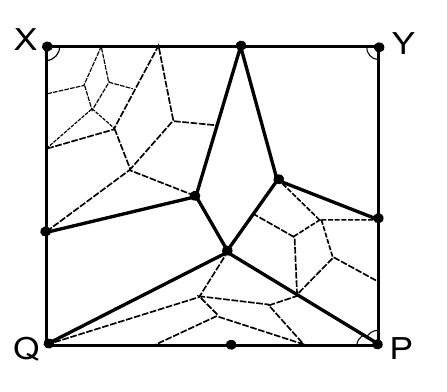}
\caption{Простые углы отмечены дугами.}
\label{fig:tileangle}
\end{figure}

{\bf Доказательство}. Сначала разберемся с операцией разбиения. Рассмотрим некоторую плитку $T$. У нее четыре угла $X$, $Y$, $P$, $Q$ (рисунок~\ref{fig:tileangle}).

При разбиении некоторые углы разбиваются ребрами, а некоторые нет. Будем называть {\it простыми} те углы, которые не разбиваются ребрами (такими будут
углы $X$ и $Y$). Из правил разбиения плиток следует, что простые углы уже никогда не будут разбиты ребрами при следующих разбиениях.

Назовем угол {\it ограниченным}, если либо при следующем разбиении он будет разбит на простые углы, либо на простые и ограниченные углы. Например, угол $P$ (правый нижний в плитке $T$) -- ограниченный угол, так как на следующем разбиении он будет подразбит на левые верхние углы новообразованных двух плиток, и дальше уже подразбиваться не будет.

Можно проверить, что каждая из создаваемых боковых и внутренних вершин создает вокруг себя только простые и ограниченные углы, и если проводить только разбиения (без подклеек) то через два этапа после создания вершины, в нее не будут добавляться входящие ребра и ее степень перестает меняться.

Заметим, что степень вершины может возрасти в двух случаях: при разбиении плитки, в которой она участвует и при подклейке плитки, где она выступает в качестве вершины $X_1$, $X_2$, $Z_1$, $Z_2$ из определения подклейки. Разбиение плиток, начиная с некоторого момента, перестает менять степень выбранной вершины. Участвовать же в подклееных плитках можно только на трех этапах (следующих после того, когда вершина создается) так как в определении подклеек участвуют вершины от $k-2$ уровня до $k$ уровня.

Таким образом, начиная со второй итерации после появления, степень вершины может возрастать только при {\it Разбиении}, а после четвертого разбиения она уже не поменяется.

\medskip

Таким образом, каждая вершина в комплексе имеет конечную степень, зависящую от уровня комплекса. При этом, число ребер не ведущих в подклееные плитки, ограничено в совокупности для всех вершин.

\medskip

\medskip

Плоские части комплекса представляют собой макроплитки некоторого уровня. Каждая такая макроплитка либо является начальной, базовой плоскостью, либо в какой-то момент была подклеена, а потом несколько раз разбита. Макроплитки состоят из плиток (простых циклов длины $4$).

\medskip

\begin{figure}[hbtp]
\centering
\includegraphics[width=0.4\textwidth]{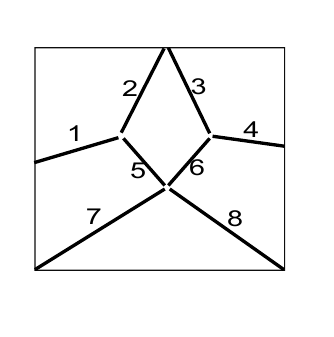}
\caption{Восемь типов внутренних ребер.}
\label{fig:inneredges}
\end{figure}

\begin{proposition}[О боковой вершине] \label{side_node}

1) Каждая боковая вершина лежит на середине стороны в какой-либо макроплитке или в двух макроплитках одного уровня, лежащих в одной плоскости;

2) Если боковая вершина не находится на границе исходного комплекса первого уровня, то она лежит на одном из восьми внутренних ребер (рисунок~\ref{fig:inneredges} ) в некоторой макроплитке, либо лежит на границе подклееной макроплитки.

\end{proposition}

{\bf Доказательство}.
В комплексе первого уровня свойство выполняется. На втором уровне также все верно для всех созданных боковых вершин.

Пусть для уровня $k$ комплекса все эти свойства выполнены. Для уже существующих боковых вершин оба свойства будут сохраняться при дальнейших разбиениях и подклейках. Если новая боковая вершина возникла при подклейке, то это она использовалась в качестве $T_2$ или $T_3$ из определения подклейки, и в качестве макроплитки можно взять как раз подклееную в тот момент макроплитку.

 Созданная при разбиении вершина, очевидно, лежит на середине стороны только что разбитой плитки. Если это подклееное ребро или граница начального комплекса первого уровня, то такая плитка одна (и тогда это граница подклееной макроплитки, либо всего комплекса), в остальных случаях таких макроплиток две.

Для уровней $1$ и $2$ комплекса заметим, что все стороны плиток, не лежащие на границе, являются частью какого-то большего ребра, классифицируемого, как одного из восьми типов из рисунка~\ref{fig:inneredges}. При дальнейших подклейках и разбиениях эти ребра не могут сменить тип. При этом создаваемые при разбиении внутренние ребра, очевидно, получаются одного из этих восьми типов.

Если ребро, при разбиении которого образовалась вершина не является границей для всего комплекса или подклееной макроплитки, то  оно принадлежит одному из восьми типов из рисунка~\ref{fig:inneredges} для какой-то макроплитки. Значит и вершина тоже принадлежит такому ребру.

\medskip

{\bf Примечание.} Для каждой боковой вершины $X$ можно выделить пару вершин $(Y,Z)$, являющихся концами ребра, на котором лежит $X$. Эту пару вершин будем называть {\it ``начальниками''} для $X$.

\medskip

\subsection{Пути на комплексе}

В этом параграфе мы установим несколько свойств путей, проходящих по комплексу.

Путь, идущий по макроплитке, может быть локально преобразован:

\begin{proposition}[О переброске пути]  \label{path_flip}

Пусть $XYZT$ -- некоторая макроплитка. Рассмотрим путь $XYZ$ (состоящий из двух соседних макроребер). Тогда, если разрешается менять подпуть из двух соседних ребер любой плитки на подпуть из двух других ребер (с общими началом и концом у этих подпутей), то путь $XYZ$ может быть преобразован в $XTZ$ -- путь по другим двум макроребрам.

\end{proposition}

\begin{figure}[hbtp]
\centering
\includegraphics[width=0.7\textwidth]{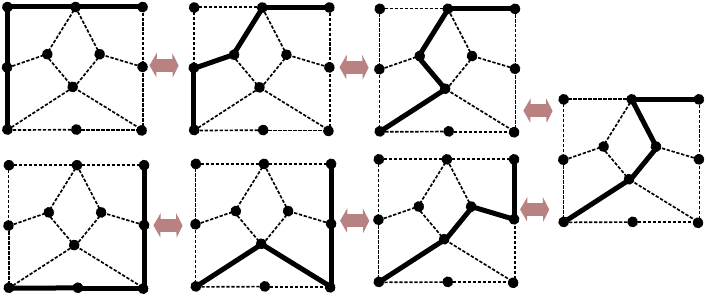}
\caption{Переброска пути на другую сторону макроплитки.}
\label{fig:transformpath}
\end{figure}

{\bf Доказательство по индукции}. Для случая, когда макроплитка -- это плитка, преобразование можно сделать сразу. Шаг индукции можно совершить, выполняя
локальные преобразования по правилам, показанным на рисунке~\ref{fig:transformpath}.

\medskip

{\bf Примечание.} Применяя аналогичные рассуждения, можно доказать то, что также перебросить можно и путь, начало и конец которого лежат в серединах противоположных ребер макроплитки.

\medskip

\begin{definition}
Замены подпутей из пары соседних ребер плитки на подпуть из другой пары будем называть {\it локальными заменами}. А возможность с помощью таких замен достичь некоторого состояния, будем называть {\it локальным преобразованием}.
\end{definition}

\medskip

\begin{definition}
 Будем считать форму пути {\it нулевой}, если она содержит подпуть длины $2$ по ребру некоторой минимальной плитки (туда-обратно).
\end{definition}

\medskip

\begin{proposition}[О ``выносе'' пути на границу]  \label{to_border}

Пусть начало и конец пути $P$, проходящего по макроплитке $T$, лежат на границе $T$. Тогда можно реализовать одну из двух возможностей:

1) $P$ может быть локально преобразован в нулевую форму.

2) $P$ может быть локально преобразован в форму $P'$ так, что $P'$ полностью лежит на границе $T$.
\end{proposition}

{\bf Доказательство по индукции}. Любой путь по минимальной плитке и так лежит на ее границе. Пусть макроплитка уровня $n$. Рассмотрим ее разбиение на $6$ (макро)плиток уровня $n-1$. Путь $P$ разбивается на несколько участков, каждый из которых имеет начало и конец на границе какой-то из этих подплиток. По предположению индукции, эти участки можно локально преобразовать в форму $P_1$ так, чтобы весь путь $P_1$ проходил только по границам подплиток.

\medskip

Обозначим буквами $L$, $U$, $R$, $D$ середины соответственно левой, верхней, правой, нижней сторон $T$, и внутренние вершины как $A$, $B$, $C$. Допустим, путь $P_1$ содержит дважды некоторую вершину  нашей макроплитки $T$. Тогда он должен полностью обходить некоторую подмакроплитку уровня $n-1$, примыкающую к этой вершине. В этом случае можно применить лемму о переброске пути и получить форму, где выполнено 1).

Пусть $P_1$ содержит $A$ один раз. Войти и выйти $P_1$ может только по трем ребрам (направленным в сторону $U$, $L$ и $C$). Если вход и выход -- по одному ребру, то условие 1) выполнено.

Пусть во входе и выходе задействованы ребра в сторону $U$ и $L$. Если, например, выйти в сторону $L$, то сойти с ребра мы не можем (так как $P_1$ проходит только по границам макроплиток $n-1$ уровня) и вернуться назад тоже (так как тогда будет выполнено 1). Заметим, что тогда $P_1$ содержит подпуть $LAU$ или $UAL$, в обоих случаях этот подпуть, согласно лемме о переброске пути, можно локально преобразовать чтобы он проходил по границе макроплитки $T$ и, таким образом, не будет содержать $A$.

Если во входе и выходе задействованы ребра в сторону $L$ и $C$, то $P_1$ содержит подпуть $LAC$ или $CAL$ и он может быть преобразован в форму соответственно $LXC$ или $CXL$, где $X$ -- левый нижний угол, и в этом случае опять форма не будет содержать $A$.

Если во входе и выходе задействованы ребра в сторону $U$ и $C$, то $P_1$ содержит подпуть $UAC$ или $CAU$ и он может быть преобразован в форму соответственно $UBC$ или $CBL$, и тоже форма не будет содержать $A$.

\medskip

Проделав аналогичные рассуждения для симметричного случая вершины $B$, мы можем заключить, что путь можно преобразовать в форму $P_2$, такую что она либо вообще не содержит $A$ и $B$, либо содержит кусок $UBC$ или $CBU$ (который может быть преобразован в $UAC$ или $CAU$ соответственно).

Если $P_2$ не содержит $A$ и $B$ и  содержит $C$, то войти и выйти в $C$ можно только по ребрам, ведущим в левый нижний и правый нижний углы. В этом случае можно применить лемму о переброске пути для участка, содержащего эти углы и $C$ и получить форму, не содержащую вершин $A$, $B$, $C$.

Пусть $P_2$ содержит кусок $UBC$ или $CBU$. Посмотрим, как $P_2$ проходит через $C$: один вход (или выход) в сторону $B$, другой не может идти в сторону $A$, так как $P_2$ не содержит $A$, значит он идет в сторону нижнего левого или нижнего правого угла. Если это нижний правый угол $Y$ то $P_2$ содержит кусок $UBCY$ или $YCBU$. Применяя лемму о переброске путей, $UBCY$ переводится в $UBRY$ и далее в $UZRY$, где $Z$ правый верхний угол. Аналогично, $YCBU$ переводится в $YRZU$.
Если это нижний правый угол $X$, то $P_2$ содержит кусок $UBCX$ или $XCBU$. Аналогично применяя лемму~\ref{path_flip} о переброске пути, переводим $UBCX$ в $UACX$ и далее в $UALX$ и потом в $UWLX$, где $W$ -- левый верхний угол. Аналогично $XCBU\Rightarrow XCAU \Rightarrow XLAU   \Rightarrow XLWU$.

Таким образом, мы получим форму не содержащую вершин $A$, $B$, $C$. В этом случае, если нельзя выполнить условие 1), то весь путь проходит по границе $T$, и выполнено условие~2).

\medskip

\begin{proposition}[О ``шевелении'' пути] \label{path_changing}

 Пусть $P$ -- путь, состоящий из макроребра $XY$ некоторой макроплитки $T$ уровня не менее $3$.
Тогда $P$ может быть локально преобразован в форму $P'$, состоящей из трех частей: первая и третья это участки по $1/4$ пути $P$, примыкающие к $X$ и $Y$ соответственно, а вторая представляет собой путь по двум соседним ребрам подклееной макроплитки, два других ребра которой являются частью $P$.
\end{proposition}

Для доказательства достаточно проследить разбиение $T$ (на три уровня в глубину).
Обозначим середину $XY$ как $H$, середины $XH$ и $HY$ как $G_1$ и $G_2$ соответственно, и середины $G_1H$ и $HG_2$ как $G_3$ и $G_4$. Так как $G_3$ и $G_4$ будут иметь уровень на $2$ больше, чем $H$, значит, была проведена подклейка для пути $G_1G_3HG_4G_2$. Пусть $F$ -- созданная при этой подклейке угловая вершина. По предыдущей лемме, $G_1HG_2$ можно локально преобразовать в $G_1FG_2$.

\medskip

\medskip

\begin{proposition}[О выделении локального участка] \label{longpath}

Если путь $P$, проходящий по макроплитке $T$ уровня $n$, имеет длину не менее $5\cdot 2^{n-2}$, он может быть локально преобразован в нулевую форму $P'$.

\smallskip

Если путь $P$, проходящий по макроплитке $T$ уровня $n$, имеет длину не менее $2^{n}+1$, и начинается в углу $T$, либо в середине стороны $T$, он может быть локально преобразован в нулевую форму $P'$.

\end{proposition}

{\bf Доказательство}. Согласно лемме о выносе на границу, можно считать, что весь путь $P$ проходит по границе $T$. Пусть при этом не образовалось нулевой формы, то есть путь покрывает какую-то часть границы $T$, совершая движение либо по часовой стрелке (либо против).  Отметим на границе $T$ восемь точек -- углы и середины сторон. В обоих случаях 1) и 2) наш путь включает в себя не менее $5$ из отмеченных точек, иначе он был бы короче, чем $5\cdot 2^{n-2}$ или $2^{n}+1$ соответственно. Тогда найдутся две отмеченные точки, диаметрально противоположные относительно $T$, и такие, что путь $P$ полностью покрывает половину границы между ними. Заметим при этом, что хотя бы одна из этих отмеченных точек не является концевой для $P$. Пользуясь леммой о переброске (с примечанием) можно перевести подпуть из половины границы на другую половину.
Ясно, что получившаяся форма будет нулевой.

\medskip

\begin{proposition}[О непродолжаемом пути] \label{bad_path}

Пусть путь $P$ лежит в некоторой макроплитке $T$, причем начало и конец $P$ лежат в серединах противоположных сторон $T$. Тогда любые плоские пути вида $WP$ или $PW$, где длина $W$ более длины $P$, могут быть приведены к нулевой форме.
\end{proposition}

{\bf Доказательство}. Случаи $WP$ и $PW$ аналогичны, рассмотрим $PW$. Пусть $a$ -- это та сторона $T$, в середине $X$ которой заканчивается $P$ и начинается $W$. Если выход из $X$ идет внутрь $P$, либо по стороне $a$, тогда рассматривая путь $P$ плюс одна вершина, мы получаем путь длины $2^{n}+1$, лежащий в макроплитке $P$ и можем применить лемму о выделении локального участка. Значит, выйти из $X$ может только в макроплитку $T'$, соседствующую с $T$ по стороне $a$. Допустим, существует вершина $Y$ из пути $W$, лежащая на границе $T'$, отличная от $X$. В этом случае, кусок пути $W$ от $X$ до $Y$ может быть локально преобразован так, что он будет проходить по границе $T'$. Но тогда он выход этого куска из $X$ будет по стороне $a$, и в этом случае опять получаем путь длины $2^{n}+1$, лежащий в макроплитке $P$.  Если такой вершины $Y$ не найдется, то все остальные вершины $W$ лежат внутри $T'$, то есть получаем опять путь длины $2^{n}+1$, на этот раз лежащий на $T'$.

\smallskip

\begin{definition}
{\it Паттерном пути} будем называть последовательность его
вершин, из которой выброшены все боковые вершины, вход и выход в которые проходят по ребрам, лежащим на границе между одной и той же парой макроплиток (главные ребра для боковой вершины). То есть в паттерн пути входят все черные (внутренние) вершины, и все боковые вершины, где происходит вход(выход) в (из) макроплитку.
Фактически паттерн пути -- это его карта, показывающая маршрут на ориентировочных точках.
\end{definition}

\begin{definition}
{\it Мертвым} будем называть такой паттерн $P$, что любой достаточно большой путь $W$, cобственным подпутем которого является путь с паттерном $P$, может быть локально преобразован к нулевой форме.
\end{definition}

\begin{proposition}[О мертвых паттернах] \label{DeadPaterns}

Рассмотрим некоторую макроплитку $T$ и обозначим в ней внутренние вершины $A$, $B$, $C$ и боковые $U$, $R$, $D$, $L$ (аналогично обозначениям при разбиениях).
Тогда паттерны $AUB$, $ACB$, $CXD$ (где $X$ -- любая вершина) являются мертвыми.
\end{proposition}

{\bf Доказательство}. Пути с паттернами $AUB$ и $ACB$ локально преобразуются друг в друга, так что достаточно рассмотреть один из этих случаев. Пусть $W$ -- достаточно большой путь, содержащий подпуть с паттерном $ACB$. Тогда можно считать, что $W$ содержит вершины на границе $T$, лежащие до и после нашего подпути. Рассмотрим разбиение макроплитки $T$ на $6$ дочерних подплиток (в соответствии с операцией разбиения) и будем локально преобразовывать $W$ так, чтобы он проходил по границам этих дочерних подплиток. Тогда из вершины $B$ путь
должен идти в по ребру к $R$, а в $A$ путь должен входить по ребру из $L$, в остальных случаях образуется обход подлитки по трем ребрам (и этот кусок преобразуется к нулевой форме). Таким образом, в $W$ можно выделить подпуть, проходящий по макроплитке $T$, начало которого в $L$ и конец в $R$. Можно воспользоваться леммой о непродолжаемом пути и получить требуемое.

Для паттерна $CXD$, пусть $T'$ -- макроплитка, соседняя с $T$ по нижней стороне. Заметим, что выход из вершины $D$ обязательно должен быть внутрь макроплитки $T'$, иначе можно применить лемму о выделении локального участка для нижней дочерней подплитки $T$. Проследим путь $W$ до ближайшей точки $X$ на границе $T'$ (после $D$) и преобразуем  кусок $DX$ чтобы он проходил по периметру $T'$.
Но тогда из $D$ путь $W$ продолжится по граничной между $T$ и $T'$ стороне.
То есть теперь, все таки, можно применить лемму о  выделении локального участка для нижней дочерней подплитки $T$ и получить требуемую нулевую форму.

\medskip

{\bf Примечание.} Ясно, что и паттерны $BUA$, $BCA$, $DXC$ тоже мертвые, доказательство полностью аналогично.

\medskip
\begin{proposition}[О мертвых путях в нижней подплитке] \label{DeadPaths}

Рассмотрим некоторую макроплитку $T$ . Пусть путь $XYZ$ лежит в $T$, причем $X$ лежит на внутреннем ребре, идущем из левого нижнего угла $T$ во внутреннюю вершину $C$, $Y$ -- нижний угол $T$, а $Z$ лежит на нижней стороне $T$. Тогда существует такое натуральное $N$,
что для любых путей $W_1, W_2$, длины которых более $N$, путь $W_1XYZW_2$ можно локально преобразовать к нулевой форме.

\end{proposition}

{\bf Доказательство}. Пусть $Y$ -- левый нижний угол. Обозначим середину нижней стороны $T$ как $D$. Пусть $T'$ нижняя подплитка $T$. Преобразуем путь $W_1XYZW_2$ так, чтобы он проходил по границам макроплиток уровня $T'$ (не заходя в более мелкие плитки). Ясно, что $W_1$ должен проходить через $C$, а $W_2$ через $D$, иначе возникает очевидный кусок с нулевой формой. Таким образом, после преобразования, путь будет содержать подпуть с паттерном $CXD$, который является мертвым, по лемме о мертвых паттернах. Следовательно, для достаточно больших $W_1$ и $W_2$ наш путь можно преобразовать к нулевой форме.

\medskip

{\bf Замечание}. Доказательство аналогично переносится на случай, когда $Y$ -- правый нижний угол.

\medskip

\begin{definition}
 {\it Некорректным участком} пути будем называть такой подпуть $XYZ$, что вершина $Y$ лежит на границе некоторой макроплитки $T$, а вершины $X$ и $Y$ лежат внутри $T$ (не на границе).
\end{definition}

\begin{proposition}[О некорректных участках] \label{uncorrect_sectors}

Пусть есть некорректный участок $XYZ$. Тогда существует константа $N$, что для любых путей $W_1, W_2$, длины более $N$, путь $W_1XYZW_2$, может быть локально преобразован к нулевой форме.
\end{proposition}

{\bf Доказательство}. Пусть $T$ -- минимальная макроплитка, внутри которой проходит путь $XYZ$. Из минимальности следует, что вершина $Y$ лежит в углу $T$, либо в середине стороны. Из иерархии разбиений следует, что из всех углов кроме нижнего левого и из середины нижней стороны выходит не более одного ребра внутрь $T$, так что $Y$ может лежать только на середине левой, правой или верхней стороны или в левом нижнем углу.

{\bf 1.} Пусть $Y$ -- левый нижний угол. Из этого угла внутрь выходит два ребра, оба к внутренней C-вершине, одна в макроплитке $T$, другая в нижней дочерней подплитке $T$. То есть, тогда путь $XYZ$ попадает под условие леммы о нижней подплитке и все доказано.

\begin{figure}[hbtp]
\centering
\includegraphics[width=0.5\textwidth]{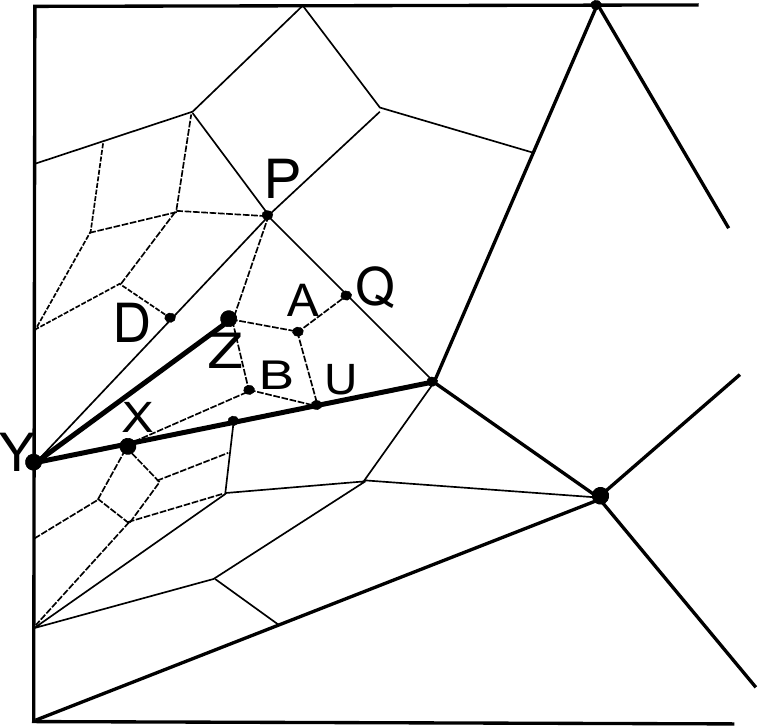}
\caption{Вершина $Y$ -- середина левой стороны.}
\label{fig:badwayleftcase}
\end{figure}

{\bf 2.} Пусть $Y$ --  середина левой стороны. Из этой вершины выходит три ребра и поэтому есть три варианта расположения пути $XYZ$. В двух случаях можно применить лемму о нижней подплитке. Оставшийся случай изображен на рисунке~\ref{fig:badwayleftcase}.

Из вершины $Z$ дальше путь может пойти:

 i) в вершину $Z_1$. В этом случае кусок $YZZ_1$ преобразуется в $YDZ_1$ и для участка $XYD$ можно применить лемму о нижней подплитке.

 ii) в вершину $A$. Посмотрим, куда путь может пойти дальше. Если это вершина $Z_2$, то участок $YZAZ_2$ преобразуется в $YZZ_1Z_2$ и далее в  $YDZ_1Z_2$, после чего опять можно применить лемму о нижней подплитке. Во втором случае, если путь идет в вершину $U$, то путь $XYZAU$ сразу можно привести к нулевой форме:
$XYZAU\rightarrow XYZBU  \rightarrow XYXBU$.

\medskip

{\bf 3.} Пусть $Y$ --  середина правой стороны. Этот случай симметричен второму, только вместо левой верхней подплитки мы имеем дело с правой нижней. Все рассуждения полностью аналогичны второму случаю.

\begin{figure}[hbtp]
\centering
\includegraphics[width=0.6\textwidth]{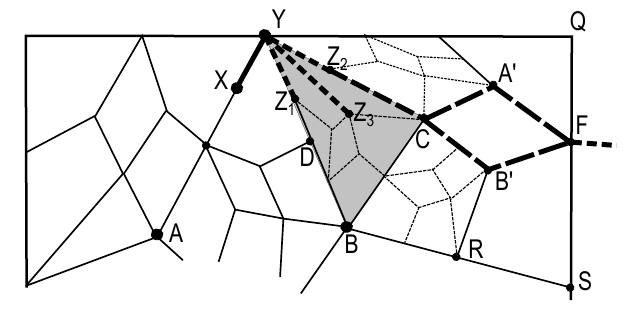}
\caption{Вершина $Y$ -- середина верхней стороны.}
\label{fig:badwayupcase}
\end{figure}

{\bf 4.}  Пусть $Y$ --  середина верхней стороны. Из этой вершины исходят четыре ребра
(рисунок~\ref{fig:badwayupcase}).

В случае, когда $X$ и $Z$ не лежат на левом ребре (идущем из $Y$ к внутренней вершине $A$ макроплитки $T$), то получается ситуация, полностью аналогичная второму и третьему случаям выше, на этот раз вместо левой верхней подплитки мы имеем дело с правой верхней. Итак, пусть верхняя правая подплитка $T$ это $T'$ и,   ради определенности, $X$ лежит на ребре, уходящем из $Y$ к внутренней вершине $A$ макроплитки $T$, а $Z$ -- на одном из трех ребер, попадающих в $T'$ (Отмечены пунктиром на рисунке~\ref{fig:badwayupcase}).

Проведем локальные преобразования пути $W_1$ так, чтобы он проходил по границам макроплиток уровня $T'$ (не заходя внутрь них). Последняя вершина $W_1$ это $X$. Таким образом, после этого преобразования либо последний участок будет $YX$ (и тогда весь наш путь получит нулевую форму) либо преобразованный $W_1$ будет проходить через внутреннюю вершину $A$ макроплитки $T$. Теперь рассмотрим путь $W_2$. Наша цель -- преобразовать его к виду, чтобы он проходил через внутреннюю вершину $B$ макроплитки ~$T$.

Пусть $T''$ это нижняя подплитка $T'$ (отмечена цветом на рисунке). Проведем локальные преобразования $W_2$ так, чтобы он проходил по границам макроплиток уровня $T''$ (не заходя внутрь них). Ясно, что ближайшая к $Z$ часть $W_2$ пройдет по одному из участков $YCA'F$ или $YCB'F$, указанных длинным пунктиром на рисунке~\ref{fig:badwayupcase} (либо $W_2$ пройдет через $B$). Таким образом, $W_2$ теперь проходит через вершину $F$.  Теперь, если из $F$ путь $W_2$ идет к вершине $Q$, то оба возможных куска $YCA'FQ$ и $YCB'FQ$ приводятся к нулевой форме. Если путь из $F$ идет в $S$, то кусок $YCA'FS$ (или $YCBFS$) приводится к $YDBRS$. Остается возможность, когда $W_2$ покидает макроплитку $T$, уходя в соседнюю макроплитку по правой стороне (короткий пунктир). В этом случае весь кусок $W_2$ после $F$ мы преобразуем так, чтобы он проходил по границам макроплиток уровня $T$. Тогда преобразованный участок будет уходить из $F$ только в $Q$ или $S$, а эти случаи мы разобрали.

Итак $W_1$ можно привести к форме, проходящей через внутреннюю вершину $A$ макроплитки $T$, а $W_2$ к форме, проходящей через внутреннюю вершину $B$.
Теперь получившийся путь имеет подпуть с мертвым паттерном $AUB$, вследствие чего, для достаточно больших путей $W_1$ и $W_2$ мы можем привести путь $W_1XYZW_2$ к нулевой форме.

\medskip

{\bf Примечание.} Рассуждения не меняются, если рассмотренная макроплитка $T$ имеет меньший уровень. В этом случае вершина $Z_1$ совпадает с $D$, а $Z_3$ с $C$.

\medskip

\begin{proposition}[О корректности путей] \label{correct_paths}

Пусть путь $P$ представляет собой проход по двум соседним ребрам некоторой макроплитки. Тогда любые локальные преобразования не могут привести $P$ к форме, содержащей некорректный участок или к нулевой форме.

\end{proposition}

{\bf Доказательство}. Для двух точек комплекса определено расстояние как наименьшая из длин соединяющих их путей. Путь, проходящий по двум соседним ребрам макроплитки является кратчайшим для вершин в его концах, что легко доказывается по индукции по уровню макроплитки. Локальные преобразования не меняют длину пути. В случае, если после очередного преобразования появляется нулевая форма, путь не может быть кратчайшим.

Теперь, допустим, что наш путь можно привести к форме, содержащей некорректный участок. Заметим, что путь после любого числа преобразований имеет форму $W_1UW_2$, где $W_1$, $W_2$ полностью лежат на границах макроплитки, а у подпути $U$ концы лежат на границе, а все остальные вершины лежат внутри. Для макроплиток первых уровней это проверяется непосредственно, а для остальных уровней можно проверить по индукции. В случае, если путь содержит некорректный участок, он не может иметь такую форму.

\medskip

\medskip

\begin{proposition}[О расстоянии от края макроплитки до выхода в подклейку] \label{pasting_distance}

Пусть вершина $X$ лежит на краю некоторой макроплитки $T$, вершина $Y$ принадлежит $T$, но не находится на ее границе, а выход из $Y$ идет в подклееную макроплитку уровня $n \geq 2$. Тогда расстояние от $X$ до $Y$ (длина кратчайшего пути по ребрам) не менее $2^{n-1}$.
\end{proposition}

{\bf Доказательство} проводится индукцией по уровню макроплитки $T$. Для уровней $2$,$3$ подклеек внутри $T$ не проводится. Для уровня $4$ ближайшая вершина с подклееным ребром лежит в середине любого ребра, выходящего на середину любого из граничных ребер. Расстояние от нее до границы $T$ будет равно $2$, при этом подклееная макроплитка будет иметь уровень $2$.

Пусть уровень $T$ равен $k>4$. Возьмем кратчайший путь $XWY$. Можно считать, что $W$ не содержит вершин на границе $T$.
Кроме того, $Y$ не может лежать внутри какой-то из подплиток $T$, иначе $X$ можно взять на границе этой подплитки.
Значит, $Y$ обязательно лежит на одном из ребер, которые разделяют макроплитку $T$ на $6$ макроплиток при разбиении. Опять используя минимальность пути, замечаем, что $Y$ лежит в углу подклееной макроплитки. Из определения подклейки получаем, что уровень $Y$ был на $1$ меньше максимального в тот момент, когда эта подклейка производилась. Если после этого разбиений не было, то расстояние $XY$ не менее $2$, а также макроплитка при $Y$ -- второго уровня. Если после этого провели $l\geq 1$ разбиений, расстояние $XY$ увеличилось в $2^l$ раз, а уровень макроплитки при $Y$ возрос на $\mathbf{l}$.

\medskip

\begin{corollary}
Расстояние между двумя выходами в подклееные области не менее $2^{n-1}$, где $n$ уровень макроплитки, куда идет второй выход.
То есть если путь $P$ имеет форму $XP_1YZ$, где  $P_1$ некоторый путь, выходы из вершин $X$ и $Y$ осуществляются по подклееным ребрам, а все остальные входы во все вершины, и выходы из них плоские. Уровень подклееной макроплитки, куда идет выход из $Y$ равен $n$. Тогда расстояние от $X$ до $Y$ не менее $2^{n-1}$.
\end{corollary}

\medskip

Таким образом, если в пути встречаются два ребра выходящие из вершин по подклееным ребрам, то они разделены, как минимум, одним плоским ребром.

\medskip

Каждая подклееная макроплитка имеет ограниченные размеры, соответственно, путь, не имеющий нулевой формы и не возвращающийся из подклееной макроплитки будет иметь ограничения на длину.

\begin{proposition}[Об ограниченности пути, уходящем в подклееную часть два раза подряд] \label{two_in_row}

Пусть путь $P$ имеет форму $XP_1YP_2$, где $P_1Y$ -- плоский участок (лежащий в какой-то макроплитке) и выходящие из вершин $X$ и $Y$ ребра идет в подклееные плитки, и больше выходов в подклееную область или входов в вершину из подклееной области нет. Пусть также $P$ не может быть приведен к нулевой форме.
Тогда длина $YP_2$ менее, чем $\frac{5}{2}$ длины $XP_1Y$.
\end{proposition}

{\bf Доказательство}. Пусть уровень макроплитки, куда уходит ребро из $Y_2$, равен $n$. Тогда длина участка $YP_2$ меньше чем $5 \cdot 2^{n-2}$ по лемме о выделении локального участка. Кроме того, длина $XP_1Y$ не менее $2^{n-1}$ по следствию леммы о расстоянию от границы макроплитки до подклееного ребра. Следовательно, длина $YP_2$ меньше, чем $\frac{5}{2}$ длины $XP_1Y$.

\medskip

\section{Строение полугруппы путей} \label{pathsemigroup}

В предыдущих разделах мы описали геометрическую структуру комплекса.  В следующей главе~\ref{coding_section} будет показано, как связать с ним специальную кодировку. Мы покажем, что с помощью конечного числа букв можно закодировать узел. Кроме того, конечным числом букв можно закодировать входящие и выходящие ребра в каждый узел, как плоские, так и выходящие в подклейку.

В этом разделе мы покажем, для чего нужна такая кодировка и как необходимая нам полугруппа строится с помощью кодировки путей на комплексе.

Пусть $X_i$ -- буквы, кодирующие входящие ребра, $Y_i$ -- буквы, кодирующие узлы (их типы, окружения и информации), $Z_i$ -- буквы, кодирующие выходящие ребра. (Детали кодирования будут показаны в последующих разделах.)

Будем говорить, что слово $W$ имеет  $\mathbf{CODE}$-{\it форму}, если в нем справа от любой (не последней в слове) буквы семейства $X$ обязательно стоит буква семейства $Y$, справа от любой (не последней в слове)  буквы семейства $Y$ стоит буква семейства $Z$, а справа от любой (не последней в слове)  буквы семейства $Z$ стоит буква семейства $X$.

Часть таких $\mathbf{CODE}$-слов соответствует путям на комплексе.
Рассмотрим полугруппу с нулем $S$ с порождающими $\{X_i,Y_i,Z_i \}$.

\medskip

\begin{proposition}[О стандартной форме слова] \label{standardform}

В полугруппе с нулем $S$ с порождающими $\{X_i,Y_i,Z_i \}$ можно ввести конечное число определяющих мономиальных соотношений так, чтобы все не $\mathbf{CODE}$-слова, можно было привести к нулевому слову.

\end{proposition}

\medskip

{\bf Доказательство}. Действительно, введем соотношения $X_iZ_j=0$, $X_iX_j=0$, $Y_iX_j=0$, $Z_iY_j=0$, $Y_iY_j=0$, $Z_iZ_j=0$, для всевозможных пар $(i,j)$. В получившейся полугруппе, любое ненулевое слово будет иметь $\mathbf{CODE}$-форму: за буквой $X_i$ может следовать только $Y_j$, за $Y_j$ -- только $Z_k$, за $Z_k$ -- только $X_m$.

\medskip

Заметим также, что можно выписать все $\mathbf{CODE}$-слова длины не более $4$, которые кодируют какой-нибудь путь на комплексе. Пусть $F_1, \dots ,F_N$ -- все оставшиеся $\mathbf{CODE}$-слова длины не более $4$. Их мы будем считать {\it запрещенными} и введем в полугруппе $S$ соотношения $F_i=0$, для всех $i=1,\dots ,N$.

Таким образом, любое слово $W$ в полугруппе $S$ либо приводится к нулю, либо имеет $\mathbf{CODE}$-форму, где каждое подслово $W$ длины не более $4$ отвечает некоторому пути на комплексе.

\medskip

В разделах~\ref{flip_section} и~\ref{pasting_section} будет показано, что для двух эквивалентных путей, не содержащих мертвых паттернов и нулевых форм, по коду одного из них можно восстановить код другого. Это дает возможность ввести конечное множество соотношений вида $W_i=W_j$, где $W_i,W_j$ -- коды двух эквивалентных путей.

Таким образом, появляется возможность преобразовывать слова, меняя подпути на им эквивалентные.
При этом если до преобразования слово являлось кодировкой некоторого пути на комплексе, то и после преобразования слово оно будет кодировкой эквивалентного ему пути. Допустим, после нескольких таких замен образуется запрещенное подслово. Это значит, что соответствующий этой кодировке путь не может быть реализован на комплексе. Но тогда и эквивалентный ему путь не мог быть реализован.
Таким образом, если, преобразовывая слова-кодировки по установленным правилам, мы получаем запрещенное слово, значит изначальное слово не соответствует пути на комплексе.

\medskip

{\bf Нулевая форма.} Если слово $W$ представляет собой кодировку пути, имеющего {\it нулевую форму}, мы вводим в полугруппе соотношение $W=0$. Рассмотрим слово-кодировку пути, концы которого лежат на периметре некоторой макроплитки $T$, и хотя бы один конец попадает в угол $T$. Заметим, что если путь не кратчайший (то есть существует путь короче по длине с теми же концами), то согласно предложению~\ref{longpath}, путь может быть преобразован в нулевую форму. Это означает, что соответствующее слово-кодировка приводится к нулю.

Фактически, обнуление нулевых форм позволяет приводить к нулю все слова, соответствующие некратчайшим путям.

\medskip

{\bf Мертвые паттерны.} Мертвые паттерны представляют собой некоторые типы слов-кодировок, которые никогда не могут встретиться в достаточно больших ненулевых словах. Нас будут интересовать только мертвые паттерны $AUB$, $ACB$, $CXD$, $BUA$, $BCA$, $DXC$, которые обсуждаются в предложении~\ref{DeadPaterns}. Выпишем все слова, длины не более $4$, содержащие хотя бы один из перечисленных мертвых паттернов. Для каждого такого слова $W$, введем соотношение $W=0$.

Заметим, что все достаточно большие ненулевые слова не могут содержать мертвый паттерн (и не могут быть приведены к форме, его содержащей).

\medskip

Будем называть слово {\it регулярным}, если оно является кодировкой некоторого пути на построенном геометрическом комплексе и при этом не может быть преобразовано в нулевую форму или форму содержащую мертвый паттерн.

\medskip

Заметим, что замена некоторого подслова регулярного слова на эквивалентное ему, не нарушает регулярность. То есть, если слово можно преобразовать так, что оно будет содержать запрещенное подслово, нулевую форму или мертвый паттерн, то значит оно приводится к нулю, и не может быть регулярным.

Нашей основной целью будет показать, что ненулевыми словами в полугруппе являются только регулярные слова, то есть, что любое нерегулярное слово может быть приведено к нулю. Фактически, это будет означать, что ненулевыми элементами полугруппы являются только слова, кодирующие кратчайшие пути, не содержащие мертвые паттерны.

\medskip

{\bf Примечание о использовании мертвых паттернов для конструкции}. Ценой некоторого усложнения кодировки можно добиться того, чтобы пути, содержащие мертвые паттерны, допускали такие же локальные преобразования, как и обычные слова. В этом случае, можно не вводить мономиальные соотношения для слов, их содержащих. И тогда ненулевыми элементами полугруппы будут просто кратчайшие пути на построенном комплексе. Но в целях упрощения конструкции, представляется разумным ввести такие соотношения.

%Фактически, множество мертвых паттернов представляет собой радикал в полугруппе путей.

\medskip

\section{Кодировка вершин и путей на комплексе} \label{coding_section}

Рассмотрим путь на комплексе. Для каждой вершины из этого пути, кроме первой и последней, есть ребро входа и ребро выхода. Ребра входа и выхода могут вести в подклееные области, а могут принадлежать той же базовой плоскости, что и сама вершина.

Каждая вершина принадлежит некоторой макроплитке. Кроме того из нее могут выходить ребра ведущие в подклееные плитки/макроплитки. То есть, каждая вершина лежит на своей базовой плоскости и еще участвует в нескольких подклееных плоскостях. Фактически, вершина играет свою роль (занимает определенное положение) в своей макроплитке на базовой плоскости, и кроме этого выполняет другие какие-то роли в своих подклееных макроплитках. При этом в подклееных макроплитках эта вершина лежит всегда на краю.

Мы хотим определить систему кодирования путей на комплексе. Каждый путь это  конечная последовательность вершин, причем любые две соседние соединены ребром плитки минимального уровня (не макроплитки). Для кодирования путей нужно сначала закодировать все вершины, которые могут встретиться на этом пути, конечным числом букв.

Каждая вершина может принимать участие во множестве подклееных макроплиток, то есть входящие в нее ребра либо относятся к плоскости вершины либо классифицируются по принадлежности к различным подклееным макроплиткам. В общем случае, путь приходит в вершину из одной подклееной макроплитки и выходит в другую. Общая кодировка вершины будет состоять из трех {\it плоских} кодов: первый представляет собой плоский код вершины (без учета подклееных плиток) второй- кодирует положение этой же вершины в плоскости подклееной макроплитки, откуда приходит входящее ребро, и третий кодирует положение этой же вершины в плоскости подклееной макроплитки, куда уходит выходящее ребро.

\medskip

Если рассматривать вершину только в рамках одной макроплитки, в которую она входит, имеет смысл говорить о {\it плоском} коде. Например, если весь путь лежит в одной плоскости и не выходит в подклееные области, то код такого пути будет плоским. Плоский код одной вершины будет состоять из следующих частей:

1) тип вершины;

2) уровень вершины;

3) окружение;

4) информация.

Мы последовательно определим, что значит каждая из этих частей.

\medskip

\subsection{Параметр ``тип'' для вершин}

Все вершины, встречающиеся на комплексе, мы разделим на следующие категории:

1) {\it Угловые} (лежащие в углах подклееных макроплиток или всего комплекса). {\it Тип угловой вершины } определим как один из четырех вариантов углов, в котором она может находиться: $\mathbb{CUL}$, $\mathbb{CUR}$, $\mathbb{CDR}$, $\mathbb{CDL}$. (Corner Up-Left и так далее.)

2) {\it Краевые} (лежащие на стороне подклееной макроплитки). Каждая такая вершина лежит в середине стороны некоторой макроплитки, прилегающей к краю. {\it Тип краевой вершины } определим в соответствии с тем, серединой какой стороны в этой макроплитке она является: $\mathbb{L}$, $\mathbb{R}$, $\mathbb{D}$, $\mathbb{U}$.

3)  {\it Регулярные} (лежащие внутри некоторой макроплитки, либо на границе между двумя макроплитками).
В этой категории определим три типа {\it внутренних} вершин: $\mathbb{A}$, $\mathbb{B}$, $\mathbb{C}$, отвечающих внутренним, черным узлам макроплиток. А также определим восемь типов {\it боковых} вершин, лежащих в середине сторон двух соседних макроплиток. Учитывая ориентацию этих макроплиток, типы будут следующие: $\mathbb{DR}$, $\mathbb{RD}$, $\mathbb{DL}$, $\mathbb{LD}$, $\mathbb{UR}$,  $\mathbb{RU}$, $\mathbb{UL}$, $\mathbb{LU}$.

\begin{figure}[hbtp]
\centering
\includegraphics[width=0.9\textwidth]{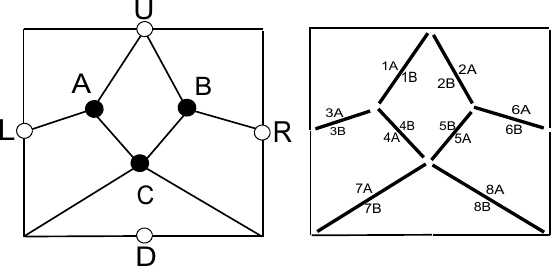}
\caption{Черные (внутренние) и белые (боковые) вершины; типы ребер}
\label{fig:blacknodes}
\end{figure}

В дальнейшем, вершину типа $\mathbb{A}$ будем, для простоты, называть $\mathbb{A}$-вершиной или $\mathbb{A}$-узлом. Аналогично для других типов. Вообще, вершины на графе мы также будем называть узлами.
В левой части рисунка~\ref{fig:blacknodes} представлены внутренние и боковые вершины. Следующее предложение показывает, что возможны только указанные сочетания сторон.

\begin{proposition}[О комбинациях сторон на границе макроплиток] \label{side_combos}

Пусть вершина $X$ лежит в середине сторон двух макроплиток, лежащих в одной плоскости. Тогда возможны только следующие сочетания ориентаций сторон: правая и нижняя, верхняя и правая,  левая и нижняя, либо верхняя и левая.
\end{proposition}

{\bf Доказательство}. По второй части леммы о боковой вершине, $X$ лежит на некотором ребре $YZ$, которое относится к одному из восьми типов. (Правая часть рисунка~\ref{fig:blacknodes}). Рассмотрим вершину $F$, лежащую в середине этого ребра.  $F$ тоже боковая вершина и сочетания ориентаций сторон у нее могут быть только право и низ, верх и право, лево и низ, верх и лево, то есть как в условии леммы. Если $F$ совпадает с $X$, то все доказано. Иначе рассмотрим разбиение наших соприкасающихся макроплиток.  Легко проверить, что, согласно правилам разбиения, середины ребер $YF$ и $FZ$ тоже могут иметь сочетание ориентаций только право и низ, верх и право, лево и низ, верх и лево ($\mathbb{RD}$, $\mathbb{UR}$, $\mathbb{LD}$, $\mathbb{UL}$ и транспозиции). Например, если макроплитки соприкасаются нижней и правой сторонами, то $F$ имеет тип $\mathbb{DR}$ или $\mathbb{RD}$, и тогда $YF$ и $FZ$ будут иметь типы $\mathbb{RU}$ и $\mathbb{UL}$ (или наоборот). Повторяя операцию разбиения, можно получить, что все боковые вершины на ребре $YZ$ будут относиться к одному из указанных четырех типов, то есть являются серединами правой и нижней, либо верхней и правой, либо левой и нижней, либо верхней и левой сторон.

\medskip

\subsection{Типы ребер и параметр ``уровень'' для вершин}

Будем считать, что вершина {\it принадлежит} макроплитке, если она образуется при первом разбиении этой макроплитки (когда она из плитки становится макроплиткой второго уровня). То есть, белые вершины принадлежат тем макроплиткам, в середине сторон которых они лежат, а черные -- макроплиткам минимального уровня, внутри которых они находятся.

\medskip

{\bf Типы ребер.} Всего существует $16$ внутренних ребер (с учетом стороны, рисунок~\ref{fig:blacknodes}), это ребра $\mathbf{1A}$, $\mathbf{1B}$, \dots ,$\mathbf{8B}$ и $4$ краевых ребра: $\mathbf{left}$, $\mathbf{top}$, $\mathbf{right}$, $\mathbf{bottom}$.

\medskip

{\bf Уровни вершин.}
Для черных, внутренних вершин уровень определять не будем (то есть, можно считать, что у всех один уровень).

%Для черных, внутренних вершин будет три возможных уровня. Первый будет у вершин, принадлежащих макроплиткам второго уровня (самым маленьким макроплиткам), второй у вершин, принадлежащих макроплиткам третьего уровня, и третий -- у всех остальных черных вершин.

%У угловых вершин уровень определим аналогично: первый у вершин, принадлежащих макроплиткам второго уровня (только что подклееная макроплитка), второй у вершин, принадлежащих макроплиткам третьего уровня, и третий у всех остальных угловых вершин (принадлежащих макроплиткам четвертого уровня и выше)

Для белых, боковых вершин уровня будет три: первый у вершин, принадлежащих макроплиткам второго уровня, второй -- у вершин, принадлежащих макроплиткам третьего уровня, третий -- у всех остальных белых вершин.

\medskip

Рассмотрим черную вершину $X$. Если $X$ имеет тип $\mathbb{A}$ и $\mathbb{B}$, то она находится на стыке трех макроплиток (разбивающих макроплитку, которой $X$ принадлежит), для $\mathbb{C}$ таких макроплиток четыре. Назовем ребра, исходящие из $X$ и при этом лежащие на границе каких-то двух из этих макроплиток, {\it главными ребрами} для $X$.

%Если $X$ -- принадлежит макроплитке второго уровня, то все выходящие из нее ребра главные, их три для $\mathbb{A}$ и $\mathbb{B}$ и четыре для $\mathbb{C}$. Если $X$ -- принадлежит макроплитке третьего уровня, для $\mathbb{A}$, $\mathbb{B}$, $\mathbb{C}$ к главным ребрам добавляются $2$, $3$, $3$ неглавных соответственно.
%Для вершины $\mathbb{С}$ принадлежащей макроплитке четвертого уровня и выше, добавляются еще три неглавных ребра. Рассмотрим два соседних исходящих из $X$ ребра. Они ограничивают некоторую плитку или макроплитку, причем $X$ является одним из четырех ее углов. Легко проверяется, что этот угол при $X$ является простым, то есть не между нашими соседними ребрами не добавится еще ребро при дальнейших разбиениях. Таким образом, при дальнейших разбиениях исходящих ребер для $X$ не добавится, и поэтому все черные вершины можно классифицировать по описанным трем уровням.

\medskip

Пусть  $X$ -- белая вершина, лежащая на середине сторон двух макроплиток $T_1$ и $T_2$. Из $X$ выходит два ребра, лежащих на границе между $T_1$ и $T_2$, их будем называть {\it главными}. Если боковая вершина лежит на краю только одной макроплитки, то главные ребра -- это лежащие на ее границе. В зависимости от типа $X$, из нее выходят еще несколько неглавных ребер. Мы можем разобрать все белые вершины третьего уровня и убедиться, что все прилежащие к $X$ углы будут простые и при дальнейших разбиениях новых исходящих ребер не добавляется.

\medskip

{\bf Ребра входа и выхода. Обозначения.}
Первым главным ребром для вершин типа $\mathbb{A}$ и $\mathbb{B}$ будем называть то главное ребро, которое уходит к середине верхней стороны той макроплитки, которой принадлежит данная вершина.
Для $\mathbb{C}$ первым главным ребром будем считать ребро, идущее в $\mathbb{A}$-узел.

\medskip

Для боковых вершин оба главных ребра являются двумя частями некоторого внутреннего ребра. Первым из них будем считать то ребро, относительно которого А-сторона (из определения типа внутреннего ребра) остается по правую сторону. Второе и третье главные ребра определяются по часовой стрелке после первого. Все главные ребра будем обозначать цифрами (от $1$ до $4$).

\medskip

Для краевых вершин первым главным ребром будем считать то ребро, которое соответствует обходу макроплитки по часовой стрелке. То есть, макроплитка остается по правую сторону от направления этого ребра. Ребро, идущее против часовой стрелке, будет вторым.

\medskip

У вершин типа $\mathbb{A}$ может быть два неглавных ребра, одно уходит в левую верхнюю подплитку, другое в левую нижнюю. Обозначим их, соответственно, $\mathbf{lu}$  и $\mathbf{ld}$.
У вершин типа $\mathbb{B}$ может быть три неглавных ребра, одно уходит в правую верхнюю подплитку, второе в правую нижнюю, третье в среднюю. Обозначим их, соответственно, $\mathbf{ru}$, $\mathbf{rd}$, $\mathbf{mid}$.
У вершин типа $\mathbb{C}$ может быть шесть неглавных ребра, два уходят в левую нижнюю подплитку, два в среднюю, и два в нижнюю. В каждой паре одно из ребер иерархически старше (появляется при разбиении раньше). Обозначим ребра, соответственно, $\mathbf{ld}_1$, $\mathbf{ld}_2$, $\mathbf{mid}_1$, $\mathbf{mid}_2$, $\mathbf{d}_1$, $\mathbf{d}_2$. Более старшему ребру даем первый номер, другому -- второй.

\begin{figure}[hbtp]
\centering
\includegraphics[width=0.5\textwidth]{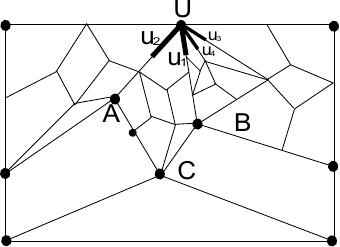}
\caption{Неглавные ребра выходящие из бокового узла}
\label{fig:upedges}
\end{figure}

Для боковых и краевых вершин будем давать неглавным ребрам имена в соответствии с типом вершины и стороны, откуда они уходят:

$\mathbb{U}$-узел или $U$-часть $\mathbb{UL}$-узла или $\mathbb{UR}$-узла -- это $\mathbf{u}_1$, $\mathbf{u}_2$, $\mathbf{u}_3$, $\mathbf{u}_4$;

$\mathbb{R}$-узел или $R$-часть $\mathbb{UR}$-узла или $\mathbb{DR}$-узла -- $\mathbf{r}$, $\mathbf{r}_2$, $\mathbf{r}_3$;

$\mathbb{L}$-узел или $L$-часть $\mathbb{UL}$-узла или $\mathbb{LD}$-узла -- $\mathbf{l}$, $\mathbf{l}_2$, $\mathbf{l}_3$.

В $\mathbb{D}$-узел неглавные ребра не приходят.

Ребра нумеруются по иерархическому старшинству, при равном старшинстве -- по часовой стрелке. Например, есть четыре неглавных ребра из узла типа $\mathbb{LU}$, входящих внутрь плитки $T$, примыкающей верхней стороной (рисунок~\ref{fig:upedges}). $\mathbf{u}_1$ будет ребром в сторону $\mathbb{B}$-узла $T$, $\mathbf{u}_2$ -- ребром в сторону $\mathbb{A}$-узла $T$. $\mathbf{u}_3$ и $\mathbf{u}_4$,  будут входить в верхнюю правую подплитку $T$.

\medskip

{\bf Типы ребер, ведущих в подклееные области.}
Рассмотрим вершину $X$. В подклееной макроплитке она может иметь один из следующих типов:
$\mathbb{L}$, $\mathbb{U}$, $\mathbb{CUR}$, $\mathbb{CUL}$, $\mathbb{CDL}$.

\begin{figure}[hbtp]
\centering
\includegraphics[width=0.9\textwidth]{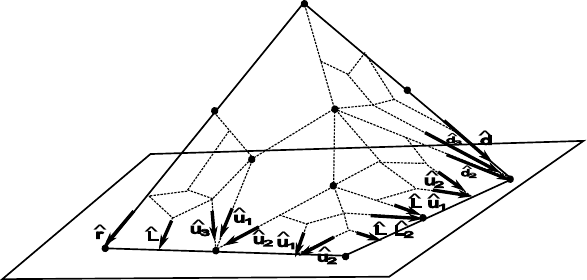}
\caption{Ребра в подклееные области}
\label{fig:pastingedges}
\end{figure}

Из вершины $\mathbb{L}$ могут выходить ребра $\mathbf{l}$, $ \mathbf{l}_2 $, $\mathbf{l}_3$ (как и в плоском случае). Чтобы показать, что имеются в виду ребра в подклееную область, будем записывать их как $\widehat{\mathbf{l}_1}$, $\widehat{\mathbf{l}_2}$, $\widehat{\mathbf{l}_3}$ (рисунок~\ref{fig:pastingedges}).

Аналогично, из вершины $\mathbb{U}$ могут выходить ребра $\widehat{\mathbf{u}_1}$, $\widehat{\mathbf{u}_2}$, $\widehat{\mathbf{u}_3}$, $\widehat{\mathbf{u}_4}$.

Из вершины $\mathbb{CUR}$ может выходить только одно ребро в подклееную область -- это ребро по правой стороне подклееной макроплитки. Будем обозначать его как $\widehat{\mathbf{r}}$.

Из вершины $\mathbb{CDL}$ может выходить два ребра в подклееную область -- по нижней стороне подклееной макроплитки, по ребру $7$ подклееной макроплитки, а также по $8$ ребру макроплитки, образующейся при разбиении нижней подплитки подклееной макроплитки. Будем обозначать их как $\widehat{\mathbf{d}}$, $\widehat{\mathbf{d}_2}$, $\widehat{\mathbf{d}_3}$ соответственно.

Из вершины $\mathbb{CUL}$ (она же -- ядро подклееной макроплитки) ребер в подклееные области не выходит.

\medskip

{\bf Цепи.}
Выберем некоторую вершину $X$. Обозначим как $T(X)$ множество всех макроплиток одного размера, для которых $X$ является левым верхним углом. ($T(X)$ может быть выбрано несколькими способами, для разных размеров макроплиток.)

\begin{definition}

Совокупность ребер $1$ и $3$ типов в макроплитках из $T(X)$ и {\bf боковых} узлов в их концах будем называть {\it цепью} для $X$. Будем называть вершину $X$ {\it центром} цепи.

\end{definition}

Из определения следует, что каждый узел может входить только в одну цепь, причем центр для этой цепи определяется однозначно.
Под {\it макроплитками из цепи} будем понимать макроплитки из соответствующего множества $T(X)$.
Заметим, что между ребрами, выходящими из $X$ и вершина цепи, а  то есть, количество узлов в любой цепи ограничено сверху (Напомним, что в цепь могут входить только боковые вершины).

\smallskip

Для вершин в одной цепи можно задать упорядоченность: первой будем считать тот узел, который лежит на первом главном ребре, выходящем из $X$, а остальные перечисляются в соответствии с выходящими из $X$ ребрами по часовой стрелке. Номер вершины в ее цепи будем называть {\it указателем}.

\medskip

Для выбранной вершины $X$ можно построить несколько цепей с центром $X$ (в каждой будут макроплитки одного размера). Каждой цепи можно присвоить уровень: для самой крупной цепи -- нулевой ($0$-цепь), далее первый и так далее. Если узел $X$ принадлежит некоторой макроплитке $T$, то ребра цепи уровня $k$ принадлежат макроплитке, получающейся после $k+1$ разбиений $T$.

\medskip

\begin{figure}[hbtp]
\centering
\includegraphics[width=1\textwidth]{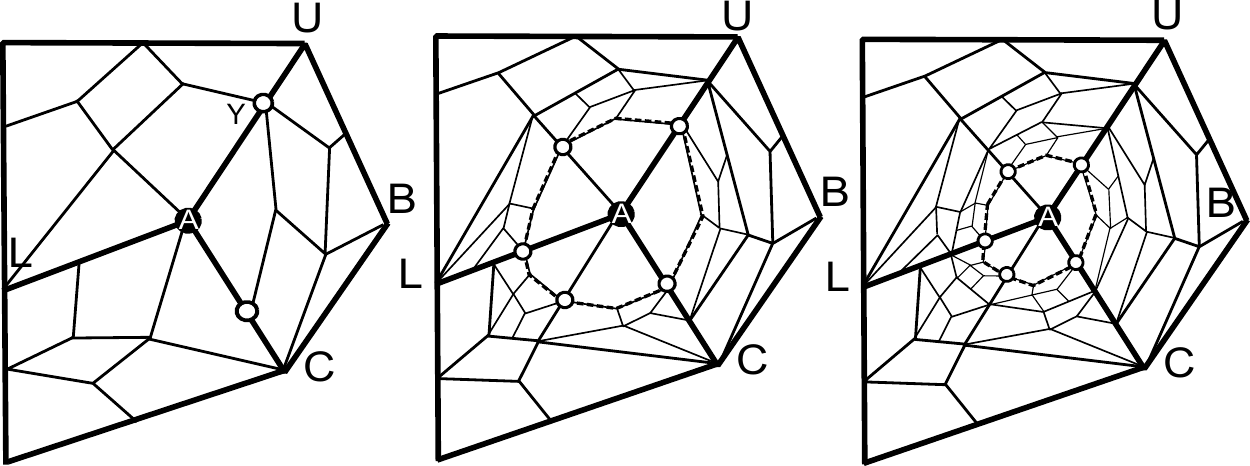}
\caption{Цепи узла $\mathbb{A}$.}
\label{fig:chainAdef}
\end{figure}

На рисунке показаны цепи нулевого, первого и второго уровней c центром типа $\mathbb{A}$.

\medskip

{\bf Параметр ``окружение'' для узлов.}

\begin{definition}

Для каждой макроплитки можно записать упорядоченную четверку типов ребер, на котором лежат ее стороны: левая, верхняя, правая, нижняя (строго в данном порядке). Такую четверку будем называть {\it окружением} данной макроплитки.

\end{definition}

\medskip

Рассмотрим некоторую вершину $X$. Она может входить в одну или несколько подклееных областей. Рассмотрим одну из таких подклееных областей $P$.

\begin{definition}

 {\it Меткой вершины относительно $P$} будем называть:

{\bf 1.} Для регулярных внутренних вершин (типы $\mathbb{A}$, $\mathbb{B}$, $\mathbb{C}$) -- окружение макроплитки, которой принадлежит данная вершина.

\smallskip
 {\bf Пример:}  {\it У узла типа $\mathbb{A}$ внутри средней подплитки метка $[\mathbf{1B},\mathbf{2B},\mathbf{4B},\mathbf{5B}]$}.

\smallskip

{\bf 2.} Для регулярных боковых вершин (типы $\mathbb{UL/LU}$, $\mathbb{UR/RU}$, $\mathbb{LD/DL}$, $\mathbb{RD/DR}$) -- упорядоченную пару окружений двух макроплиток, в середине сторон которых лежит данная вершина.

\smallskip

    {\bf Пример: } {\it У узла $Y$ на левой части рисунка~\ref{fig:chainAdef} метка $[x,y,\mathbf{1A},\mathbf{3A}]-[\mathbf{1B},\mathbf{2B},\mathbf{4B},\mathbf{5B}]$}, где $x$, $y$ -- {\it типы} левой и верхней сторон всей макроплитки рисунка.

\smallskip

{\bf 3.} Для угловых вершин -- метка будет ($\mathbf{left}$,$\mathbf{top}$,$\mathbf{right}$,$\mathbf{bottom}$), то есть как окружение подклееной макроплитки.

{\bf 4.} Для краевых вершин -- метка соответствует окружению макроплитки, в середине стороны которой лежит вершина.

\end{definition}

\medskip

{\bf Замечание}. У вершины может быть несколько меток -- по числу подклееных областей, в которые она входит. Ясно, что только в своей базовой плоскости вершина может быть регулярной. В подклееных областях она будет обязательно угловой или краевой.

\medskip

\begin{definition}

{\it Окружением цепи} будем называть упорядоченный набор меток всех вершин цепи.

\end{definition}

\smallskip

  {\bf Пример.} {\it У цепи в средней части рисунка~\ref{fig:chainAdef} окружение

   $$ [\mathbf{1A},\mathbf{8A},\mathbf{5A},\mathbf{6B}]:[\mathbf{2B},\mathbf{4B},\mathbf{1A},\mathbf{3A}]\times [\mathbf{2B},\mathbf{4B},\mathbf{1A},\mathbf{3A}]:[\mathbf{8B},\mathbf{4A},\mathbf{4A},\mathbf{7B}]\times $$

   $$ \times [\mathbf{8B},\mathbf{4A},\mathbf{4A},\mathbf{7B}]:[\mathbf{3B},\mathbf{8A},\mathbf{5A},\mathbf{6B}]\times [\mathbf{3B},\mathbf{8A},\mathbf{5A},\mathbf{6B}]:[\mathbf{8B},\mathbf{3A},\mathbf{3A},\mathbf{7B}] \times  $$

   $$\times [\mathbf{8B},\mathbf{3A},\mathbf{3A},\mathbf{7B}]:[\mathbf{1A},\mathbf{8A},\mathbf{5A},\mathbf{6B}]$$}

\smallskip

Ниже будет показано, что тип узла в центре цепи уровня выше нулевого полностью определяет окружение этой цепи.

\medskip

\begin{definition}

{\it Окружением вершины} относительно плоскости $P$ будем называть:

{\bf 1.} Для регулярных боковых вершин типов $\mathbb{UL/LU}$, $\mathbb{UR/RU}$, $\mathbb{LD/DL}$, а также для краевых вершин типов $\mathbb{U}$ и $\mathbb{L}$  -- упорядоченную пару $(x,y)$, где $x$ --  окружение цепи, в которую входит данная вершина, $y$ -- указатель, номер данной вершины в ее цепи.

\smallskip

{\bf 2.} Для регулярных боковых вершин типов $\mathbb{RD/DR}$ -- тип ребра, на которой находится данная вершина (это ребро $2$, $3$, $5$ или $6$).

{\bf 3.} Для всех остальных вершин (угловые, краевые, регулярные типов $\mathbb{A}$, $\mathbb{B}$, $\mathbb{C}$) -- метку этой вершины относительно плоскости $P$.

\end{definition}

Таким образом, если вершина входит в цепь, то зная окружение вершины, мы знаем окружения всех вершин в этой цепи.
Окружение вершины $X$ относительно плоскости $P$ будем обозначать как $\mathbf{SurrP}(X)$. Если из контекста ясно, о какой плоскости идет речь, будем использовать обозначение $\mathbf{Surr}(X)$.

\medskip

Пусть $X$ -- узел, входящий в некоторую цепь. Если $Y$ -- центр этой цепи, то $X$ соответствует одному из выходящих из $Y$ ребер. Номер этого ребра, (начиная с первого выходящего из $X$) мы называем {\it указателем} для вершины $X$. Очевидно, что по указателю и окружению цепи можно установить окружение самого узла.

\medskip

\begin{definition}

И наконец, определим {\it расширенное окружение} относительно области~$P$. Это упорядоченная пара окружения относительно $P$ и окружения относительно базовой плоскости вершины.
\end{definition}

Таким образом, расширенное окружение несет информацию не только об окрестности вершины в подклейке, но и об окрестности в базовой плоскости.

\medskip

\subsection{Параметр ``информация'' для узлов}

{\bf Начальники.}

Определим понятие {\it начальников} для регулярных вершин.
Пусть вершина $X$ -- регулярная ( то есть не угловая, не боковая). Тогда можно выбрать такую минимальную макроплитку $T$, что $X$  находится внутри (не на границе) $T$. В этом случае, $X$ либо является внутренней вершиной ($\mathbb{A}$, $\mathbb{B}$ или $\mathbb{C}$), либо лежит на одном из восьми внутренних ребер $T$.

\begin{definition}

{\it Первым начальником} вершины $X$ будем называть узел в середине верхней стороны макроплитки $T$.

Для боковых вершин, лежащих на ребрах $2$, $5$, $6$ макроплитки $T$, {\it вторым начальником} будем называть узел в правом нижнем углу макроплитки $T$.

Для боковых вершин, лежащих на ребрах $7$ и $8$, а также для для вершин типа $\mathbb{C}$ определим {\it второго} и {\it третьего начальников} как узлы в левом нижнем и правом нижнем углах $T$ соответственно.

\end{definition}

У краевых вершин начальников не будет. Из угловых вершин определим начальника только для $\mathbb{CDR}$-узла, и им будет $\mathbb{CDL}$-узел в той же макроплитке.

\medskip

Теперь дадим определение информации. Оно будет зависеть от расположения вершины внутри макроплитки.

\begin{definition}

Для регулярных боковых вершин, лежащих на внутренних ребрах $1$, $3$, $4$, а также для $\mathbb{A}$-узлов {\it информацией вершины} будем называть упорядоченную тройку состоящую из типа, уровня и расширенного окружения ее единственного начальника.

Для регулярных боковых вершин, лежащих на внутренних ребрах $2$, $5$, $6$, а также для $\mathbb{B}$-узлов {\it информацией вершины} будем называть упорядоченную четверку состоящую из типа, уровня и расширенного окружения ее первого начальника, а также типа ее второго начальника.

Для регулярных боковых вершин, лежащих на внутренних ребрах $7$, $8$, а также для $\mathbb{C}$-узлов {\it информацией вершины} будем называть упорядоченную девятку состоящую из типов, уровней и расширенных окружения всех трех ее начальников.

У краевых вершин информации не будет.

Из угловых вершин информацию определим только для $\mathbb{CDR}$-узла, и она будет равна упорядоченной тройке состоящей из типа, уровня и расширенного окружения ее единственного начальника, $\mathbb{CDL}$-узла в той же макроплитке.

\end{definition}

\medskip

{\bf Обозначения.} Информацию о первом начальнике узла $X$ будем обозначать как $\mathbf{FBoss}(X)$. Иногда это обозначение будет употребляться для типа или окружения начальника, тогда это будет понятно по контексту. Аналогичны обозначения для второго и третьего начальников:
$\mathbf{SBoss}(X)$, $\mathbf{TBoss}(X)$.

\medskip

Таким образом, информация -- это данные о начальниках, а для боковых вершин на ребрах $2$, $5$, $6$ -- это еще знание, какой тип у узла в правом нижнем углу.

\medskip

\subsection{Ядро подклейки и флаг макроплитки}

Рассмотрим некоторую макроплитку $T$, подклееную к ребрам $e_1$ и $e_2$ некоторой вершины $X$.
Вершину $X$ будем называть {\it ядром подклейки}. Она лежит в левом верхнем углу $T$. Сочетание значений параметров типа, базового окружения, информации $X$, а также упорядоченную пару типов выходящих ребер $e_1$ и $e_2$ будем называть {\it флагом подклейки} для макроплитки $T$.

Для каждой вершины $Y$, лежащей внутри подклееной макроплитки $T$ (не на ее краю) определим параметр {\it флаг подклейки}, со значением, равным значению этого параметра для макроплитки. Таким образом, для всех таких вершин $Y$, лежащих внутри подклееной макроплитки значение этого параметра будет одинаковым.

\subsection{Кодировка путей}

\smallskip

Рассмотрим путь $X_1 e_1 X_2 e_2 X_3 e_3 \dots e_{n-1} X_n$, где $X_i$ -- вершины, $e_i$ -- ребра между соседними вершинами в пути.
Обозначим $\mathbf{In}(e_i)$ -- тип $e_i$ как выходящего ребра из вершины $X_i$. Также пусть  $\mathbf{Out}(e_i)$ -- тип $e_i$ как входящего ребра в вершину $X_{i+1}$.

\medskip

Рассмотрим некоторую вершину $X_i$. Если и входящее в $X_i$ ребро $e_{i-1}$, и выходящее ребро $e_i$ являются плоскими, то кодом $X_i$ будем считать упорядоченную четверку, состоящую из типа $X_i$, ее уровня, ее окружения и ее информации.

\smallskip

Если одно из ребер, например входящее ребро $e_{i-1}$ приходит из подклееной области $P_{i-1}$, а выходящее ребро $e_i$ -- плоское, то кодом $X_i$ будем считать упорядоченный набор из четырех элементов:
\begin{enumerate}

\item {\it первый элемент} -- упорядоченная пара из типа $X_i$ в области $P_{i-1}$ и типа $X_i$ в базовой плоскости;

\item {\it второй элемент} -- упорядоченная пара из уровня $X_i$ в области $P_{i-1}$ и уровня $X_i$ в базовой плоскости;

\item {\it третий элемент} -- упорядоченная пара из окружения $X_i$ в области $P_{i-1}$ и окружения $X_i$ в базовой плоскости;

\item {\it четвертый элемент} -- информация у $X_i$.

\end{enumerate}

\smallskip

Таким образом кодировка учитывает и положение вершины в макроплитке, откуда пришло ребро и в макроплитке базовой плоскости. Для случая, когда входящее ребро -- плоское, а выходящее ведет в подклейку, определение аналогично, только в упорядоченной паре (типов, уровней, окружений) сначала будет идти базовая плоскость, а потом подклееная.

\medskip

Если оба ребра, и входящее $e_{i-1}$, и выходящее $e_i$, ведут в подклееные области $P_{i-1}$ и $P_i$, то кодом $X_i$ будем считать упорядоченный набор из четырех элементов:

\begin{enumerate}

\item {\it первый элемент} -- упорядоченная тройка -- тип $X_i$ в области $P_{i-1}$, тип $X_i$ в базовой плоскости, тип $X_i$ в области $P_i$;

\item {\it второй элемент} -- упорядоченная тройка -- уровень $X_i$ в области $P_{i-1}$, уровень $X_i$ в базовой плоскости, уровень $X_i$ в области $P_i$;

\item {\it третий элемент } -- упорядоченная тройка -- окружение $X_i$ в области $P_{i-1}$, окружение $X_i$ в базовой плоскости, окружение $X_i$ в области $P_i$;

\item {\it четвертый элемент} -- информация у $X_i$.

\end{enumerate}

\medskip

{\bf Замечание}.  Как видно, информация не подчиняется общей логике. Это происходит по причине того, что любая вершина в своих подклееных макроплитках всегда занимает место на верхней или левой сторонах, либо в лежит в любом углу, кроме правого нижнего. Во всех этих случаях информация вершины в рамках такой подклееной макроплитки, по определению, пустая.

\smallskip

Код вершины $X$ будем обозначать как $\mathbf{Code}(X)$. Теперь определим код всего пути.

\begin{definition}
{\it Кодом пути} $X_1 e_1 X_2 e_2 X_3 e_3 \dots e_{n-1} X_n$ будем называть упорядоченный набор

\smallskip

$\mathbf{Code}(X_1)$, $\mathbf{Out}(e_1)$, $\mathbf{In}(e_1)$,  $\mathbf{Code}(X_2)$, $\mathbf{Out}(e_2)$, \dots , $\mathbf{In}(e_{n-1})$, $\mathbf{Code}(X_n)$.

\smallskip

Здесь $\mathbf{Code}(X_i)$ -- код вершины $X_i$, $\mathbf{Out}(e_i)$ -- тип $e_i$ как выходящего ребра из вершины $X_i$, $\mathbf{In}(e_i)$ -- тип $e_i$ как входящего ребра в вершину $X_{i+1}$.

\end{definition}

\medskip

Нашей основной целью будет определить локальные преобразования кодов путей так, чтобы переход от одного кода к другому соответствовал переходу к эквивалентному пути на комплексе.

При этом элементарным преобразованиям путей будут соответствовать замены одних слов-кодов путей на другие. То есть элементарные преобразования путей являются аналогами определяющих соотношений в полугруппе слов-кодов.

Для того, чтобы задать такие определяющие соотношения, нужно рассмотреть все локальные преобразования путей и для каждого из них указать, какой паре равных кодов соответствует данное локальное преобразование.

Часто знание типа, окружения или информации некоторой вершины помогает нам понять устройство ближайшей окрестности. Для удобства, мы зададим набор функций на узлах, аргументом которых являются типы, окружения или информации некоторых специально расположенных узлы, а значениями -- окружения или типы узлов в окрестности узла-аргумента. В следующем параграфе мы определим набор таких функций.

\medskip

\section{Функции на узлах и структура цепей} \label{functions}

Ниже приведены функции, которые облегчают нам работу по восстановлению параметров вершин.

\medskip

\begin{enumerate}

\item $\mathbf{TopFromCorner}$ Значение: окружение узла в середине верхней стороны $T$.

\item $\mathbf{RightCorner}$  Значение: окружение узла в правом нижнем углу $T$.

\item $\mathbf{TopRightType}$  Значение: тип узла в правом верхнем углу $T$, а также тип ребра для боковых узлов.

\item $\mathbf{BottomLeftType}$  Значение: тип узла в левом нижнем углу $T$,
  а также тип ребра для боковых узлов.

\item $\mathbf{LevelPlus}$  Значение: окружение цепи с уровнем на один больше, чем у $X$, с тем же центром и указателем.

\item $\mathbf{BottomRightTypeFromRight}$ Значение: тип узла в правом нижнем углу $T$, а также тип ребра для боковых узлов.

\item $\mathbf{TopFromRight}$ Значение: окружение узла в середине верхней стороны $T$.

\end{enumerate}

Функции 1 и 2 содержат два аргумента: первый -- это окружение узла $X$, являющимся левым нижним углом в макроплитке $T$, второй -- тип ребра выхода из $X$, соответствующего нижней стороне $T$.

Функции 3, 4, 5 содержат один аргумент: окружение узла $X$, являющегося серединой верхней стороны в макроплитке $T$

Функции 6 и 7 содержат два аргумента: окружение {\bf и информацию} узла $X$, являющимся серединой правой стороны в макроплитке~$T$.

\smallskip

Еще раз отметим, что под окружением узла в цепи мы понимаем окружение всей этой цепи вместе с указателем -- типом ребра входа-выхода, на котором лежит узел из цепи.

\medskip

\subsection{Свойства цепей}

В этом параграфе мы опишем свойства цепей для различных типов узлов.
В частности, мы установим, что зная тип центра, можно узнать окружение цепи, и наоборот. Помимо этого, мы используем свойства цепей для доказательства того, что по аргументу функции можно узнать ее значение.

\medskip

Мы разберем цепи с различными центрами и в каждом случае будем проверять, что если мы знаем аргумент функции, то можем узнать ее значение.

\medskip

\begin{figure}[hbtp]
\centering
\includegraphics[width=1\textwidth]{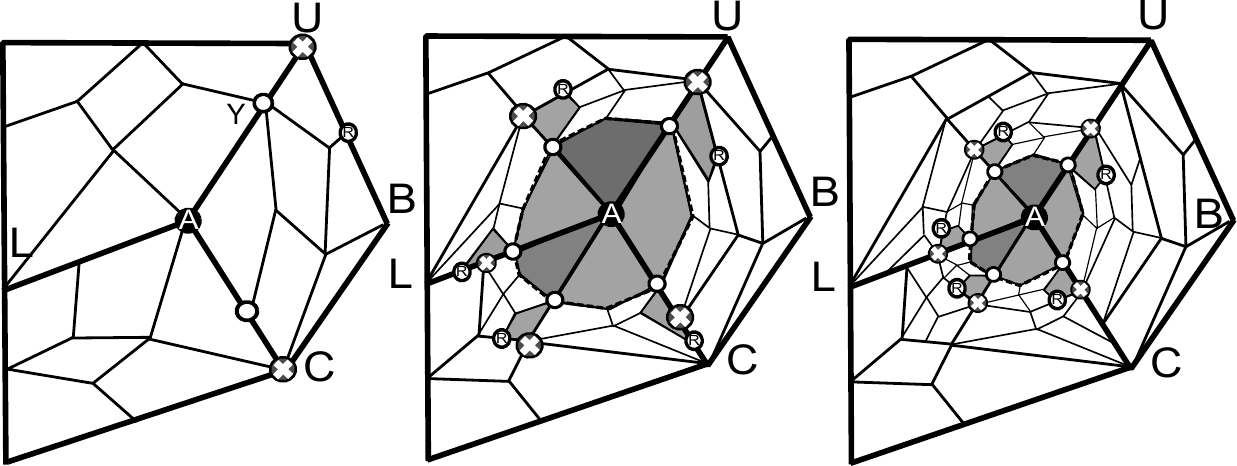}
\caption{Цепи узла $\mathbb{A}$.}
\label{fig:chainA}
\end{figure}

{\bf Цепи с центром в узле типа $\mathbb{A}$.}
На рисунке~\ref{fig:chainA} белыми точками отмечены цепи узла $\mathbb{A}$ нулевого, первого и второго уровней. Заметим, что на втором уровне все пять макроплиток, левый верхний угол которых попадает в узел $\mathbb{A}$, занимают левое верхнее положение в своих родительских макроплитках. Значит, если взять следующий уровень подразбиения, типы ребер-границ макроплиток не изменятся. То есть окружения вершин в цепях будут те же, то есть, окружение цепи третьего и последующих уровней совпадает с окружением цепи второго уровня.

Таким образом, существует только три возможных окружения цепи с центром в узле типа $\mathbb{A}$. Заметим, что по окружению цепи, мы можем установить, является ли центр цепи узлом типа $\mathbb{A}$, и какого уровня цепь. Действительно, одна из макроплиток в цепи с центром в $\mathbb{A}$ обязательно имеет окружение $(\mathbf{4B},\mathbf{1B},\mathbf{2B},\mathbf{5B})$ или $(\mathbf{4B},\mathbf{1B},\mathbf{1A},\mathbf{3A})$, причем, макроплитки с такими окружениями не могут входить в цепи с другими центрами. Цепь нулевого уровня содержит макроплитку с окружением $(\mathbf{4B},\mathbf{1B},\mathbf{2B},\mathbf{5B})$ (остальные -- нет). Для цепи второго уровня все четверки типов ребер имеют вид $(x,y,\mathbf{1A},\mathbf{3A})$, а для первого уровня это не так.

Забегая вперед, можно сказать, что возможность установить тип центра по окружению цепи будет и для центров других типов (не только $\mathbb{A}$).

Заметим также, что мы можем выписать полностью все окружение $1$-цепи и $2$-цепи с центром в узле типа $\mathbb{A}$, а зная окружение узла $\mathbb{A}$, можем выписать и окружение $0$-цепи.

\medskip

\medskip

Из сказанного выше следует, что для $\mathbb{A}$-узлов значение функции $\mathbf{LevelPlus}$ может быть установлено по ее аргументу.

На рисунке~\ref{fig:chainA} знаком ``$\otimes$'' отмечены узлы, являющиеся верхними правыми или левыми нижними углами в макроплитках, где середина верхней стороны попадает в узел цепи $X$ . Можно убедиться в том, что для цепей первого и второго уровней все типы узлов с крестами, а также типы ребер, на которых они лежат, мы можем выписать, для каждого заданного $X$. Для $0$-цепи мы можем выписать типы узлов со знаком ``$\otimes$'' (а также типы ребер), если знаем окружение центра (узла $\mathbb{A}$).
Таким образом, мы можем вычислить функции $\mathbf{TopRightType}$ и $\mathbf{BottomLeftType}$.

%Черными точками на рисунке~\ref{fig:chainA} помечены узлы являющиеся правыми нижними углами в макроплитках, где середина верхней стороны попадает в узел цепи $X$. Заметим, что зная окружение и информацию $X$, мы можем установить тип этого правого нижнего угла. Действительно, для всех узлов кроме одного (на $3$ выходящем из $\mathbb{A}$ ребре в $1$-цепи), мы сразу можем установить их тип. В оставшемся случае, этот угол может быть записан как $\mathbf{Next.FBoss}(X)$. Таким образом, мы можем вычислить функцию $\mathbf{BottomRightTypeFromTop}$.

\medskip

Серым цветом на рисунке выделены плитки, где левый нижний угол (при дальнейшем разбиении) будет попадать в вершину цепи. Можно заметить, что зная эту вершину (ее место в цепи) мы можем установить и окружение правого нижнего угла в соответствующей макроплитке. Кроме того, мы также можем установить тип левого верхнего угла, а значит и всю цепь, содержащую середину верхней стороны в соответствующей макроплитке. Все это значит, что мы можем вычислить функции
$\mathbf{TopFromCorner}$, $\mathbf{RightCorner}$.

 %Рассмотрим все макроплитки, где середина верхней стороны $X$ попадает в цепь. При этом середина
 %правой стороны попадет в один из узлов, помеченных $<<R>>$. Заметим, что зная окружение и информацию
 %$X$, мы можем вычислить окружение этого $R$-узла. (Во всех случаях, кроме одного, хватает просто
 %окружения $X$, для $X$ лежащей на $3$ ребре выхода из $\mathbb{A}$, окружение $R$-вершины будет
 %$\mathbf{Next.FBoss}(X)$). Таким образом, функция $\mathbf{RightFromTop}$ также вычисляется.

Функции $\mathbf{TopFromRight}$ и $\mathbf{BottomRightTypeFromRight}$ могут быть применимы только в одном случае: для $0$-цепи, если аргументом является узел в середине ребра типа $1$ (обозначим этот узел как $Y$, а макроплитку, где $Y$ является серединой правой стороны, как $T$). $\mathbf{BottomRightTypeFromRight}$ в этом случае это просто наш узел $\mathbb{A}$.

Для вычисления $\mathbf{TopFromRight}$ нам надо найти окружение узла в середине верхней стороны макроплитки $T$. Заметим, что этот узел входит в цепь с центром в левом верхнем углу $T$, причем уровень этой цепи на один выше, чем уровень цепи, куда входит правый верхний угол $T$. Окружение этой цепи мы можем установить по информации $Y$, так как правый верхний угол $T$ является серединой верхней стороны в макроплитке $T'$, родительской к $T$ и является первым начальником узла $Y$. По окружению этой цепи можно установить тип левого верхнего угла $T$ и окружение цепи на один уровень выше, что и требуется.

\medskip

Таким образом, мы вычислили значения всех функций для аргументов, входящих в $\mathbb{A}$-цепи.

\medskip

\begin{figure}[hbtp]
\centering
%\leftskip=-1cm
\includegraphics[width=1\textwidth]{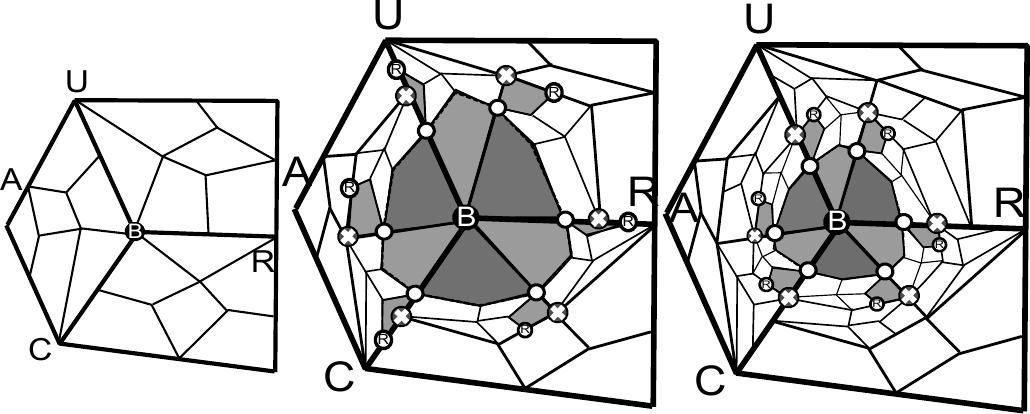}
\caption{Цепи узла $\mathbb{B}$.}
\label{fig:chainB}
\end{figure}

{\bf Цепи с центром в узле типа $\mathbb{B}$.}
На рисунке~\ref{fig:chainB} белыми точками отмечены цепи узла $\mathbb{B}$. $0$-цепи с центром в $\mathbb{B}$ не существует (все боковые узлы в соответствующей окрестности имеют тип $\mathbb{DR}$). Заметим, что в $2$-цепи все шесть макроплиток, левый верхний угол которых попадает в узел $\mathbb{B}$, занимают левое верхнее положение в своих родительских макроплитках. Значит, если взять следующий уровень подразбиения, типы ребер-границ макроплиток не изменятся. То есть окружения вершин в цепях будут те же, то есть, окружение цепи третьего и последующих уровней совпадает с окружением цепи второго уровня.

Таким образом, существует только два возможных окружения цепи с центром в узле типа $\mathbb{B}$. Заметим, что по окружению цепи, мы можем установить, является ли центр цепи узлом типа $\mathbb{B}$, и какого уровня цепь. Действительно, $1$-цепь с центром в $\mathbb{B}$ содержит узел с окружением $(\mathbf{8B},\mathbf{5B},\mathbf{5B},\mathbf{7B})$, который не может содержаться ни в какой другой цепи. $2$-цепь с центром в $\mathbb{B}$ содержит узел с окружением $(\mathbf{8B},\mathbf{5B},\mathbf{1A},\mathbf{3A})$, который также ни в какой другой цепи не встречается.

Заметим также, что мы можем выписать полностью все окружение $1$-цепи и $2$-цепи с центром в узле типа $\mathbb{B}$.

\medskip

\medskip

Из сказанного выше следует, что для $\mathbb{B}$-узлов значение функции $\mathbf{LevelPlus}$ может быть установлено по ее аргументу.
На рисунке~\ref{fig:chainB} знаком ``$\otimes$'' отмечены узлы, являющиеся верхними правыми или левыми нижними углами в макроплитках, где середина верхней стороны попадает в узел цепи $X$. Можно убедиться в том, что для цепей первого и второго уровней все типы узлов с крестами, а также типы ребер, на которых они лежат,  мы можем выписать, для каждого заданного $X$.
Таким образом, мы можем вычислить функции $\mathbf{TopRightType}$ и $\mathbf{BottomLeftType}$.

\medskip

Серым цветом на рисунке выделены плитки, где левый нижний угол (при дальнейшем разбиении) будет попадать в вершину цепи. Можно заметить, что зная эту вершину (ее место в цепи мы знаем, благодаря второму аргументу -- нижнему ребру) мы можем установить и окружение правого нижнего угла в соответствующей макроплитке. Кроме того, мы также можем установить тип левого верхнего угла, а значит и всю цепь, содержащую середину верхней стороны в соответствующей макроплитке. Все это значит, что мы можем вычислить функции
$\mathbf{TopFromCorner}$, $\mathbf{RightCorner}$.

Функции $\mathbf{TopFromRight}$ и $\mathbf{BottomRightTypeFromRight}$ не могут быть применимы, то есть аргументом этих функций никогда не может быть узел из B-цепи так как ни один узел из цепи с центром в $\mathbb{B}$ не является серединой правой стороны ни в какой макроплитке.

%Рассмотрим все макроплитки, где середина верхней стороны $X$ попадает в цепь. При этом середина правой стороны попадет в один из узлов, помеченных $<<R>>$. Заметим, что для $2$-цепи, зная окружение и информацию $X$, мы можем вычислить окружение этого R-узла.

%Для $1$-цепи возможно шесть случаев расположения $X$, в соответствии с выходящими из $\mathbb{B}$ ребрами. Для $\mathbf{ru}$, $\mathbf{rd}$, $\mathbf{mid}$ -- ребер $R$-вершина лежит в середине $5$ ребра, окружение легко вычисляется. Для $1$-ребра окружение $R$-вершины это $1$-цепь вокруг $\mathbf{FBoss}(X)$, с указателем $\mathbf{u}_2$. Для $3$-ребра окружение $R$-вершины это $1$-цепь вокруг $\mathbb{C}$-узла, с указателем $2$. Для $2$-ребра $R$-вершины входит в $1$-цепь вокруг узла в середине правой стороны макроплитки $\mathbb{B}$-узла с указателем $\mathbf{r}$. Тип этого узла может быть установлен по окружению макроплитки, то есть $U$-части окружения узла в середине верхней стороны: $\mathbf{U.FBoss}(X)$.  Для левого верхнего положения тип будет $\mathbb{UR}$, для левого нижнего $\mathbb{DR}$, для среднего $\mathbb{DR}$, правого верхнего $\mathbb{RD}$, для правого нижнего $\mathbb{RD}$. Для нижнего ребра тип будет $\mathbb{RU}$ или $\mathbb{UR}$ в соответствии с ориентацией нижнего ребра ($А$ или $B$). Таким образом, функция $\mathbf{RightFromTop}$ также вычисляется.

\medskip

Таким образом, мы вычислили значения всех функций для аргументов, входящих в $\mathbb{B}$-цепи.

\medskip

\begin{figure}[hbtp]
\centering
%\leftskip=-1cm
\includegraphics[width=1\textwidth]{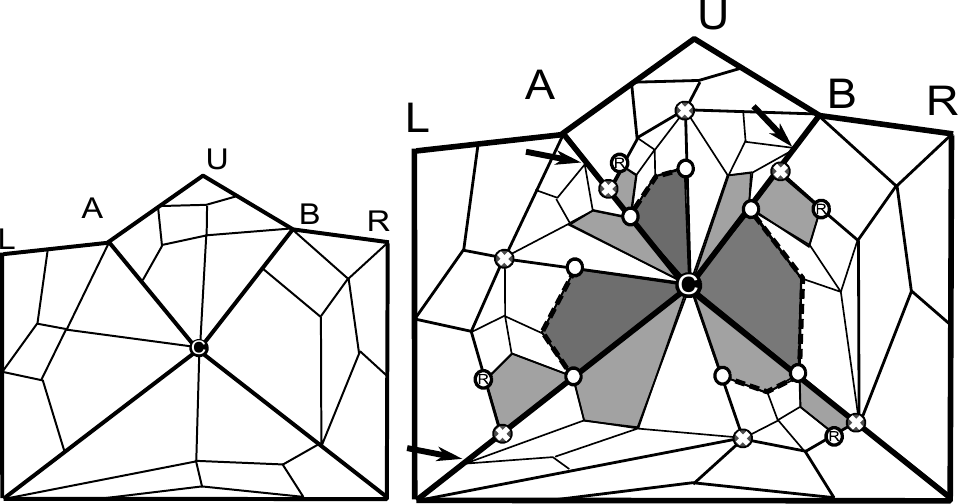}
\caption{$1$-цепь узла $\mathbb{C}$.}
\label{fig:chainC1}
\end{figure}

\begin{figure}[hbtp]
\centering
\includegraphics[width=0.8\textwidth]{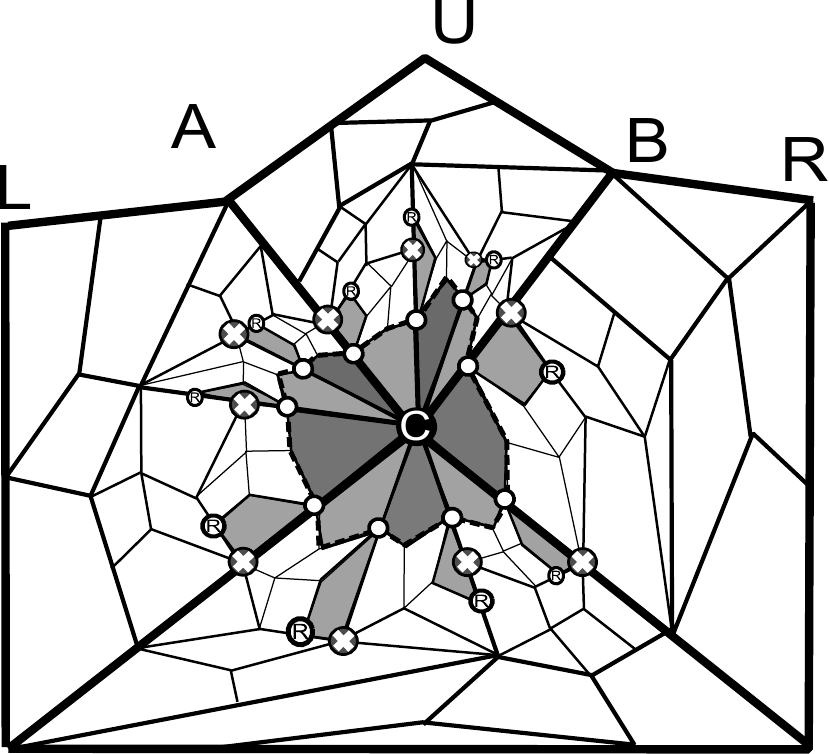}
\caption{$2$-цепь узла $\mathbb{C}$.}
\label{fig:chainC2}
\end{figure}

{\bf Цепи с центром в узле типа $\mathbb{C}$.}
На рисунках~\ref{fig:chainC1} и~\ref{fig:chainC2} белыми точками отмечены $1$-цепи и $2$-цепи узла $\mathbb{C}$. $0$-цепи с центром в $\mathbb{C}$ не существует.

Кроме указанных на рисунках цепей первого и второго уровней существует также $3$-цепь, получаемая применением операции разбиения к макроплиткам на рисунке~\ref{fig:chainC2}.

На получающемся третьем уровне все  макроплитки, левый верхний угол которых попадает в узел $\mathbb{C}$, будут занимать левое верхнее положение в своих родительских макроплитках. Значит, на следующих уровнях подразбиения, типы ребер-границ макроплиток не изменятся. То есть окружения вершин в цепях будут те же, то есть, окружение цепи четвертого и последующих уровней совпадает с окружением цепи третьего уровня.

Таким образом, существует только три возможных окружения цепи с центром в узле типа $\mathbb{C}$ (для цепей первого, второго и третьего уровней). Заметим, что по окружению цепи, мы можем установить, является ли центр цепи узлом типа $\mathbb{C}$, и какого уровня цепь. Действительно, $1$-цепь с центром в $\mathbb{C}$ содержит узел с окружением $(\mathbf{8A},\mathbf{5A},\mathbf{1A},\mathbf{3A})-(\mathbf{8B},\mathbf{5B},\mathbf{5B},\mathbf{7B})$ (на ребре $5$), который не может содержаться ни в какой другой цепи. $2$-цепь с центром в $\mathbb{C}$ содержит узел с окружением $(\mathbf{8A},\mathbf{5A},\mathbf{1A},\mathbf{3A})-(\mathbf{5B},\mathbf{8A},\mathbf{5A},\mathbf{6B})$, который также ни в какой другой цепи не встречается.
Для $3$-цепи таким узлом будет $(\mathbf{8A},\mathbf{5A},\mathbf{1A},\mathbf{3A})-(\mathbf{5B},\mathbf{8A},\mathbf{1A},\mathbf{3A})$.

Заметим также, что мы можем выписать полностью все окружение $1$-цепей, $2$-цепей, $3$-цепей с центром в узле типа $\mathbb{C}$.

\medskip

\leftskip=0cm
\medskip

Из сказанного выше следует, что для $\mathbb{C}$-узлов значение функции $\mathbf{LevelPlus}$ может быть установлено по ее аргументу.
На рисунках~\ref{fig:chainC1} и~\ref{fig:chainC2} знаком ``$\otimes$'' отмечены узлы, являющиеся верхними правыми или левыми нижними углами в макроплитках, где середина верхней стороны попадает в узел цепи $X$ . Можно убедиться в том, что для цепей первого, второго и третьего уровней все типы узлов с крестами, а также типы ребер, на которых они лежат,  мы можем выписать, для каждого заданного $X$.
Таким образом, мы можем вычислить функции $\mathbf{TopRightType}$ и $\mathbf{BottomLeftType}$.

Серым цветом на рисунке выделены плитки, где левый нижний угол (при дальнейшем разбиении) будет попадать в вершину цепи. Можно заметить, что зная эту вершину (ее место в цепи) мы можем установить и окружение правого нижнего угла в соответствующей макроплитке. Кроме того, мы также можем установить тип левого верхнего угла, а значит и всю цепь, содержащую середину верхней стороны в соответствующей макроплитке. (Все это также очевидно проверяется для $3$-цепи.) Все это значит, что мы можем вычислить функции $\mathbf{TopFromCorner}$, $\mathbf{RightCorner}$.

\medskip

\leftskip=0cm

%Рассмотрим все макроплитки, где середина верхней стороны $X$ попадает в цепь. При этом середина правой стороны попадет в один из узлов, помеченных <<$R$>>. Заметим, что зная окружение и информацию $X$, мы можем вычислить окружение этого $R$-узла.
%Таким образом, функция $\mathbf{RightFromTop}$ также вычисляется.

Функции $\mathbf{TopFromRight}$ и $\mathbf{BottomRightTypeFromRight}$ могут быть применимы только к $1$-цепи (в остальные входят только узлы типов $\mathbb{UL}$ и $\mathbb{LU}$). В случае $1$-цепи аргументом могут быть узлы, лежащие на ребрах $1$, $2$ и $4$ (относительно $\mathbb{C}$). Пусть $X$ такой узел.  Во всех случаях
$\mathbf{BottomRightTypeFromRight}$ будет $\mathbb{C}$. Узлы, на которые указывает $\mathbf{TopFromRight}$ отмечены стрелками на рисунке~\ref{fig:chainC1}. Таким образом, $\mathbf{TopFromRight}$ принимает значения: $1$-цепи вокруг $\mathbb{A}$ для $1$ ребра (с указателем $2$), $1$-цепи вокруг $\mathbb{B}$ для ребра $2$ (с указателем $3$),
$2$-цепь вокруг $\mathbf{SBoss}(X)$ для ребра $4$ (с указателем, соответствующим входу ребра $7$ в левый нижний угол (см параграф~\ref{pointers} ``Указатели'').

\medskip

Таким образом, мы вычислили значения всех функций для аргументов, входящих в $\mathbb{C}$-цепи.

\medskip

\begin{figure}[hbtp]
\centering
\includegraphics[width=1\textwidth]{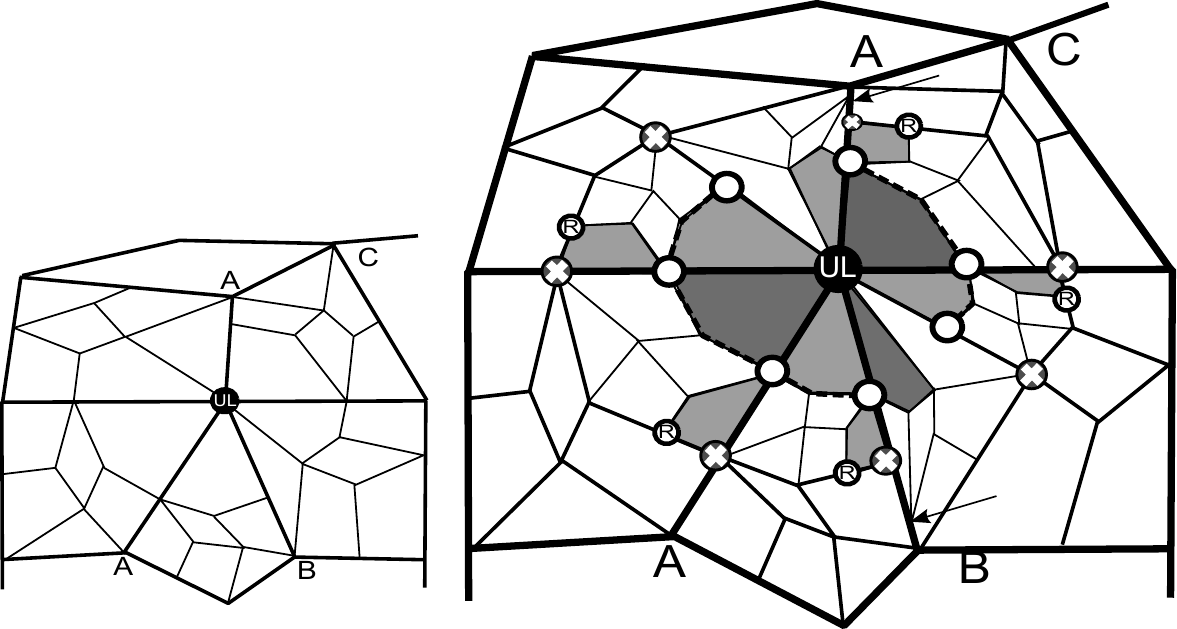}
\caption{$1$-цепь узла $\mathbb{UL}$.}
\label{fig:chainUL1}
\end{figure}

\begin{figure}[hbtp]
\centering
\includegraphics[width=0.8\textwidth]{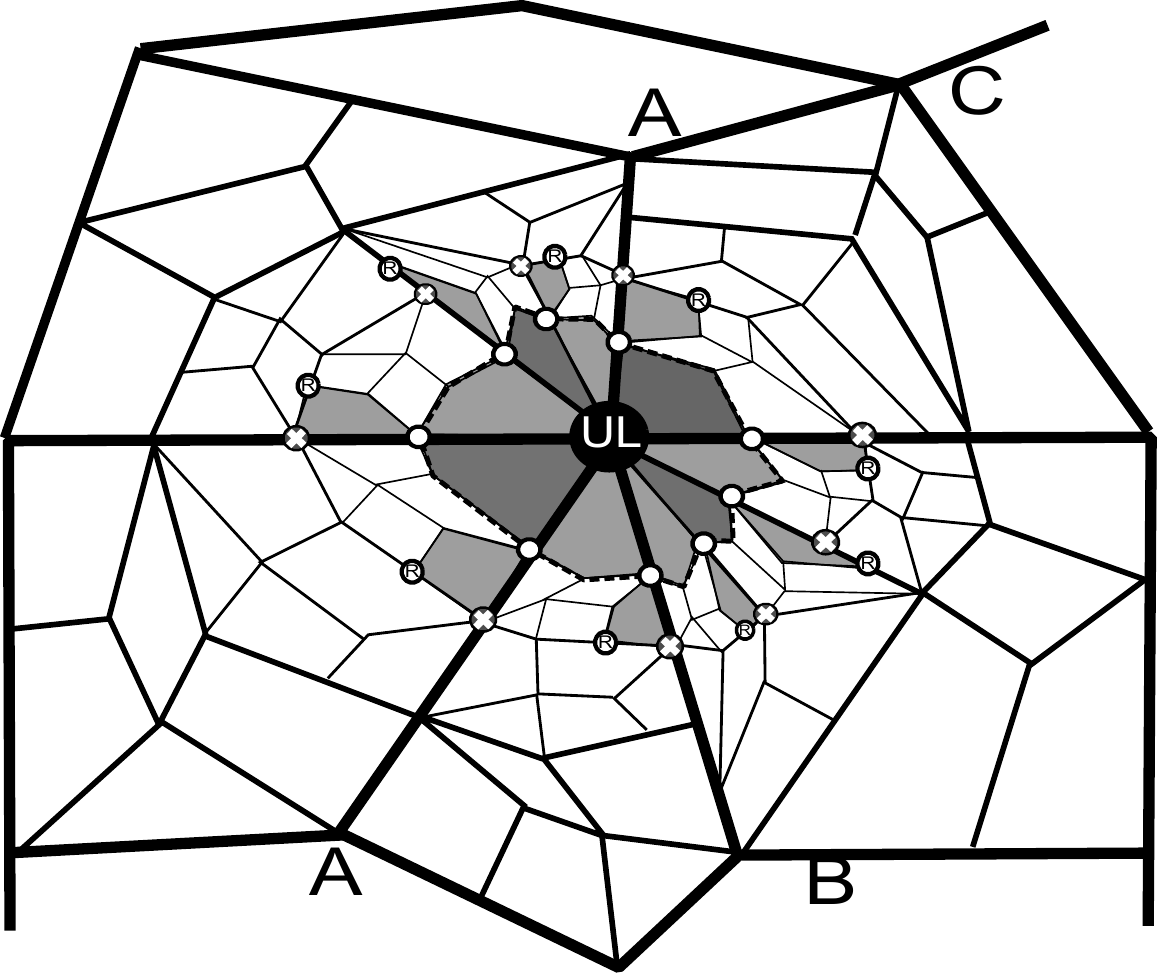}
\caption{$2$-цепь узла $\mathbb{UL}$.}
\label{fig:chainUL2}
\end{figure}

{\bf Цепи с центром в узлах типа $\mathbb{UL}$ и $\mathbb{LU}$.}
Случаи $\mathbb{UL}$ и $\mathbb{LU}$ узлов симметричны, для определенности будем далее разбирать случай $\mathbb{UL}$ узла. Будем считать, что главное ребро имеет тип $t$, то есть с $U$-стороны это $tA$, а с $L$-стороны $tB$.
На рисунках~\ref{fig:chainUL1} и~\ref{fig:chainUL2} белыми точками отмечены $1$-цепи и $2$-цепи узла $\mathbb{UL}$. $0$-цепи с центром в $\mathbb{UL}$ не существует.
Кроме указанных на рисунках цепей первого и второго уровней существует также $3$-цепь, получаемая применением операции разбиения к макроплиткам на рисунке~\ref{fig:chainUL2}.

На получающемся третьем уровне все  макроплитки, левый верхний угол которых попадает в узел $\mathbb{UL}$, будут занимать левое верхнее положение в своих родительских макроплитках. Значит, на следующих уровнях подразбиения, типы ребер-границ макроплиток не изменятся. То есть окружения вершин в цепях будут те же, то есть, окружение цепи четвертого и последующих уровней совпадает с окружением цепи третьего уровня.

Таким образом, существует только три возможных конфигурации окружения цепи с центром в узле типа $\mathbb{UL}$ (для цепей первого, второго и третьего уровней, при этом еще может быть выбрано разное главное ребро $t$). Заметим, что по окружению цепи, мы можем установить, является ли центр цепи узлом типа $\mathbb{UL}$, какого уровня цепь, а также сам параметр $t$. Действительно, $1$-цепь с центром в $\mathbb{UL}$ содержит узлы с окружением $(tA,\mathbf{1A},\mathbf{6A},\mathbf{2A})-(\mathbf{1B},\mathbf{2B},\mathbf{6A},\mathbf{2A})$ и $(\mathbf{8B},\mathbf{3A},\mathbf{3A},\mathbf{7B})-(tB,\mathbf{3B},\mathbf{6A},\mathbf{2A})$, которые вместе ни в какой другой цепи не встречаются. $2$-цепь с центром в $\mathbb{UL}$ содержит одновременно узлы с окружением $(tA,\mathbf{1A},\mathbf{1A},\mathbf{3A})-(\mathbf{1B},\mathbf{2B},\mathbf{1A},\mathbf{3A})$ и $(\mathbf{3A},\mathbf{8A},\mathbf{5A},\mathbf{6B})-(tB,\mathbf{3B},\mathbf{1A},\mathbf{3A})$, которые вместе также ни в какой другой цепи не встречаются.
Для $3$-цепи такой парой будет $(tA,\mathbf{1A},\mathbf{1A},\mathbf{3A})-(\mathbf{1B},\mathbf{2B},\mathbf{1A},\mathbf{3A})$ и $(\mathbf{3A},\mathbf{8A},\mathbf{1A},\mathbf{3A})-(tB,\mathbf{3B},\mathbf{1A},\mathbf{3A})$

Заметим также, что зная тип ребра $t$, мы можем выписать полностью все окружение $1$-цепей, $2$-цепей, $3$-цепей с центром в узле типа $\mathbb{UL}$.

\medskip

\leftskip=0cm

Из сказанного выше следует, что для $\mathbb{UL}$-узлов значение функции $\mathbf{LevelPlus}$ может быть установлено по ее аргументу.

На рисунках~\ref{fig:chainUL1} и~\ref{fig:chainUL2} знаком ``$\otimes$'' отмечены узлы, являющиеся верхними правыми или левыми нижними углами в макроплитках, где середина верхней стороны попадает в узел цепи $X$ . Можно убедиться в том, что для цепей первого, второго и третьего уровней все типы узлов с крестами, а также типы ребер на которых они лежат, мы можем выписать, для каждого заданного $X$.
Таким образом, мы можем вычислить функции $\mathbf{TopRightType}$ и $\mathbf{BottomLeftType}$.

Серым цветом на рисунке выделены плитки, где левый нижний угол (при дальнейшем разбиении) будет попадать в вершину цепи. Можно заметить, что зная эту вершину (ее место в цепи) мы можем установить и окружение правого нижнего угла в соответствующей макроплитке. Кроме того, мы также можем установить тип левого верхнего угла, а значит и всю цепь, содержащую середину верхней стороны в соответствующей макроплитке. (Все это также очевидно проверяется для $3$-цепи.) Все это значит, что мы можем вычислить функции $\mathbf{TopFromCorner}$, $\mathbf{RightCorner}$.

\medskip

\leftskip=0cm

%Рассмотрим все макроплитки, где середина верхней стороны $X$ попадает в цепь. При этом середина правой стороны попадет в один из узлов, помеченных <<$R$>>. Заметим, что зная окружение и информацию $X$, мы можем вычислить окружение этого $R$-узла.
%Таким образом, функция $\mathbf{RightFromTop}$ также вычисляется.

Функции $\mathbf{TopFromRight}$ и $\mathbf{BottomRightTypeFromRight}$ могут быть применимы только к $1$-цепи (в остальные входят только узлы типов $\mathbb{UL}$ и $\mathbb{LU}$). В случае $1$-цепи аргументом могут быть узлы, лежащие на ребрах $\mathbf{u}_1$ и $\mathbf{l}$ (относительно центрального $\mathbb{UL}$-узла). Пусть $X$ такой узел.  Во всех случаях
$\mathbf{BottomRightTypeFromRight}$ будет наш же узел $\mathbb{UL}$. Узлы, на которые указывает $\mathbf{TopFromRight}$ отмечены стрелками на рисунке~\ref{fig:chainUL1}. Таким образом, $\mathbf{TopFromRight}$ принимает значения: $1$-цепи вокруг $\mathbb{B}$ (указатель $1$) для ребра $\mathbf{u}_1$, $1$-цепи вокруг $\mathbb{A}$ (указатель $3$) для ребра $\mathbf{l}$.

\medskip

Таким образом, мы вычислили значения всех функций для аргументов, входящих в $\mathbb{UL}$-цепи.

\medskip

{\bf Цепи с центром в узлах типа $\mathbb{UR}$ и $\mathbb{RU}$.}
Случаи узлов $\mathbb{UR}$ и $\mathbb{RU}$ симметричны, для определенности будем далее разбирать случай $\mathbb{UR}$ узла. Будем считать, что главное ребро имеет тип $t$, то есть с $U$-стороны это $tA$, а с $R$-стороны $tB$.

\begin{figure}[hbtp]
\centering
\includegraphics[width=0.9\textwidth]{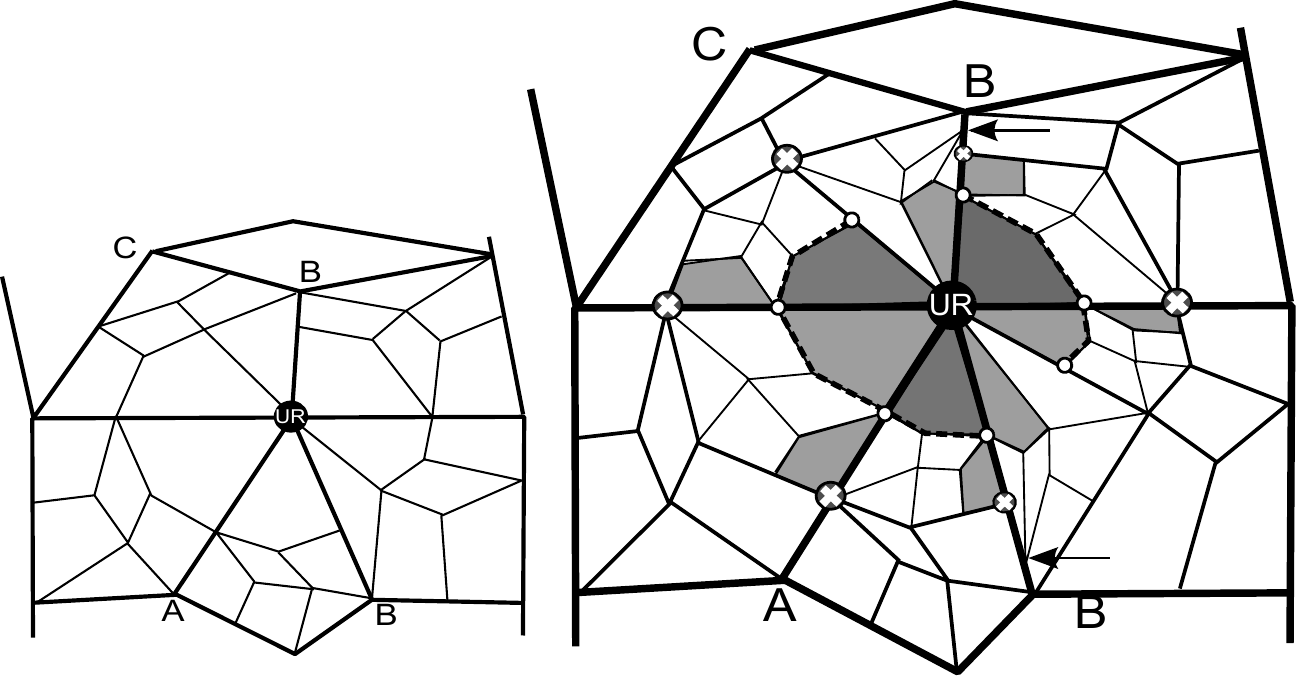}
\caption{$1$-цепь узла $\mathbb{UR}$.}
\label{fig:chainUR1}
\end{figure}

\begin{figure}[hbtp]
\centering
\includegraphics[width=0.8\textwidth]{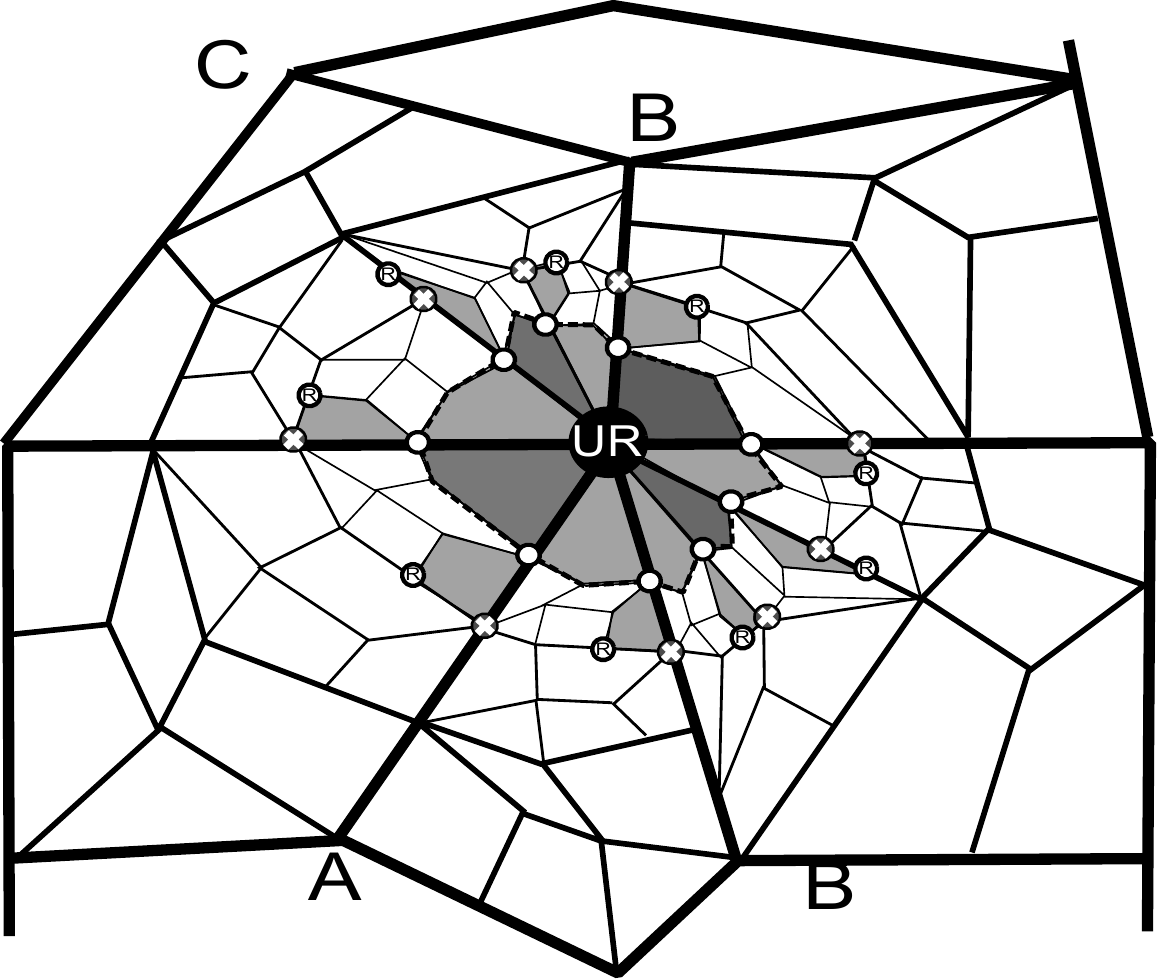}
\caption{$2$-цепь узла $\mathbb{UR}$.}
\label{fig:chainUR2}
\end{figure}

На рисунках~\ref{fig:chainUR1} и~\ref{fig:chainUR2} белыми точками отмечены $1$-цепи и $2$-цепи узла $\mathbb{UR}$. $0$-цепи с центром в $\mathbb{UR}$ не существует.
Кроме указанных на рисунках цепей первого и второго уровней существует также $3$-цепь, получаемая применением операции разбиения к макроплиткам на рисунке~\ref{fig:chainUR2}.

На получающемся третьем уровне все  макроплитки, левый верхний угол которых попадает в узел $\mathbb{UR}$, будут занимать левое верхнее положение в своих родительских макроплитках. Значит, на следующих уровнях подразбиения, типы ребер-границ макроплиток не изменятся. То есть окружения вершин в цепях будут те же, то есть, окружение цепи четвертого и последующих уровней совпадает с окружением цепи третьего уровня.

Таким образом, существует только три возможных конфигурации окружения цепи с центром в узле типа $\mathbb{UR}$ (для цепей первого, второго и третьего уровней, при этом еще может быть выбрано разное главное ребро $t$). Заметим, что по окружению цепи, мы можем установить, является ли центр цепи узлом типа $\mathbb{UR}$, какого уровня цепь, а также сам параметр $t$. Действительно, $1$-цепь с центром в $\mathbb{UR}$ содержит одновременно узлы с окружением $(tA,\mathbf{1A},\mathbf{6A},\mathbf{2A})-(\mathbf{1B},\mathbf{2B},\mathbf{6A},\mathbf{2A})$ и $(tB,\mathbf{6A},\mathbf{6A},\mathbf{2A})-(\mathbf{8B},\mathbf{6B},\mathbf{6B},\mathbf{7B})$, которые вместе ни в какой другой цепи не встречаются. $2$-цепь с центром в $\mathbb{UR}$ содержит одновременно узлы с окружением $(tA,\mathbf{1A},\mathbf{1A},\mathbf{3A})-(\mathbf{1B},\mathbf{2B},\mathbf{1A},\mathbf{3A})$ и $(tB,\mathbf{6A},\mathbf{1A},\mathbf{3A})-(\mathbf{6B},\mathbf{8A},\mathbf{5A},\mathbf{6B})$, которые вместе также ни в какой другой цепи не встречаются.
Для $3$-цепи такой парой будет $(tA,\mathbf{1A},\mathbf{1A},\mathbf{3A})-(\mathbf{1B},\mathbf{2B},\mathbf{1A},\mathbf{3A})$ и $(tB,\mathbf{6A},\mathbf{1A},\mathbf{3A})-(\mathbf{6B},\mathbf{8A},\mathbf{1A},\mathbf{3A})$

Заметим также, что зная тип ребра $t$, мы можем выписать полностью все окружение $1$-цепей, $2$-цепей, $3$-цепей с центром в узле типа $\mathbb{UR}$.

\medskip

\leftskip=0cm

\medskip
Из сказанного выше следует, что для $\mathbb{UR}$-узлов значение функции $\mathbf{LevelPlus}$ может быть установлено по ее аргументу.

\medskip

\leftskip=0cm

На рисунках~\ref{fig:chainUR1} и~\ref{fig:chainUR2} знаком ``$\otimes$'' отмечены узлы, являющиеся верхними правыми или левыми нижними углами в макроплитках, где середина верхней стороны попадает в узел цепи $X$. Можно убедиться в том, что для цепей первого, второго и третьего уровней все типы узлов с крестами, а также типы ребер на которых они лежат, мы можем выписать для каждого заданного $X$.
Таким образом, мы можем вычислить функции $\mathbf{TopRightType}$ и $\mathbf{BottomLeftType}$.

Серым цветом на рисунке выделены плитки, где левый нижний угол (при дальнейшем разбиении) будет попадать в вершину цепи. Можно заметить, что зная эту вершину (ее место в цепи) мы можем установить и окружение правого нижнего угла в соответствующей макроплитке. Кроме того, мы также можем установить тип левого верхнего угла, а значит и всю цепь, содержащую середину верхней стороны в соответствующей макроплитке. (Все это также очевидно проверяется для $3$-цепи.) Все это значит, что мы можем вычислить функции $\mathbf{TopFromCorner}$, $\mathbf{RightCorner}$.

\medskip

Функции $\mathbf{TopFromRight}$ и $\mathbf{BottomRightTypeFromRight}$ могут быть применимы только к $1$-цепи (в остальные входят только узлы типов $\mathbb{UL}$ и $\mathbb{LU}$). В случае $1$-цепи аргументом могут быть узлы, лежащие на ребрах $\mathbf{u}_1$ и $\mathbf{r}$ (относительно $\mathbb{UR}$). Пусть $X$ такой узел.  Во всех случаях
$\mathbf{BottomRightType}$ будет наш же узел $\mathbb{UR}$. Узлы, на которые указывает $\mathbf{TopFromRight}$ отмечены стрелками на рисунке~\ref{fig:chainUR1}. Таким образом, $\mathbf{TopFromRight}$ принимает значения: $1$-цепи вокруг $\mathbb{B}$ (указатель $1$) для $\mathbf{u}_1$ ребра, $1$-цепи вокруг $\mathbb{B}$ (указатель $2$) для $\mathbf{l}$ ребра.

\medskip

Таким образом, мы вычислили значения всех функций для аргументов, входящих в $\mathbb{UR}$-цепи.

\medskip

\medskip

\begin{figure}[hbtp]
\centering
\includegraphics[width=0.9\textwidth]{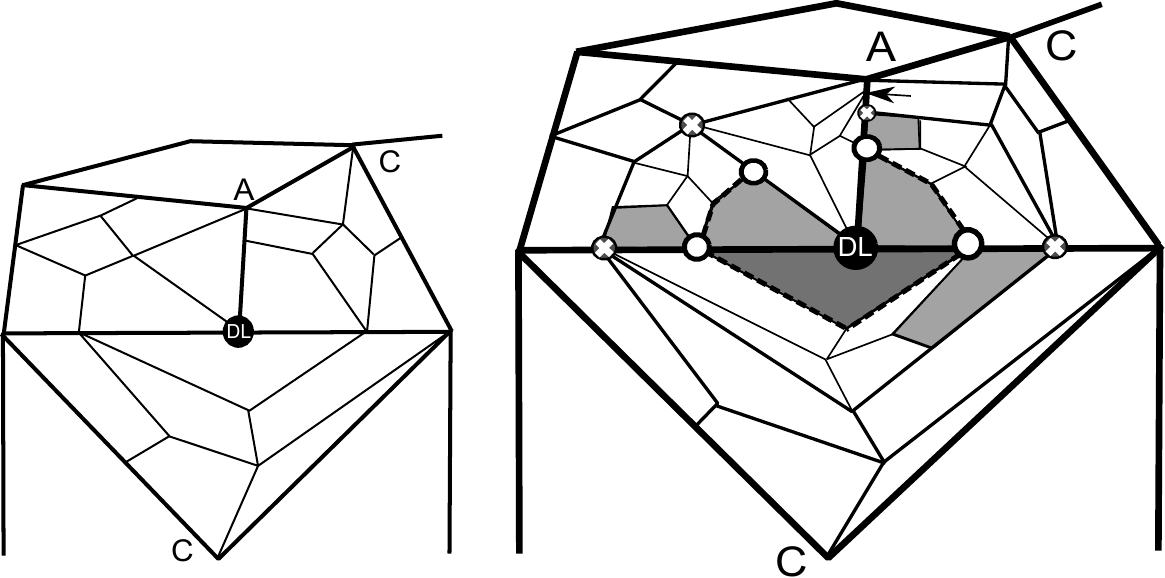}
\caption{$1$-цепь узла $\mathbb{DL}$.}
\label{fig:chainDL1}
\end{figure}

\begin{figure}[hbtp]
\centering
\includegraphics[width=0.8\textwidth]{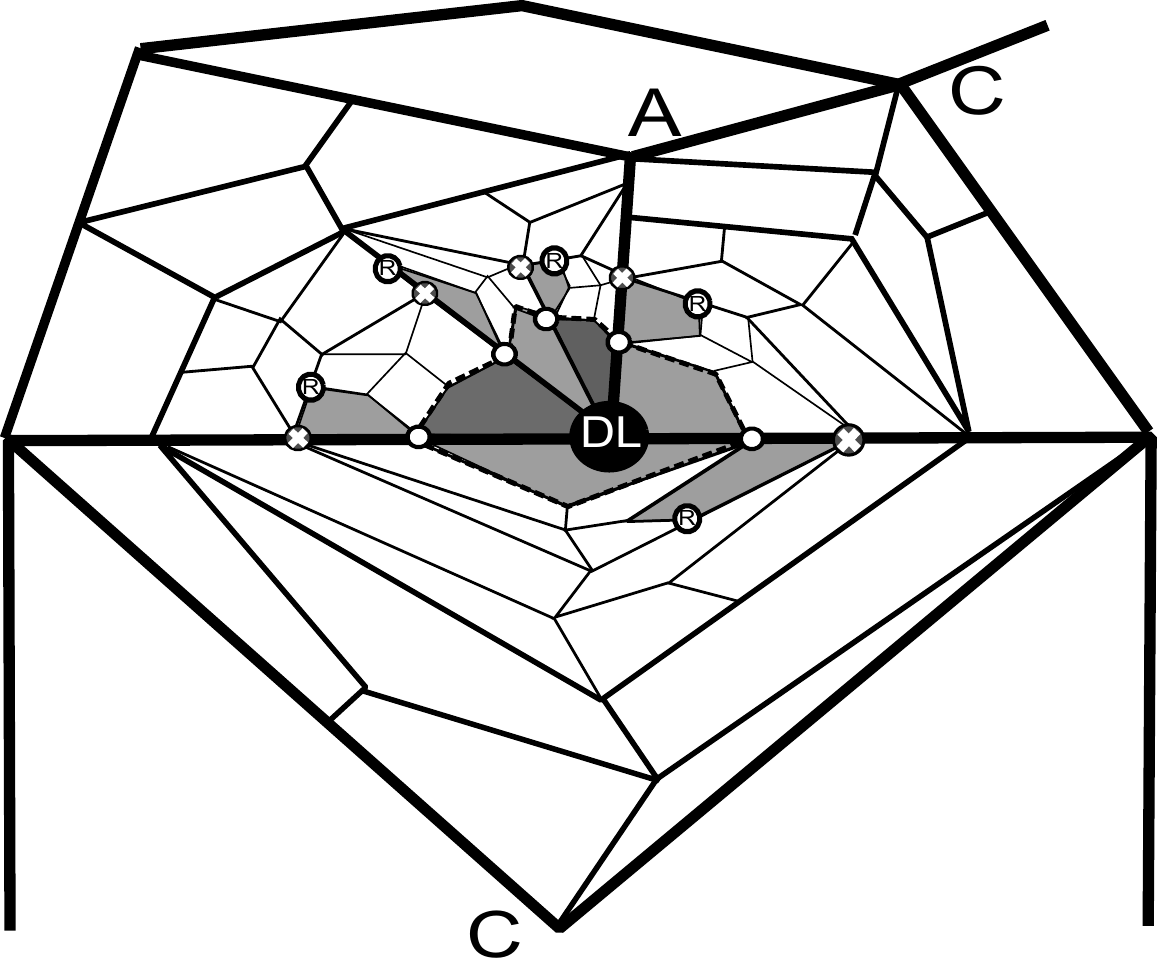}
\caption{$2$-цепь узла $\mathbb{DL}$.}
\label{fig:chainDL2}
\end{figure}

{\bf Цепи с центром в узлах типа $\mathbb{DL}$ и $\mathbb{LD}$.}
Случаи узлов $\mathbb{DL}$ и $\mathbb{LD}$ симметричны, для определенности будем далее разбирать случай $\mathbb{DL}$ узла. Будем считать, что главное ребро имеет тип $t$, то есть с $D$-стороны это $tA$, а с $L$-стороны $tB$.
На рисунках~\ref{fig:chainDL1} и~\ref{fig:chainDL2} белыми точками отмечены $1$-цепи и $2$-цепи узла $\mathbb{DL}$. $0$-цепи с центром в $\mathbb{DL}$ не существует.
Кроме указанных на рисунках цепей первого и второго уровней существует также $3$-цепь, получаемая применением операции разбиения к макроплиткам на рисунке~\ref{fig:chainDL2}.

\medskip

\leftskip=0cm

\medskip

На получающемся третьем уровне все  макроплитки, левый верхний угол которых попадает в узел $\mathbb{DL}$, будут занимать левое верхнее положение в своих родительских макроплитках. Значит, на следующих уровнях подразбиения, типы ребер-границ макроплиток не изменятся. То есть окружения вершин в цепях будут те же, то есть, окружение цепи четвертого и последующих уровней совпадает с окружением цепи третьего уровня.

Таким образом, существует только три возможных конфигурации окружения цепи с центром в узле типа $\mathbb{DL}$ (для цепей первого, второго и третьего уровней, при этом еще может быть выбрано разное главное ребро $t$). Заметим, что по окружению цепи, мы можем установить, является ли центр цепи узлом типа $\mathbb{DL}$, какого уровня цепь, а также сам параметр $t$. Действительно, $1$-цепь с центром в $\mathbb{DL}$ содержит одновременно узлы с окружением $(tA,tA,\mathbf{6A},\mathbf{2A})-(tB,\mathbf{3B},\mathbf{6A},\mathbf{2A})$ и $(tA,tA,\mathbf{6A},\mathbf{2A})-(\mathbf{7A},tB,\mathbf{3B},\mathbf{4A})$, которые вместе ни в какой другой цепи не встречаются. $2$-цепь с центром в $\mathbb{DL}$ содержит одновременно узлы с окружением $(tA,tA,\mathbf{1A},\mathbf{3A})-(tB,\mathbf{3B},\mathbf{1A},\mathbf{3A})$ и $(\mathbf{3A},\mathbf{8A},\mathbf{5A},\mathbf{6B})-(\mathbf{8B},\mathbf{7B},\mathbf{7B},\mathbf{7B})$, которые вместе также ни в какой другой цепи не встречаются.
Для $3$-цепи такой парой будет $(tA,tA,\mathbf{1A},\mathbf{3A})-(tB,\mathbf{3B},\mathbf{1A},\mathbf{3A})$ и $(\mathbf{3A},\mathbf{8A},\mathbf{1A},\mathbf{3A})-(\mathbf{8B},\mathbf{7B},\mathbf{1A},\mathbf{3A})$

Заметим также, что зная тип ребра $t$, мы можем выписать полностью все окружение $1$-цепей, $2$-цепи, $3$-цепей с центром в узле типа $\mathbb{DL}$.

\medskip

\leftskip=0cm

Из сказанного выше следует, что для $\mathbb{DL}$-узлов значение функции $\mathbf{LevelPlus}$ может быть установлено по ее аргументу.
На рисунках~\ref{fig:chainDL1} и~\ref{fig:chainDL2} знаком ``$\otimes$'' отмечены узлы, являющиеся верхними правыми или левыми нижними углами в макроплитках, где середина верхней стороны попадает в узел цепи $X$ . Можно убедиться в том, что для цепей первого, второго и третьего уровней все типы узлов с крестами, а также типы ребер на которых они лежат, мы можем выписать для каждого заданного $X$. Таким образом, мы можем вычислить функции $\mathbf{TopRightType}$ и $\mathbf{BottomLeftType}$.

Серым цветом на рисунке выделены плитки, где левый нижний угол (при дальнейшем разбиении) будет попадать в вершину цепи. Можно заметить, что зная эту вершину (ее место в цепи) мы можем установить и окружение правого нижнего угла в соответствующей макроплитке. Кроме того, мы также можем установить тип левого верхнего угла, а значит и всю цепь, содержащую середину верхней стороны в соответствующей макроплитке. (Все это также очевидно проверяется для $3$-цепи.) Все это значит, что мы можем вычислить функции $\mathbf{TopFromCorner}$, $\mathbf{RightCorner}$.

Функции $\mathbf{TopFromRight}$ и $\mathbf{BottomRightTypeFromRight}$ могут быть применимы только к $1$-цепи (в остальные входят только узлы типов $\mathbb{UL}$ и $\mathbb{LU}$). В случае $1$-цепи аргументом могут быть только узел, лежащий на $\mathbf{l}$ ребре (относительно $\mathbb{DL}$ узла). Пусть $X$ этот узел.  Тогда
$\mathbf{BottomRightTypeFromRight}$ будет наш же узел $\mathbb{DL}$. Узел, на которые указывает $\mathbf{TopFromRight}$ отмечен стрелкой на рисунке~\ref{fig:chainDL1}. Таким образом, $\mathbf{TopFromRight}$ это $1$-цепь вокруг $\mathbb{A}$ (указатель $3$).

\medskip

\begin{figure}[hbtp]
\centering
\includegraphics[width=0.9\textwidth]{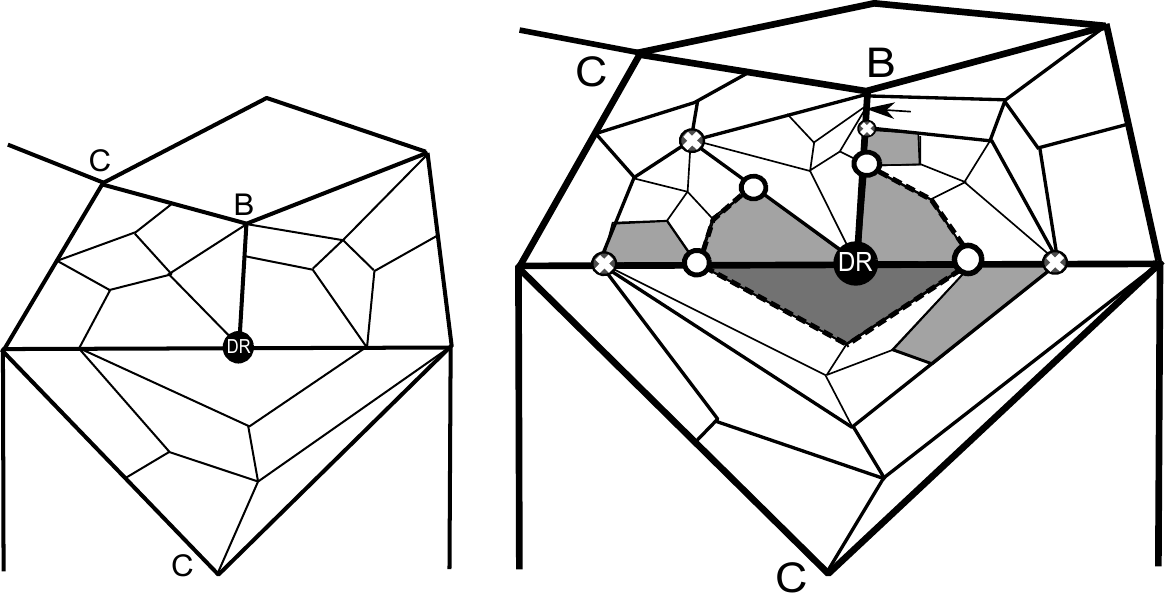}
\caption{$1$-цепь узла $\mathbb{DR}$.}
\label{fig:chainDR1}
\end{figure}

\begin{figure}[hbtp]
\centering
\includegraphics[width=0.8\textwidth]{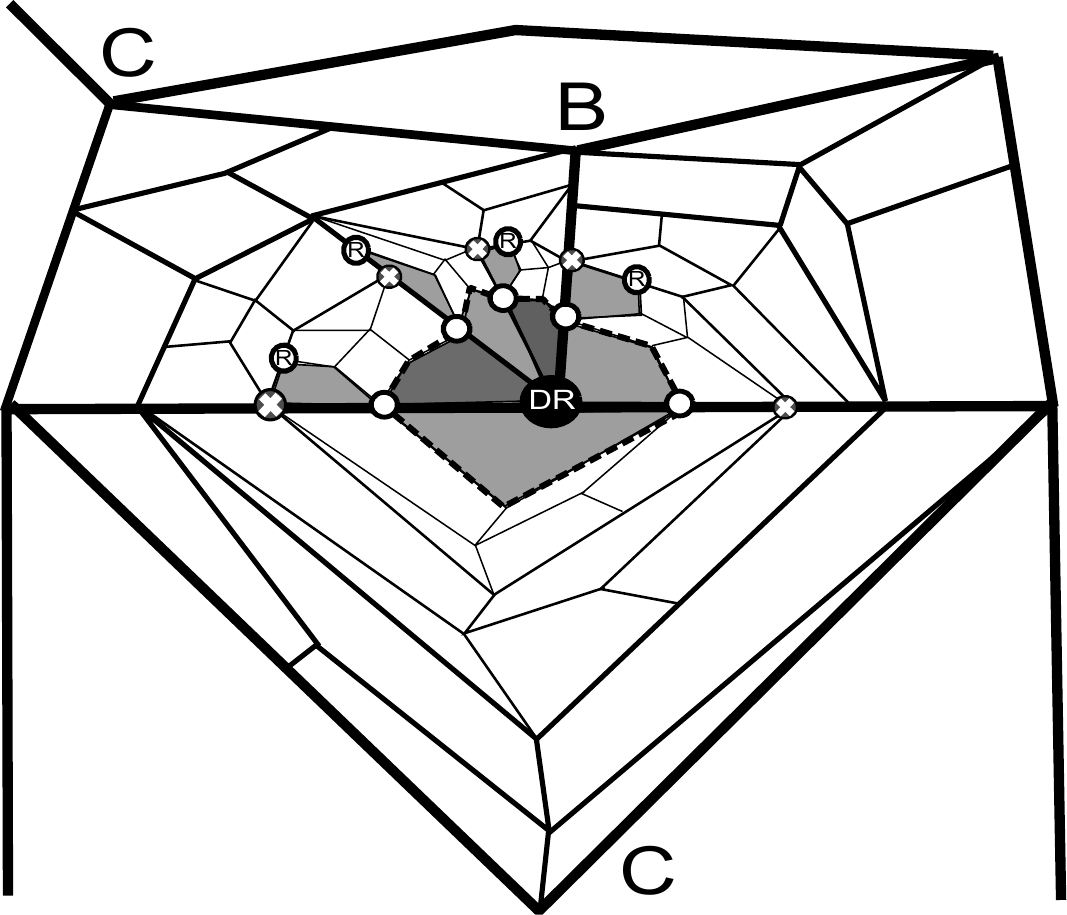}
\caption{$2$-цепь узла $\mathbb{DR}$.}
\label{fig:chainDR2}
\end{figure}

{\bf Цепи с центром в узлах типа $\mathbb{DR}$ и $\mathbb{RD}$.}
На рисунках~\ref{fig:chainDR1} и~\ref{fig:chainDR2} белыми точками отмечены $1$-цепи и $2$-цепи узла $\mathbb{DR}$. $0$-цепи с центром в $\mathbb{DR}$ не существует.

Кроме указанных на рисунках цепей первого и второго уровней существует также $3$-цепь, получаемая применением операции разбиения к макроплиткам на рисунке~\ref{fig:chainDL2}.
Полностью аналогично цепи $\mathbb{DL}$, можно показать, что существует только три возможных конфигурации окружения цепи с центром в узле типа $\mathbb{DR}$.
Кроме того, зная тип ребра $t$, мы можем выписать полностью все окружение $1$-цепи, $2$-цепи, $3$-цепи с центром в узле типа $\mathbb{DR}$.

\medskip

\leftskip=0cm

\medskip

Таким образом, что для $\mathbb{DR}$-узлов значение функции $\mathbf{LevelPlus}$ может быть установлено по ее аргументу.

На рисунках~\ref{fig:chainDR1} и~\ref{fig:chainDR2} знаком ``$\otimes$'' отмечены узлы, являющиеся верхними правыми или левыми нижними углами в макроплитках, где середина верхней стороны попадает в узел цепи $X$.
Вычисление функций $\mathbf{TopRightType}$, $\mathbf{BottomLeftType}$, $\mathbf{TopFromCorner}$, $\mathbf{RightCorner}$ полностью аналогично случаю $\mathbb{DL}$-цепи.

\medskip

\leftskip=0cm

\medskip

Функции $\mathbf{TopFromRight}$ и $\mathbf{BottomRightTypeFromRight}$ могут быть применимы только к $1$-цепи (в остальные входят только узлы типов $\mathbb{UL}$ и $\mathbb{LU}$). В случае $1$-цепи аргументом могут быть только узел, лежащий на $\mathbf{r}$ ребре (относительно $\mathbb{DR}$ узла). Пусть $X$ этот узел.  Тогда
$\mathbf{BottomRightType}$ будет наш же узел $\mathbb{DR}$. Узел, на которые указывает $\mathbf{TopFromRight}$ отмечен стрелкой на рисунке~\ref{fig:chainDR1}. Таким образом, $\mathbf{TopFromRight}$ это $1$-цепь вокруг $\mathbb{B}$ (указатель $2$).

\medskip

Нам остается разобрать цепи с центрами в краевых и угловых вершинах. Важное замечание: при разборе цепи с центром в краевой вершине часть узлов в попадает на край макроплитки. Под их окружением, в данном случае, мы понимаем базовое окружение в рамках нашей макроплитки (эта макроплитка может быть подклееной и тогда вершина на ее краю одновременно имеет какое-то другое окружение в основной своей области, вот это окружение мы не рассматриваем).

С учетом этого замечания, для рассмотрения цепей с центрами в вершинах типов $\mathbb{U}$, $\mathbb{L}$, $\mathbb{R}$, $\mathbb{D}$ достаточно применить те же рассуждения как и для цепей $\mathbb{UL}$, $\mathbb{LD}$, $\mathbb{RD}$, $\mathbb{DL}$, но рассматривать только половину картинки для каждой цепи. В силу полной аналогичности этого рассмотрения, мы не будем тут приводить этот разбор.

Цепи с центрами в угловых вершинах рассмотрим отдельно.

\medskip

\begin{figure}[hbtp]
\centering
\includegraphics[width=0.9\textwidth]{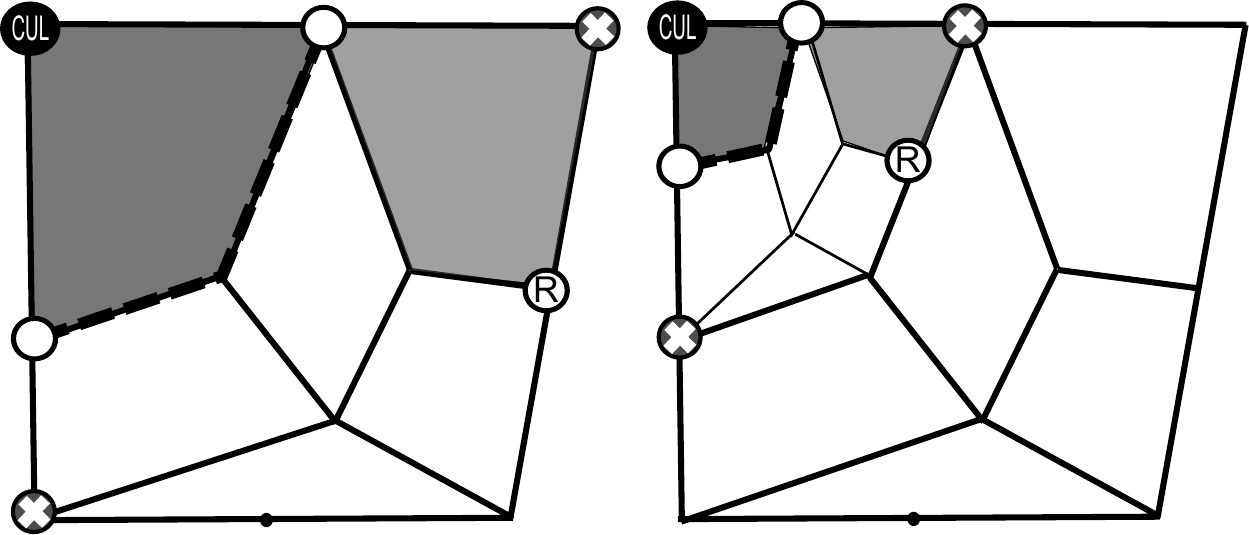}
\caption{$0$-цепь и $1$-цепь узла $\mathbb{CUL}$.}
\label{chainCUL01}
\end{figure}

\begin{figure}[hbtp]
\centering
\includegraphics[width=0.9\textwidth]{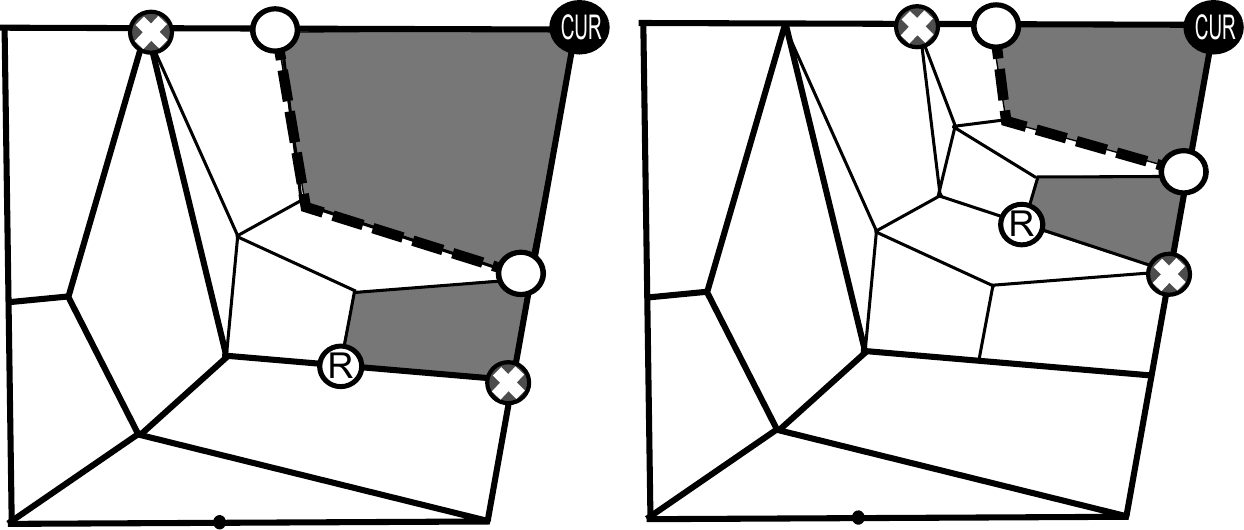}
\caption{$1$-цепь и $2$-цепь узла $\mathbb{CUR}$.}
\label{chainCUR12}
\end{figure}

\begin{figure}[hbtp]
\centering
\includegraphics[width=0.9\textwidth]{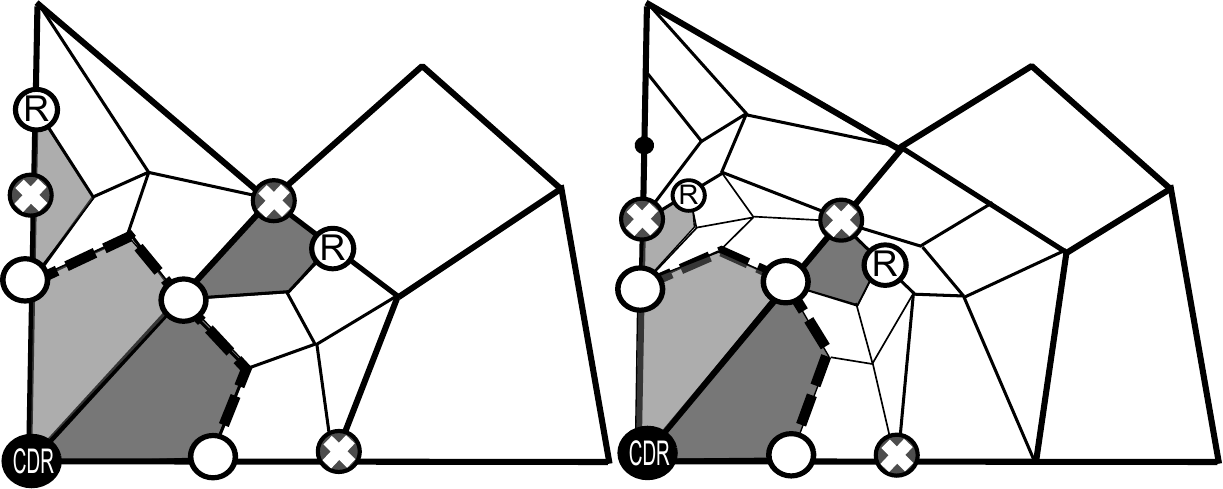}
\caption{$1$-цепь и $2$-цепь узла $\mathbb{CDR}$.}
\label{chainCDR12}
\end{figure}

\begin{figure}[hbtp]
\centering
\includegraphics[width=0.9\textwidth]{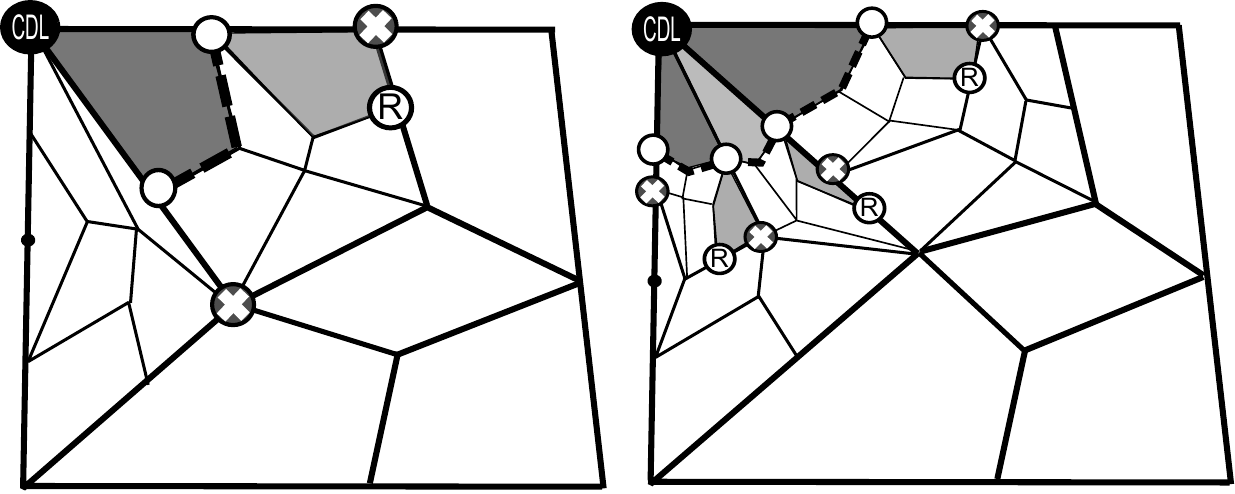}
\caption{$1$-цепь и $2$-цепь узла $\mathbb{CDL}$.}
\label{chainCDL12}
\end{figure}

{\bf Цепи с центром в узлах типа $\mathbb{CUL}$, $\mathbb{CUR}$, $\mathbb{CDR}$, $\mathbb{CDL}$.}
На рисунках~\ref{chainCUL01}, ~\ref{chainCUR12},~\ref{chainCDR12},~\ref{chainCDL12}    белыми точками отмечены цепи c вершинами в узлах типов  $\mathbb{CUL}$, $\mathbb{CUR}$, $\mathbb{CDR}$, $\mathbb{CDL}$.

Аналогично предыдущим цепям, можно показать, что существует только две возможных конфигурации окружения цепи в каждом из этих четырех случаев.

\medskip

\leftskip=0cm

\medskip

\leftskip=0cm

\medskip

\leftskip=0cm

\medskip

\leftskip=0cm

\medskip

Мы используем те же обозначения, что и в предыдущих случаях, крестами отмечены вершины, являющиеся правыми верхними или левыми нижними.

Вычисление функций $\mathbf{LevelPlus}$, $\mathbf{TopRightType}$, $\mathbf{BottomLeftType}$, $\mathbf{TopFromCorner}$, $\mathbf{RightCorner}$ очевидно, а функции
 $\mathbf{TopFromRight}$ и $\mathbf{BottomRightTypeFromRight}$ не могут быть применимы к вершинам данным цепям.

\medskip

{\bf Узлы типа $\mathbb{DR}$ и $\mathbb{RD}$.}
Для некоторых функций нужно также рассмотреть ситуации, когда аргументом является узел $\mathbb{DR}$ (или симметричный случай $\mathbb{RD}$).

\begin{figure}[hbtp]
\centering
\includegraphics[width=0.5\textwidth]{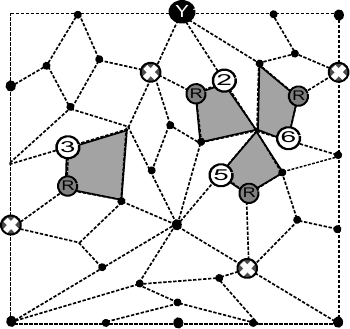}
\caption{Случаи расположения $\mathbb{DR}$-узла.}
\label{RDnode}
\end{figure}

Окружение $\mathbb{DR}$-узла может быть четырех видов, по числу внутренних ребер, на которых он может располагаться. На рисунке~\ref{RDnode} отмечены эти четыре возможные ситуации. Во всех случаях начальником нашего $\mathbb{DR}$-узла является вершина $Y$ в середине верхней стороны. На рисунке отмечены серым цветом макроплитки, левый нижний угол которых попадает в $\mathbb{DR}$-узел. Таким образом, значение функции $\mathbf{RightCorner}$ попадает в узлы отмеченные как ``R''. Это $\mathbb{B}$-узел, и его окружение во всех случаях можно вычислить, учитывая, что мы знаем окружение узла в середине верхней стороны.

\medskip

\leftskip=0cm

\medskip

Значение $\mathbf{TopFromCorner}$  указывает на $1$ цепь вокруг $\mathbb{A}$-узла с указателем $\mathbf{ld}$ для одного из случаев расположения, и $1$-цепь вокруг $\mathbb{B}$-узла для трех остальных случаев, указатели ясны из рисунка. Значение $\mathbf{BottomRightTypeFromRight}$ также очевидно во всех четырех случаях.

Теперь рассмотрим функцию $\mathbf{TopFromRight}$. Вершины, являющиеся ее значениями, отмечены знаком ``$\otimes$'' на рисунке. Поскольку мы знаем окружение вершины $Y$ в середине верхней стороны, мы можем вычислить и все требуемые окружения.

\medskip

\subsection{Указатели} \label{pointers}
Пусть узел $X$ принадлежит некоторой цепи, а узел $Y$ является ее центром, причем известно окружение $X$ и тип внутреннего ребра, на котором лежит $X$ (один из восьми типов внутренних ребер), а также известно окружение макроплитки $T$, которой принадлежит это ребро. Покажем, как восстановить {\it указатель} $X$, то есть тип ребра входа-выхода для $Y$.

Почти для всех ребер {\it указатель}, являющийся ребром выхода из $Y$ устанавливается сразу по типу ребра:
Для  ребра типа $1$, {\it указатель} будет очевидно $\mathbf{u}_2$, для ребра типа $2$ будет $\mathbf{u}_1$, для ребра типа $3$ будет $\mathbf{l}$, для ребра типа $4$ -- в зависимости от типа $Y$: для $\mathbb{A}$ указатель $2$, для $\mathbb{C}$ будет $1$. Для $5$ ребра опять в зависимости от типа $Y$ -- для $\mathbb{B}$ -- указатель $3$, для $\mathbb{C}$ будет $2$. Для $6$ ребра указатель $\mathbf{r}$.

Пусть ребро имеет тип $7$. Зная окружение макроплитки $T$, мы можем установить ее положение в родительской макроплитке. Далее вычисляем значение указателя: для левого-верхнего положения $T$ это $\mathbf{l}_2$, для левого нижнего $\mathbf{ld}$, для среднего $\mathbf{mid}$, для правого верхнего $\mathbf{u}_3$, для правого нижнего -- $\mathbf{rd}$, для нижнего -- $\mathbf{d}$.

Пусть ребро имеет тип $8$. Зная окружение макроплитки $T$, мы можем установить ее положение в родительской макроплитке. Далее вычисляем значение указателя: для левого-верхнего $\mathbf{lu}$ , для левого нижнего $\mathbf{ld}$, для среднего $\mathbf{mid}$, для правого верхнего $\mathbf{ru}$, для правого нижнего -- $\mathbf{rd}$. В случае нижнего положения $T$,  по типу верхней (и правой) стороны $T$ можно узнать, какое положение уже родительская макроплитка для $T$  занимает в своей родительской макроплитке  и в зависимости от этого установить значение указателя:  для левого-верхнего положения  -- $\mathbf{l}_3$, для левого нижнего $\mathbf{ld}_2$, для среднего $\mathbf{mid}_2$, для правого верхнего $\mathbf{u}_4$, для правого нижнего -- $\mathbf{r}_3$, для нижнего -- $\mathbf{d}_2$.

\medskip

\subsection{Дополнительные функции} \label{addfunc}

\begin{figure}[hbtp]
\centering
\includegraphics[width=0.4\textwidth]{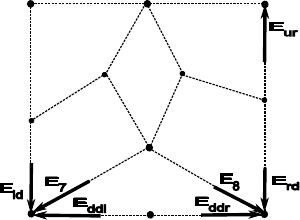}
\caption{Функции $E_7$, $E_8$, $E_{\mathbf{ld}}$, $E_{\mathbf{rd}}$, $E_{\mathbf{ur}}$, $E_{\mathbf{ddl}}$, $E_{\mathbf{ddr}}$ обозначают входящие или выходящие ребра в углы.}
\label{fig:e7e8}
\end{figure}

Символами $E_7$ и $E_8$ будем обозначать код соответствующего входящего и выходящего ребра. Очевидно, что зная по окружению узла в левом нижнем углу можно установить $E_7$, а по окружению узла в правом нижнем углу можно установить $E_8$ -- (рисунок~\ref{fig:e7e8}).

\medskip

Аналогично, обозначим символами  $E_{\mathbf{ld}}$,  $E_{\mathbf{rd}}$, $E_{\mathbf{ur}}$ коды входящих или выходящих ребер в нижние углы по левой или правой сторонам и в верхний угол по правой стороне.  Также символами $E_{\mathbf{ddl}}$,  $E_{\mathbf{ddr}}$ обозначим коды ребер в нижние углы по нижней стороне, слева и справа. Очевидно, что зная окружение соответствующего узла, можно установить букву, обозначающую нужное ребро (рисунок~\ref{fig:e7e8}).

\medskip

Пусть узел $X$ принадлежит некоторой цепи. Поскольку нам известны окружения всех узлов в цепи, а также типы указателей и их порядок следования, мы можем рассмотреть следующий по часовой стрелке узел в цепи, после $X$. Будем обозначать его как $\mathbf{Next}(X)$. Аналогично можно определить следующий после $X$ узел против часовой стрелки: $\mathbf{Prev}(X)$. Например, если в некоторой макроплитке известно окружение $X$ узла в середине верхней стороны, то окружение узла в середине левой стороны будет $\mathbf{Next}(X)$.

\medskip

Также отдельно определим функцию $\mathbf{BottomRightType}$, аргументами которой являются окружение и информация узла $X$, являющимся серединой верхней стороны в некоторой макроплитке $T$,
а значением --  тип узла в правом нижнем углу $T$ и ребро, на котором он расположен, если это боковой узел.

Значение этой функции вычисляется следующим образом. Мы можем установить положение $T$ в родительской макроплитке $T'$. Далее, для левого верхнего и левого нижнего положений тип правого нижнего угла будет $\mathbb{A}$, для среднего, правого верхнего и правого нижнего положений тип правого нижнего угла будет $\mathbb{B}$.

В случае нижнего положения $T$, правый нижний угол попадает в левый нижний угол $T'$. По типу правого и верхнего ребер $T$ можно установить, какое положение занимает уже $T'$ в своей родительской макроплитке $T''$.

Для левого верхнего положения правый нижний угол $T$ попадет в середину левой стороны $T''$ и его тип будет вычисляться функцией $\mathbf{Next.FBoss}(X)$.

Для левого нижнего, среднего и нижнего положений правый нижний угол $T$ попадет в $\mathbb{C}$-узел $T''$.

Для правого верхнего положения правый нижний угол $T$ попадет в середину верхней стороны $T''$ и его тип будет вычисляться как $\mathbf{FBoss}(X)$.

Для правого нижнего положения правый нижний угол $T$ попадет в $\mathbb{B}$ узел $T''$.

\medskip

\begin{figure}[hbtp]
\centering
\includegraphics[width=0.7\textwidth]{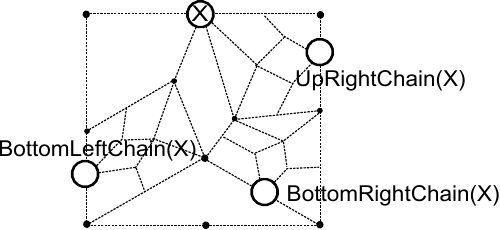}
\caption{Узлы, на которые указывают функции $\mathbf{BottomLeftChain}$, $\mathbf{BottomRightChain}$, $\mathbf{UpRightChain}$.}
\label{fig:addfunctions}
\end{figure}

Если узел $X$ в середине верхней стороны макроплитки известен, то функции $\mathbf{BottomRightType}$, $\mathbf{BottomLeftType}$, $\mathbf{TopRightType}$ дают нам типы узлов в трех соответствующих углах. Учитывая, что мы можем установить ребра входа в углы (символы $E_{\mathbf{ld}}$, $E_{\mathbf{rd}}$, $E_{\mathbf{ur}}$), мы можем установить окружения отмеченных на рисунке~\ref{fig:addfunctions} трех узлов. Обозначим соответствующие функции как $\mathbf{BottomLeftChain}$, $\mathbf{BottomRightChain}$, $\mathbf{UpRightChain}$.

\medskip

\medskip

Также определим функцию $\mathbf{RightFromB}(X)$, которая вычисляет окружение узла в середине правой стороны в макроплитке, по известному $\mathbb{B}$-узлу $X$. Значение этой функции вычисляется следующим образом.

Мы можем установить положение $T$ в родительской макроплитке $T'$.
Для левого верхнего положения, значением будет $\mathbb{UR}$-узел c окружением, соответствующим A0-цепи, с указателем $1$. Типы левой и верхней стороны в $T'$ те же что и в $T$, то есть мы их знаем по окружению $X$.
Для среднего, левого нижнего, правого верхнего, правого нижнего положений значением будет $\mathbb{RD}$-узел с окружением соответственно $3$, $2$, $6$, $5$ внутреннем ребрам.
Для нижнего положения значением будет $1$-цепь вокруг второго начальника $X$ (мы знаем его тип), с указателем $E_{\mathbf{ddl}}$ (вход по нижнему ребру в левый нижний угол).

\medskip

\medskip

\section{Оценка количества букв} \label{count}

Алфавит состоит из вершинных букв и реберных букв. Реберные буквы кодируют всевозможные ребра входа и выхода, включая возможные выходы в подклееные области. Плоских ребер, выходящих из вершины -- не более $10$ (больше всего в $\mathbb{C}$-вершине), всех подклееных ребер  не более $11$ ($\widehat{\mathbf{r}}$, $\widehat{\mathbf{u}_i}$; где $i=1,\dots, 4$, $\widehat{\mathbf{l}}$, $\widehat{\mathbf{l}_2}$, $\widehat{\mathbf{l}_3}$, $\widehat{\mathbf{d}}$, $\widehat{\mathbf{d}_2}$, $\widehat{\mathbf{d}_3}$). То есть реберных букв всего $21$.
Каждая вершинная буква кодирует одну комбинацию значений параметров вершины, то есть тип, уровень, расширенное окружение, информация и флаг подклейки.

\medskip

Посчитаем количество вершинных букв, используемых при кодировании. Четверок типов ребер всего существует $210$, подробно этот подсчет произведен в приложении~\ref{Appendix}. То есть, существует $210$ вариантов базовых окружений для типов $\mathbb{A}$, $\mathbb{B}$ и $\mathbb{C}$.
В цепях с центром в $\mathbb{UL}$-вершине существует $25$ возможных вершин ($7$ в $1$-цепи, и по $9$ в $2$-цепи и $3$-цепи). Каждая может быть одного из трех уровней, то есть всего $75$ возможных сочетаний базового окружения и уровня.

В таблице~\ref{chains} ниже представлено, сколько различных базовых окружений может быть у вершины заданного типа, входящей в цепь указанного уровня.

  \begin{table}[hbtp]
  \caption{Количество вершин с различными окружениями, сгруппированных по цепям}
\centering
 \begin{tabular}{|c|c|c|c|c|c|}   \hline
  тип центра цепи & $0$-цепь & $1$-цепь & $2$-цепь & $3$-цепь & всего \cr \hline
$\mathbb{UL}$/$\mathbb{LU}$ & -- & $7$ & $9$ & $9$ & $50$ \cr \hline
$\mathbb{UR}$/$\mathbb{RU}$ & -- & $7$ & $9$ & $9$ & $50$ \cr \hline
$\mathbb{DL}$/$\mathbb{LD}$ & -- & $4$ & $5$ & $5$ & $28$ \cr \hline
$\mathbb{DR}$/$\mathbb{RD}$ & -- & $4$ & $5$ & $5$ & $28$ \cr \hline
$\mathbb{A}$ & $420$ & $5$ & $5$ & -- & $430$ \cr \hline
$\mathbb{B}$ & -- & $6$ & $6$ & -- & $12$ \cr \hline
$\mathbb{C}$ & -- & $7$ & $10$ & $10$ & $27$ \cr \hline
$\mathbb{D}$ & -- & $2$ & $2$ & -- & $4$ \cr \hline
$\mathbb{L}$ & -- & $4$ & $5$ & $5$ & $14$ \cr \hline
$\mathbb{R}$ & -- & $4$ & $5$ & $5$ & $14$ \cr \hline
$\mathbb{U}$ & -- & $5$ & $6$ & $6$ & $17$ \cr \hline
$\mathbb{CUL}$ & -- & $2$ & $2$ & -- & $4$ \cr \hline
$\mathbb{CUR}$ & -- & $2$ & $2$ & -- & $4$ \cr \hline
$\mathbb{CDL}$ & -- & $2$ & $4$ & $4$ & $10$ \cr \hline
$\mathbb{CDR}$ & -- & $3$ & $3$ & -- & $6$ \cr \hline
  \end{tabular}
\label{chains}
\end{table}

Таким образом, всего возможно не более $698$ окружений для вершин, входящих в некоторую цепь, причем для каждой из этих вершин также известен ее тип. Каждая из них может быть одного из трех уровней, итого $2094$ сочетаний {\it тип-уровень-окружение}.

\medskip

Теперь посчитаем количество вершин, не входящих в цепи. В таблице~\ref{nonchains} ниже подсчитаны все возможные окружения вершин, не входящих в цепи.

  \begin{table}[hbtp]
  \caption{Количество вершин с разными окружениями, сгруппированных по типам}
\centering
 \begin{tabular}{|c|c|c|c|}   \hline
  тип вершины & число окружений & число уровней & всего  \cr \hline
$\mathbb{DR}$/$\mathbb{RD}$ & $4$ & $3$ & $12$ \cr \hline
$\mathbb{A}$ & $210$ & -- & $210$  \cr \hline
$\mathbb{B}$ & $210$ & -- & $210$  \cr \hline
$\mathbb{C}$ & $210$ & -- & $210$   \cr \hline
$\mathbb{D}$ & $1$ & $3$ & $3$  \cr \hline
$\mathbb{R}$ & $2$ & $3$ & $6$ \cr \hline
$\mathbb{CUL}$ & $1$ & -- & $1$  \cr \hline
$\mathbb{CUR}$ & $1$ & -- & $1$ \cr \hline
$\mathbb{CDL}$ & $1$ & -- & $1$  \cr \hline
$\mathbb{CDR}$ & $1$ & -- & $1$  \cr \hline
  \end{tabular}
\label{nonchains}
\end{table}

Заметим, что у $\mathbb{D}$-вершины может быть только окружение $(\mathbf{left},\mathbf{top},\mathbf{right},\mathbf{bottom})$, так как нижней стороной макроплитка примыкает к краю только являясь при этом подклееной плиткой.
Аналогично, у $\mathbb{R}$-вершины может быть два окружения $(\mathbf{left},\mathbf{top},\mathbf{right},\mathbf{bottom})$ и $(\mathbf{8B},\mathbf{bottom},\mathbf{bottom},\mathbf{7B})$, так как она макроплитка может примыкать правой стороной только к правому или нижнему краям и только в двух этих случаях.

То есть всего существует не более $655$ различных окружений для вершин не из цепей. Итого всех сочетаний тип-уровень-базовое окружение- не более $698\times 3 + 655 = 2749$.

\medskip

Рассмотрим теперь тип в подклееной области. Он может быть один из следующих: $\mathbb{CUL}$, $\mathbb{CUR}$, $\mathbb{CDL}$, $\mathbb{U}$, $\mathbb{L}$.  У первых трех по одному возможному подклееному окружению. У $\mathbb{U}$-вершины всего $6$ вариантов (по одному узлу в $\mathbb{U}1$, $\mathbb{U}2$, $\mathbb{U}3$, $\mathbb{L}1$, $\mathbb{L}2$, $\mathbb{L}3$, цепях). Аналогично у $\mathbb{L}$-вершины тоже $6$ вариантов подклееных окружений. С учетом трех возможных уровней, получается $3+6\times 3+ 6\times 3=39$ сочетаний.

Из определения операции подклейки следует, что к вершине типов $\mathbb{A}$, $\mathbb{B}$, $\mathbb{C}$ и к угловым вершинм подклееные макроплитки не могут примыкать левой или верхней стороной (только левым-верхним углом). Значит, вершина с нетривиальным подклееным окружением, не являющаяся ядром, в базовой плоскости либо имеет тип $\mathbb{DR}$/$\mathbb{RD}$, $\mathbb{D}$, $\mathbb{R}$, либо входит в некоторую цепь.
То есть для таких вершин базовых окружений всего $698\times 3 + 12 +3 +6 =2115$. То есть расширенных окружений (сочетаний базового и подклееного) имеется не более $2115\times 39 =82485$.

\medskip

{\bf Число значений параметра ``информация''}.
Для простоты будем говорить не ``число значений параметра информация'', а просто ``число информаций''.
Оценим это число. Пусть сначала все начальники не лежат на левом или верхнем краю подклееных плиток (то есть у них только базовое окружение).

Тогда первый начальник обязательно является вершиной, входящей в некоторую цепь, то есть $2094$ варианта. Второй начальник (левый нижний угол) либо попадает в вершину из цепи, либо $\mathbb{DR}$/$\mathbb{RD}$-вершину, либо $\mathbb{C}$-вершину либо $\mathbb{CDL}$-вершину. Всего $2094+12+210+1=2317$ вариантов. Окружение третьего начальника (правый нижний угол) может быть вычислено по окружению второго (функция $\mathbf{RightCorner}$).

Итак, если начальник один, есть $2094$ вариантов информации. Можно проверить, что тип правого нижнего угла можно восстановить по окружению узла в середине верхней стороны, для всех случаев, кроме двух -- цепи $\mathbb{A}0$ с указателем $3$ и цепи  $\mathbb{B}1$ c указателем $2$. То есть, всего возможных сочетаний окружения первого начальника и типа второго существует не более $2092+19+19=2130$.

Если начальника три, то окружение первого и третьего восстанавливаются по окружению второго (функции $\mathbf{TopFromCorner}$ и $\mathbf{RightCorner}$), то есть всего вариантов информации не более $2317$. Итого вариантов информации для плоского случая не более $2094+2317+2130=6541$.

\medskip

Теперь рассмотрим случай, когда хотя бы один начальник лежит на левом или верхнем краю подклееной макроплитки. Тогда в случае одного начальника вариантов информации будет не более $2115\times 39 =82485$ (любое сочетание базового и подклееного окружения). Для вершин с двумя начальниками заметим, что когда первый начальник имеет краевой тип, мы можем вычислить подклееный тип правого нижнего угла по подклееному окружению первого начальника. То есть тут столько же вариантов информаций, $82485$.

Для вершин с тремя начальниками, третий начальник не может попасть на левый или верхней край подклееной макроплитки, то есть у него будет только базовое окружение, которое можно вычислить по функции $\mathbf{RightCorner}$.
Возможны три случая:

{\bf 1.} И первый и второй начальники имеют и подклееное и базовое окружение (то есть они оба попали на левую и верхнюю стороны подклееной макроплитки). Тогда
расширенное окружение второго начальника можно выбрать не более чем $82485$ способами. Подклееное окружение первого начальника можно вычислить по функции $\mathbf{TopFromCorner}$. Если считать, что базовое окружение первого начальника -- любое,  то всего
вариантов информации $82485\cdot 2115$.

{\bf 2.} Только первый начальник попадает на край. Тогда число вариантов информации не превосходит $2749 \cdot 2115$. (Произвольное базовое окружение первого начальника и произвольное базовое окружение второго).

{\bf 3.} Только второй начальник попадает на край.  Окружение первого начальника можно вычислить по функции $\mathbf{TopFromCorner}$. Тогда число вариантов информации не превосходит $82485$.
Таким образом, общее число всех вариантов не более $82485\cdot 2115+ 2115\cdot 2749+82485=180352395$.

Любых возможных информаций не более $181\cdot 10^6$. Всего сочетаний (тип, уровень, расширенное окружение, информация) не более чем $181\cdot 10^6 \cdot (82485+2749)$ (то есть либо произвольное расширенное окружение, либо базовое окружение, без подклееного). Это число не превосходит $16 \cdot 10^{12}$.

\medskip

{\bf Число флагов макроплиток}.
Ядро, то есть вершина в левом верхнем углу подклееной макроплитки, может иметь произвольное базовое окружение и произвольную информацию. Флагом макроплитки называется сочетание типа, уровня, базового окружения и информации ядра этой макроплитки и упорядоченной пары двух выходящих ребер, которые определяют стороны макроплитки. Одно из ребер может вести в подклееные области.
Плоских выходящих ребер у любой вершины не более $10$. Различных подклееных ребер существует не более $11$. То есть возможных флагов существует не более $2749\cdot 11\cdot 10 \cdot 16 \cdot 10^{12}<5\cdot 10^{18}$.

\medskip

Таким образом, число всех возможных букв в алфавите не превосходит $82485 \cdot 16\cdot 10^{12} \cdot 5\cdot 10^{18}< 7\cdot 10^{36}$.

\section{Разбор случаев расположения путей} \label{flip_section}

\begin{figure}[hbtp]
\centering
\includegraphics[width=0.6\textwidth]{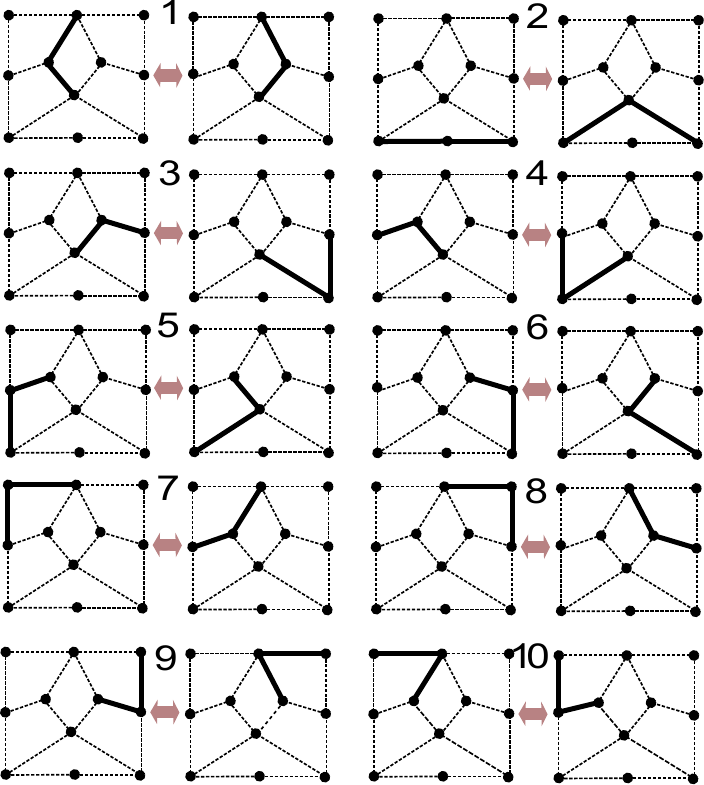}
\caption{По известному коду одного пути из пары можно восстановить код другого.}
\label{fig:flippaths}
\end{figure}

На рисунке~\ref{fig:flippaths} изображены пары путей. Далее мы покажем, как, зная код одного пути из пары, можно определить код другого пути либо установить, что код является мертвым.
В целом для этого требуется уметь вычислять код любой вершины в любой макроплитки по известным кодам трех остальных вершин. Рассмотрим произвольный путь $X_1e_1e_2X_2e_3e_4X_3$, где $X_1$, $X_2$, $X_3$ -- буквы, отвечающие кодам вершин, а $e_1$, $e_2$, $e_3$, $e_4$ -- буквы, отвечающие ребрам входа и выхода. Мы должны установить, к какой из десяти конфигураций относится наш путь, а также показать, как провести локальное преобразование, то есть получить код другого пути из соответствующей пары.

\medskip

В этой главе мы рассмотрим случай, когда все буквы $e_1$, $e_2$, $e_3$, $e_4$ отвечают плоским ребрам, то есть наш путь не выходит в подклееные области.
Случай выхода в подклееные области, когда хотя бы одна буква отвечает ребру, выходящему в подклейку, мы разберем в следующей главе.
Фактически, в этой и следующей главах мы вводим определяющие соотношения в полугруппе.

\medskip

\subsection{Обзор перебора случаев: сочетания типов сторон макроплитки}

Чтобы провести восстановление вершин для случаев $(1)-(10)$, требуется знать что расположено вокруг плитки, в которой идет наш путь.

\medskip

Итак, мы рассматриваем путь из двух последовательных ребер, проходящий внутри либо по границе некоторой макроплитки $T$. На рисунке~\ref{fig:flippaths} изображены макроплитки $T$ c проходящими внутри них нашими путями. Основная задача состоит в том, чтобы в каждом случае {\it восстановить путь}, то есть по известным значениям параметров трех вершин и входящим ребрам вычислить параметры четвертой вершины и код парного пути.

Для плоского пути, когда среди ребер вдоль пути не встречаются выходы в подклейки, параметр {\it флаг подклейки} восстанавливается тривиально -- для четвертой вершины он будет такой же, как и для остальных вершин (среди которых можно выбрать гарантированно лежащую не на границе подклееной макроплитки).
То есть нужно вычислить остальные параметры --  тип, уровень, окружение и информацию.

\medskip

Локальные преобразования с $1$ по $6$ мы рассмотрим каждое отдельно. При этом мы исследуем различные варианты расположения макроплитки, в которой проходит путь, и покажем, как вводятся определяющие соотношения в различных случаях расположения.

Локальные преобразования с $7$ по $10$ мы будем рассматривать все вместе для каждого случая расположения макроплитки.

\medskip

{\bf Обозначения.}
Буквами $e_1$, $e_2$, $e_3$, $e_4$, $e_{\mathbf{u}_1}$, $e_{\mathbf{u}_2}$, $e_{\mathbf{u}_3}$, $e_{\mathbf{ld}}$, $e_{\mathbf{r}}$, $e_{\mathbf{r}_2}$, $e_{\mathbf{r}_3}$ (и все подобные) будем обозначать входящие и выходящие ребра, соответствующие названию.  Символами $E_{\mathbf{l}}$, $E_{\mathbf{ld}}$, $E_{\mathbf{ddl}}$, (и подобными, указанными в разделе~\ref{addfunc} ``дополнительные функции'') будем обозначать ребра, получаемые в результате применения указанных функций. Все эти буквы мы будем использовать эти буквы при введении определяющих соотношений.

\smallskip

Иногда нам нужно будет зафиксировать часть параметров некоторой вершины $X$, чтобы потом использовать их для вычисления параметров других вершин. В этом случае мы будем использовать слово {\it назначение} $X$.

\smallskip

Пусть часть параметров вершины $X$ зафиксированы и известны. Пусть также оперируя известными параметрами $X$, мы вычисляем некоторые параметры вершин $X_1$, $X_2$, $X_3$. {\it Разрешенной} комбинацией будем называть такую упорядоченную четверку букв ($x$, $x_1$, $x_2$, $x_3$), что $x$ кодирует $X$ c зафиксированными параметрами, $x_i$ кодирует $X_i$, причем параметры совпадают с вычисленными. Заметим, что возможно несколько комбинаций разрешенных четверок, например, при разных параметрах флага подклейки. Часто кодирующие буквы мы будем обозначать так же, как сами вершины, например, $Z$, $J$, $Y$, $F$.

\medskip

{\bf Замечание.}
Рассматриваемые нами пути могут быть кусками более длинных путей. В том числе возможна ситуация, когда после прохождения нашего подпути, далее путь уходит по подклееному ребру, например, после вершины $J$. В этой ситуации, согласно определению, окружением $J$ является расширенное окружение, состоящее из базового окружения (в плоскости нашего пути) и подклееного окружения (в плоскости той подклееной плитки, куда уходит наш путь). При этом подклееное окружение часто может быть произвольным. Это отражается в том, что даже если тип, базовое окружение, флаг подклейки и информация у $J$ вычислены нами из условия расположения пути и остальных вершин, подклееное окружение $J$ часто может быть произвольным. В этом случае в качестве {\it разрешенных} букв, кодирующих $J$, можно выбрать буквы из множества с фиксированными параметрами типа, базового окружения, флаг подклейки и информации, но при этом с разными подклееными окружениями.

В этом случае мы будем говорить, что $J$ выбрана с точностью до подклееного окружения. Указанная ситуация может быть применима и к другим параметрам, например, информации.

\medskip

\subsection{Локальное преобразование 1}

\begin{figure}[hbtp]
\centering
\includegraphics[width=0.7\textwidth]{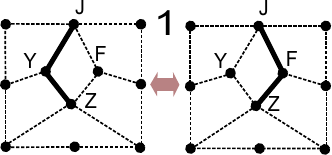}
\caption{Локальное преобразование $1$}
\label{flip1}
\end{figure}

Рассмотрим пару путей на рисунке~\ref{flip1}.

\medskip

\medskip

\medskip

{\bf Определяющие соотношения}.

\medskip

{\bf Назначение $Z$}. Пусть $Z$ -- вершина типа $\mathbb{C}$. Зафиксируем ее флаг подклейки, базовое окружение, подклееное окружение и информацию.

Теперь можно вычислить базовые окружения вершин $J$, $Y$, $F$, так как $J$ -- это просто первый начальник $Z$, а базовые окружения $Y$ и $F$ совпадают с окружением макроплитки, которое содержится в базовом окружении $J$. Зафиксируем некоторую информацию $J$. Информацию $Y$ и $F$ известна, так как известны окружения $J$ и вершины в правом нижнем углу (третий начальник $Z$).

Таким образом, выбрав произвольно базовое и подклееное окружение $Z$, подклееное окружение и информацию $J$, а также флаг макроплитки, мы можем вычислить значение остальных параметров. То есть, мы можем выписать множество разрешенных четверок букв ($Z$, $J$, $Y$, $F$), кодирующих описанное положение вершин.

Заметим, что если буква $Z$ зафиксирована, то буквы $Y$, $F$  -- точно определены, а буква $J$ -- определена с точностью до информации, то есть $J$ могут кодировать разные буквы, отличающиеся информацией.

\medskip

Для каждой разрешенной комбинации четырех букв $J$, $Z$, $Y$, $F$, введем следующие определяющие соотношения:

$Ze_1e_2Ye_1e_{\mathbf{u}_2}J=Ze_2e_3Fe_1e_{\mathbf{u}_1}J$

$Je_{\mathbf{u}_2}e_1Ye_2e_{1}J=Ze_{\mathbf{u}_1}e_1Fe_3e_{2}J$

\medskip

Два соотношения вводятся, потом что путь может быть пройден как в прямом порядке, так и в обратном.

Итак, $Z$ -- это произвольная буква типа $\mathbb{C}$ (произвольные расширенное окружение, флаг подклейки и информация), а у $J$ может быть произвольная информация и произвольное расширенное окружение. Остальные буквы в этом соотношении заданы однозначно.

\medskip

\subsection{Локальное преобразование 2}

\begin{figure}[hbtp]
\centering
\includegraphics[width=0.7\textwidth]{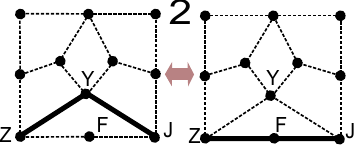}
\caption{Локальное преобразование $2$}
\label{flip2}
\end{figure}

Рассмотрим пару путей на рисунке~\ref{flip2}.

\medskip

\medskip

{\bf Определяющие соотношения}.
Мы рассмотрим по очереди шесть случаев расположения. Для каждого случая мы определим разрешенные множества букв, которыми могут быть вершины $Z$, $J$, $F$, $Y$. После чего для каждой комбинации букв, введем два соотношения, одно отвечает прямому обходу данного пути (по стрелкам), а другое -- обратному (против стрелок).

\medskip

\begin{figure}[hbtp]
\centering
\includegraphics[width=0.9\textwidth]{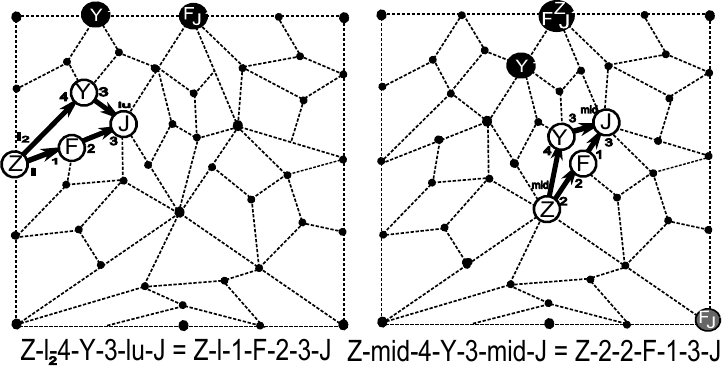}
\caption{Левое-верхнее и среднее положения. При входе и выходе в вершины обозначены коды входящих и выходящих ребер}
\label{pathplace2-1}
\end{figure}

{\bf 1. Левое-верхнее расположение}  (левая часть рисунка~\ref{pathplace2-1}).

\medskip

{\bf Назначение $Z$}. Пусть $Z$ -- вершина типа $\mathbb{LU/UL}$, $\mathbb{LD/DL}$ или $\mathbb{L}$  c произвольным расширенным окружением и информацией.

Начальники вершин $J$, $Y$, $F$ отмечены на рисунке~\ref{pathplace2-1} (левая часть) черными кругами.
Заметим, что мы можем вычислить их окружения. Одна из них входит в ту же цепь что и $Z$, а другая вычисляется с помощью функции $\mathbf{LevelPlus}$.  Значит, мы можем определить базовые окружения и информации вершин $J$, $Y$, $F$. Флаг подклейки может быть назначен произвольно (одинаковое значение для всех четырех вершин).

Таким образом, можно сформировать множество разрешенных четверок букв $Z$, $J$, $Y$, $F$.

\medskip

Для каждой разрешенной комбинации заданных букв, введем следующее определяющие соотношения:

\smallskip

$Ze_{\mathbf{l}_2}e_4Ye_3e_{\mathbf{lu}}J=Ze_{\mathbf{l}}e_1Fe_2e_3J$

$Je_{\mathbf{lu}}e_3Ye_4e_{\mathbf{l}_2}Z=Je_3e_2Fe_1e_{\mathbf{l}}Z$

\smallskip

Два соотношения отвечают двум направлениям обхода заданного участка, по стрелкам или против стрелок.

\medskip

{\bf Восстановление кода.} Зная коды $J$ и $Z$ можно вычислить код как $Y$, так и $F$, так как их начальники и базовые окружения легко вычисляются по известному окружению и начальникам $Z$, в частности, у $F$ общий набор начальников с $J$.

\medskip

\medskip

\medskip

{\bf 2. Среднее расположение} (правая часть рисунка~\ref{pathplace2-1}).

\medskip

{\bf Назначение $Z$}. Пусть $Z$ -- вершина с типом $\mathbb{C}$ и произвольным расширенным окружением, флагом подклейки и информацией.

Начальники вершин $J$, $Y$, $F$ отмечены черными кругами. Зная код $Z$, можно вычислить их окружения. То есть можно вычислить коды $J$, $Y$, $F$.

\medskip

Заметим, что если буква $Z$ зафиксирована, то $Y$, $F$ -- точно определены, а $J$ -- с точностью до произвольного подклееного окружения, при заданном базовом. То есть можно выделить множество разрешенных четверок букв, кодирующих наши вершины.

\medskip

Для каждой разрешенной комбинации заданных букв, введем следующие определяющие соотношения:

\smallskip

$Ze_{\mathbf{mid}}e_4Ye_3e_{\mathbf{mid}}J=Ze_2e_2Fe_1e_3J$

$Je_{\mathbf{mid}}e_3Ye_4e_{\mathbf{mid}}Z=Je_3e_1Fe_2e_2Z$

\smallskip

Два соотношения отвечают двум направлениям обхода заданного участка, по стрелкам или против стрелок.

\medskip

{\bf Восстановление кода.} Аналогично левому-верхнему расположению.

\medskip

\begin{figure}[hbtp]
\centering
\includegraphics[width=0.9\textwidth]{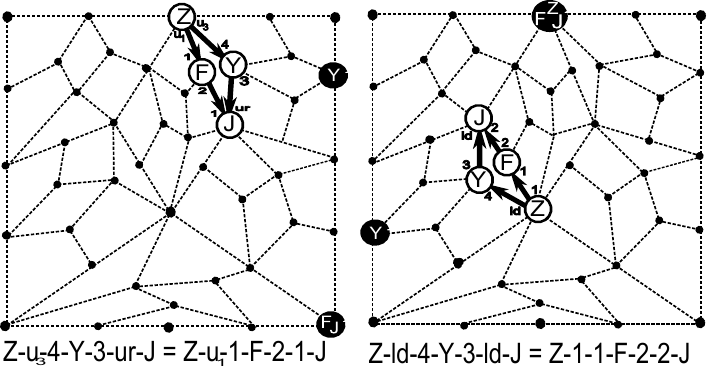}
\caption{Правое-верхнее и левое нижнее положения. При входе и выходе в вершины обозначены коды входящих и выходящих ребер}
\label{pathplace2-2}
\end{figure}

{\bf 3. Правое-верхнее расположение} (рисунок~\ref{pathplace2-2}, левая часть).

\medskip

{\bf Назначение $Z$}. Пусть $Z$ --  вершина типа $\mathbb{UL/LU}$ или $\mathbb{UR/RU}$  c произвольным расширенным окружением, флагом подклейки и информацией.

Окружение первого начальника $Y$ мы полуаем, используя функцию $\mathbf{UpRightChain}(Z)$.
Тип узла в правом нижнем углу -- $\mathbf{BottomRightType}(Z)$. То есть, коды всех остальных вершин можно вычислить.

Заметим, что если буква $Z$ зафиксирована, то буквы $F$, $Y$  -- точно определены, а буква $J$ -- с точностью до произвольного подклееного окружения, при заданном базовом.

\medskip

Для каждой разрешенной комбинации заданных букв, введем следующие определяющие соотношения:

\smallskip

$Ze_{\mathbf{u}_3}e_4Ye_3e_{\mathbf{ur}}J=Ze_{u_1}e_1Fe_2e_1J$

$Je_{\mathbf{ur}}e_3Ye_4e_{\mathbf{u}_3}Z=Je_1e_2Fe_1e_{\mathbf{u}_1}Z$

\smallskip

{\bf Восстановление кода.} Аналогично левому-верхнему расположению.

\medskip

\medskip

{\bf 4. Левое-нижнее расположение}  (рисунок~\ref{pathplace2-2}, правая часть).

\medskip

{\bf Назначение $Z$}. Пусть $Z$ вершина типа $\mathbb{C}$ c произвольным расширенным окружением, флагом подклейки и информацией.

Цепь вокруг левого нижнего угла $T$ можно получить как $\mathbf{BottomLeftChain}(Z)$. То есть, окружения вершин, отмеченных кругами можно вычислить, и коды всех вершин $J$, $Y$, $F$ можно вычислить, зная код $Z$.

Заметим, что если буква $Z$ зафиксирована, то буквы $F$, $Y$  -- точно определены, а буква $J$ -- с точностью до произвольного подклееного окружения, при заданном базовом.

\medskip

Для каждой разрешенной комбинации заданных букв, введем следующие определяющие соотношения:

\smallskip
$Ze_{\mathbf{ld}}e_4Ye_3e_{\mathbf{ld}}J=Ze_1e_1Fe_2e_2J$

$Je_{\mathbf{ld}}e_3Ye_4e_{\mathbf{ld}}Z=Je_2e_2Fe_1e_1Z$
\smallskip

{\bf Восстановление кода.} Аналогично левому-верхнему расположению.

\medskip

\begin{figure}[hbtp]
\centering
\includegraphics[width=0.9\textwidth]{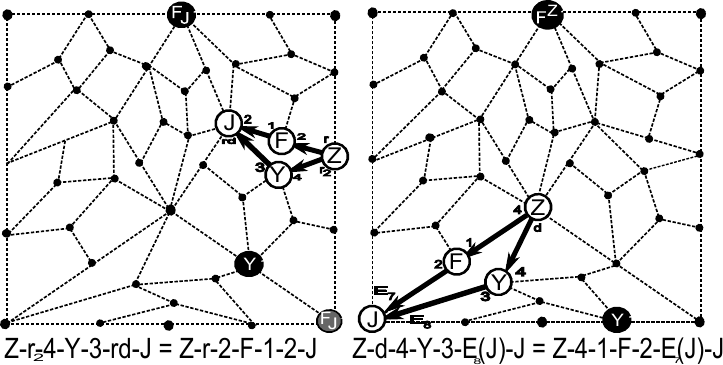}
\caption{Правое-нижнее и нижнее положения. При входе и выходе в вершины обозначены коды входящих и выходящих ребер}
\label{pathplace2-3}
\end{figure}

{\bf 5. Правое-нижнее расположение}  (рисунок~\ref{pathplace2-3}, левая часть).

\medskip

{\bf Назначение $Z$}. Пусть $Z$ -- вершина типа $\mathbb{RD/DR}$, $\mathbb{RU/UR}$  или $\mathbb{R}$  c произвольным расширенным окружением, флагом подклейки и информацией.

Цепь вокруг правого нижнего угла $T$ можно получить как $\mathbf{BottomRightChain}(Z)$. То есть, можно вычислить окружения всех начальников, значит и все коды вершин $J$, $Y$, $F$.

Заметим, что если буква $Z$ зафиксирована, то $F$, $Y$  -- точно определены, а $J$ -- с точностью до одного начальника и произвольного подклееного окружения, при заданном базовом.

\medskip

Для каждой разрешенной комбинации заданных букв, введем следующие определяющие соотношения:

\smallskip

$Ze_{\mathbf{r}_2}e_4Ye_3e_{\mathbf{rd}}J=Ze_{\mathbf{r}}e_2Fe_1e_2J$

$Je_{\mathbf{rd}}e_3Ye_4e_{\mathbf{r}_2}Z=Je_2e_1Fe_2e_{\mathbf{r}}Z$

\smallskip

{\bf Восстановление кода.} Аналогично левому-верхнему расположению.

\medskip

\medskip

{\bf 6. Нижнее расположение} (рисунок~\ref{pathplace2-3}, правая часть).

\medskip

{\bf Назначение $Z$}. Пусть $Z$ -- вершина типа $\mathbb{C}$ c произвольным расширенным окружением, флагом подклейки и информацией.

Заметим, что $J$ -- второй начальник $Z$. Окружения и информации $F$ и $Y$  мы также можем вычислить.

\medskip

{\bf Назначение $J$}. Тип, окружение и флаг подклейки для $J$ уже определены. Назначим произвольную информацию.

Теперь буквы $F$, $Y$  -- точно определены, а буква $J$ -- с точностью до произвольной информации.

\medskip

$E_7(J)$ и  $E_8(J)$ --  реберные буквы, полученные с помощью применения функций $E_7$ и $E_8$ к узлу $J$, то есть это тоже буквы, кодирующие ребра входов и выходов, но в зависимости от узла $J$.

\medskip

Для каждой разрешенной комбинации заданных букв, введем следующие определяющие соотношения:

\smallskip
$Ze_{\mathbf{d}}e_4Ye_3E_8(J)J=Ze_{4}e_1Fe_2E_7(J)J$

$JE_8(J)e_3Ye_4e_{\mathbf{d}}Z=JE_7(J)e_2Fe_1e_4Z$
\smallskip

{\bf Восстановление кода.} Аналогично левому-верхнему расположению.

\medskip

\subsection{Локальное преобразование 3}

\begin{figure}[hbtp]
\centering
\includegraphics[width=0.7\textwidth]{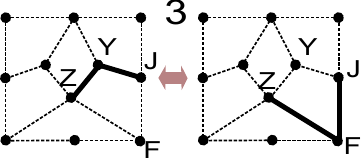}
\caption{Локальное преобразование $3$}
\label{flip3}
\end{figure}

Рассмотрим пару путей на рисунке~\ref{flip3}.

\medskip

\medskip

\begin{figure}[hbtp]
\centering
\includegraphics[width=0.9\textwidth]{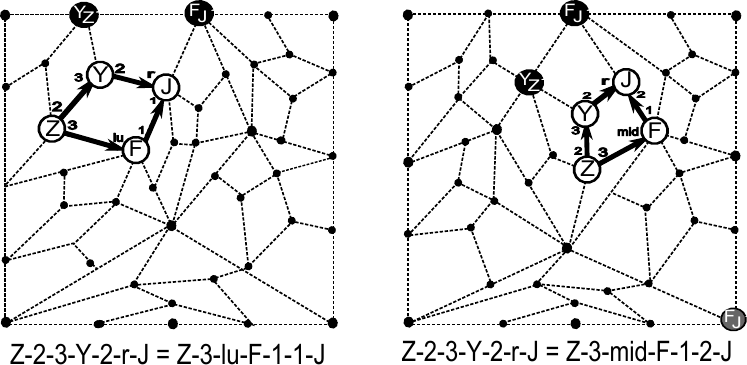}
\caption{Левое верхнее и среднее положения. При входе и выходе в вершины обозначены коды входящих и выходящих ребер}
\label{pathplace3-1}
\end{figure}

{\bf 1. Левое-верхнее расположение} (рисунок~\ref{pathplace3-1}, левая часть).

\medskip

{\bf Назначение $F$}. Пусть $F$ -- вершина типа $\mathbb{A}$ c произвольным базовым окружением, флагом подклейки и информацией.

Аналогично предыдущим случаям, этого достаточно, чтобы вычислить коды других вершин.

Заметим, что буквы $F$, $Y$ -- точно определены, а множества разрешенных букв $J$ и $Z$ -- различаются произвольным подклееным окружением, при заданном базовом.

\medskip

Для каждой разрешенной комбинации заданных букв, введем следующие определяющие соотношения:

\smallskip
$Ze_{2}e_3Ye_2e_{\mathbf{r}}J=Ze_{3}e_{\mathbf{lu}}Fe_1e_1J$

$Je_{\mathbf{r}}e_2Ye_3e_{2}Z=Je_1e_1Fe_{\mathbf{lu}}e_{3}Z$
\smallskip

\medskip

{\bf Восстановление кода.} Зная коды $J$ и $Z$ можно вычислить код как $Y$, так и $F$: у $Y$ базовое окружение и первый начальник как у $Z$, а $F$ является третьим начальником $Z$ и имеет общий набор начальников с $J$.

\medskip

\medskip

{\bf 2. Среднее расположение}  (рисунок~\ref{pathplace3-1}, правая часть).

\medskip

{\bf Назначение $F$}. Пусть $F$ --  вершина типа $\mathbb{B}$ c произвольным базовым окружением, флагом подклейки и информацией.

Аналогично предыдущим случаям, мы можем вычислить коды всех вершин ($J$ и $Z$ -- с точностью до произвольного подклееного окружения, при заданном базовом).

\medskip

Для каждой разрешенной комбинации заданных букв, введем следующие определяющие соотношения:

\smallskip
$Ze_{2}e_3Ye_2e_{\mathbf{r}}J=Ze_{3}e_{\mathbf{mid}}Fe_1e_2J$

$Je_{\mathbf{r}}e_2Ye_3e_{2}Z=Je_2e_1Fe_{\mathbf{mid}}e_{3}Z$
\smallskip

{\bf Восстановление кода.} Аналогично левому-верхнему расположению.

\medskip

\begin{figure}[hbtp]
\centering
\includegraphics[width=0.9\textwidth]{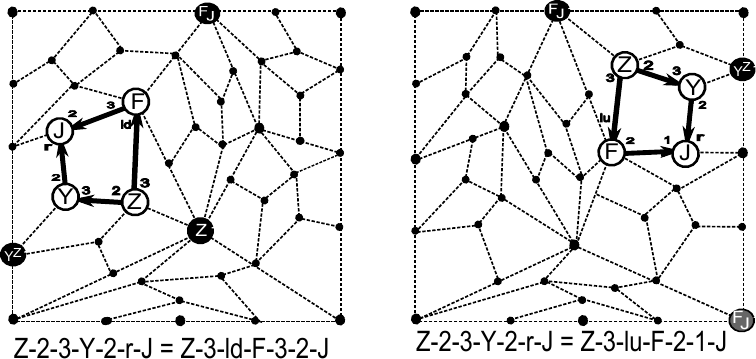}
\caption{Левое-нижнее и правое-верхнее положения. При входе и выходе в вершины обозначены коды входящих и выходящих ребер}
\label{pathplace3-2}
\end{figure}

{\bf 3. Левое нижнее расположение}  (рисунок~\ref{pathplace3-2}, левая часть).

\medskip

{\bf Назначение $F$}. Пусть $F$ -- вершина типа $\mathbb{A}$ c произвольным базовым окружением, флагом подклейки и информацией.

Цепь вокруг левого нижнего угла $T$ можно получить как $\mathbf{BottomLeftChain.FBoss}(F)$. Учитывая это, можно вычислить коды всех вершин ($J$ и $Z$ -- с точностью до произвольного подклееного окружения, при заданном базовом).

\medskip

Для каждой разрешенной комбинации заданных букв, введем следующие определяющие соотношения:

$Ze_{2}e_3Ye_2e_{\mathbf{r}}J=Ze_{3}e_{\mathbf{ld}}Fe_3e_2J$

$Je_{\mathbf{r}}e_2Ye_3e_{2}Z=Je_2e_3Fe_{\mathbf{ld}}e_{3}Z$

\medskip

{\bf Восстановление кода.} Аналогично левому-верхнему расположению.

\medskip

\medskip

{\bf 4. Правое-верхнее расположение} (рисунок~\ref{pathplace3-2}, правая часть).

\medskip

{\bf Назначение $F$}. Пусть $F$ -- вершина типа $\mathbb{B}$ c произвольным базовым окружением, флагом подклейки и информацией.

Цепь вокруг правого верхнего угла можно получить как  $\mathbf{UpRightChain.FBoss}(F)$.
Учитывая это, можно вычислить коды всех вершин ($J$ и $Z$ -- с точностью до произвольного подклееного окружения, при заданном базовом).

\medskip

Для каждой разрешенной комбинации заданных букв, введем следующие определяющие соотношения:

\smallskip
$Ze_{2}e_3Ye_2e_{\mathbf{r}}J=Ze_{3}e_{\mathbf{lu}}Fe_2e_1J$

$Je_{\mathbf{r}}e_2Ye_3e_{2}Z=Je_1e_2Fe_{\mathbf{lu}}e_{3}Z$

\medskip

{\bf Восстановление кода.} Аналогично левому-верхнему расположению.
\medskip

\begin{figure}[hbtp]
\centering
\includegraphics[width=0.9\textwidth]{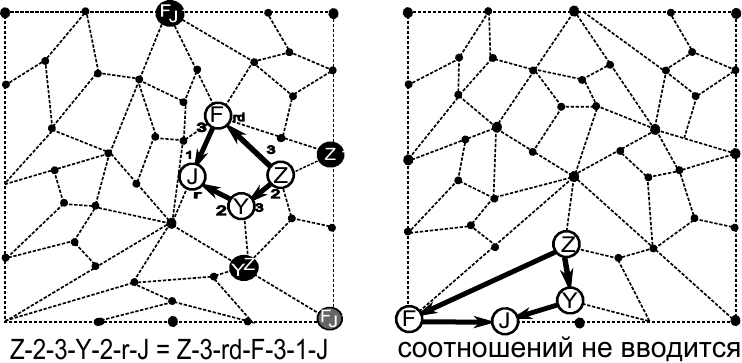}
\caption{Правое-нижнее и нижнее положения. При входе и выходе в вершины обозначены коды входящих и выходящих ребер}
\label{pathplace3-3}
\end{figure}

{\bf 5. Правое-нижнее расположение} (рисунок~\ref{pathplace3-3}, левая часть).

\medskip

{\bf Назначение $F$}. Пусть $F$ -- вершина типа $\mathbb{A}$ c произвольным базовым окружением, флагом подклейки и информацией.

Окружение вершины в середине правой стороны можно вычислить с помощью функции $\mathbf{RightFromB}$. Учитывая это, мы можем вычислить коды всех вершин ($J$ и $Z$ -- с точностью до произвольного подклееного окружения, при заданном базовом).

\medskip

Для каждой разрешенной комбинации заданных букв, введем следующие определяющие соотношения:

\smallskip
$Ze_{2}e_3Ye_2e_{\mathbf{r}}J=Ze_{3}e_{\mathbf{rd}}Fe_3e_1J$

$Je_{\mathbf{r}}e_2Ye_3e_{2}Z=Je_1e_3Fe_{\mathbf{rd}}e_{3}Z$

\medskip

{\bf Восстановление кода.} Аналогично левому-верхнему расположению.

\medskip

\medskip

{\bf 6. Нижнее расположение} (рисунок~\ref{pathplace3-3}, правая часть).

В этом случае соотношений мы не вводим, так как расположение пути удовлетворяет признакам мертвого паттерна. То есть участок пути с данным кодом не может являться подпутем достаточно длинного ненулевого пути.

\medskip

\subsection{Локальное преобразование 4}

\begin{figure}[hbtp]
\centering
\includegraphics[width=0.7\textwidth]{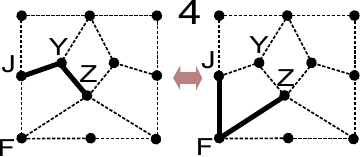}
\caption{Локальное преобразование $4$}
\label{flip4}
\end{figure}

Рассмотрим пару путей на рисунке~\ref{flip4}.

\medskip

\begin{figure}[hbtp]
\centering
\includegraphics[width=0.9\textwidth]{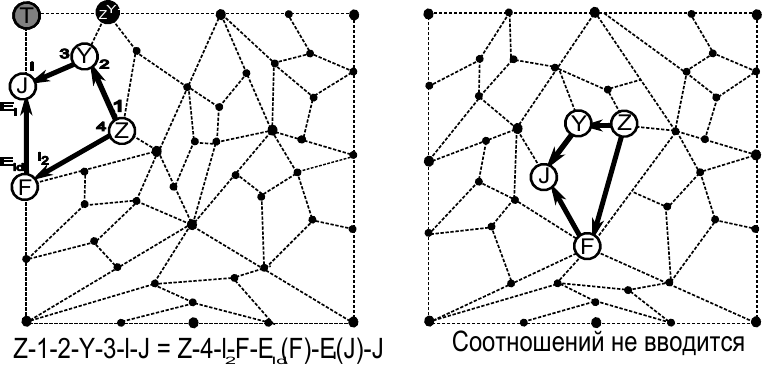}
\caption{Левое верхнее и среднее положения. При входе и выходе в вершины обозначены коды входящих и выходящих ребер}
\label{pathplace4-1}
\end{figure}

{\bf 1. Левое-верхнее расположение} (рисунок~\ref{pathplace4-1}, левая часть).

В этом случае сначала определим тип узла $T$ в левом верхнем углу. При этом $J$ и $F$ будут лежать на ребре выхода из $T$.

\medskip

{\bf Назначение $F$}. Пусть $F$ -- буква, кодирующая вершину типа $\mathbb{L}$, $\mathbb{LD}$ или $\mathbb{LU}$, лежащую в цепи вокруг $T$.

Заметим, учитывая расположение на рисунке, по известному коду $F$, мы можем вычислить окружения и начальников остальных трех вершин. То есть мы можем выписать все комбинации четырех букв $Z$, $Y$, $J$, $F$ такие, что их окружения и начальники соответствуют окружению и начальникам изображенных вершин. Будем называть такие комбинации {\it разрешенными}. При этом, буквы $Y$ и $F$ будут определены однозначно, а $Z$ и $J$ могут быть любыми буквами с заданным базовым окружением и начальниками, но с различными подклееными окружениями.

\medskip

Кроме того, буквы $E_\mathbf{l}$ и $E_{\mathbf{ld}}$ (кодирующие соответствующие на рисунке ребра входа и выхода) легко вычисляются по известным кодам $J$ и $F$.

\medskip

Для каждой разрешенной комбинации заданных букв, введем следующие определяющие соотношения:

\smallskip
$Ze_{1}e_2Ye_3e_{\mathbf{l}}J=Ze_{4}e_{\mathbf{l}_2}FE_{\mathbf{ld}}E_{\mathbf{l}}J$

$Je_{\mathbf{l}}e_{3}Ye_2e_{1}Z=JE_{\mathbf{l}}E_{\mathbf{ld}}Fe_{\mathbf{l}_2}e_{4}Z$

\medskip

{\bf Восстановление кода.} Зная коды $J$ и $Z$ можно вычислить код как $Y$, так и $F$: для $Y$ это очевидно, а $F$ является вторым начальником $Z$ и имеет общий набор начальников с $J$.

\medskip

\medskip

{\bf 2. Среднее расположение}, рисунок~\ref{pathplace4-1}, правая часть.

В этом случае соотношений мы не вводим, так как расположение пути удовлетворяет признакам мертвого паттерна. То есть участок пути с данным кодом не может являться подпутем достаточно длинного ненулевого пути.

\medskip

\begin{figure}[hbtp]
\centering
\includegraphics[width=0.9\textwidth]{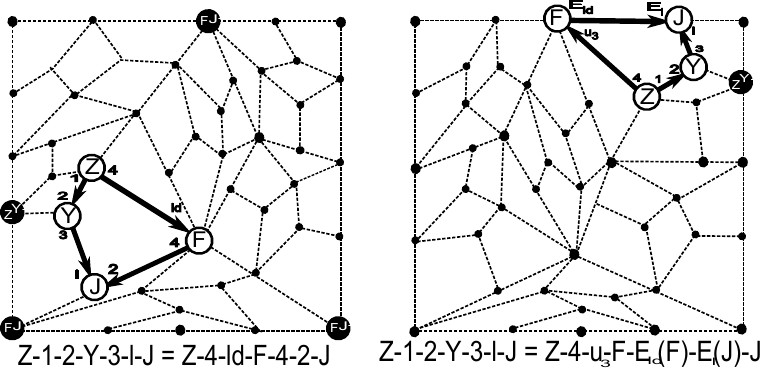}
\caption{Левое-нижнее и правое-верхнее положения. При входе и выходе в вершины обозначены коды входящих и выходящих ребер}
\label{pathplace4-2}
\end{figure}

{\bf 3. Левое нижнее расположение} (рисунок~\ref{pathplace4-2}, левая часть).

\medskip

{\bf Назначение $F$}. Пусть $F$ -- буква, кодирующая вершину типа $\mathbb{C}$, с произвольным окружением и начальниками.

Тогда базовые окружения вершин $Z$, $Y$, $J$ легко вычисляются. Начальники всех вершин отмечены черными кругами, и очевидно все они вычисляются, если мы знаем начальников $F$. Таким образом, мы можем аналогично определить все разрешенные четверки букв $Z$, $Y$, $J$, $F$, соответствующих вычисленным значениям окружений и начальников. У вершин $Z$ и $J$ могут быть разные подклееные окружения.

\medskip

Для каждой разрешенной комбинации заданных букв, введем следующие определяющие соотношения:

\smallskip
$Ze_{1}e_2Ye_3e_{\mathbf{l}}J=Ze_{4}e_{\mathbf{ld}}Fe_4e_2J$

$Je_{\mathbf{l}}e_3Ye_2e_{1}Z=Je_2e_4Fe_{\mathbf{ld}}e_{4}Z$

\medskip

{\bf Восстановление кода.} Аналогично левому-верхнему расположению.

\medskip

\medskip

{\bf 4. Правое верхнее расположение} (рисунок~\ref{pathplace4-2}, правая часть).

\medskip

{\bf Назначение $F$}. Пусть $F$ -- буква, кодирующая вершину типа $\mathbb{U}$, $\mathbb{UL}$ или $\mathbb{UR}$.

Тогда базовые окружения вершин $Z$, $Y$, $J$ легко вычисляются, учитывая, что мы можем применить функцию $\mathbf{UpRightChain}(F)$ и узнать тип вершины в правом верхнем углу.

Ясно также что начальники всех вершин также вычисляются, например у $J$ они такие же, как у $F$.  Таким образом, мы можем аналогично определить все разрешенные четверки букв $Z$, $Y$, $J$, $F$, соответствующих вычисленным значениям окружений и начальников. У вершин $Z$ и $J$ могут быть разные подклееные окружения.

\medskip

Кроме того, буквы $E_{\mathbf{l}}$ и $E_{\mathbf{ld}}$ (кодирующие соответствующие на рисунке ребра входа и выхода) легко вычисляются по известным кодам $J$ и $F$.

\medskip

Для каждой разрешенной комбинации заданных букв, введем следующие определяющие соотношения:

\smallskip
$Ze_{1}e_2Ye_3e_{\mathbf{l}}J=Ze_{4}e_{\mathbf{u}_3}FE_{\mathbf{ld}}E_{\mathbf{l}}J$

$Je_{\mathbf{l}}e_3Ye_2e_{1}Z=JE_1E_{\mathbf{ld}}Fe_{\mathbf{u}_3}e_{4}Z$

\medskip

{\bf Восстановление кода.} Аналогично левому-верхнему расположению.

\medskip

\begin{figure}[hbtp]
\centering
\includegraphics[width=0.9\textwidth]{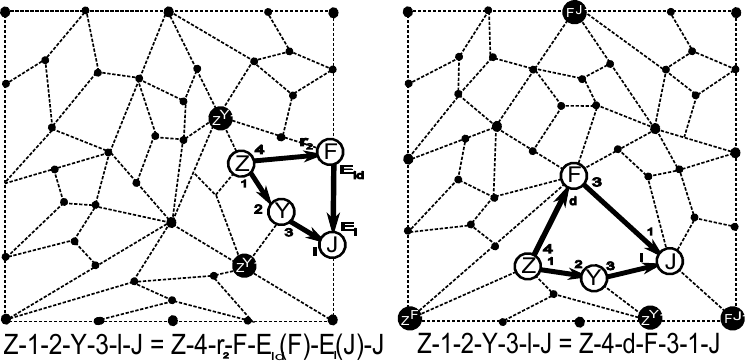}
\caption{Правое-нижнее и нижнее положения. При входе и выходе в вершины обозначены коды входящих и выходящих ребер}
\label{pathplace4-3}
\end{figure}

{\bf 5. Правое-нижнее расположение} (рисунок~\ref{pathplace4-3}, левая часть).

\medskip

{\bf Назначение $F$}. Пусть $F$ -- буква, кодирующая вершину типа $\mathbb{R}$, $\mathbb{RD}$ или $\mathbb{RU}$  c произвольным окружением и начальниками.

Тогда базовые окружения вершин $Z$, $Y$, $J$ легко вычисляются. Начальники всех вершин отмечены черными кругами, и очевидно все они вычисляются, если мы знаем начальников $F$. Таким образом, мы можем аналогично определить все разрешенные четверки букв $Z$, $Y$, $J$, $F$, соответствующих вычисленным значениям окружений и начальников. У вершин $Z$ и $J$ могут быть разные подклееные окружения.

\medskip
Кроме того, буквы $E_{\mathbf{l}}$ и $E_{\mathbf{ld}}$ (кодирующие соответствующие на рисунке ребра входа и выхода) легко вычисляются по известным кодам $J$ и $F$.

\medskip

Для каждой разрешенной комбинации заданных букв, введем следующие определяющие соотношения:

\smallskip

$Ze_{1}e_2Ye_3e_{\mathbf{l}}J=Ze_{4}e_{\mathbf{r}_2}FE_{\mathbf{ld}}E_{\mathbf{l}}J$

$Je_{\mathbf{l}}e_3Ye_2e_{1}Z=JE_{\mathbf{l}}E_{\mathbf{ld}}Fe_{\mathbf{r}_2}e_{4}Z$

\medskip

{\bf Восстановление кода.} Аналогично левому-верхнему расположению.

\medskip

\medskip

{\bf 6. Нижнее расположение} (рисунок~\ref{pathplace4-3}, правая часть).

\medskip

{\bf Назначение $F$}. Пусть $F$ -- буква, кодирующая вершину типа $\mathbb{C}$, c произвольным окружением и начальниками.

Тогда базовые окружения вершин $Z$, $Y$, $J$ легко вычисляются. Начальники всех вершин отмечены черными кругами, и очевидно все они вычисляются, если мы знаем начальников $F$. Таким образом, мы можем аналогично определить все разрешенные четверки букв $Z$, $Y$, $J$, $F$, соответствующих вычисленным значениям окружений и начальников. У вершин $Z$ и $J$ могут быть разные подклееные окружения.

\medskip

Для каждой разрешенной комбинации заданных букв, введем следующие определяющие соотношения:

\smallskip
$Ze_{1}e_2Ye_3e_{\mathbf{l}}J=Ze_{4}e_{\mathbf{d}}Fe_3e_1J$

$Je_{\mathbf{l}}e_3Ye_2e_{1}Z=Je_1e_3Fe_{\mathbf{d}}e_{4}Z$

\medskip

{\bf Восстановление кода.} Аналогично левому-верхнему расположению.

\medskip

\subsection{Локальное преобразование 5}

\begin{figure}[hbtp]
\centering
\includegraphics[width=0.7\textwidth]{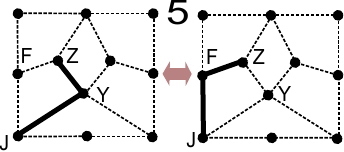}
\caption{Локальное преобразование $5$}
\label{flip5}
\end{figure}

Рассмотрим пару путей на рисунке~\ref{flip5}.

\medskip

\medskip

\begin{figure}[hbtp]
\centering
\includegraphics[width=0.9\textwidth]{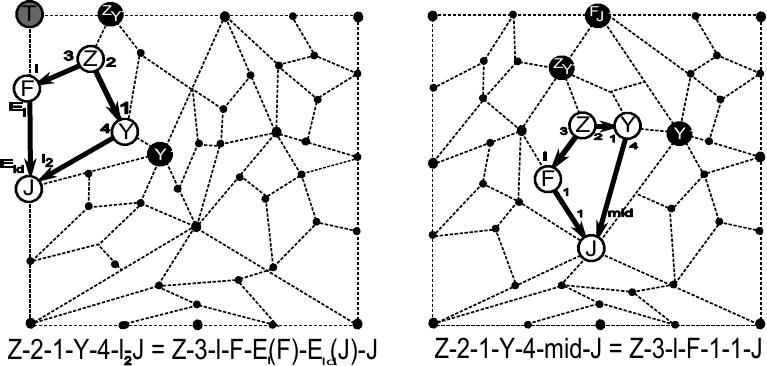}
\caption{Левое верхнее и среднее положения. При входе и выходе в вершины обозначены коды входящих и выходящих ребер}
\label{pathplace5-1}
\end{figure}

{\bf 1. Левое-верхнее расположение} (рисунок~\ref{pathplace5-1}, левая часть).

В этом случае сначала определим тип узла $T$ в левом верхнем углу. При этом $J$ и $F$ будут лежать на ребре выхода из $T$.

\medskip

{\bf Назначение $J$}. Пусть $J$ -- буква, кодирующая вершину типа $\mathbb{L}$, $\mathbb{LD}$ или $\mathbb{LU}$, лежащую в цепи вокруг $T$.

Аналогично предыдущему случаю, по известному коду $J$, мы можем вычислить окружения и начальников остальных трех вершин. То есть мы можем выписать все {\it разрешенные} комбинации четырех букв $Z$, $Y$, $J$, $F$ то есть такие, что их окружения и начальники соответствуют окружению и начальникам изображенных вершин. При этом, буквы $Y$ и $F$ будут определены однозначно, а $Z$ и $J$ могут быть любыми буквами с заданным базовым окружением и начальниками, но с различными подклееными окружениями.

\medskip

Кроме того, буквы $E_{\mathbf{l}}$ и $E_{\mathbf{ld}}$ (кодирующие соответствующие на рисунке ребра входа и выхода) легко вычисляются по известным кодам $J$ и $F$.

\medskip

Для каждой разрешенной комбинации заданных букв, введем следующие определяющие соотношения:

\smallskip
$Ze_{2}e_1Ye_4e_{\mathbf{l}_2}J=Ze_{3}e_{\mathbf{l}}FE_{\mathbf{l}}E_{\mathbf{ld}}J$

$Je_{\mathbf{l}_2}e_{4}Ye_1e_{2}Z=JE_{\mathbf{ld}}E_{\mathbf{l}}Fe_{\mathbf{l}}e_{3}Z$

\medskip

{\bf Восстановление кода.} Зная коды $J$ и $Z$ можно вычислить код $Y$ и $F$: для $Y$ это очевидно, а $F$ соответствует $\mathbf{Next.FBoss}(Z)$ и имеет общий набор начальников с $J$.

\medskip

\medskip

{\bf 2. Среднее расположение} (рисунок~\ref{pathplace5-1}, правая часть).

\medskip

{\bf Назначение $J$}. Пусть $J$ -- буква, кодирующая вершину типа $\mathbb{C}$, с произвольным окружением и начальниками.

Тогда базовые окружения вершин $Z$, $Y$, $F$ легко вычисляются. Начальники всех вершин отмечены черными кругами, и очевидно все они вычисляются, если мы знаем начальников $J$. Таким образом, мы можем аналогично определить все разрешенные четверки букв $Z$, $Y$, $J$, $F$, соответствующих вычисленным значениям окружений и начальников. У вершин $Z$ и $J$ могут быть разные подклееные окружения.

\medskip

Для каждой разрешенной комбинации заданных букв, введем следующие определяющие соотношения:

\smallskip
$Ze_{2}e_1Ye_4e_{\mathbf{mid}}J=Ze_{3}e_{\mathbf{l}}Fe_1e_1J$

$Je_{\mathbf{mid}}e_4Ye_1e_{2}Z=Je_1e_1Fe_{\mathbf{l}}e_{3}Z$

\medskip

{\bf Восстановление кода.} Аналогично левому-верхнему расположению.

\medskip

\begin{figure}[hbtp]
\centering
\includegraphics[width=0.9\textwidth]{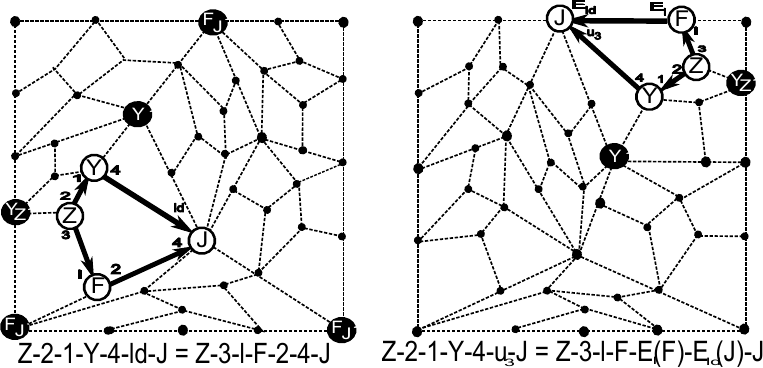}
\caption{Левое-нижнее и правое-верхнее положения. При входе и выходе в вершины обозначены коды входящих и выходящих ребер}
\label{pathplace5-2}
\end{figure}

{\bf 3. Левое нижнее расположение} (рисунок~\ref{pathplace5-2}, левая часть).

\medskip

{\bf Назначение $J$}. Пусть $J$ -- буква, кодирующая вершину типа $\mathbb{C}$, с произвольным окружением и начальниками.

Тогда базовые окружения вершин $Z$, $Y$, $F$ легко вычисляются. Начальники всех вершин отмечены черными кругами, и очевидно все они вычисляются, если мы знаем начальников $J$. Таким образом, мы можем аналогично определить все разрешенные четверки букв $Z$, $Y$, $J$, $F$, соответствующих вычисленным значениям окружений и начальников. У вершин $Z$ и $J$ могут быть разные подклееные окружения.

\medskip

Для каждой разрешенной комбинации заданных букв, введем следующие определяющие соотношения:

\smallskip
$Ze_{2}e_1Ye_4e_{\mathbf{ld}}J=Ze_{3}e_{\mathbf{l}}Fe_2e_4J$

$Je_{\mathbf{ld}}e_4Ye_1e_{2}Z=Je_4e_2Fe_{\mathbf{l}}e_{3}Z$

\medskip

{\bf Восстановление кода.} Аналогично левому-верхнему расположению.

\medskip

\medskip

{\bf 4. Правое верхнее расположение} (рисунок~\ref{pathplace5-2}, правая часть).

\medskip

{\bf Назначение $J$}. Пусть $J$ -- буква, кодирующая вершину типа $\mathbb{U}$, $\mathbb{UL}$ или $\mathbb{UR}$.

Тогда базовые окружения вершин $Z$, $Y$, $J$ легко вычисляются, учитывая, что мы можем применить функцию $\mathbf{UpRightChain}(F)$ и узнать тип вершины в правом верхнем углу.

Начальники всех вершин также вычисляются, например у $F$ они такие же, как у $J$.  Таким образом, мы можем аналогично определить все разрешенные четверки букв $Z$, $Y$, $J$, $F$, соответствующих вычисленным значениям окружений и начальников. У вершин $Z$ и $J$ могут быть разные подклееные окружения.

\medskip

Кроме того, буквы $E_{\mathbf{l}}$ и $E_{\mathbf{ld}}$ (кодирующие соответствующие на рисунке ребра входа и выхода) легко вычисляются по известным кодам $J$ и $F$.

\medskip

Для каждой разрешенной комбинации заданных букв, введем следующие определяющие соотношения:

\smallskip
$Ze_{2}e_1Ye_4e_{\mathbf{u}_3}J=Ze_{3}e_{\mathbf{l}}FE_{\mathbf{l}}E_{\mathbf{ld}}J$

$Je_{\mathbf{u}_3}e_4Ye_1e_{2}Z=JE_{\mathbf{ld}}E_{\mathbf{l}}Fe_{\mathbf{l}}e_{3}Z$

\medskip

{\bf Восстановление кода.} Аналогично левому-верхнему расположению.

\medskip

\begin{figure}[hbtp]
\centering
\includegraphics[width=0.9\textwidth]{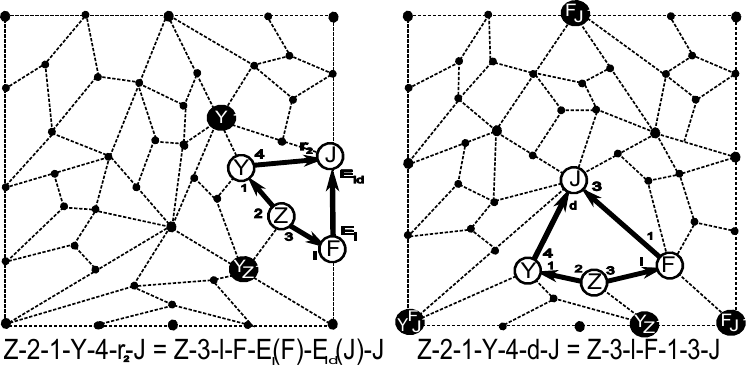}
\caption{Правое-нижнее и нижнее положения. При входе и выходе в вершины обозначены коды входящих и выходящих ребер}
\label{pathplace5-3}
\end{figure}

{\bf 5. Правое-нижнее расположение} (рисунок~\ref{pathplace5-3}, левая часть).

\medskip

{\bf Назначение $J$}. Пусть $J$ -- буква, кодирующая вершину типа $\mathbb{R}$, $\mathbb{RD}$ или $\mathbb{RU}$  c произвольным окружением и начальниками.

Тогда базовые окружения вершин $Z$, $Y$, $F$ (рисунок~\ref{pathplace5-3}, левая часть), очевидны. Начальники всех вершин отмечены черными кругами, и очевидно все они вычисляются, если мы знаем начальников $J$. Таким образом, мы можем аналогично определить все разрешенные четверки букв $Z$, $Y$, $J$, $F$, соответствующих вычисленным значениям окружений и начальников. У вершин $Z$ и $J$ могут быть разные подклееные окружения.

\medskip
Кроме того, буквы $E_{\mathbf{l}}$ и $E_{\mathbf{ld}}$ (кодирующие соответствующие на рисунке ребра входа и выхода) легко вычисляются по известным кодам $J$ и $F$.

\medskip

Для каждой разрешенной комбинации заданных букв, введем следующие определяющие соотношения:

\smallskip
$Ze_{2}e_1Ye_4e_{\mathbf{r}_2}J=Ze_{3}e_{\mathbf{l}}FE_{\mathbf{l}}E_{\mathbf{ld}}J$

$Je_{\mathbf{r}_2}e_4Ye_1e_{2}Z=JE_{\mathbf{ld}}E_{\mathbf{l}}Fe_{\mathbf{l}}e_{3}Z$

\medskip

{\bf Восстановление кода.} Аналогично левому-верхнему расположению.

\medskip

\medskip

{\bf 6. Нижнее расположение} (рисунок~\ref{pathplace5-3}, правая часть).

\medskip

{\bf Назначение $J$}. Пусть $J$ -- буква, кодирующая вершину типа $\mathbb{C}$, c произвольным окружением и начальниками.

Тогда базовые окружения вершин $Z$, $Y$, $F$ на вышеуказанном рисунке очевидны. Начальники всех вершин отмечены черными кругами, и очевидно все они вычисляются, если мы знаем начальников $J$. Таким образом, мы можем аналогично определить все разрешенные четверки букв $Z$, $Y$, $J$, $F$, соответствующих вычисленным значениям окружений и начальников. У вершин $Z$ и $J$ могут быть разные подклееные окружения.

\medskip

Для каждой разрешенной комбинации заданных букв, введем следующие определяющие соотношения:

\smallskip
$Ze_{2}e_1Ye_4e_{\mathbf{d}}J=Ze_{3}e_{\mathbf{l}}Fe_1e_3J$

$Je_{\mathbf{d}}e_4Ye_1e_{2}Z=Je_3e_1Fe_{\mathbf{l}}e_{3}Z$

\medskip

{\bf Восстановление кода.} Аналогично левому-верхнему расположению.

\medskip

\subsection{Локальное преобразование 6}

\begin{figure}[hbtp]
\centering
\includegraphics[width=0.7\textwidth]{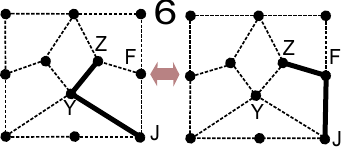}
\caption{Локальное преобразование $6$}
\label{flip6}
\end{figure}

Рассмотрим пару путей на рисунке~\ref{flip6}.

\medskip

\medskip

\begin{figure}[hbtp]
\centering
\includegraphics[width=0.9\textwidth]{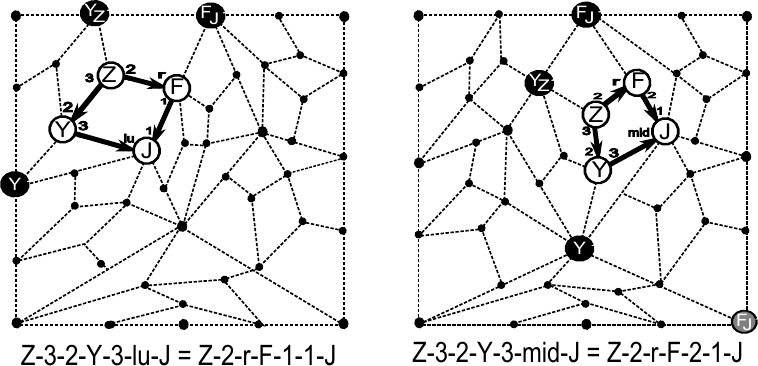}
\caption{Левое верхнее и среднее положения. При входе и выходе в вершины обозначены коды входящих и выходящих ребер}
\label{pathplace6-1}
\end{figure}

{\bf 1. Левое-верхнее расположение} (рисунок~\ref{pathplace6-1}, левая часть).

\medskip

{\bf Назначение $J$}. Пусть $J$ -- буква, кодирующая вершину типа $\mathbb{A}$, c произвольным окружением и произвольным начальником.

По известному коду $J$, мы можем выписать все {\it разрешенные} комбинации четырех букв $Z$, $Y$, $J$, $F$, то есть такие, что их окружения и начальники соответствуют окружению и начальникам изображенных на рисунке~\ref{pathplace6-1} вершин. При этом, буквы $Y$ и $F$ будут определены однозначно, а $Z$ и $J$ могут быть любыми буквами с заданным базовым окружением и начальниками, но с различными подклееными окружениями.

\medskip

Для каждой разрешенной комбинации заданных букв, введем следующие определяющие соотношения:

\smallskip
$Ze_{3}e_2Ye_3e_{\mathbf{lu}}J=Ze_{2}e_{\mathbf{r}}Fe_{1}e_{1}J$

$Je_{\mathbf{lu}}e_{3}Ye_2e_{3}Z=Je_{1}e_{1}Fe_{\mathbf{r}}e_{2}Z$

\medskip

{\bf Восстановление кода.} Зная коды $J$ и $Z$ можно вычислить код как $Y$, так и $F$. Окружение $Y$ и первый начальник такие же как у $Z$, второй начальник соответствует $\mathbf{Next.FBoss}(J)$, а третий -- вершина $J$. Окружение $F$ соответствует $0$-цепи вокруг $J$ с указателем $1$, начальник тот же что и у $J$.

\medskip

\medskip

{\bf 2. Среднее расположение} (рисунок~\ref{pathplace6-1}, правая часть).

\medskip

{\bf Назначение $J$}. Пусть $J$ -- буква, кодирующая вершину типа $\mathbb{B}$, с произвольным окружением и начальниками.

Тогда базовые окружения вершин $Z$, $Y$, $F$ очевидны. Начальники всех вершин отмечены черными кругами, и очевидно все они вычисляются, если мы знаем начальников $J$. Таким образом, мы можем аналогично определить все разрешенные четверки букв $Z$, $Y$, $J$, $F$, соответствующих вычисленным значениям окружений и начальников. У вершин $Z$ и $J$ могут быть разные подклееные окружения.

\medskip

Для каждой разрешенной комбинации заданных букв, введем следующие определяющие соотношения:

\smallskip

$Ze_{3}e_2Ye_3e_{\mathbf{mid}}J=Ze_{2}e_{\mathbf{r}}Fe_2e_1J$

$Je_{\mathbf{mid}}e_3Ye_2e_{3}Z=Je_1e_2Fe_{\mathbf{r}}e_{2}Z$

\medskip

{\bf Восстановление кода.} Зная коды $J$ и $Z$ можно вычислить код как $Y$, так и $F$. Окружение $Y$ и первый начальник такие же как у $Z$, второй начальник это вершина типа $\mathbb{C}$ с окружением как у $J$, а третий -- вершина $J$. Окружение $F$ вычисляется по известному окружению $J$, начальник тот же что и у $J$.

\medskip

\begin{figure}[hbtp]
\centering
\includegraphics[width=0.9\textwidth]{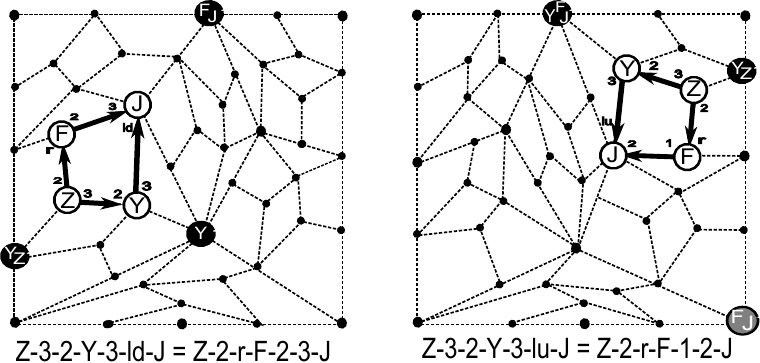}
\caption{Левое-нижнее и правое-верхнее положения. При входе и выходе в вершины обозначены коды входящих и выходящих ребер}
\label{pathplace6-2}
\end{figure}

{\bf 3. Левое нижнее расположение}  (рисунок~\ref{pathplace6-2}, левая часть).

\medskip

{\bf Назначение $J$}. Пусть $J$ -- буква, кодирующая вершину типа $\mathbb{A}$, с произвольным окружением и начальником.

Тогда базовые окружения вершин $Z$, $Y$, $F$ очевидны. Начальники всех вершин отмечены черными кругами, и очевидно все они вычисляются, если мы знаем код $J$. Таким образом, мы можем аналогично определить все разрешенные четверки букв $Z$, $Y$, $J$, $F$, соответствующих вычисленным значениям окружений и начальников. У вершин $Z$ и $J$ могут быть разные подклееные окружения.

\medskip

Для каждой разрешенной комбинации заданных букв, введем следующие определяющие соотношения:

\smallskip
$Ze_{3}e_2Ye_3e_{\mathbf{ld}}J=Ze_{2}e_{\mathbf{r}}Fe_2e_3J$

$Je_{\mathbf{ld}}e_3Ye_2e_{3}Z=Je_3e_2Fe_{\mathbf{r}}e_{2}Z$

\medskip

{\bf Восстановление кода.} Зная коды $J$ и $Z$ можно вычислить код как $Y$, так и $F$. Окружение $Y$ и первый начальник такие же как у $Z$, второй начальник соответствует узлу $\mathbb{C}$ c окружением как у $J$, а третий -- сама вершина $J$. Окружение $F$ вычисляется по известному окружению $J$, начальник тот же что и у $J$.

\medskip

\medskip

{\bf 4. Правое верхнее расположение} (рисунок~\ref{pathplace6-2}, правая часть).

\medskip

{\bf Назначение $J$}. Пусть $J$ -- буква, кодирующая вершину типа $\mathbb{B}$, c произвольным окружением и начальниками.

Тогда базовые окружения вершин $Z$, $Y$, $J$, очевидны.
Заметим, что мы можем применить функцию $\mathbf{UpRightChain}(F)$ и узнать тип вершины в правом верхнем углу. То есть, начальники всех вершин, также вычисляются. Таким образом, мы можем аналогично определить все разрешенные четверки букв $Z$, $Y$, $J$, $F$, соответствующих вычисленным значениям окружений и начальников. У вершин $Z$ и $J$ могут быть разные подклееные окружения.

\medskip

Для каждой разрешенной комбинации заданных букв, введем следующие определяющие соотношения:

\smallskip

$Ze_{3}e_2Ye_3e_{\mathbf{lu}}J=Ze_{2}e_{\mathbf{r}}Fe_1e_{2}J$

$Je_{\mathbf{lu}}e_3Ye_2e_{3}Z=Je_{2}e_{1}Fe_{\mathbf{r}}e_{2}Z$

\medskip

{\bf Восстановление кода.} Зная коды $J$ и $Z$ можно вычислить код как $Y$, так и $F$. Окружение $Y$ и первый начальник такие же как у $Z$, второй начальник соответствует $\mathbf{FBoss}(J)$, а третий -- вершина $J$. Окружение $F$ вычисляется по известному окружению $J$, начальники те же что и у $J$.

\medskip

\begin{figure}[hbtp]
\centering
\includegraphics[width=0.9\textwidth]{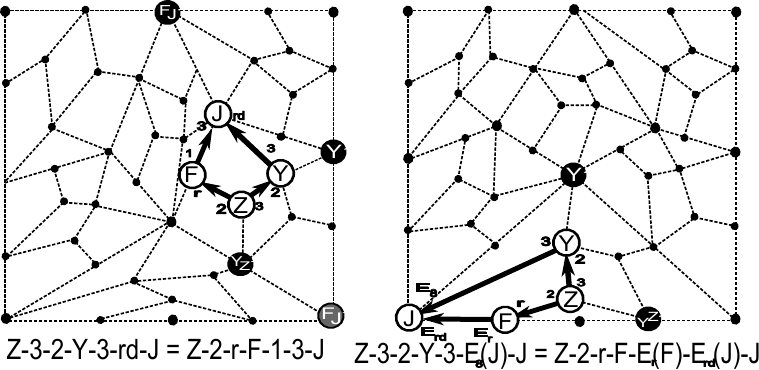}
\caption{Правое-нижнее и нижнее положения. При входе и выходе в вершины обозначены коды входящих и выходящих ребер}
\label{pathplace6-3}
\end{figure}

{\bf 5. Правое-нижнее расположение} (рисунок~\ref{pathplace6-3}, левая часть).

\medskip

{\bf Назначение $J$}. Пусть $J$ -- буква, кодирующая вершину типа $\mathbb{B}$, c произвольным окружением и начальниками.

Тогда базовые окружения вершин $Z$, $Y$, $F$ (рисунок~\ref{pathplace6-3}, левая часть), очевидны. Начальники всех вершин отмечены черными кругами, и очевидно все они вычисляются, если мы знаем начальников $J$. Таким образом, мы можем аналогично определить все разрешенные четверки букв $Z$, $Y$, $J$, $F$, соответствующих вычисленным значениям окружений и начальников. У вершин $Z$ и $J$ могут быть разные подклееные окружения.

\medskip

Для каждой разрешенной комбинации заданных букв, введем следующие определяющие соотношения:

\smallskip
$Ze_{3}e_2Ye_3e_{\mathbf{rd}}J=Ze_{2}e_{\mathbf{r}}Fe_{1}e_{3}J$

$Je_{\mathbf{rd}}e_3Ye_2e_{3}Z=Je_{3}e_{1}Fe_{\mathbf{r}}e_{2}Z$

\medskip

{\bf Восстановление кода.} Зная коды $J$ и $Z$ можно вычислить код как $Y$, так и $F$. Окружение $Y$ и первый начальник такие же как у $Z$, второй начальник вычисляется с помощью процедуры $\mathbf{RightFromB}(J)$, а третий -- вершина $J$. Окружение $F$ вычисляется по известному окружению $J$, начальники те же что и у $J$.

\medskip

\medskip

{\bf 6. Нижнее расположение} (рисунок~\ref{pathplace6-3}, правая часть).

\medskip

{\bf Назначение $J$}. Пусть $J$ -- буква, кодирующая вершину типа $\mathbb{C}$, c произвольным окружением и начальниками.

Тогда базовые окружения вершин $Z$, $Y$, $F$ на вышеуказанном рисунке очевидны. Начальники всех вершин отмечены черными кругами, и очевидно все они вычисляются, если мы знаем начальников $J$. Таким образом, мы можем аналогично определить все разрешенные четверки букв $Z$, $Y$, $J$, $F$, соответствующих вычисленным значениям окружений и начальников. У вершин $Z$ и $J$ могут быть разные подклееные окружения.

\medskip

Для каждой разрешенной комбинации заданных букв, введем следующие определяющие соотношения:

\smallskip
$Ze_{2}e_1Ye_4e_{\mathbf{d}}J=Ze_{3}e_{\mathbf{l}}Fe_1e_3J$

$Je_{\mathbf{d}}e_4Ye_1e_{2}Z=Je_3e_1Fe_{\mathbf{l}}e_{3}Z$

\medskip

{\bf Восстановление кода.} Зная коды $J$ и $Z$ можно вычислить код как $Y$, так и $F$. Окружение $Y$ и первый начальник такие же как у $Z$, второй начальник -- вершина $J$, третий соответствует $\mathbf{SBoss}(J)$. Окружение $F$ соответствует $1$-цепи вокруг третьего начальника $J$ с указателем, соответствующим входу в узел по нижней стороне, функция $E_{ddr}(\mathbf{TBoss}(J))$. Начальники те же что и у $J$.

\bigskip

\medskip

{\bf Локальные преобразования 7, 8, 9, 10.}

Рассмотрим узел, лежащий в середине верхней стороны макроплитки, где расположен наш путь. Этот узел должен входить в некоторую цепь. То есть может быть три варианта цепи вокруг центра, который может иметь следующий тип:

\smallskip

$\mathbb{C}$, $\mathbb{B}$, $\mathbb{A}$, $\mathbb{UL}/ \mathbb{UL}$,  $\mathbb{UR}/ \mathbb{RU}$, $\mathbb{DL}/ \mathbb{LD}$,  $\mathbb{DR}/ \mathbb{RD}$,  $\mathbb{U}$, $\mathbb{L}$, $\mathbb{D}$, $\mathbb{R}$, $\mathbb{CUL}$, $\mathbb{CUR}$, $\mathbb{CDL}$, $\mathbb{CDR}.$

\smallskip

 Мы разберем все эти случаи и введем соотношения в соответствии с устройством всех возможных цепей.

\medskip

{\bf Замечание.} Ниже мы не будем выписывать симметричное к введенному соотношение, отвечающее проходу пути в обратном порядке. Просто будем считать, что соотношений вводится в два раза больше.

\medskip

\subsection{Случай цепи $\mathbb{C}1$; преобразования 8 и 9}

\begin{figure}[hbtp]
\centering
\leftskip=-1.0cm
\includegraphics[width=1\textwidth]{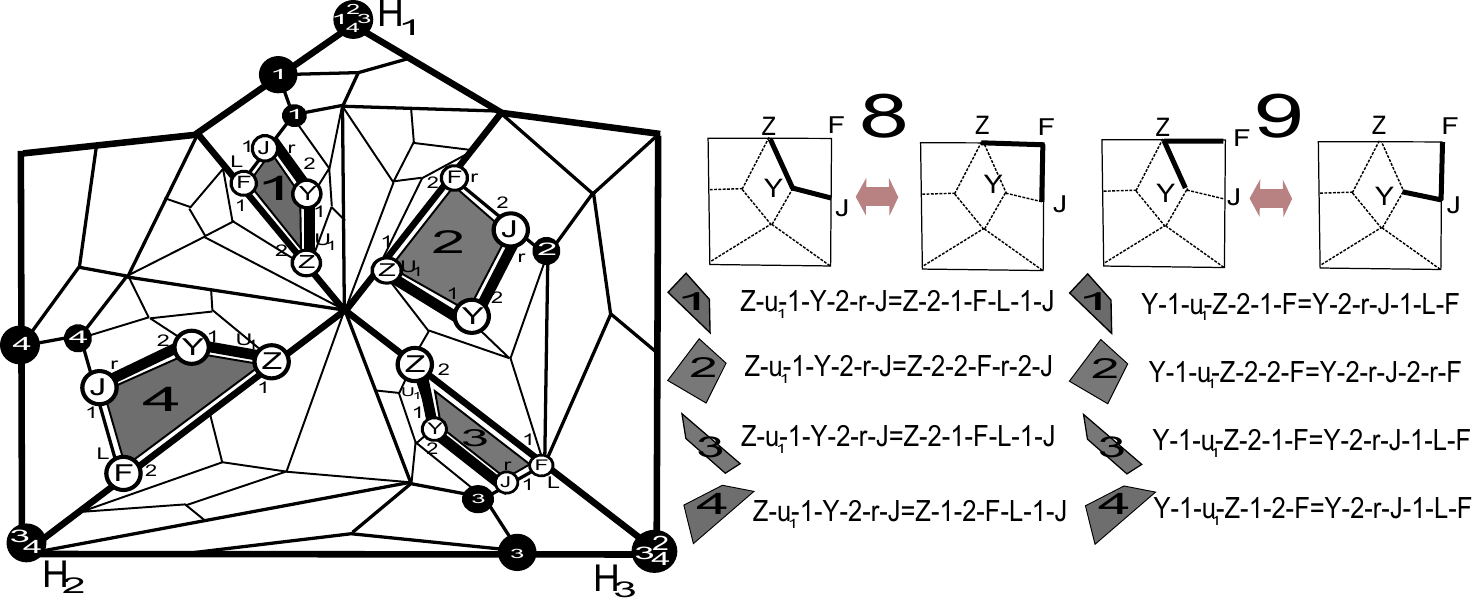}
\caption{Случаи расположения пути вокруг C1-цепи и соответствующие им определяющие соотношения для локальных преобразований 8 и 9}
\label{C1b}
\end{figure}

В правой части рисунка~\ref{C1b} изображены локальные преобразования $8$ и $9$.
Обозначим макроплитку, в которой проходят пути как $T$. Вершина $Z$ попадает в середину верхней стороны $T$.

\medskip

Итак, пусть середина верхней стороны $T$ входит в $\mathbb{C}1$-цепь. Различные случаи расположения интересующих нас путей показаны в левой части рисунка~\ref{C1b}.

\medskip

Сначала введем определяющие соотношения для данного случая.
Зафиксируем $\mathbb{C}$-вершину с некоторым базовым окружением и начальниками $\mathbf{FBoss}(Z)$, $\mathbf{SBoss}(Z)$, $\mathbf{TBoss}(Z)$. Исходя из этих данных, мы можем определить коды вершин $Z$, $Y$, $F$, $J$ во всех четырех случаях расположения и выписать определяющие соотношения для локальных преобразований $8$ и $9$.

\medskip

Покажем, как можно определить эти коды.
Вершины $H_1$, $H_2$, $H_3$  являются начальниками центральной $\mathbb{C}$-вершины. Зная типы, окружения и уровни вершин $H_1$, $H_2$, $H_3$, можно выписать типы, окружения и уровни $Z$, $Y$, $F$, $J$ во всех четырех случаях расположения, что видно на рисунке~\ref{C1b}. Например, окружение $F$ в случае $4$ соответствует $1$-цепи вокруг $H_2$ с указателем $E_7(H_2)$, и уровень $F$ равен $2$.

На рисунке~\ref{C1b} черными кругами отмечены вершины, являющиеся начальниками хотя бы одной из вершин $Z$, $Y$, $F$, $J$ в каждом из четырех случаев, число в круге обозначает номер случая. Окружение всех вершин, отмеченных черными кругами, также можно выписать, зная начальников центральной $\mathbb{C}$-вершины. То есть, мы можем вычислить начальников каждой из четырех вершин во всех случаях.

Уровни вершин $F$ во всех случаях равны $2$, а остальных боковых вершин ($J$ и $Z$) равны~$1$.

\medskip

Таким образом, мы можем определить множество букв $Z$, кодирующих вершины с заданным базовым окружением и начальниками, аналогично с остальными вершинами. Для локального преобразования $8$ буквы $F$ и $Y$ будут выбраны однозначно, а  $Z$ и $J$ -- с точностью до подклееного окружения. Для локального преобразования $9$ все наоборот.

\medskip

{\bf Соотношения для преобразования $8$.}\

\medskip

{\bf Назначение $Z$}. Для каждого из четырех случаев расположения, обозначим символом $Z$ буквы в алфавите, соответствующие кодам вершин с заданными типом, уровнем, базовым окружением и информацией. То есть эти буквы отличаются друг от друга только различными подклееными окружениями.

\medskip

{\bf Назначение $J$}. Аналогично определим буквы $J$, как буквы в алфавите, соответствующие кодам вершин $J$ с заданными типом, уровнем, базовым окружением и информацией, и произвольным подклееным окружением.

\medskip

{\bf Назначение $Y$ и $F$}. Буквы $Y$ и $F$ мы определим как конкретные буквы в алфавите с заданными типом, уровнем, базовым окружением и информацией, при пустом подклееном окружении. Мы это делаем, так как путь проходит через $F$ и $Y$ по плоским ребрам.

Теперь для каждого случая и для каждого разрешенного набора букв $Z$, $F$, $Y$, $J$ выпишем соотношение, представленное в правой части рисунка~\ref{C1b}.
Буквы $Z$, $F$, $Y$, $J$ мы только что определили, а остальные символы в соотношении обозначают входящие и выходящие ребра.

\medskip

{\bf Соотношения для преобразования $9$.}
Тут уже $Z$ и $J$ определены однозначно, а $Y$ и $F$ с точностью для подклееного окружения.
В остальном, все аналогично.

\medskip

{\bf Оценка числа соотношений.}
Так как в каждом соотношении две буквы определены однозначно, а другие две -- с точностью для подклееных окружений, то для каждого случая мы выписываем число соотношений, равное квадрату числа возможных подклееных окружений. Кроме того, все соотношения могут быть выписаны для любого значения параметра {\it флаг макроплитки}.

Мы вводим не более $8FP^2\mathbf{Num}(\mathbb{C})$ соотношений, где $F$ -- число различных флагов макроплиток, $P$ -- число различных подклееных окружений, $\mathbf{Num}(\mathbb{C})$ -- число вершин типа $\mathbb{C}$, то есть число сочетаний ``тип-уровень-окружение-информация''.

\medskip

{\bf Характеризация.}
Пусть есть слово $W$ представляющее собой код пути $X_1 e_1 e_2 X_2 e_3 e_4 X_3$. Покажем, как по нему установить, имеем ли мы дело с локальным преобразованием $8$ или $9$, а также как его провести.

На рисунке~\ref{C1b} отмечены входящие и выходящие ребра в каждом из четырех случаев. Исходя из этого легко выписать свойства пути, позволяющие нам установить, с каким случаем мы имеем дело. Все показано в таблицах~\ref{tableC1ba} и~\ref{tableC1bb}.

\medskip

\begin{table}[hbtp]
\caption{Характеристические условия на принадлежность пути к случаю $\mathbb{C}1$ цепи. }
\centering
\begin{tabular}{|c|c|c|c|}   \hline
Условие на буквы  & \x{Симметричное условие \cr (проход в обратном порядке)} & \x{ случай}  & \x{Локальное \cr преобр} \cr \hline
$\mathbf{Surr}(X_1)=\mathbb{C}11$, $e_1 e_2 e_3 e_4 = u_1 1 2 r$  & $\mathbf{Surr}(X_3)=\mathbb{C}11$, $e_1 e_2 e_3 e_4 = r 2 1 u_1$ & 1  & $8$ (левая) \cr \hline
$\mathbf{Surr}(X_1)=\mathbb{C}11$, $e_1 e_2 e_3 e_4 = 2 1 l 1$ & $\mathbf{Surr}(X_3)=\mathbb{C}11$, $e_1 e_2 e_3 e_4 = 1 l 1 2$  & 1 & $8$ (правая) \cr \hline
$\mathbf{Surr}(X_1)=\mathbb{C}12$, $e_1 e_2 e_3 e_4 = u_1 1 2 r$ & $\mathbf{Surr}(X_3)=\mathbb{C}12$, $e_1 e_2 e_3 e_4 = r 2 1 u_1$  & 2  & $8$ (левая) \cr \hline
$\mathbf{Surr}(X_1)=\mathbb{C}12$, $e_1 e_2 e_3 e_4 = 1 2 r 2$ & $\mathbf{Surr}(X_3)=\mathbb{C}12$, $e_1 e_2 e_3 e_4 = 2 r 2 1$  & 2 & $8$ (правая)\cr \hline
$\mathbf{Surr}(X_1)=\mathbb{C}13$, $e_1 e_2 e_3 e_4 = u_1 1 2 r$ & $\mathbf{Surr}(X_3)=\mathbb{C}13$, $e_1 e_2 e_3 e_4 = r 2 1 u_1$  & 3  & $8$ (левая) \cr \hline
$\mathbf{Surr}(X_1)=\mathbb{C}13$, $e_1 e_2 e_3 e_4 = 2 1 l 1$ & $\mathbf{Surr}(X_3)=\mathbb{C}13$, $e_1 e_2 e_3 e_4 = 1 l 1 2$  & 3 & $8$ (правая)\cr \hline
$\mathbf{Surr}(X_1)=\mathbb{C}14$, $e_1 e_2 e_3 e_4 = u_1 1 2 r$ & $\mathbf{Surr}(X_3)=\mathbb{C}14$, $e_1 e_2 e_3 e_4 = r 2 1 u_1$  & 4  & $8$ (левая) \cr \hline
$\mathbf{Surr}(X_1)=\mathbb{C}14$, $e_1 e_2 e_3 e_4 = 1 2 l 1$ & $\mathbf{Surr}(X_3)=\mathbb{C}14$, $e_1 e_2 e_3 e_4 = 1 l 2 1$  & 4 & $8$ (правая)\cr \hline
  \end{tabular}
\label{tableC1ba}
\end{table}

\medskip

\begin{table}[hbtp]
\caption{Характеристические условия на принадлежность пути к случаю $\mathbb{C}1$ цепи. }
\centering
\begin{tabular}{|c|c|c|c|}   \hline
Условие на буквы  & \x{Симметричное условие \cr (проход в обратном порядке)} & \x{ случай}  & \x{Локальное \cr преобр} \cr \hline
\x{ $\mathbf{Type}(X_1)=\mathbb{B}$, $\mathbf{FBoss}(X_1)=\mathbb{C}11$, \cr $e_1 e_2 e_3 e_4 = 1 u_1 2 1$}  & \x{$\mathbf{Type}(X_3)=\mathbb{B}$,  $\mathbf{FBoss}(X_3)=\mathbb{C}11$, \cr $e_1 e_2 e_3 e_4 = 1 2 u_1 1$} & 1  & $9$ (левая) \cr \hline
\x{ $\mathbf{Type}(X_1)=\mathbb{B}$, $\mathbf{FBoss}(X_1)=\mathbb{C}11$, \cr $e_1 e_2 e_3 e_4 = 2 r 1 l$} & \x{ $\mathbf{Type}(X_3)=\mathbb{B}$, $\mathbf{FBoss}(X_3)=\mathbb{C}11$, \cr $e_1 e_2 e_3 e_4 = l 1 r 2$}  & 1 & $9$ (правая) \cr \hline
\x{ $\mathbf{Type}(X_1)=\mathbb{B}$, $\mathbf{FBoss}(X_1)=\mathbb{C}12$, \cr $e_1 e_2 e_3 e_4 = 1 u_1 1 2$}  & \x{$\mathbf{Type}(X_3)=\mathbb{B}$, $\mathbf{FBoss}(X_3)=\mathbb{C}12$, \cr $e_1 e_2 e_3 e_4 = 2 1 u_1 1$} & 2  & $9$ (левая) \cr \hline
\x{ $\mathbf{Type}(X_1)=\mathbb{B}$, $\mathbf{FBoss}(X_1)=\mathbb{C}12$, \cr $e_1 e_2 e_3 e_4 = 2 r 2 r$} & \x{$\mathbf{Type}(X_3)=\mathbb{B}$, $\mathbf{FBoss}(X_3)=\mathbb{C}12$, \cr $e_1 e_2 e_3 e_4 = r 2 r 2$}  & 2 & $9$ (правая) \cr \hline
\x{ $\mathbf{Type}(X_1)=\mathbb{B}$, $\mathbf{FBoss}(X_1)=\mathbb{C}13$, \cr $e_1 e_2 e_3 e_4 = 1 u_1 2 1$}  & \x{$\mathbf{Type}(X_3)=\mathbb{B}$, $\mathbf{FBoss}(X_3)=\mathbb{C}13$, \cr $e_1 e_2 e_3 e_4 = 1 2 u_1 1$} & 3  & $9$ (левая) \cr \hline
\x{ $\mathbf{Type}(X_1)=\mathbb{B}$, $\mathbf{FBoss}(X_1)=\mathbb{C}13$, \cr $e_1 e_2 e_3 e_4 = 2 r 1 l$} & \x{$\mathbf{Type}(X_3)=\mathbb{B}$, $\mathbf{FBoss}(X_3)=\mathbb{C}13$, \cr $e_1 e_2 e_3 e_4 = l 1 r 2$}  & 3 & $9$ (правая) \cr \hline
\x{$\mathbf{Type}(X_1)=\mathbb{B}$, $\mathbf{FBoss}(X_1)=\mathbb{C}14$, \cr $e_1 e_2 e_3 e_4 = 1 u_1 1 2$}  & \x{$\mathbf{Type}(X_3)=\mathbb{B}$, $\mathbf{FBoss}(X_3)=\mathbb{C}14$, \cr $e_1 e_2 e_3 e_4 = 2 1 u_1 1$} & 4  & $9$ (левая) \cr \hline
\x{$\mathbf{Type}(X_1)=\mathbb{B}$, $\mathbf{FBoss}(X_1)=\mathbb{C}14$, \cr $e_1 e_2 e_3 e_4 = 2 r 1 l$} & \x{$\mathbf{Type}(X_3)=\mathbb{B}$, $\mathbf{FBoss}(X_3)=\mathbb{C}14$, \cr $e_1 e_2 e_3 e_4 = l 1 r 2$}  & 4 & $9$ (правая) \cr \hline
  \end{tabular}
\label{tableC1bb}
\end{table}

\medskip

Ясно, что указанные условия на буквы полностью определяют конфигурацию пути, а также его код, с точностью до подклееных окружений. Иначе говоря, не бывает никакого другого пути с заданными условиями на буквы, кроме пути, указанного нами.

То есть, для заданной конфигурации есть конечное множество слов $Q$, которые могли бы кодировать такой путь. Два слова в этом множестве отличаются только кодами подклееных окружений для крайних вершинных букв.

\medskip

{\bf Восстановление кода.}
Зафиксируем некоторый случай расположения $n$, один из четырех. Допустим, мы знаем коды трех из четырех вершин (из числа $Z$, $Y$, $J$, $F$). Тогда, учитывая характеризацию, мы знаем с каким из случаев расположения мы имеем дело. Далее, окружение вершины $H_1$ мы вычисляем, так как она первый начальник $Z$ и $F$. В первом и втором случаях расположения знания окружения $H_1$ достаточно, чтобы определить окружение оставшейся вершины. Для третьего и четвертого случаев, мы знаем также окружение $H_2$ и $H_3$ -- второй и третий начальник $F$ и $Z$. Этого достаточно, чтобы определить оставшееся окружение.

Теперь покажем, как определить начальников.
У $F$ и $Z$ общее множество начальников. Первым начальником $Y$ всегда является $Z$, а тип второго очевиден из расположения. В  случаях расположения $1$, $3$, $4$ начальники $J$ лежат в той же цепи, что и $F$.

В случае $2$, первый начальник $J$ соответствует $1$-цепи вокруг $H_3$. Тип $H_3$ мы знаем, так как это второй начальник $F$ и $Z$. Указатель можно установить по окружению макроплитки, $\mathbf{FBoss}(Z)$. Тип второго начальника $J$ -- это $\mathbb{B}$.

Таким образом, зная код пути $ZYJ$, можно вычислить код пути $ZFJ$, и наоборот. А также зная код пути $FZY$, можно вычислить код пути $FJY$, и наоборот.

\medskip

Пусть слово $W$ подходит под один из перечисленных случаев. Тогда мы
знаем с каким случаем мы имеем дело, а также какую конфигурацию имеет путь. Слово $W$ должно кодировать этот путь, то есть $W$ должно входить в $Q$.  Если слово $W$ не удовлетворяет этому условию, например, если код средней вершины не такой, какой должен быть исходя из конфигурации пути, то $W$ кодирует невозможный путь. В этом случае $W=0$, так как мы ввели обнуляющие соотношения для всех достаточно коротких невозможных путей.

\medskip
Если же $W$ входит в $Q$, то к $W$ можно применить соотношение введенное нами, из числа указанных справа на рисунке~\ref{C1b}, и получить код для другой части локального преобразования.
Таким образом, мы можем осуществить локальное преобразование пути через операцию с его кодом.

\medskip

\subsection{Случай цепи $\mathbb{C}1$; преобразования 7 и 10}

\begin{figure}[hbtp]
\centering
\leftskip=-1.0cm
\includegraphics[width=1\textwidth]{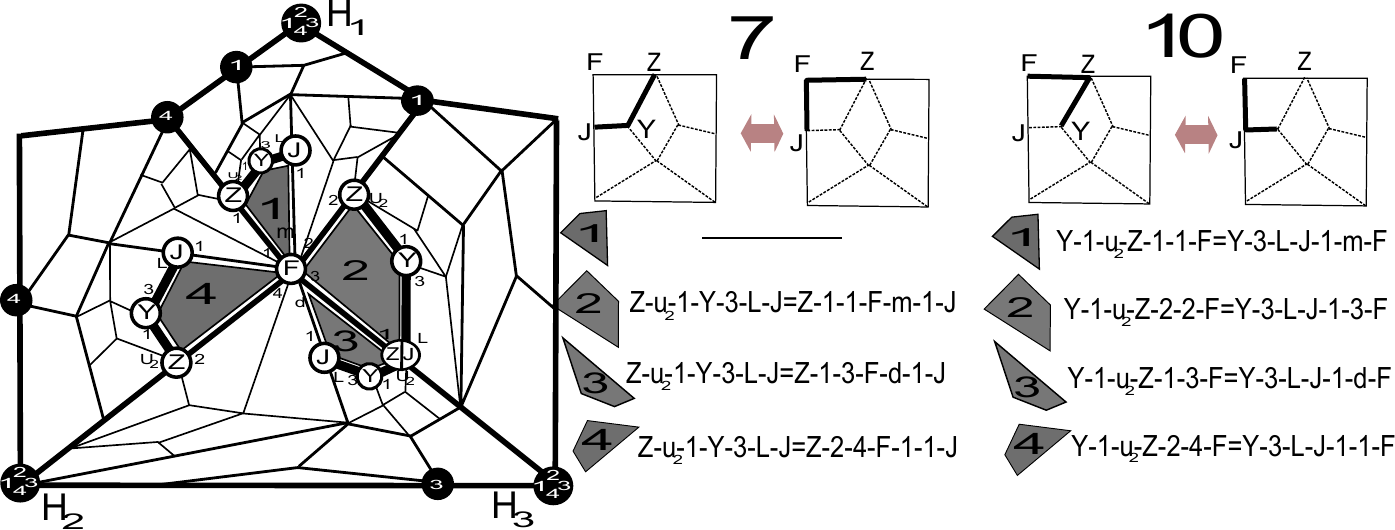}
\caption{Случаи расположения пути вокруг C1-цепи и соответствующие им определяющие соотношения для локальных преобразований 7 и 10}
\label{C1a}
\end{figure}

В правой части рисунка~\ref{C1a} изображены локальные преобразования $7$ и $10$.

\medskip

Итак, зафиксируем $\mathbb{C}$-узел с некоторым окружением и начальниками. Различные случаи расположения интересующих нас путей показаны в левой части рисунка~\ref{C1a}.

\medskip

Зная окружение и начальников центрального $\mathbb{C}$-узла, можно вычислить базовые окружения и начальников вершин $Z$, $Y$, $F$, $J$ во всех случаях расположения.
То есть мы можем определить коды вершин $Z$, $Y$, $F$, $J$ во всех четырех случаях расположения и выписать определяющие соотношения для локальных преобразований $7$ и $10$.

Уровни боковых вершин (то есть $J$ и $Z$) во всех случаях равны $1$.

\medskip

{\bf Соотношения для  преобразования 7.}

\medskip

{\bf Назначение $Z$}. Для каждого из четырех случаев расположения, обозначим символом $Z$ буквы в алфавите, соответствующие кодам вершин $Z$ с заданными типом, уровнем, базовым окружением и информацией. То есть эти буквы отличаются друг от друга только различными подклееными окружениями.

\medskip

{\bf Назначение $J$}. Аналогично определим буквы $J$, как буквы в алфавите, соответствующие кодам вершин $J$ с заданными типом, уровнем, базовым окружением и информацией, и произвольным подклееным окружением.

\medskip

{\bf Назначение $Y$ и $F$}. Буквы $Y$ и $F$ мы определим как конкретные буквы в алфавите с заданными типом, уровнем, базовым окружением и информацией, при пустом подклееном окружении. Мы это делаем, так как путь проходит через $F$ и $Y$ по плоским ребрам.

Теперь для каждого случая и для каждого разрешенного набора букв $Z$, $F$, $Y$, $J$ выпишем соотношение, представленное в правой части рисунка~\ref{C1a}.
Буквы $Z$, $F$, $Y$, $J$ мы только что определили, а остальные символы в соотношении обозначают входящие и выходящие ребра.

\medskip

{\bf Соотношения для преобразования 10.}
Тут уже $Z$ и $J$ определены однозначно, а $Y$ и $F$ с точностью для подклееного окружения.
В остальном, все аналогично.

\medskip

{\bf Оценка числа соотношений.}
Мы вводим не более $7FP^2\mathbf{Num}(\mathbb{C})$ соотношений, где $F$ -- число различных флагов макроплиток, $P$ -- число различных подклееных окружений, $\mathbf{Num}(\mathbb{C})$ -- число вершин типа $\mathbb{C}$.

В последующих случаях отношения считаются аналогично.

\medskip

{\bf Характеризация.}
Пусть есть слово $W$ представляющее собой код пути $X_1 e_1 e_2 X_2 e_3 e_4 X_3$. Покажем, как по нему установить, имеем ли мы дело с локальным преобразованием $7$ или $10$, а также как его провести.

На рисунке~\ref{C1a} отмечены входящие и выходящие ребра в каждом из четырех случаев. Исходя из этого легко выписать свойства пути, позволяющие нам установить, с каким случаем мы имеем дело. Все показано в таблицах~\ref{tableC1aa} и~\ref{tableC1ab}.

\begin{table}[hbtp]
\caption{Характеристические условия на принадлежность пути к случаю $\mathbb{C}1$ цепи. }
\centering
\begin{tabular}{|c|c|c|c|}   \hline
Условие на буквы  & \x{Симметричное условие \cr (проход в обратном порядке)} & \x{ случай}  & \x{Локальное \cr преобр} \cr \hline
$\mathbf{Surr}(X_1)=\mathbb{C}11$, $e_1 e_2 e_3 e_4 = u_1 1 3 l$  & $\mathbf{Surr}(X_3)=\mathbb{C}11$, $e_1 e_2 e_3 e_4 = l 3 1 u_1$ & 1  & $7$ (левая) \cr \hline
$\mathbf{Surr}(X_1)=\mathbb{C}11$, $e_1 e_2 e_3 e_4 = 1 1 m 1$ & $\mathbf{Surr}(X_3)=\mathbb{C}11$, $e_1 e_2 e_3 e_4 = 1 m 1 1$  & 1 & $7$ (правая) \cr \hline
$\mathbf{Surr}(X_1)=\mathbb{C}12$, $e_1 e_2 e_3 e_4 = u_2 1 3 l$ & $\mathbf{Surr}(X_3)=\mathbb{C}12$, $e_1 e_2 e_3 e_4 = l 3 1 u_2$  & 2  & $7$ (левая) \cr \hline
$\mathbf{Surr}(X_1)=\mathbb{C}12$, $e_1 e_2 e_3 e_4 = 2 2 3 1$ & $\mathbf{Surr}(X_3)=\mathbb{C}12$, $e_1 e_2 e_3 e_4 = 1 3 2 2$  & 2 & $7$ (правая)\cr \hline
$\mathbf{Surr}(X_1)=\mathbb{C}13$, $e_1 e_2 e_3 e_4 = u_2 1 3 l$ & $\mathbf{Surr}(X_3)=\mathbb{C}13$, $e_1 e_2 e_3 e_4 = l 3 1 u_2$  & 3  & $7$ (левая) \cr \hline
$\mathbf{Surr}(X_1)=\mathbb{C}13$, $e_1 e_2 e_3 e_4 = 1 3 d 1$ & $\mathbf{Surr}(X_3)=\mathbb{C}13$, $e_1 e_2 e_3 e_4 = 1 d 3 1$  & 3 & $7$ (правая)\cr \hline
$\mathbf{Surr}(X_1)=\mathbb{C}14$, $e_1 e_2 e_3 e_4 = u_2 1 3 l$ & $\mathbf{Surr}(X_3)=\mathbb{C}14$, $e_1 e_2 e_3 e_4 = l 3 1 u_2$  & 4  & $7$ (левая) \cr \hline
$\mathbf{Surr}(X_1)=\mathbb{C}14$, $e_1 e_2 e_3 e_4 = 2 4 1 1$ & $\mathbf{Surr}(X_3)=\mathbb{C}14$, $e_1 e_2 e_3 e_4 = 1 1 4 2$  & 4 & $7$ (правая)\cr \hline
  \end{tabular}
\label{tableC1aa}
\end{table}

\begin{table}[hbtp]
\caption{Характеристические условия на принадлежность пути к случаю $\mathbb{C}1$ цепи.}
\centering
\begin{tabular}{|c|c|c|c|}   \hline
Условие на буквы  & \x{Симметричное условие \cr (проход в обратном порядке)} & \x{ случай}  & \x{Локальное \cr преобр} \cr \hline
\x{ $\mathbf{Type}(X_1)=\mathbb{B}$, $\mathbf{FBoss}(X_1)=\mathbb{C}11$, \cr $e_1 e_2 e_3 e_4 = 1 u_1 1 1$}  & \x{$\mathbf{Type}(X_3)=\mathbb{A}$,  $\mathbf{FBoss}(X_3)=\mathbb{C}11$, \cr $e_1 e_2 e_3 e_4 = 1 1 u_1 1$} & 1  & $10$ (левая) \cr \hline
\x{ $\mathbf{Type}(X_1)=\mathbb{A}$, $\mathbf{FBoss}(X_1)=\mathbb{C}11$, \cr $e_1 e_2 e_3 e_4 = 3 l 1 m$} & \x{ $\mathbf{Type}(X_3)=\mathbb{A}$, $\mathbf{FBoss}(X_3)=\mathbb{C}11$, \cr $e_1 e_2 e_3 e_4 = m 1 l 3$}  & 1 & $10$ (правая) \cr \hline
\x{ $\mathbf{Type}(X_1)=\mathbb{A}$, $\mathbf{FBoss}(X_1)=\mathbb{C}12$, \cr $e_1 e_2 e_3 e_4 = 1 u_2 2 2$}  & \x{$\mathbf{Type}(X_3)=\mathbb{A}$, $\mathbf{FBoss}(X_3)=\mathbb{C}12$, \cr $e_1 e_2 e_3 e_4 = 2 2 u_2 1$} & 2  & $10$ (левая) \cr \hline
\x{ $\mathbf{Type}(X_1)=\mathbb{A}$, $\mathbf{FBoss}(X_1)=\mathbb{C}12$, \cr $e_1 e_2 e_3 e_4 = 3 l 1 3$} & \x{$\mathbf{Type}(X_3)=\mathbb{A}$, $\mathbf{FBoss}(X_3)=\mathbb{C}12$, \cr $e_1 e_2 e_3 e_4 = 3 1 l 3$}  & 2 & $10$ (правая) \cr \hline
\x{ $\mathbf{Type}(X_1)=\mathbb{A}$, $\mathbf{FBoss}(X_1)=\mathbb{C}13$, \cr $e_1 e_2 e_3 e_4 = 1 u_2 1 3$}  & \x{$\mathbf{Type}(X_3)=\mathbb{A}$, $\mathbf{FBoss}(X_3)=\mathbb{C}13$, \cr $e_1 e_2 e_3 e_4 = 3 1 u_2 1$} & 3  & $10$ (левая) \cr \hline
\x{ $\mathbf{Type}(X_1)=\mathbb{A}$, $\mathbf{FBoss}(X_1)=\mathbb{C}13$, \cr $e_1 e_2 e_3 e_4 = 3 l 1 d$} & \x{$\mathbf{Type}(X_3)=\mathbb{A}$, $\mathbf{FBoss}(X_3)=\mathbb{C}13$, \cr $e_1 e_2 e_3 e_4 = d 1 l 3$}  & 3 & $10$ (правая) \cr \hline
\x{$\mathbf{Type}(X_1)=\mathbb{A}$, $\mathbf{FBoss}(X_1)=\mathbb{C}14$, \cr $e_1 e_2 e_3 e_4 = 1 u_2 2 4$}  & \x{$\mathbf{Type}(X_3)=\mathbb{A}$, $\mathbf{FBoss}(X_3)=\mathbb{C}14$, \cr $e_1 e_2 e_3 e_4 = 4 2 u_2 1$} & 4  & $10$ (левая) \cr \hline
\x{$\mathbf{Type}(X_1)=\mathbb{B}$, $\mathbf{FBoss}(X_1)=\mathbb{C}14$, \cr $e_1 e_2 e_3 e_4 = 3 l 1 ld$} & \x{$\mathbf{Type}(X_3)=\mathbb{A}$, $\mathbf{FBoss}(X_3)=\mathbb{C}14$, \cr $e_1 e_2 e_3 e_4 = ld 1 l 3$}  & 4 & $10$ (правая) \cr \hline
  \end{tabular}
 \label{tableC1ab}
\end{table}

Ясно, что указанные условия на буквы полностью определяют конфигурацию пути, а также его код, с точностью до подклееных окружений. Иначе говоря, не бывает никакого другого пути с заданными условиями на буквы, кроме пути, указанного нами.

То есть, для заданной конфигурации есть конечное множество слов $Q$, которые могли бы кодировать такой путь. Два слова в этом множестве отличаются только кодами подклееных окружений для крайних вершинных букв.

\medskip

{\bf Восстановление кода.}
Зная код $F$, все остальные коды легко вычисляются. То есть остается разобраться, как определить код $F$ по известным кодам остальных вершин, в $7$ локальном преобразовании,  случаи расположения $2$, $3$, $4$.
Окружение $F$ совпадает с $U$-частью окружения $H_1$ первого начальника $Z$. Заметим, что $Z$ в случае $3$ и $4$, и $J$ в случае $2$ имеют тех же начальников, что и $F$.

Таким образом, зная код пути $ZYJ$, можно вычислить код пути $ZFJ$, и наоборот. А также зная код пути $FZY$, можно вычислить код пути $FJY$, и наоборот.

\medskip

Пусть слово $W$ подходит под один из перечисленных случаев.
Аналогично предыдущему параграфу,  мы можем определить с каким случаем мы имеем дело, а также какую конфигурацию имеет путь. Слово $W$ должно кодировать этот путь, то есть $W$ должно входить в $Q$. В противном случае $W$ приводится к нулю.

\medskip
Если же $W$ входит в $Q$, то к $W$ можно применить соотношение введенное нами, из числа указанных справа на рисунке~\ref{C1a}, и получить код для другой части локального преобразования.

Таким образом, мы можем осуществить локальное преобразование пути через операцию с его кодом.

\medskip

\subsection{Случай цепи $\mathbb{C}2$; преобразования 8 и 9}

\begin{figure}[hbtp]
\centering
\includegraphics[width=1\textwidth]{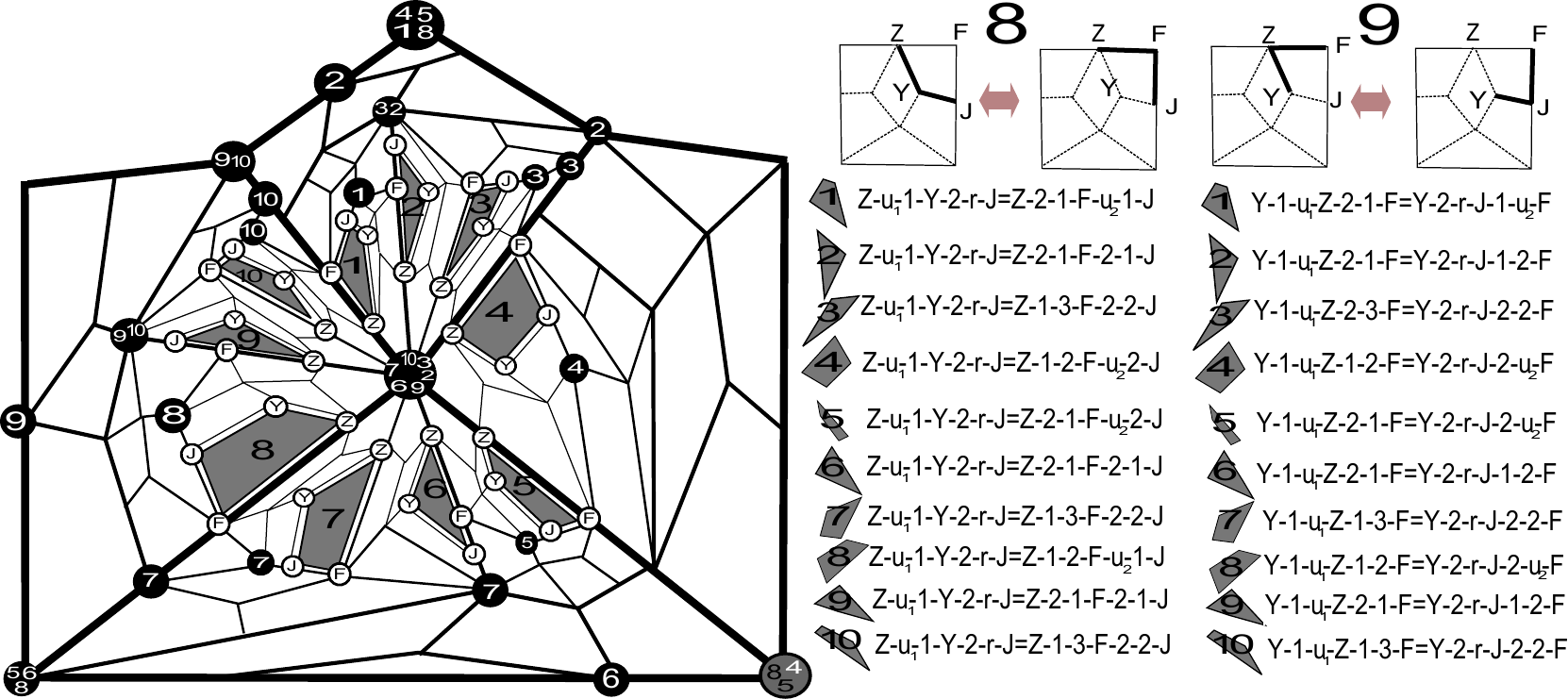}
\caption{Случаи расположения пути вокруг C2-цепи и соответствующие им определяющие соотношения для локальных преобразований 8 и 9}
\label{C2b}
\end{figure}

В правой части рисунка~\ref{C2b} изображены локальные преобразования $8$ и $9$.

\medskip

Аналогично предыдущим случаям, мы можем ввести определяющие соотношения. Фиксируем $\mathbb{C}$-узел с базовым окружением и тремя начальниками.
Ясно, что зная код центрального узла, типы, уровни и окружения вершин $Z$, $Y$, $F$, $J$ можно вычислить (с точностью до подклееными окружениями крайних вершин в пути).

Черными кругами отметим вершины, являющиеся начальниками вершин  $Z$, $Y$, $F$, $J$. Заметим, что зная код центрального $\mathbb{C}$-узла мы можем найти их типы, уровни и окружения. То есть коды вершин $Z$, $Y$, $F$, $J$ во всех десяти случаях мы можем назвать явно, с точностью до подклееных окружений.

Это позволяет ввести определяющие отношения, записанные в правой части рисунка~\ref{C2b}.

\medskip

{\bf Характеризация.}
Пусть есть слово $W$ представляющее собой код пути $X_1 e_1 e_2 X_2 e_3 e_4 X_3$.
Аналогично предыдущим случаям, мы можем выписать характеризующую последовательность кодов вершин и входящих-выходящих ребер. Это позволит для любого слова определить, является ли он кодом локального преобразования $8$ или $9$.

Мы не будем приводить здесь таблицу как для случая $\mathbb{C}1$-цепи, ее можно составить полностью аналогично, опираясь на рисунок~\ref{C2b}.

\medskip

{\bf Восстановление кода.} Учитывая характеризацию, мы знаем, с каким случаем расположения имеем дело. Кроме того, ясно, что окружение каждой из четырех вершин можно определить во всех десяти случаях. Заметим, что у $Z$ и $F$ одно и то же множество начальников во всех случаях. Первым начальником $Y$ всегда является $Z$, а тип второго очевиден из расположения. Остается определить начальников $J$. В  случаях $1$, $4$, $5$, $8$ это просто $F$, а в остальных случаях у $J$ либо те же начальники что и у $F$, либо один начальник, совпадающий с первым начальником $F$.

Таким образом, зная код пути $ZYJ$, можно вычислить код пути $ZFJ$, и наоборот. А также зная код пути $FZY$, можно вычислить код пути $FJY$, и наоборот.

\medskip

Таким образом, если слово $W$ подходит под один из перечисленных случаев то к нему можно применить соотношение введенное нами, из числа указанных справа на рисунке~\ref{C2b}, и получить код для другой части локального преобразования.

\medskip

\subsection{Случай цепи $\mathbb{C}2$; преобразования 7 и 10}

\begin{figure}[hbtp]
\centering
\includegraphics[width=1\textwidth]{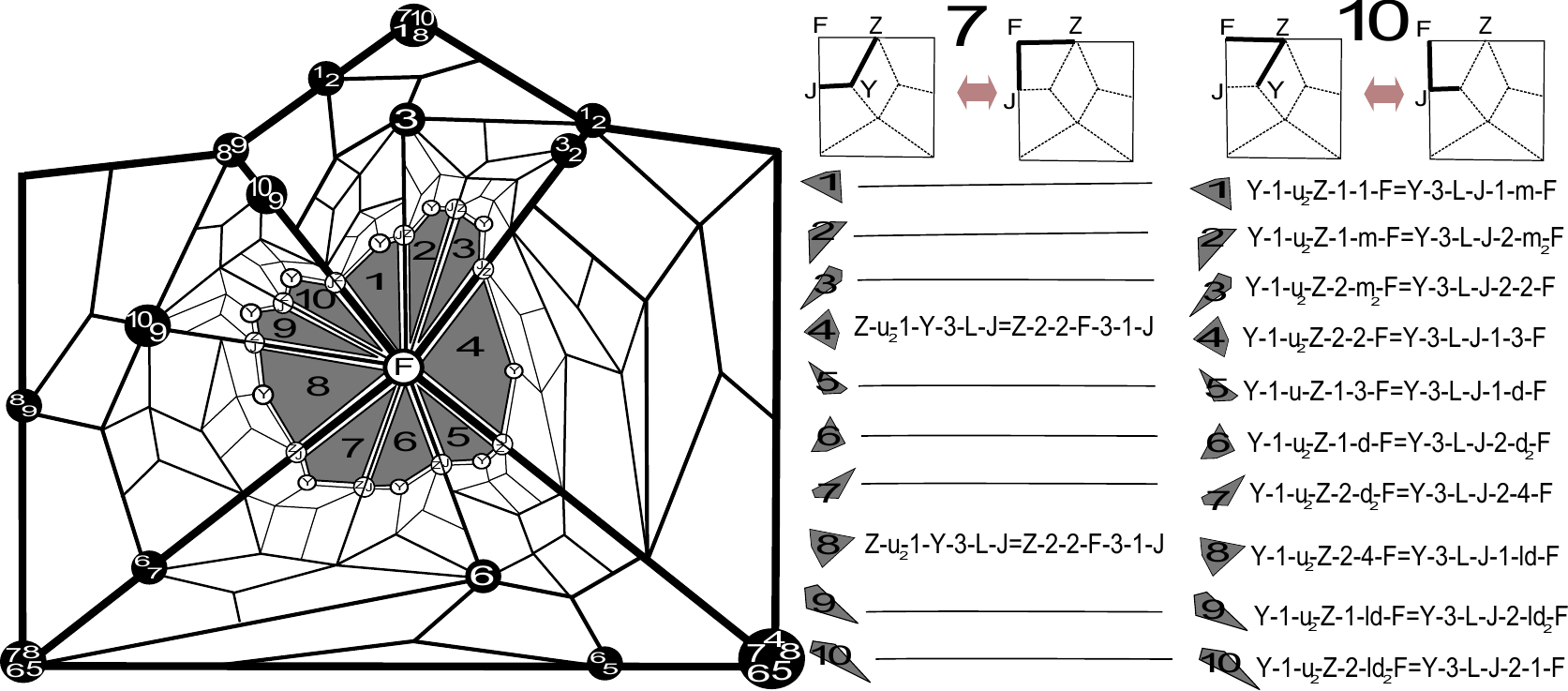}
\caption{Случаи расположения пути вокруг C2-цепи и соответствующие им определяющие соотношения для локальных преобразований 7 и 10}
\label{C2a}
\end{figure}

В правой части рисунка~\ref{C2a} изображены локальные преобразования $7$ и $10$.

\medskip

Черными кругами отметим вершины, являющиеся начальниками вершин  $Z$, $Y$, $F$, $J$. Заметим, что зная код центрального $\mathbb{C}$-узла мы можем найти их окружения. То есть коды вершин $Z$, $Y$, $F$, $J$ (включающие типы, уровни, окружения и информацию), во всех десяти случаях мы можем назвать явно, с точностью до подклееных окружений.
Это позволяет ввести определяющие отношения, записанные в правой части рисунка~\ref{C2a}.

\medskip

{\bf Характеризация.} Составление таблиц полностью аналогично случаю $\mathbb{C}1$-цепи.

\medskip

{\bf Восстановление кода.}
Рассмотрим случаи расположения $4$ и $8$ для преобразования $7$. Заметим, что $Z$ в  случае $8$ и $J$ в случае $4$ имеют тех же начальников, что и $F$. Окружение $F$ совпадает с $U$-частью его первого начальника. Значит, зная коды $Z$, $J$, $Y$, можно вычислить код $F$. Зная же код $F$, очевидно, все остальные вершины тоже вычисляются.
Для преобразования $10$, зная коды $F$, $J$, $Y$, можно вычислить $Z$, а зная $F$, $Z$, $Y$, можно вычислить $J$.

Таким образом, зная код пути $ZYJ$, можно вычислить код пути $ZFJ$, и наоборот. А также зная код пути $FZY$, можно вычислить код пути $FJY$, и наоборот.

\medskip

Таким образом, к $W$ можно применить соотношение введенное нами, из числа указанных справа на рисунке~\ref{C2a}, и получить код для другой части локального преобразования.

\medskip

\subsection{Случай цепи $\mathbb{C}3$; преобразования 7, 8, 9, 10}

Случай $\mathbb{C}3$-цепи полностью аналогичен случаю $\mathbb{C}2$-цепи, соотношения выглядят идентично, только кодировки вершин $J$, $F$, $Z$, $Y$ отвечают $\mathbb{C}3$-цепи. Все рассуждения о вычислении путей полностью аналогичны. Соотношений вводится столько же, сколько для $\mathbb{C}2$ случая, то есть $32FP^2\mathbf{Num}(\mathbb{C})$.

\medskip

\subsection{Случай цепи $\mathbb{B}1$; преобразования 8 и 9}

\begin{figure}[hbtp]
\centering
\includegraphics[width=1\textwidth]{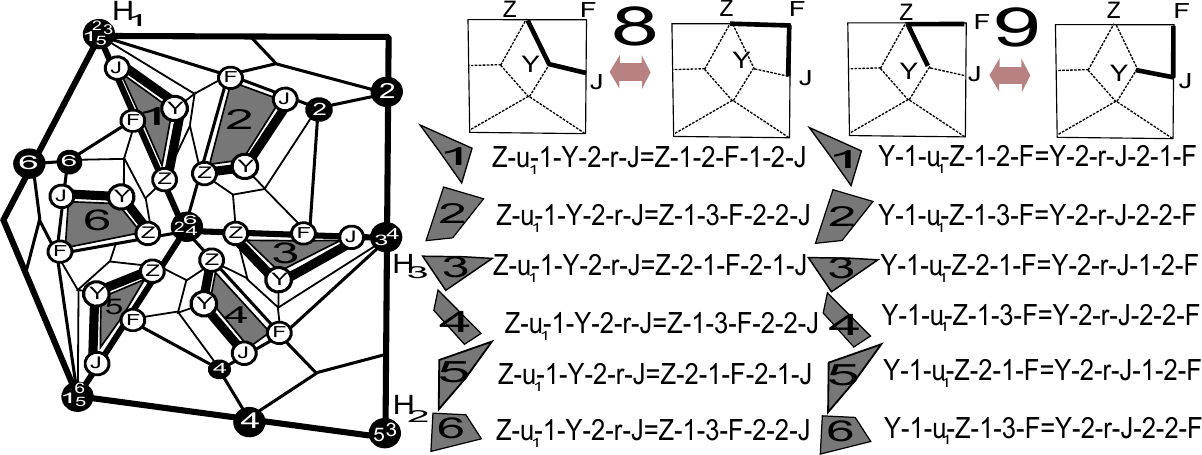}
\caption{Случаи расположения пути вокруг B1-цепи и соответствующие им определяющие соотношения для локальных преобразований 8 и 9}
\label{B1b}
\end{figure}

В правой части рисунка~\ref{B1b} изображены локальные преобразования $8$ и $9$.

\medskip

Аналогично предыдущим случаям, мы можем ввести определяющие соотношения. Фиксируем вершину типа $\mathbb{B}$ и ее начальников, вершины $H_1$ и $H_2$, причем от $H_2$ мы фиксируем только тип.
Черными кругами отметим вершины, являющиеся начальниками вершин  $Z$, $Y$, $F$, $J$.
Числа в них обозначают номер случая расположения. Заметим, что мы можем найти их окружения (кроме вершины $H_2$). Например, $H_3$ вычисляется с помощью $\mathbf{RightFromB}$. Окружения для вершин $Z$, $Y$, $F$, $J$ в каждом случае вычисляются с учетом известных нам окружения $H_1$ и типа $H_2$. То есть что коды вершин $Z$, $Y$, $F$, $J$ во всех случаях мы можем назвать явно, с точностью до подклееных окружений.
Это позволяет ввести определяющие отношения, записанные в правой части рисунка~\ref{B1b}.

\medskip

{\bf Характеризация.} Составление таблиц полностью аналогично случаю $\mathbb{C}1$-цепи.

\medskip

{\bf Восстановление кода.}
Пусть мы знаем коды трех вершин из $Z$, $J$, $Y$, $F$. Поскольку мы можем установить, с каким именно случаем расположения мы имеем дело, то и окружение оставшейся вершины мы можем вычислить.

$F$ и $Z$ имеют общее множество начальников во всех случаях. Первым начальником $Y$ всегда является $Z$, тип второго очевиден из расположения в каждом случае. В случаях~$1$, $3$, $5$ у $Y$ те же начальники что и у $F$, в случаях $2$, $4$, $6$ у $J$ первый начальник как у $F$, а второй -- как третий у $F$.

Таким образом, зная код пути $ZYJ$, можно вычислить код пути $ZFJ$, и наоборот. А также зная код пути $FZY$, можно вычислить код пути $FJY$, и наоборот.

\medskip

\subsection{Случай цепи $\mathbb{B}1$; преобразования 7 и 10}

\begin{figure}[hbtp]
\centering
\includegraphics[width=1\textwidth]{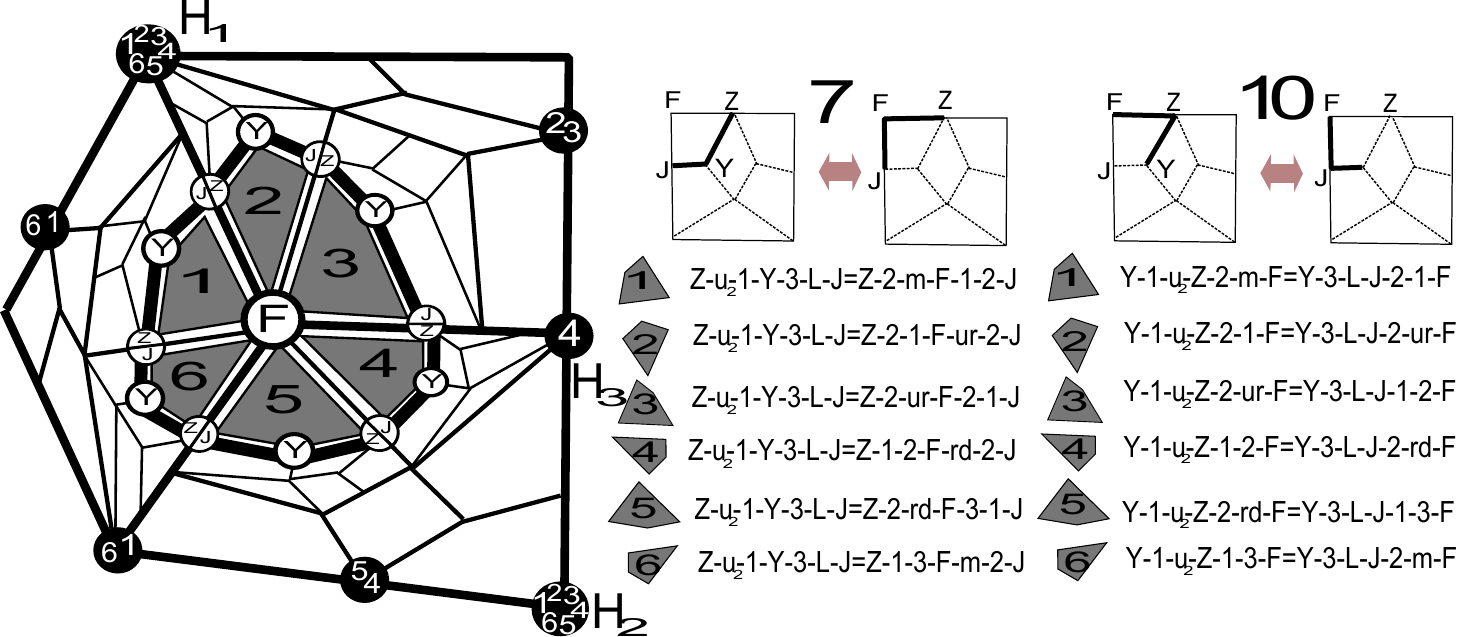}
\caption{Случаи расположения пути вокруг B1-цепи и соответствующие им определяющие соотношения для локальных преобразований 7 и 10}
\label{B1a}
\end{figure}

В правой части рисунка~\ref{B1a} изображены локальные преобразования $7$ и $10$.

\medskip

Аналогично предыдущим случаям, мы можем ввести определяющие соотношения. Фиксируем вершину типа $\mathbb{B}$ и ее начальников, вершины $H_1$ и $H_2$, причем от $H_2$ мы фиксируем только тип.
Черными кругами отметим вершины, являющиеся начальниками вершин  $Z$, $Y$, $F$, $J$.
Числа в них обозначают номер случая расположения. Заметим, что мы можем найти их типы, уровни и окружения (кроме вершины $H_2$, для нее только тип). Типы, уровни и окружения для вершин $Z$, $Y$, $F$, $J$, в каждом случае вычисляются с учетом известных нам окружения $H_1$ и типа $H_2$. То есть что коды вершин $Z$, $Y$, $F$, $J$ во всех случаях мы можем назвать явно, с точностью до подклееных окружений.
Это позволяет ввести определяющие отношения, записанные в правой части рисунка~\ref{B1a}.

\medskip

{\bf Характеризация.} Составление таблиц полностью аналогично случаю $\mathbb{C}1$-цепи.

\medskip

{\bf Восстановление кода.}
Для локального преобразования $7$, если нам нужно вычислить код $F$ по известным кодам $Z$, $J$, $Y$, можно заметить, что хотя бы у одной вершины из  $Z$, $J$, $Y$ начальником будет $F$, а у какой-либо другой начальники будут совпадать с $F$. Если же $F$ нам известно, все остальные вершины, очевидно, вычисляются.

Таким образом, зная код пути $ZYJ$, можно вычислить код пути $ZFJ$, и наоборот. А также зная код пути $FZY$, можно вычислить код пути $FJY$, и наоборот.

\medskip

\subsection{Случай цепи $\mathbb{B}2$; преобразования 8 и 9}

\begin{figure}[hbtp]
\centering
\includegraphics[width=1\textwidth]{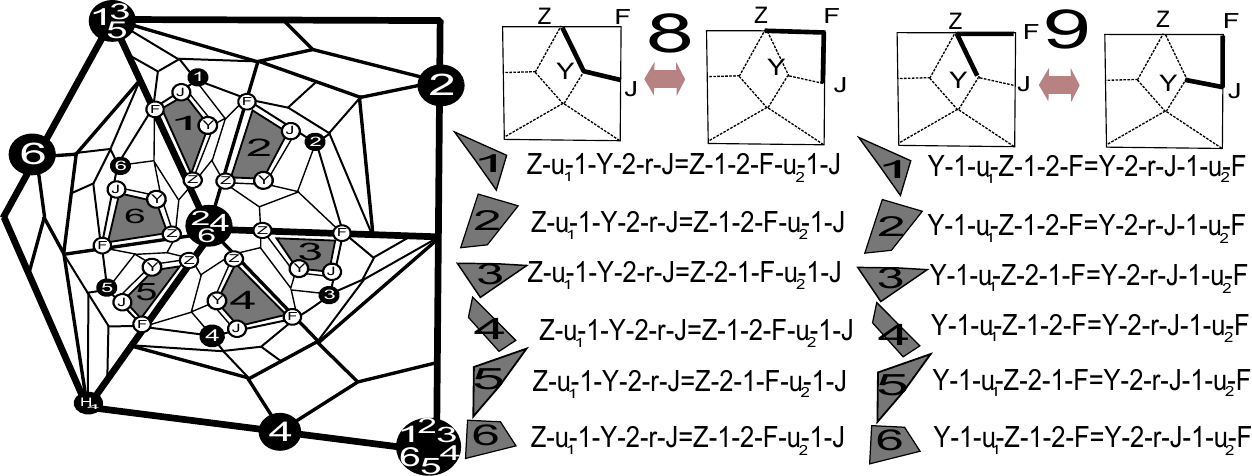}
\caption{Случаи расположения пути вокруг B2-цепи и соответствующие им определяющие соотношения для локальных преобразований 8 и 9}
\label{B2b}
\end{figure}

В правой части рисунка~\ref{B2b} изображены локальные преобразования $8$ и $9$.

\medskip

Введение определяющих соотношений полностью аналогично предыдущим случаям, можно легко проверить, что зная центральный $\mathbb{B}$-узел, коды всех четырех вершин в каждом из случаев расположения легко можно выписать. Это позволяет ввести определяющие отношения, записанные в правой части рисунка~\ref{B2b}.

\medskip

{\bf Характеризация.} Составление таблиц полностью аналогично случаю $\mathbb{C}1$-цепи.

\medskip

{\bf Восстановление кода.}
Поскольку мы можем установить, с каким именно случаем расположения мы имеем дело, то и окружение каждой из четырех вершин мы можем вычислить.
Кроме того, заметим, что каждая из четырех вершин  $Z$, $J$, $Y$, $F$ в каждом из случаев либо имеет начальников среди других трех вершин, либо у нее общие начальники с какой-то из трех других вершин. Значит, в каждом из случаев, зная коды трех вершин, можно вычислить код четвертой.

Таким образом, зная код пути $ZYJ$, можно вычислить код пути $ZFJ$, и наоборот. А также зная код пути $FZY$, можно вычислить код пути $FJY$, и наоборот.

\medskip

\subsection{Случай цепи $\mathbb{B}2$; преобразования 7 и 10}

\begin{figure}[hbtp]
\centering
\includegraphics[width=1\textwidth]{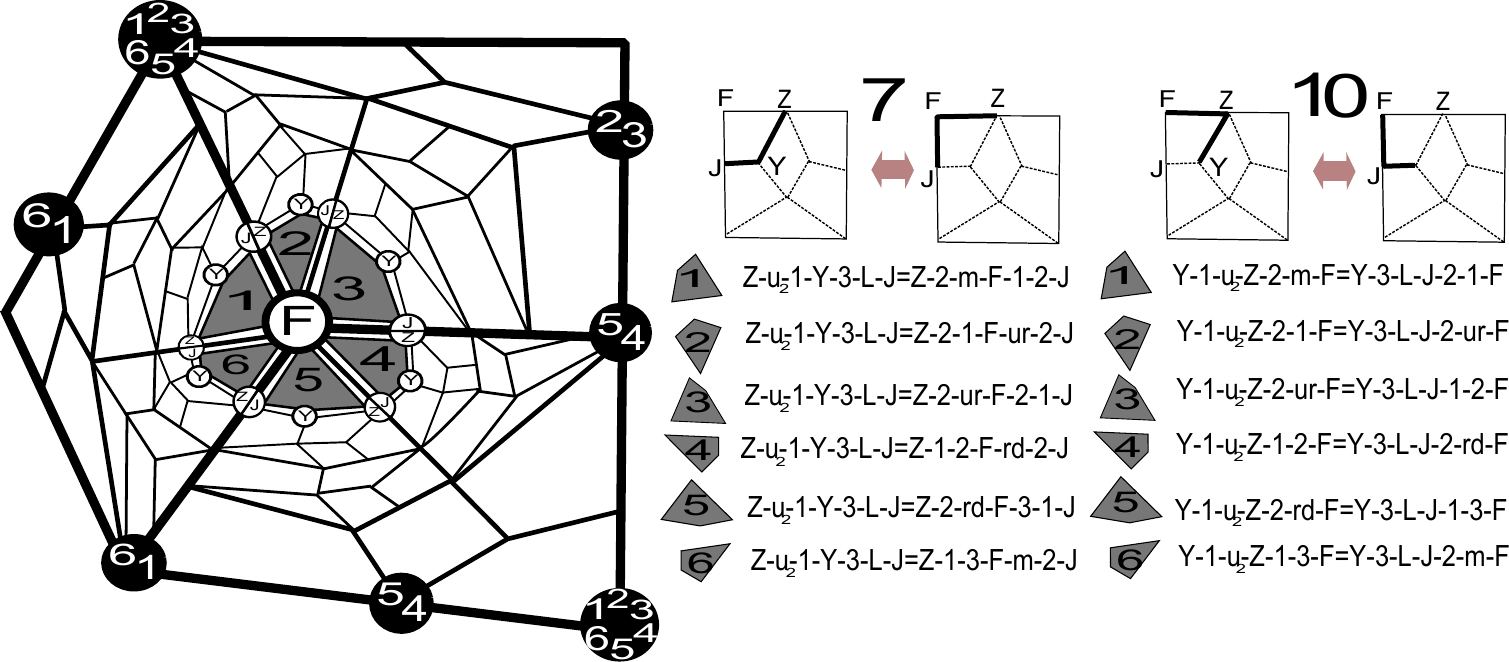}
\caption{Случаи расположения пути вокруг B2-цепи и соответствующие им определяющие соотношения для локальных преобразований 7 и 10}
\label{B2a}
\end{figure}

В правой части рисунка~\ref{B2a} изображены локальные преобразования $7$ и $10$.

\medskip

Зафиксируем вершину типа $\mathbb{B}$ и ее начальников. Аналогично предыдущим случаям, мы можем вычислить коды всех вершин и ввести определяющие соотношения.

\medskip

{\bf Характеризация.} Составление таблиц полностью аналогично случаю $\mathbb{C}1$-цепи.

\medskip

{\bf Восстановление кода.}
Для локального преобразования $7$, если нам нужно вычислить код $F$ по известным кодам $Z$, $J$, $Y$ можно заметить, что хотя бы у одной вершины из  $Z$, $J$, $Y$ начальники будут совпадать с $F$. Если же $F$ нам известно, все остальные вершины, очевидно, вычисляются.

Таким образом, зная код пути $ZYJ$, можно вычислить код пути $ZFJ$, и наоборот. А также зная код пути $FZY$, можно вычислить код пути $FJY$, и наоборот.

\medskip

\subsection{Случай цепи $\mathbb{A}0$; преобразования 8 и 9}

\begin{figure}[hbtp]
\centering
\includegraphics[width=0.8\textwidth]{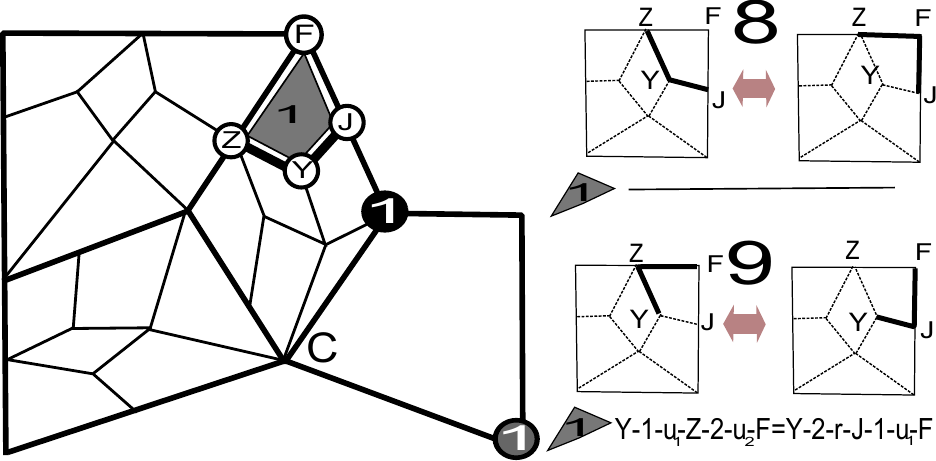}
\caption{Случаи расположения пути вокруг A0-цепи и соответствующие им определяющие соотношения для локальных преобразований 8 и 9}
\label{A1b}
\end{figure}

В правой части рисунка~\ref{A1b} изображены локальные преобразования $8$ и $9$.

\medskip

Зафиксируем вершину $X$ второго уровня с типом $\mathbb{UL}$, $\mathbb{UR}$ или $\mathbb{U}$. Мы будем вводить соотношения только для преобразования $9$. Путь в преобразовании $8$ подпадает под мертвый паттерн, для него определяющие соотношения вводить не требуется. По коду $X$ вычисляются коды остальных трех вершин, например тип второго начальника $J$ -- вершины в правом нижнем углу может быть найден с помощью функции $\mathbf{BottomRightType}(X)$.

\medskip

{\bf Характеризация.} Составление таблиц полностью аналогично случаю $\mathbb{C}1$-цепи.

\medskip

{\bf Восстановление кода.}
В преобразовании $9$, зная код $F$, мы можем вычислить и $J$ и $Z$.

\medskip

\subsection{Случай цепи $\mathbb{A}0$; преобразования 7 и 10}

\begin{figure}[hbtp]
\centering
\includegraphics[width=0.8\textwidth]{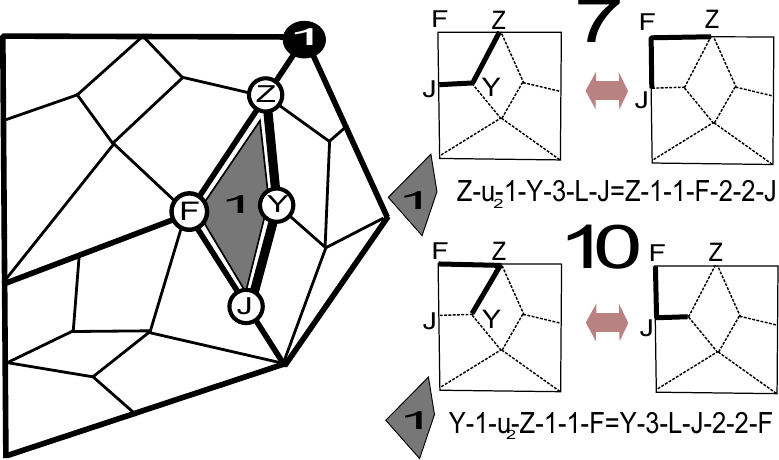}
\caption{Случаи расположения пути вокруг A0-цепи и соответствующие им определяющие соотношения для локальных преобразований 7 и 10}
\label{A1a}
\end{figure}

В правой части рисунка~\ref{A1a} изображены локальные преобразования $7$ и $10$.

\medskip

Зафиксируем вершину типа $\mathbb{A}$ и ее начальника (вершина, отмеченная черным кругом с ``$1$'', уровень которой равен $2$).
Заметим, что зная окружение этой вершины, мы можем вычислить коды вершин $Z$, $Y$, $F$, $J$ с точностью до подклееных окружений. Это позволяет ввести определяющие отношения, записанные в правой части рисунка~\ref{A1a}.

\medskip

{\bf Характеризация.} Составление таблиц полностью аналогично случаю $\mathbb{C}1$-цепи.

\medskip

{\bf Восстановление кода.}
Вершина, отмеченная черным кругом с $1$, является единственным начальником  $Z$, $F$, $J$.
Таким образом, зная код пути $ZYJ$, можно вычислить код пути $ZFJ$, и наоборот. А также зная код пути $FZY$, можно вычислить код пути $FJY$, и наоборот.

\medskip

\subsection{Случай цепи $\mathbb{A}1$; преобразования 8 и 9}

\begin{figure}[hbtp]
\centering
\includegraphics[width=1\textwidth]{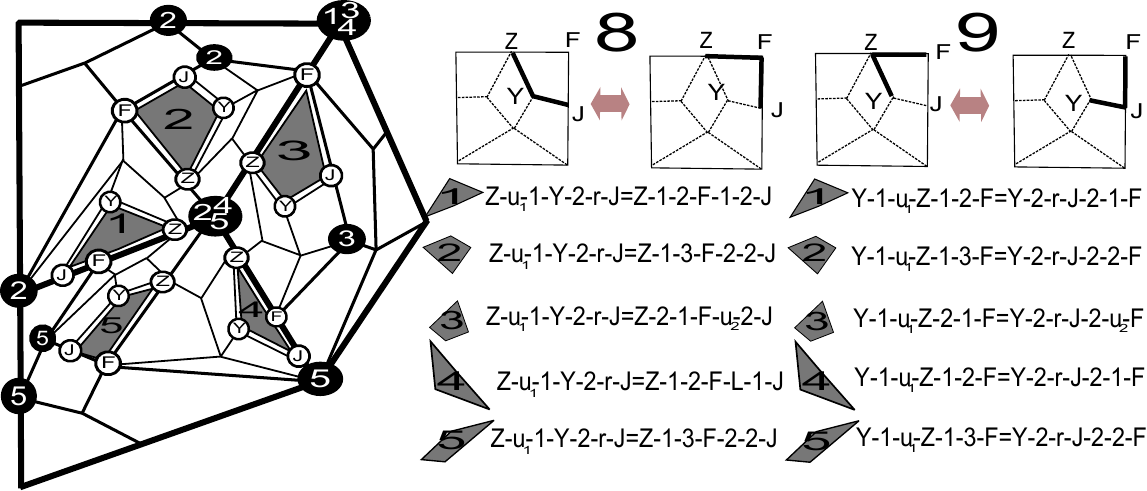}
\caption{Случаи расположения пути вокруг A1-цепи  и соответствующие им определяющие соотношения для локальных преобразований 8 и 9}
\label{A2b}
\end{figure}

В правой части рисунка~\ref{A2b} изображены локальные преобразования $8$ и $9$.

\medskip

Зафиксируем вершину типа $\mathbb{A}$ и ее начальника (третьего уровня).
Черными кругами отметим вершины, являющиеся начальниками вершин  $Z$, $Y$, $F$, $J$.
Заметим, что зная окружение центрального $\mathbb{A}$-узла, мы можем найти их типы, уровни и окружения. Так как начальники вершин  $Z$, $Y$, $F$, $J$ содержатся среди вершин, отмеченных черными кругами, то коды вершин $Z$, $Y$, $F$, $J$ во всех случаях мы можем назвать явно, с точностью до подклееных окружений. Это позволяет ввести определяющие отношения, записанные в правой части рисунка~\ref{A2b}.

\medskip

{\bf Характеризация.} Составление таблиц полностью аналогично случаю $\mathbb{C}1$-цепи.

\medskip

{\bf Восстановление кода.}
Поскольку мы можем установить, с каким именно случаем расположения мы имеем дело, то и окружение каждой из вершин мы можем вычислить, зная остальные три. Кроме того, заметим, что каждая из четырех вершин  $Z$, $J$, $Y$, $F$ в каждом из случаев либо имеет начальников среди других трех вершин, либо у нее общие начальники с какой-то из трех других вершин. Значит, в каждом из случаев, зная коды трех вершин, можно вычислить код четвертой.

Таким образом, зная код пути $ZYJ$, можно вычислить код пути $ZFJ$, и наоборот. А также зная код пути $FZY$, можно вычислить код пути $FJY$, и наоборот.

\medskip

\subsection{Случай цепи $\mathbb{A}1$; преобразования 7 и 10}

\begin{figure}[hbtp]
\centering
\includegraphics[width=1\textwidth]{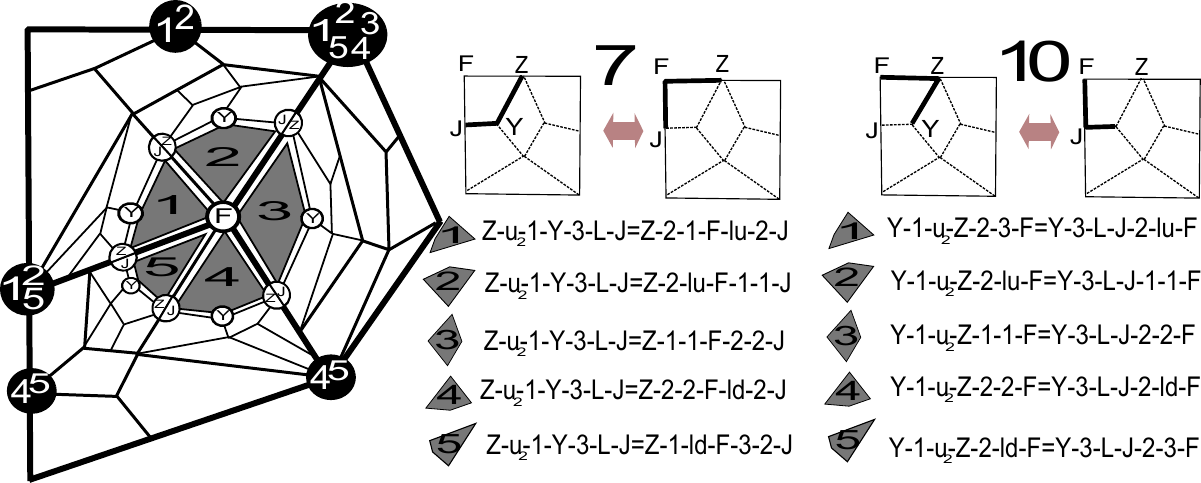}
\caption{Случаи расположения пути вокруг A1-цепи и соответствующие им определяющие соотношения для локальных преобразований 7 и 10}
\label{A2a}
\end{figure}

В правой части рисунка~\ref{A2a} изображены локальные преобразования $7$ и $10$.

\medskip

Фиксируем некоторую вершину типа $\mathbb{A}$ и ее начальника.
Черными  кругами отметим вершины, являющиеся начальниками вершин  $Z$, $Y$, $F$, $J$.
Заметим, что зная окружение этой центральной $\mathbb{A}$-вершины, мы можем вычислить коды вершин $Z$, $Y$, $F$, $J$ с точностью до подклееных окружений. Это позволяет ввести определяющие отношения, записанные в правой части рисунка~\ref{A2a}.

\medskip

{\bf Характеризация.} Составление таблиц полностью аналогично случаю $\mathbb{C}1$-цепи.

\medskip

{\bf Восстановление кода.}
Вершина $F$ имеет общего начальника с одной из вершин $J$ или $Z$. То есть, зная код пути $ZYJ$, можно вычислить код пути $ZFJ$. В обратную сторону, а также для локального преобразования 10: очевидно, что зная $F$, можно вычислить коды остальных вершин. Таким образом, мы можем осуществить локальное преобразование пути через операцию с его кодом.

\medskip

\subsection{Случай цепи $\mathbb{A}2$; преобразования 8 и 9}

\begin{figure}[hbtp]
\centering
\includegraphics[width=1\textwidth]{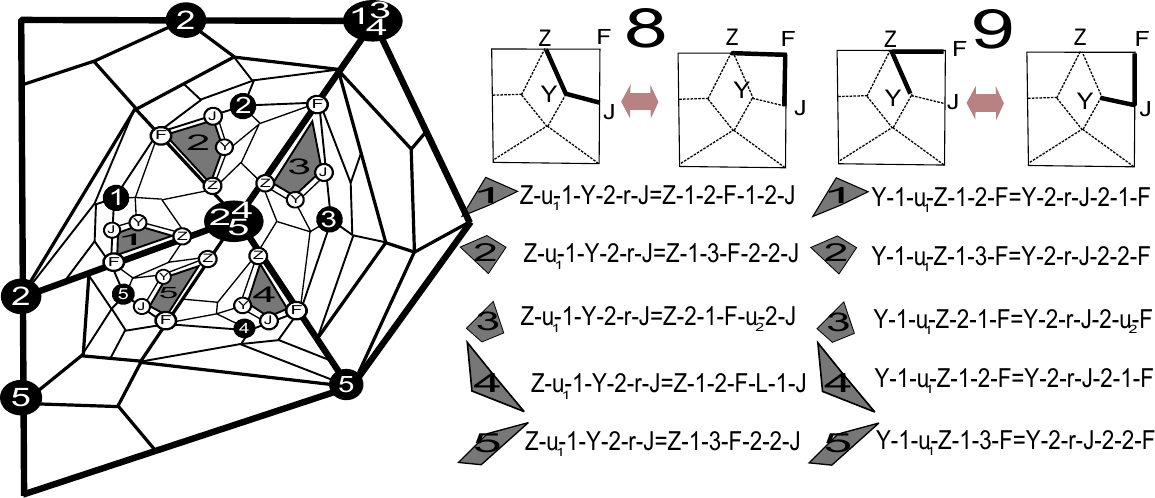}
\caption{Случаи расположения пути вокруг A2-цепи  и соответствующие им определяющие соотношения для локальных преобразований 8 и 9}
\label{A3b}
\end{figure}

В правой части рисунка~\ref{A3b} изображены локальные преобразования $8$ и $9$.

\medskip

Введение определяющих соотношений полностью аналогично случаю с $\mathbb{A}1$-цепью. Также аналогично можно установить, что зная код пути $ZYJ$, можно вычислить код пути $ZFJ$, и наоборот. А также зная код пути $FZY$, можно вычислить код пути $FJY$, и наоборот. Аналогично можно установить последовательности букв, характеризующие данную цепь и локальные преобразования.

\medskip

\subsection{Случай цепи $\mathbb{A}2$; преобразования 7 и 10}

\begin{figure}[hbtp]
\centering
\includegraphics[width=1\textwidth]{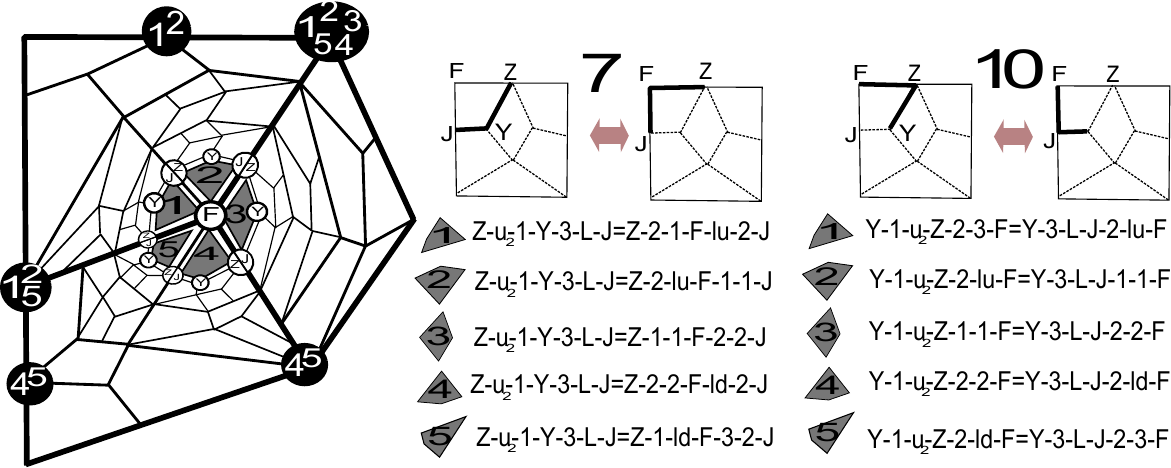}
\caption{Случаи расположения пути вокруг A2-цепи  и соответствующие им определяющие соотношения для локальных преобразований 7 и 10}
\label{A3a}
\end{figure}

В правой части рисунка~\ref{A3a} изображены локальные преобразования $7$ и $10$.

\medskip

Введение определяющих соотношений полностью аналогично случаю с $\mathbb{A}1$-цепью. Также аналогично можно установить, что зная код пути $ZYJ$, можно вычислить код пути $ZFJ$, и наоборот. А также зная код пути $FZY$, можно вычислить код пути $FJY$, и наоборот. Аналогично можно установить последовательности букв, характеризующие данную цепь и локальные преобразования.

\medskip

\subsection{Случай цепи $\mathbb{UL}1$; преобразования 8 и 9}

\begin{figure}[hbtp]
\centering
\includegraphics[width=1\textwidth]{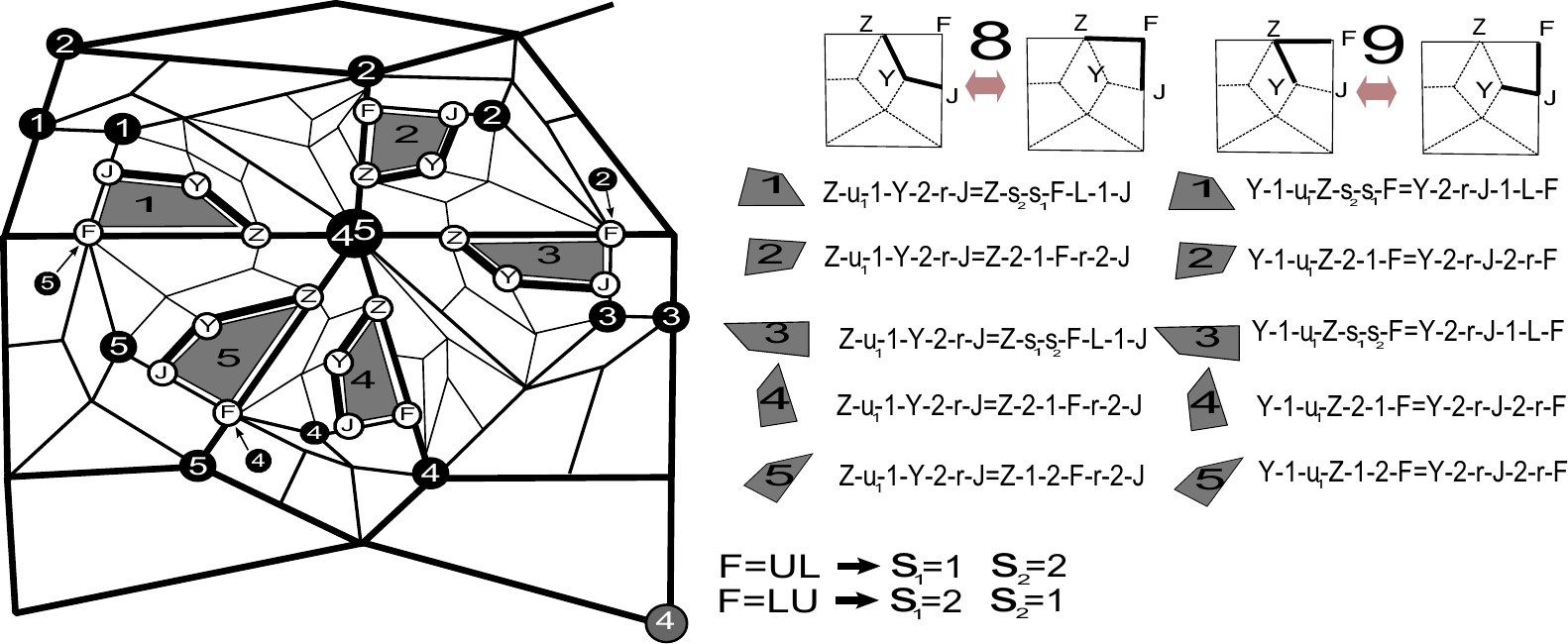}
\caption{Случаи расположения пути вокруг UL1-цепи и определяющие соотношения для  локальных преобразований 8 и 9}
\label{UL1b}
\end{figure}

В правой части рисунка~\ref{UL1b} изображены локальные преобразования $8$ и $9$.

\medskip

Зафиксируем вершину типа $\mathbb{UL}$ третьего уровня и ее начальников.
Черными кругами отметим вершины, являющиеся начальниками вершин  $Z$, $Y$, $F$, $J$. Три круга попадают в $F$-узлы других расположений, это отмечено стрелками.
Заметим, что зная окружение центрального $\mathbb{UL}$-узла, мы можем найти типы, уровни и окружения вершин с черными кругами (кроме вершины в правом нижнем углу, для нее только тип). Так как начальники вершин  $Z$, $Y$, $F$, $J$ содержатся среди вершин, отмеченных черными кругами, то коды вершин $Z$, $Y$, $F$, $J$ во всех случаях мы можем назвать явно, с точностью до подклееных окружений. Это позволяет ввести определяющие отношения, записанные в правой части рисунка~\ref{UL1b}.

\medskip

{\bf Характеризация.} Составление таблиц полностью аналогично случаю $\mathbb{C}1$-цепи.

\medskip

{\bf Восстановление кода.}
Поскольку мы можем установить, с каким именно случаем расположения мы имеем дело, то окружение оставшейся вершины мы можем вычислить, зная коды остальных трех. Кроме того, заметим, что у $Z$ и $F$ в каждом из случаев общие начальники, то есть код каждой из этих вершин может быть вычислен, исходя из кодов остальных трех вершин. У $Y$ во всех случаях первым начальником является $Z$, а тип второго в каждом случае ясен из расположения.

Начальники $J$ вычисляется следующим образом:  случай $1$ -- это $\mathbf{Prev}(F)$;  случай $2$ -- это $1$-цепь вокруг $\mathbf{TypeBottomLeft.FBoss}(Z)$, указатель $\mathbf{l}$ (первый начальник), и $\mathbb{A}$-тип второго начальника; случай $3$ -- это $\mathbf{Prev}(F)$;
 случай $4$ -- это $0$-цепь с указателем $1$ вокруг $\mathbb{A}$, окружения которого как $U$-часть $\mathbf{Fboss}(Z)$ (первый начальник), и тип второго начальника $\mathbb{B}$;
 случай $5$ -- это первый $\mathbf{Plus.FBoss}(Z)$, тип второго $\mathbb{A}$.

Таким образом, зная код пути $ZYJ$, можно вычислить код пути $ZFJ$, и наоборот. А также зная код пути $FZY$, можно вычислить код пути $FJY$, и наоборот. Таким образом, мы можем осуществить локальное преобразование пути через операцию с его кодом.

\medskip

\subsection{Случай цепи $\mathbb{UL}1$; преобразования 7 и 10}

\begin{figure}[hbtp]
\centering
\includegraphics[width=1\textwidth]{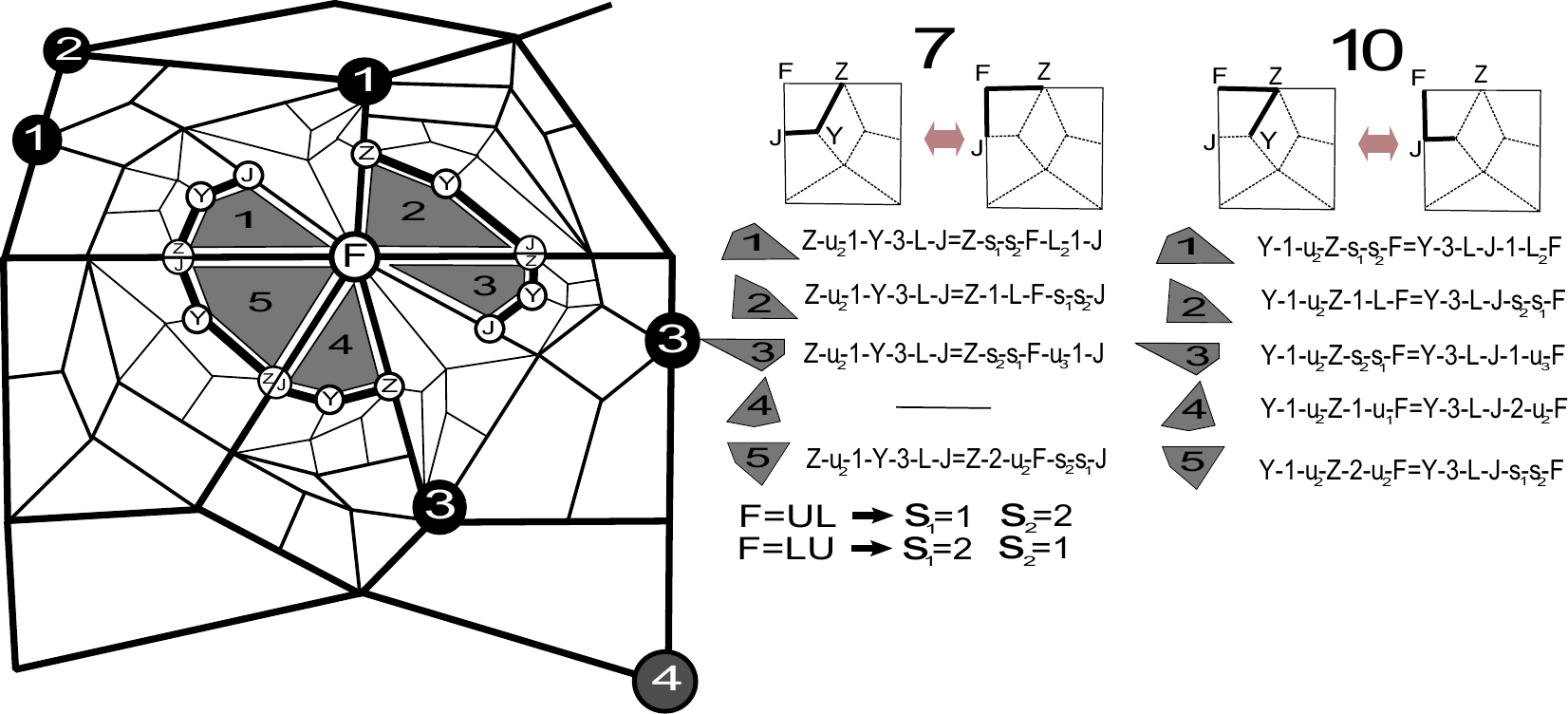}
\caption{Случаи расположения пути вокруг UL1-цепи и определяющие соотношения для локальных преобразований 7 и 10}
\label{UL1a}
\end{figure}

В правой части рисунка~\ref{UL1a} изображены локальные преобразования $7$ и $10$.

\medskip

Фиксируем некоторую вершину типа $\mathbb{UL}$ третьего уровня и ее начальников.
Черными кругами отметим вершины, являющиеся начальниками вершин  $Z$, $Y$, $F$, $J$.
Заметим, что зная окружение центральной вершины (которая совпадает с $F$), мы можем вычислить коды вершин $Z$, $Y$, $F$, $J$ с точностью до подклееных окружений во всех случаях. Это позволяет ввести определяющие отношения, записанные в правой части рисунка~\ref{UL1a}.

\medskip

{\bf Характеризация.} Составление таблиц полностью аналогично случаю $\mathbb{C}1$-цепи.

\medskip

{\bf Восстановление кода.}
Вершина $F$ имеет общего начальника с одной из вершин $J$ или $Z$, в каждом из случаев, кроме случая $4$ (а в  случае $4$ нам ее не надо вычислять).  То есть, зная код пути $ZYJ$, можно вычислить код пути $ZFJ$. В обратную сторону, а также для локального преобразования $10$: очевидно, что зная $F$, можно вычислить коды остальных вершин. Таким образом, мы можем осуществить локальное преобразование пути через операцию с его кодом.

\medskip

\subsection{Случай цепи $\mathbb{UL}2$; преобразования 8 и 9}

\begin{figure}[hbtp]
\centering
\includegraphics[width=0.9\textwidth]{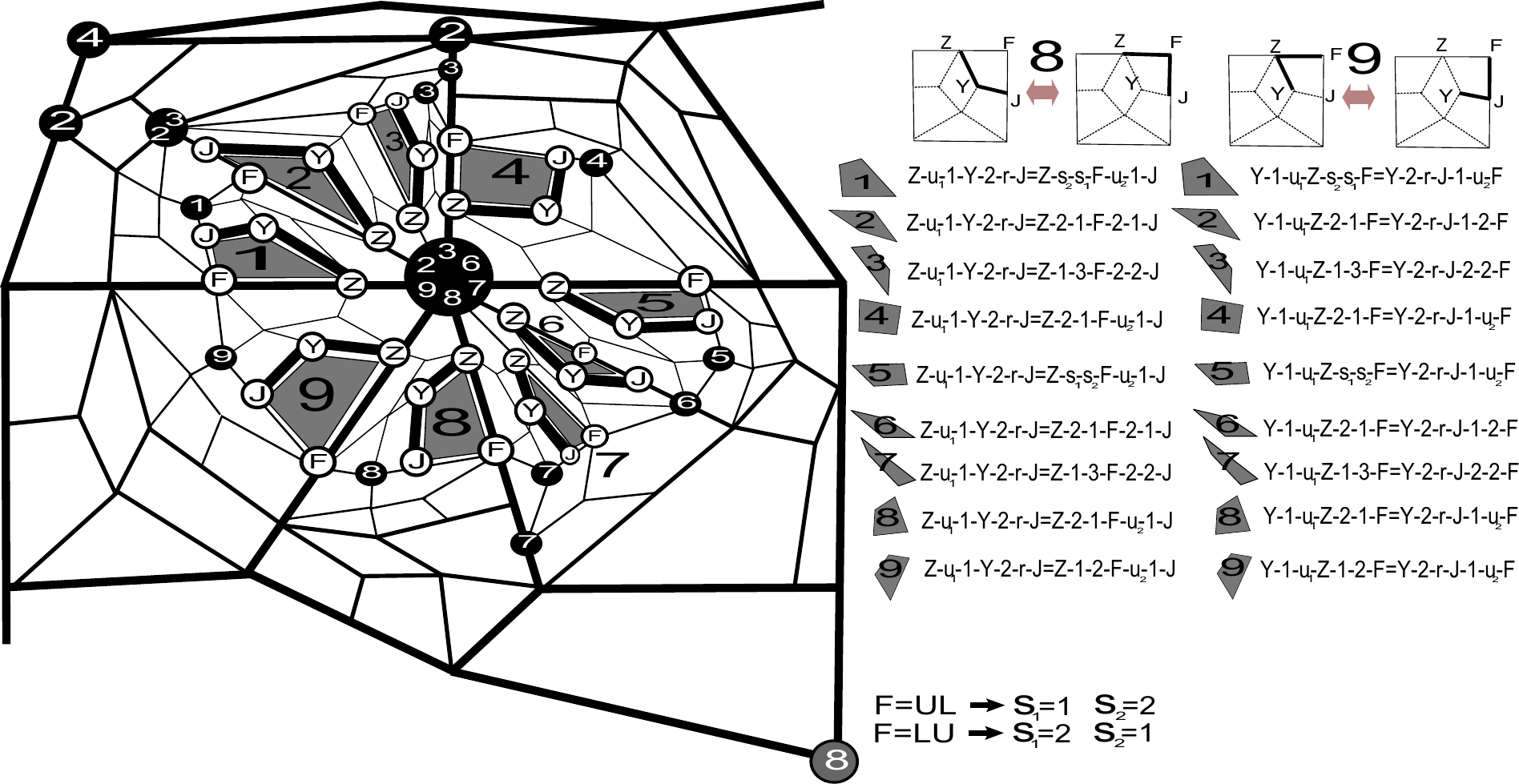}
\caption{Случаи расположения пути вокруг UL2-цепи и определяющие соотношения для локальных преобразований 8 и 9}
\label{UL2b}
\end{figure}

В правой части рисунка~\ref{UL2b} изображены локальные преобразования $8$ и $9$.

\medskip

\leftskip=-0.0cm

Зафиксируем вершину типа $\mathbb{UL}$ третьего уровня и ее начальников.
Черными кругами отметим вершины, являющиеся начальниками вершин  $Z$, $Y$, $F$, $J$.
Заметим, что зная окружение центрального $\mathbb{UL}$-узла, мы можем найти их типы, уровни и окружения. Так как начальники вершин  $Z$, $Y$, $F$, $J$ содержатся среди вершин, отмеченных черными кругами, то коды вершин $Z$, $Y$, $F$, $J$, во всех случаях мы можем назвать явно, с точностью до подклееных окружений. Это позволяет ввести определяющие отношения, записанные в правой части рисунка~\ref{UL2b}.

\medskip

{\bf Характеризация.} Составление таблиц полностью аналогично случаю $\mathbb{C}1$-цепи.

\medskip

{\bf Восстановление кода.}
Заметим, что у $Z$ и $F$ в каждом из случаев общие начальники, то есть код каждой из этих вершин может быть вычислен, исходя из кодов остальных трех вершин. $Y$ вычисляется также легко, так как во всех случаях первым начальником является $Z$, а тип второго в каждом случае легко виден из рисунка.

Начальники $J$ вычисляется следующим образом: в  случаях $1$, $4$, $5$, $8$, $9$ -- единственный начальник, узел $F$; в случаях $2$, $6$ -- начальники те же что и у $F$; в случаях  $3$ и $7$ -- первый и второй начальники как первый и третий у $F$.

Таким образом, зная код пути $ZYJ$, можно вычислить код пути $ZFJ$, и наоборот. А также зная код пути $FZY$, можно вычислить код пути $FJY$, и наоборот. Таким образом, мы можем осуществить локальное преобразование пути через операцию с его кодом.

\medskip

\subsection{Случай цепи $\mathbb{UL}2$; преобразования 7 и 10}

\begin{figure}[hbtp]
\centering
\includegraphics[width=1\textwidth]{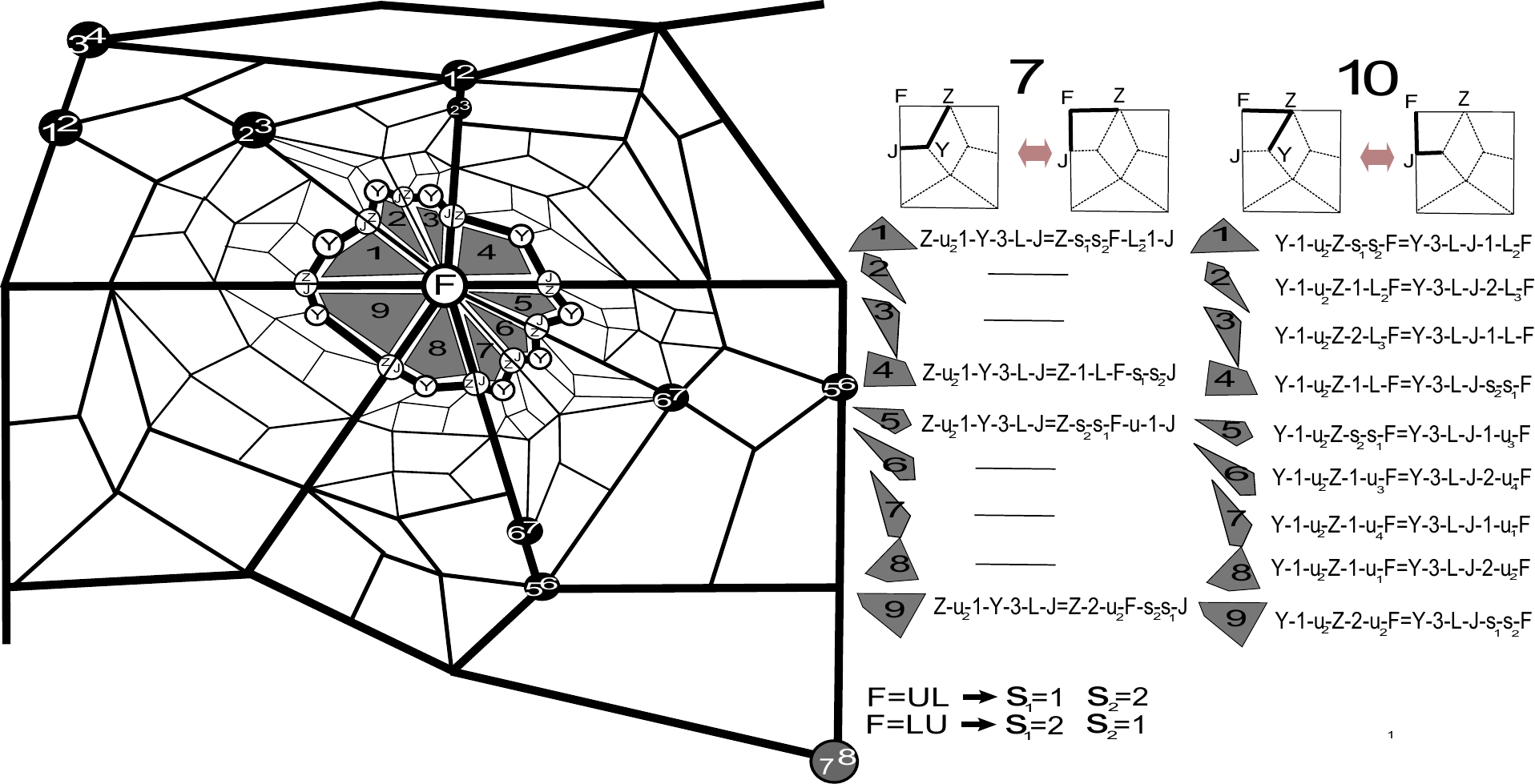}
\caption{Случаи расположения пути вокруг UL2-цепи и определяющие соотношения для локальных преобразований 7 и 10}
\label{UL2a}
\end{figure}

В правой части рисунка~\ref{UL2a} изображены локальные преобразования $7$ и $10$.

\medskip

Фиксируем некоторую вершину типа $\mathbb{UL}$ третьего уровня и ее начальников.
Черными и серыми кругами отметим вершины, являющиеся начальниками вершин  $Z$, $Y$, $F$, $J$.
Заметим, что зная окружение центральной вершины (совпадающей с $F$), мы можем вычислить коды вершин $Z$, $Y$, $F$, $J$, с точностью до подклееных окружений. Это позволяет ввести определяющие отношения, записанные в правой части рисунка~\ref{UL2a}.

Для  локального преобразования $7$ пути в случаях $2,3,6,7,8$ удовлетворяют условиям мертвого паттерна, и в этих случаях мы соотношения не вводим.

\medskip

{\bf Характеризация.} Составление таблиц полностью аналогично случаю $\mathbb{C}1$-цепи.

\medskip

{\bf Восстановление кода.}
В случаях $1,4,5,9$ вершина $F$ имеет общего начальника с одной из вершин $J$ или $Z$. То есть, зная код пути $ZYJ$, можно вычислить код пути $ZFJ$. В обратную сторону, а также для локального преобразования $10$: очевидно, что зная $F$, можно вычислить коды остальных вершин. Таким образом, мы можем осуществить локальное преобразование пути через операцию с его кодом.

\medskip

\subsection{Случай цепи $\mathbb{UL}3$; преобразования 7, 8, 9, 10}

Случай $\mathbb{UL}3$-цепи полностью аналогичен $\mathbb{UL}2$ случаю, соотношения выглядят идентично, только кодировки вершин $J$, $F$, $Z$, $Y$ отвечают $\mathbb{UL}3$-цепи. Все рассуждения о вычислении путей полностью аналогичны. Соотношений вводится столько же, сколько для $\mathbb{UL}2$ случая.

\medskip

\subsection{Случай цепи $\mathbb{UR}1$; преобразования 8 и 9}

\begin{figure}[hbtp]
\centering
\includegraphics[width=1\textwidth]{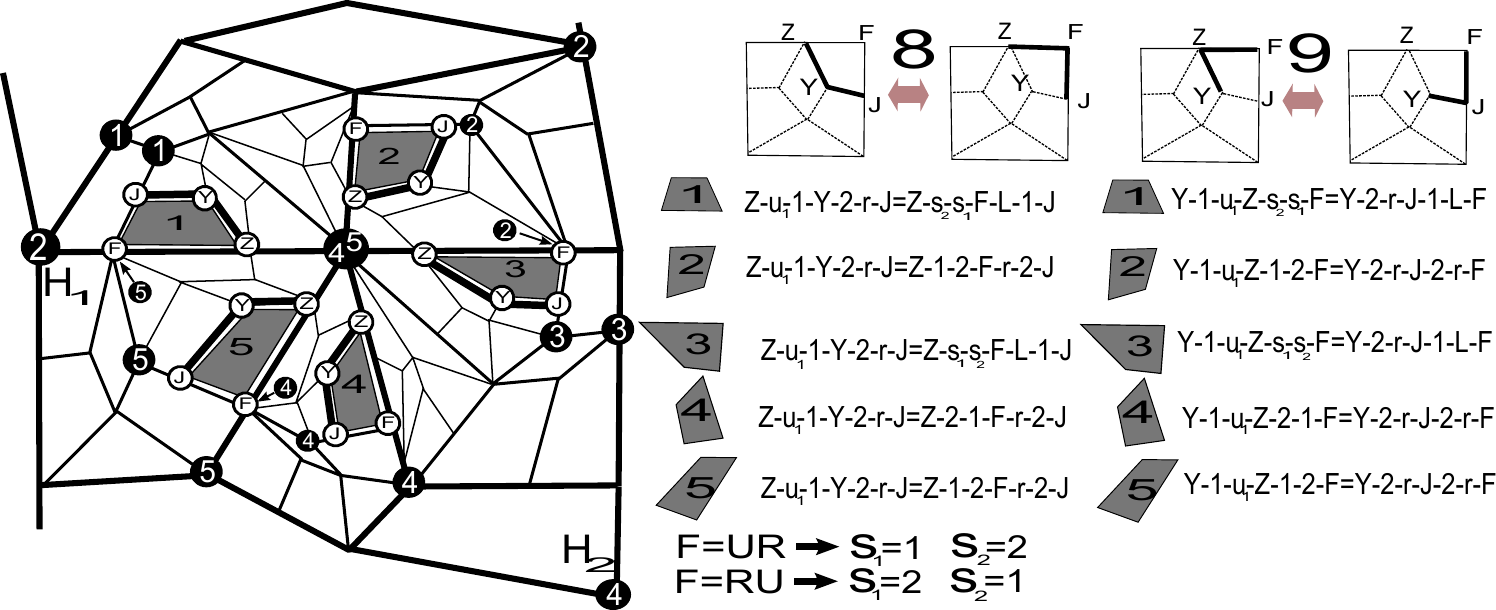}
\caption{Случаи расположения пути вокруг UR1-цепи и определяющие соотношения для локальных преобразований 8 и 9}
\label{UR1b}
\end{figure}

В правой части рисунка~\ref{UR1b} изображены локальные преобразования $8$ и $9$.

\medskip

Зафиксируем вершину типа $\mathbb{UR}$ третьего уровня и ее начальников.
Черными  кругами отметим вершины, являющиеся начальниками вершин  $Z$, $Y$, $F$, $J$. Три круга попадают в $F$-узлы других расположений, это отмечено стрелками.
Заметим, что зная окружение центрального $\mathbb{UR}$-узла, мы можем найти их типы, уровни и окружения (Кроме вершин $H_1$ и $H_2$, для которых можем найти тип). Так как начальники вершин  $Z$, $Y$, $F$, $J$ содержатся среди вершин, отмеченных черными кругами, то коды вершин $Z$, $Y$, $F$, $J$, во всех случаях мы можем назвать явно, с точностью до подклееных окружений. Это позволяет ввести определяющие отношения, записанные в правой части рисунка~\ref{UR1b}.

\medskip

{\bf Характеризация.} Составление таблиц полностью аналогично случаю $\mathbb{C}1$-цепи.

\medskip

{\bf Восстановление кода.}
Заметим, что у $Z$ и $F$ в каждом из случаев общие начальники, то есть код каждой из этих вершин может быть вычислен, исходя из кодов остальных трех вершин. $Y$ вычисляется также легко, так как во всех случаях первым начальником является $Z$, а тип второго в каждом случае виден из рисунка.

Начальники $J$ вычисляется следующим образом: случай $1$ -- это $\mathbf{Prev}(F)$; случай $2$ -- это $\mathbf{UpRightChain.FBoss}(Z)$ и $\mathbb{B}$-тип второго начальника; случай $3$ -- это $\mathbf{Prev}(F)$; случай $4$ -- это $0$-цепь с указателем $1$ вокруг $\mathbb{A}$, окружения которого как $U$-часть $\mathbf{Fboss}(Z)$, и тип второго начальника $\mathbb{B}$;
случай $5$ -- это $\mathbf{Plus.FBoss}(Z)$ и тип второго начальника $\mathbb{B}$.

Таким образом, зная код пути $ZYJ$, можно вычислить код пути $ZFJ$, и наоборот. А также зная код пути $FZY$, можно вычислить код пути $FJY$, и наоборот. Таким образом, мы можем осуществить локальное преобразование пути через операцию с его кодом.

\medskip

\subsection{Случай цепи $\mathbb{UR}1$; преобразования 7 и 10}

\begin{figure}[hbtp]
\centering
\includegraphics[width=1\textwidth]{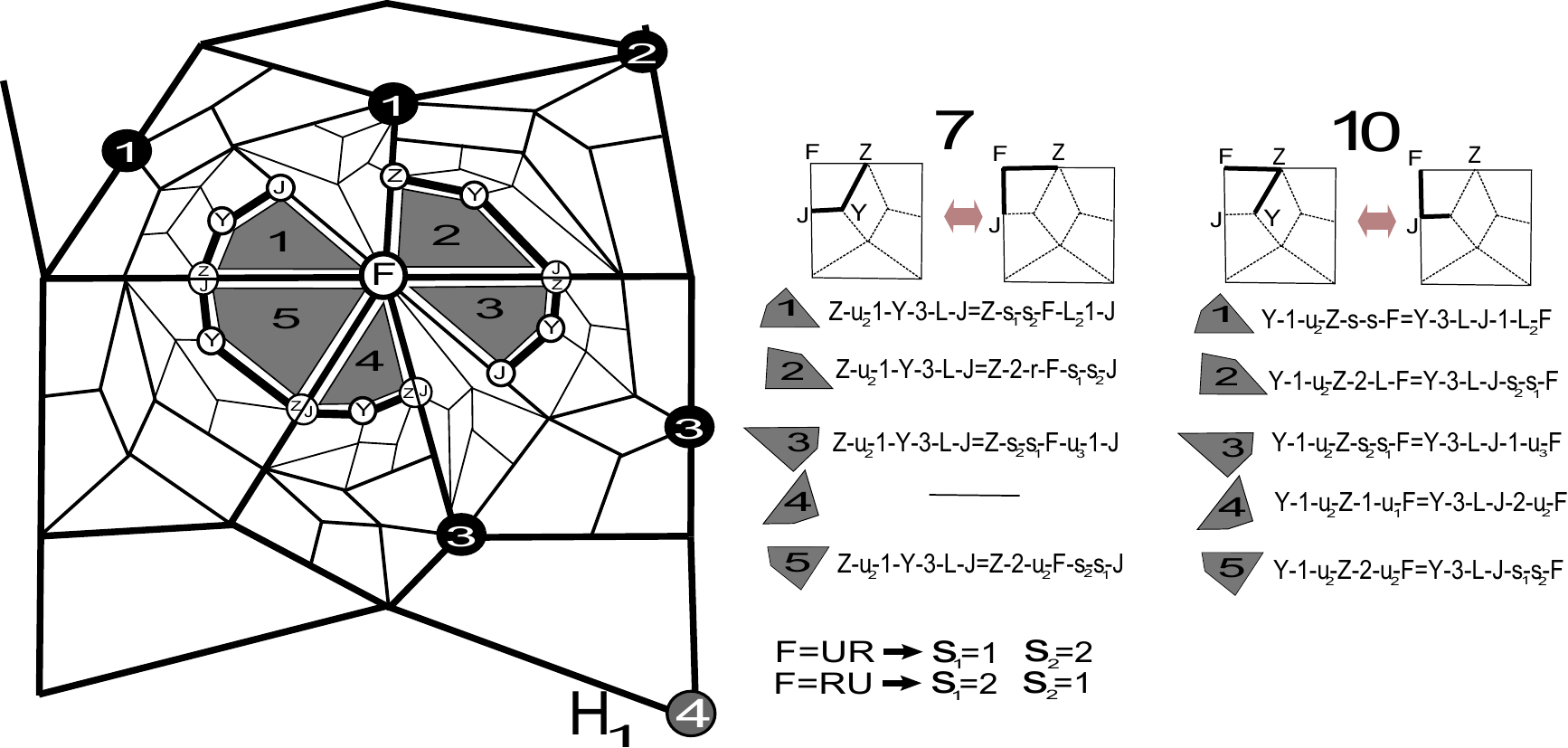}
\caption{Случаи расположения пути вокруг UR1-цепи и определяющие соотношения для локальных преобразований 7 и 10}
\label{UR1a}
\end{figure}

В правой части рисунка~\ref{UR1a} изображены локальные преобразования $7$ и $10$.

\medskip

Фиксируем некоторую вершину типа $\mathbb{UR}$ третьего уровня и ее начальников.
Черными кругами отметим вершины, являющиеся начальниками вершин  $Z$, $Y$, $F$, $J$.
Заметим, что зная окружение центральной вершины, мы можем вычислить коды вершин $Z$, $Y$, $F$, $J$, с точностью до подклееных окружений. Это позволяет ввести определяющие отношения, записанные в правой части рисунка~\ref{UR1a}.

\medskip

{\bf Характеризация.} Составление таблиц полностью аналогично случаю $\mathbb{C}1$-цепи.

\medskip

{\bf Восстановление кода.}
Вершина $F$ имеет общего начальника с одной из вершин $J$ или $Z$, в каждом из случаев, кроме случая $4$ (а в случае $4$ нам ее не надо вычислять). Таким образом, зная код пути $ZYJ$, можно вычислить код пути $ZFJ$. В обратную сторону, а также для  локального преобразования $10$: очевидно, что зная $F$, можно вычислить коды остальных вершин. То есть, мы можем осуществить локальное преобразование пути через операцию с его кодом.

\medskip

\subsection{Случай цепи $\mathbb{UR}2$; преобразования 8 и 9}

\begin{figure}[hbtp]
\centering
\includegraphics[width=1\textwidth]{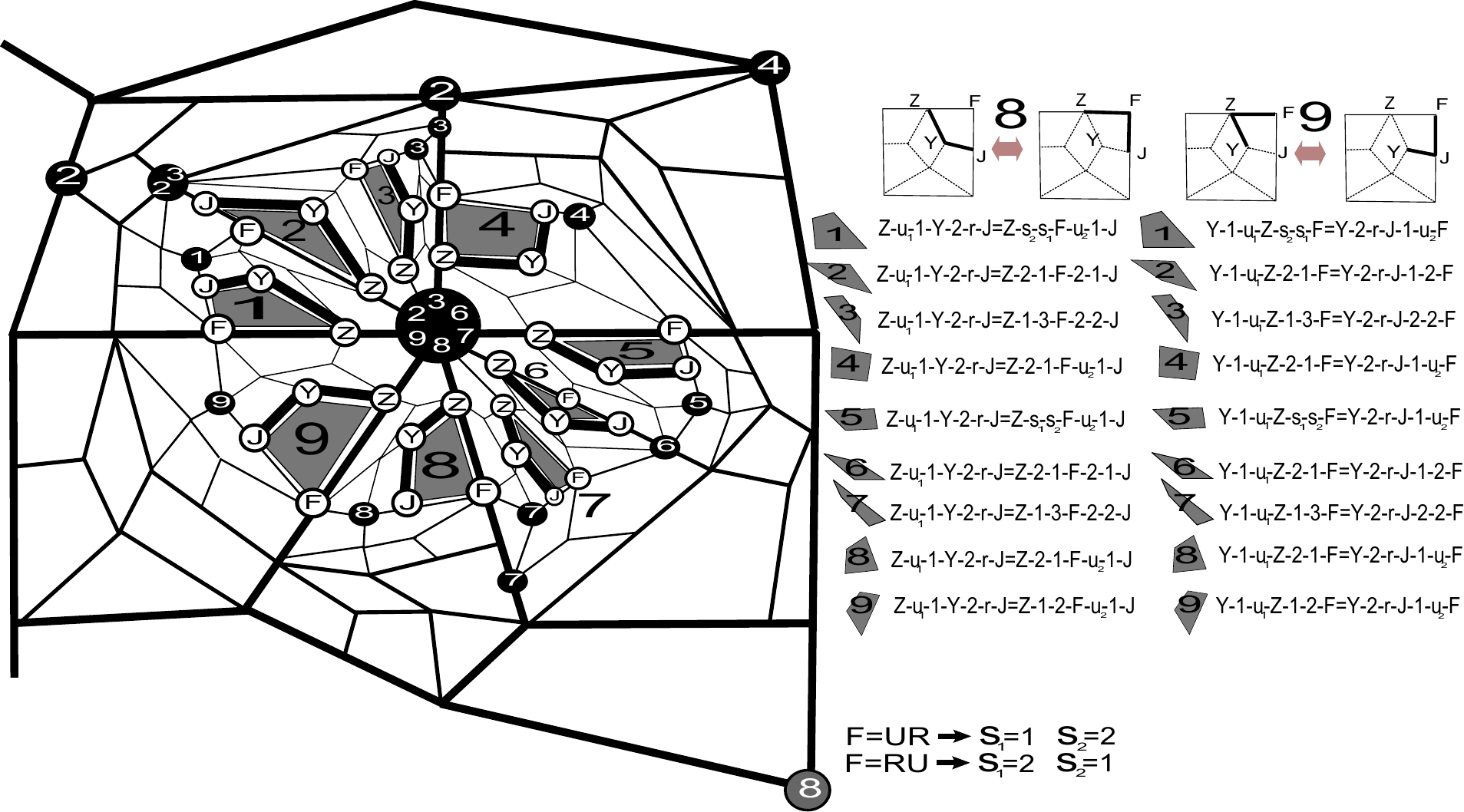}
\caption{Случаи расположения пути вокруг UR2-цепи и определяющие соотношения для локальных преобразований 8 и 9}
\label{UR2b}
\end{figure}

В правой части рисунка~\ref{UR2b} изображены локальные преобразования $8$ и $9$.

\medskip

Фиксируем вершину типа $\mathbb{UR}$ третьего уровня и ее начальников.
Черными кругами отметим вершины, являющиеся начальниками вершин  $Z$, $Y$, $F$, $J$.
Заметим, что зная окружение центрального $\mathbb{UR}$-узла, мы можем найти их типы, уровни и окружения (для вершины в правом нижнем углу -- только тип). Так как начальники вершин  $Z$, $Y$, $F$, $J$ содержатся среди вершин, отмеченных черными кругами, то коды вершин $Z$, $Y$, $F$, $J$, во всех случаях мы можем назвать явно, с точностью до подклееных окружений. Это позволяет ввести определяющие отношения, записанные в правой части рисунка~\ref{UR2b}.

\medskip

{\bf Характеризация.} Составление таблиц полностью аналогично случаю $\mathbb{C}1$-цепи.

\medskip

{\bf Восстановление кода.}
Заметим, что у $Z$ и $F$ в каждом из случаев общие начальники, то есть код каждой из этих вершин может быть вычислен, исходя из кодов остальных трех вершин. $Y$ вычисляется также легко, так как во всех случаях первым начальником является $Z$, а тип второго в каждом случае ясен.

Начальники $J$ вычисляется следующим образом: в  случаях $1$, $4$, $5$, $8$, $9$ -- единственный начальник, узел $F$; в  случаях $2$, $6$ -- начальники те же что и у $F$; в случаях $3$ и $7$ -- первый и второй начальники как первый и третий у $F$.

Таким образом, зная код пути $ZYJ$, можно вычислить код пути $ZFJ$, и наоборот. А также зная код пути $FZY$, можно вычислить код пути $FJY$, и наоборот. Таким образом, мы можем осуществить локальное преобразование пути через операцию с его кодом.

\medskip

\subsection{Случай цепи $\mathbb{UR}2$; преобразования 7 и 10}

\begin{figure}[hbtp]
\centering
\includegraphics[width=1\textwidth]{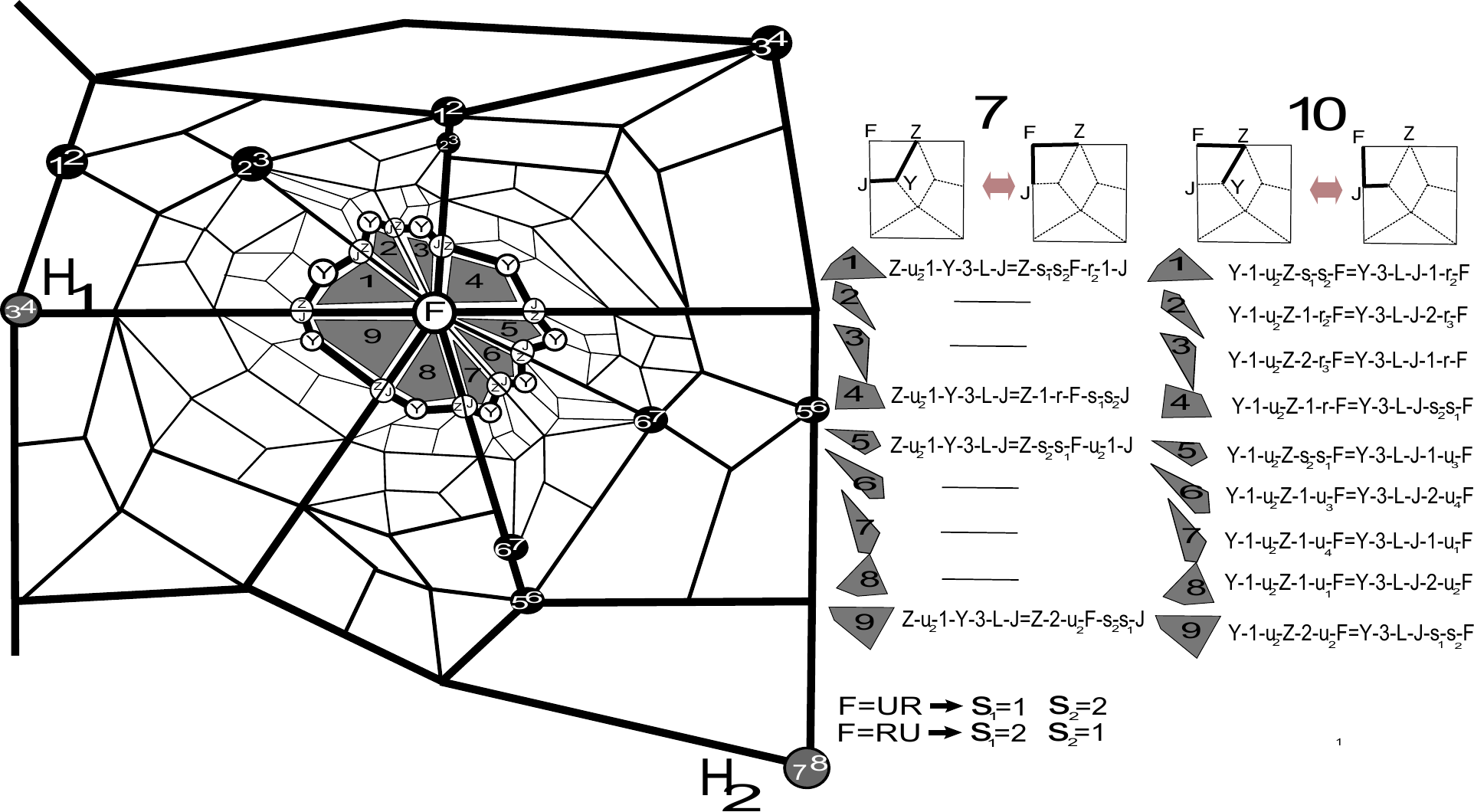}
\caption{Случаи расположения пути вокруг UR2-цепи и определяющие соотношения для локальных преобразований 7 и 10}
\label{UR2a}
\end{figure}

В правой части рисунка~\ref{UR2a} изображены локальные преобразования $7$ и $10$.

\medskip

Фиксируем некоторую вершину типа $\mathbb{UR}$ и ее начальников.
Черными кругами отметим вершины, являющиеся начальниками вершин  $Z$, $Y$, $F$, $J$.
Заметим, что зная окружение этой вершины, мы можем вычислить коды вершин $Z$, $Y$, $F$, $J$, с точностью до подклееных окружений (кроме вершин $H_1$ и $H_2$, для которых можно найти тип). Это позволяет ввести определяющие отношения, записанные в правой части рисунка~\ref{UR2a}.

Для локального преобразования $7$ в случаях $2$, $3$, $6$, $7$, $8$  пути удовлетворяют условиям мертвого паттерна, и в этих случаях мы соотношения не вводим.

\medskip

{\bf Характеризация.} Составление таблиц полностью аналогично случаю $\mathbb{C}1$-цепи.

\medskip

{\bf Восстановление кода.}
В случаях $1,4,5,9$ вершина $F$ имеет общего начальника с одной из вершин $J$ или $Z$. То есть, зная код пути $ZYJ$, можно вычислить код пути $ZFJ$. В обратную сторону, а также для  локального преобразования $10$: очевидно, что зная $F$, можно вычислить коды остальных вершин. Таким образом, мы можем осуществить локальное преобразование пути через операцию с его кодом.

\medskip

\subsection{Случай цепи $\mathbb{UR}3$; преобразования 7, 8, 9, 10}

Случай $\mathbb{UR}3$-цепи полностью аналогичен $\mathbb{UR}2$ случаю, соотношения выглядят идентично, только кодировки вершин $J$, $F$, $Z$, $Y$ отвечают $\mathbb{UR}3$-цепи. Все рассуждения о вычислении путей полностью аналогичны. Соотношений вводится столько же, сколько для $\mathbb{UR}2$ случая.

\medskip

\subsection{Случай цепи $\mathbb{DR}1$; преобразования 8 и 9}

\begin{figure}[hbtp]
\centering
\includegraphics[width=0.5\textwidth]{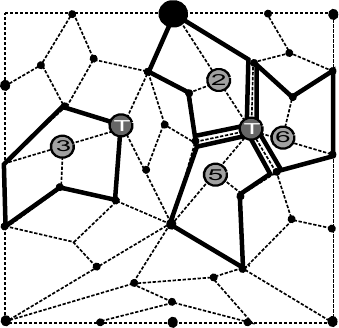}
\caption{Случаи расположения $\mathbb{DR}$-вершины}
\label{RDplacement}
\end{figure}

Аналогично предыдущим случаям, мы можем ввести определяющие соотношения. Сначала зафиксируем вершину типа $\mathbb{DR}$ третьего уровня и ее начальников. Всего возможно четыре случая расположения $\mathbb{DR}$-вершины, на четырех внутренних ребрах (типов $2$, $3$, $5$, $6$). Все эти случаи расположения показаны на рисунке~\ref{RDplacement}.

\medskip

\begin{figure}[hbtp]
\centering
\includegraphics[width=1\textwidth]{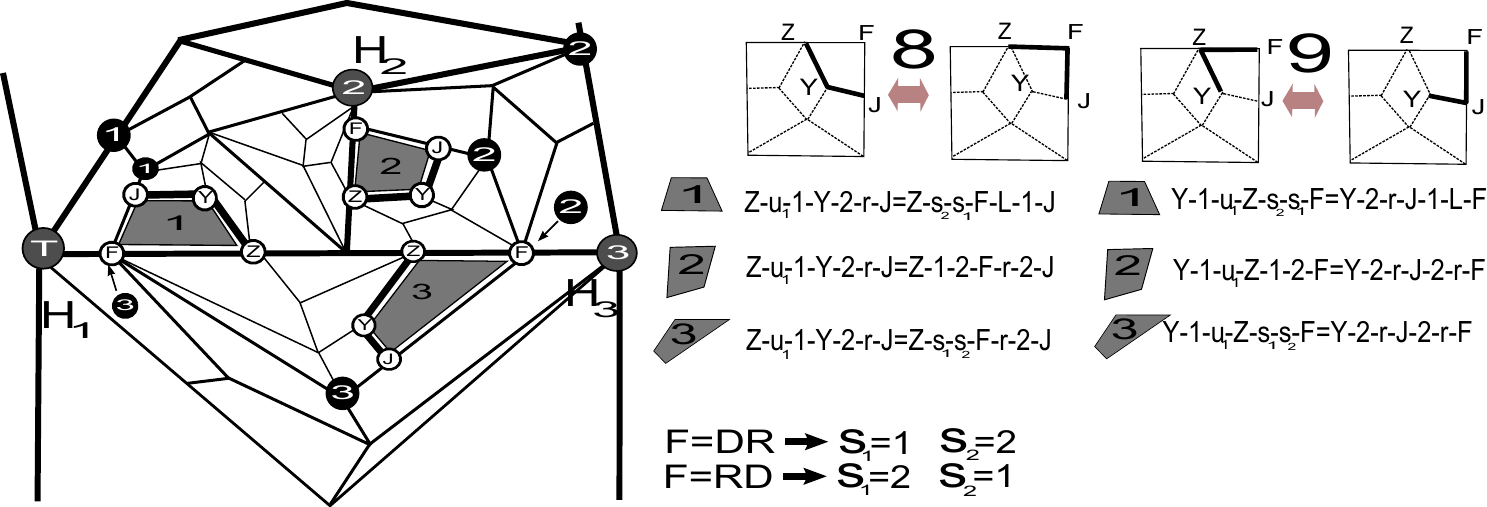}
\caption{Случаи расположения пути вокруг $\mathbb{DR}$1-цепи  и определяющие соотношения для локальных преобразований 8 и 9}
\label{DR1b}
\end{figure}

Первым начальником нашей $\mathbb{DR}$-вершины во всех случаях будет вершина в середине верхней стороны. Теперь непосредственно рассмотрим цепь с центром в этой $\mathbb{DR}$-вершине.
В правой части рисунка~\ref{DR1b} изображены локальные преобразования $8$ и $9$.

\medskip

Черными кругами отметим вершины, являющиеся начальниками вершин  $Z$, $Y$, $F$, $J$.  Два круга попадают в $F$-узлы других расположений, это отмечено стрелками.
Зная окружение центрального $\mathbb{DR}$-узла, мы можем найти их типы, уровни и окружения (для вершин $H_1$ $H_2$, $H_3$ можем найти только тип).

Тип узла, отмеченного как $T$, мы можем определить исходя из типа внутреннего ребра, на котором лежит центральный $\mathbb{DR}$-узел. Если это ребро $3$ -- то $\mathbb{A}$, в остальных случаях ($2$, $5$ или $6$) -- $\mathbb{B}$.

Так как начальники вершин  $Z$, $Y$, $F$, $J$ содержатся среди вершин, отмеченных черными кругами, то коды вершин $Z$, $Y$, $F$, $J$, во всех случаях мы можем назвать явно, с точностью до подклееных окружений. Это позволяет ввести определяющие отношения, записанные в правой части рисунка~\ref{DR1b}.

\medskip

{\bf Характеризация.} Составление таблиц полностью аналогично случаю $\mathbb{C}1$-цепи.

\medskip

{\bf Восстановление кода.}
Поскольку мы можем установить, с каким именно случаем расположения мы имеем дело, то окружение каждой вершины мы можем вычислить, зная код остальных трех.
Кроме того, заметим, что у $Z$ и $F$ в каждом из случаев общие начальники, то есть код каждой из этих вершин может быть вычислен, исходя из кодов остальных трех вершин. $Y$ вычисляется также легко, так как во всех случаях первым начальником является $Z$, а тип правого нижнего угла в каждом случае ясен.

Начальники $J$ вычисляется следующим образом:  случай $1$ -- это $\mathbf{Prev}(F)$;  случай $2$ -- первый начальник -- $\mathbf{UpRightChain.FBoss}(Z)$ и $\mathbb{B}$-тип второго;  случай $3$ -- это $1$-цепь вокруг узла $T$ c указателями в зависимости от типа ребра центрального $\mathbb{DR}$-узла: $3$ для $3$ типа, $1$ для $2$ типа, $2$ для $6$ типа, $3$ для $5$ типа.

Таким образом, зная код пути $ZYJ$, можно вычислить код пути $ZFJ$, и наоборот. А также зная код пути $FZY$, можно вычислить код пути $FJY$, и наоборот. Таким образом, мы можем осуществить локальное преобразование пути через операцию с его кодом.

\medskip

\subsection{Случай цепи $\mathbb{DR}1$; преобразования 7 и 10}

\begin{figure}[hbtp]
\centering
\includegraphics[width=1\textwidth]{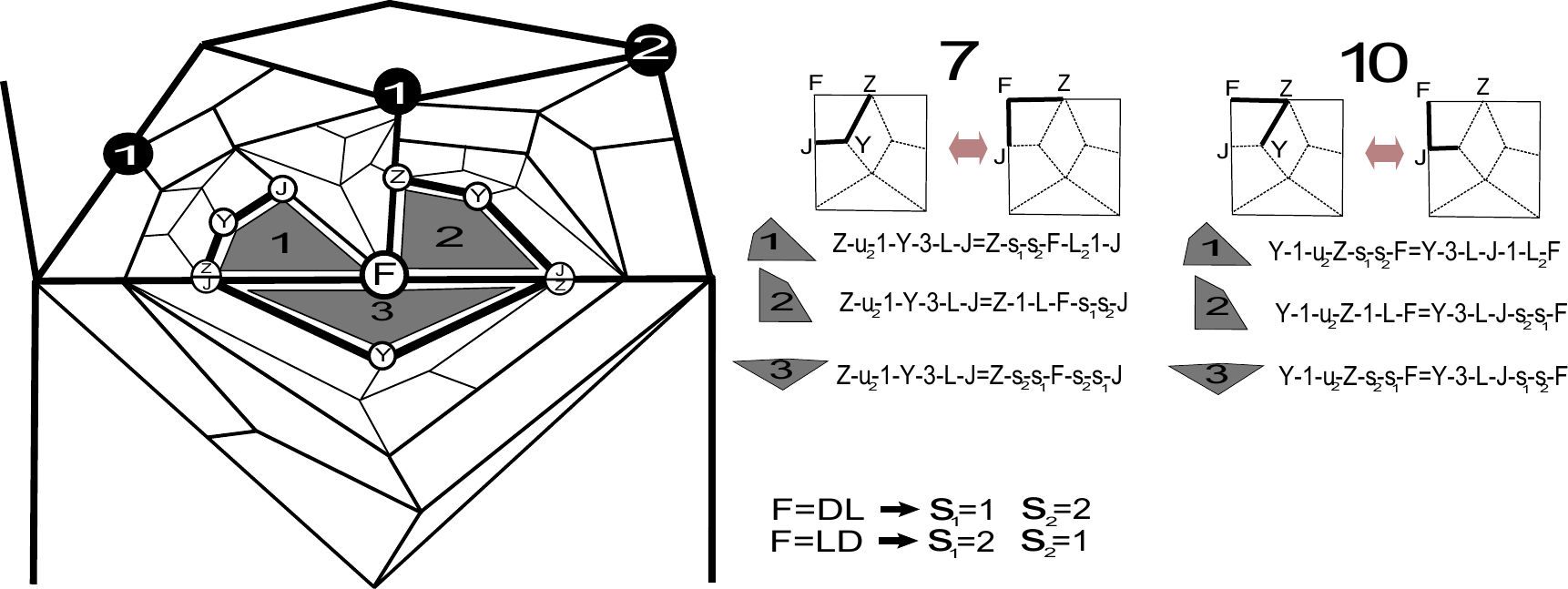}
\caption{Случаи расположения пути вокруг DR1-цепи и определяющие соотношения для локальных преобразований ~7 и 10}
\label{DR1a}
\end{figure}

В правой части рисунка~\ref{DR1a} изображены локальные преобразования $7$ и $10$.

\medskip

Фиксируем некоторую вершину типа $\mathbb{DR}$ третьего уровня и ее начальников.
Черными  кругами отметим вершины, являющиеся начальниками вершин  $Z$, $Y$, $F$, $J$.
Заметим, что зная окружение центральной $\mathbb{DR}$-вершины, мы можем в каждом случае вычислить коды вершин $Z$, $Y$, $F$, $J$ с точностью до подклееных окружений. Это позволяет ввести определяющие отношения, записанные в правой части рисунка~\ref{DR1a}.

\medskip

{\bf Характеризация.} Составление таблиц полностью аналогично случаю $\mathbb{C}1$-цепи.

\medskip

{\bf Восстановление кода.}
Вершина $F$ имеет общего начальника с одной из вершин $J$ или $Z$, в каждом из случаев.  То есть, зная код пути $ZYJ$, можно вычислить код пути $ZFJ$. В обратную сторону, а также для $10$ локального преобразования: очевидно, что зная $F$, можно вычислить коды остальных вершин. Таким образом, мы можем осуществить локальное преобразование пути через операцию с его кодом.

\medskip

\subsection{Случай цепи $\mathbb{DR}2$; преобразования 8 и 9}

\begin{figure}[hbtp]
\centering
\includegraphics[width=0.9\textwidth]{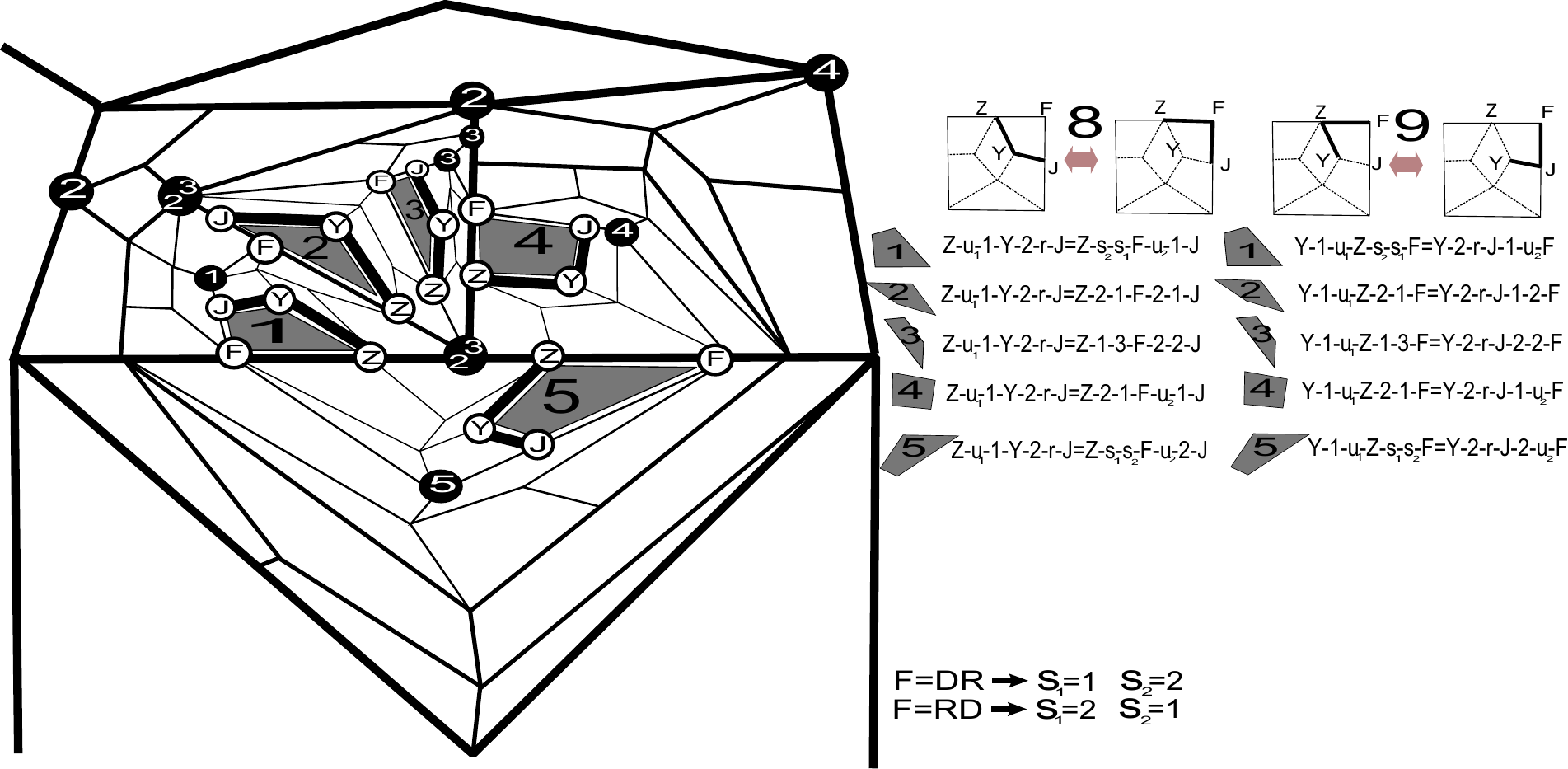}
\caption{Случаи расположения пути вокруг DR2-цепи  и определяющие соотношения для локальных преобразований 8 и 9}
\label{DR2b}
\end{figure}

В правой части рисунка~\ref{DR2b} изображены локальные преобразования $8$ и $9$.

\medskip

Зафиксируем вершину типа $\mathbb{DR}$ третьего уровня и ее начальников.
Черными кругами отметим вершины, являющиеся начальниками вершин  $Z$, $Y$, $F$, $J$.
Заметим, что зная окружение центрального $\mathbb{DR}$-узла, мы можем найти их типы, уровни и окружения. Так как начальники вершин  $Z$, $Y$, $F$, $J$ содержатся среди вершин, отмеченных черными кругами, то коды вершин $Z$, $Y$, $F$, $J$, во всех случаях мы можем назвать явно, с точностью до подклееных окружений. Это позволяет ввести определяющие отношения, записанные в правой части рисунка~\ref{DR2b}.

\medskip

{\bf Характеризация.} Составление таблиц полностью аналогично случаю $\mathbb{C}1$-цепи.

\medskip

{\bf Восстановление кода.}
Заметим, что у $Z$ и $F$ в каждом из случаев общие начальники, то есть код каждой из этих вершин может быть вычислен, исходя из кодов остальных трех вершин. $Y$ вычисляется также легко, так как во всех случаях первым начальником является $Z$, а тип второго в каждом случае ясен.

Начальники $J$ вычисляется следующим образом: в  случаях $1$, $4$ и $5$ -- единственный начальник, узел $F$; в случае $2$ -- начальники те же что и у $F$; в  случае $3$ -- первый и второй начальники как первый и третий у $F$.

Таким образом, зная код пути $ZYJ$, можно вычислить код пути $ZFJ$, и наоборот. А также зная код пути $FZY$, можно вычислить код пути $FJY$, и наоборот. Таким образом, мы можем осуществить локальное преобразование пути через операцию с его кодом.

\medskip

\subsection{Случай цепи $\mathbb{DR}2$; преобразования 7 и 10}

\begin{figure}[hbtp]
\centering
\includegraphics[width=1\textwidth]{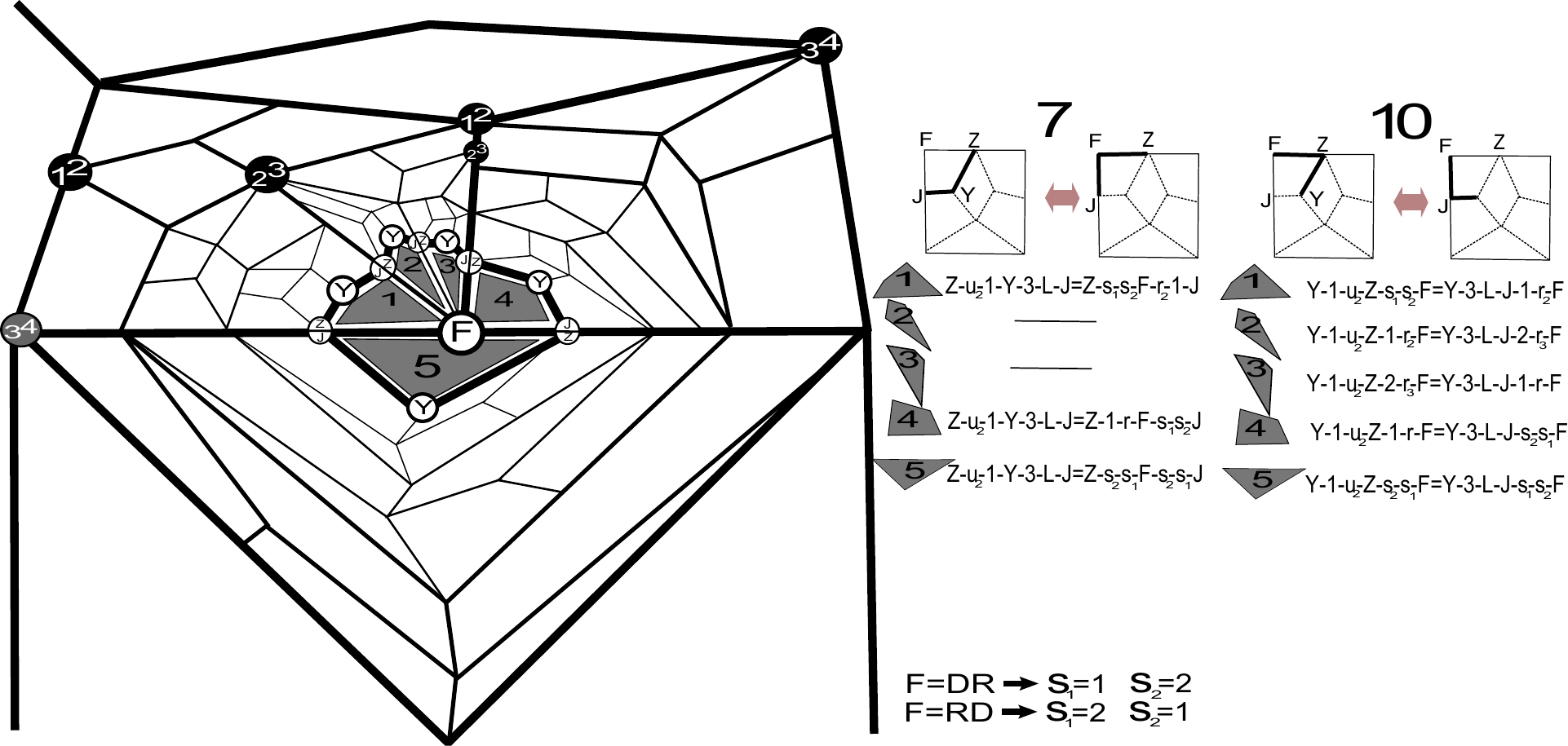}
\caption{Случаи расположения пути вокруг DR2-цепи и определяющие соотношения для локальных преобразований 7 и 10}
\label{DR2a}
\end{figure}

В правой части рисунка~\ref{DR2a} изображены локальные преобразования $7$ и $10$.

\medskip

Зафиксируем некоторую вершину типа $\mathbb{DR}$ третьего уровня и ее начальников.
Черными кругами отметим вершины, являющиеся начальниками вершин  $Z$, $Y$, $F$, $J$. Заметим, что зная окружение этой вершины, мы можем вычислить коды вершин $Z$, $Y$, $F$, $J$, с точностью до подклееных окружений. Это позволяет ввести определяющие отношения, записанные в правой части рисунка~\ref{DR2a}.

Для  локального преобразования $7$ в случаях $2$, $3$, $6$, $7$, $8$ пути удовлетворяют условиям мертвого паттерна, и в этих случаях мы соотношения не вводим.

\medskip

{\bf Характеризация.} Составление таблиц полностью аналогично случаю $\mathbb{C}1$-цепи.

\medskip

{\bf Восстановление кода.}
В случаях $1,4,5$ вершина $F$ имеет общего начальника с одной из вершин $J$ или $Z$. То есть, $ZYJ$, можно вычислить код пути $ZFJ$. В обратную сторону, а также для  локального преобразования $10$: очевидно, что зная $F$, можно вычислить коды остальных вершин. Таким образом, мы можем осуществить локальное преобразование пути через операцию с его кодом.

\medskip

\subsection{Случай цепи $\mathbb{DR}3$; преобразования 7, 8, 9, 10}

Случай $\mathbb{DR}3$-цепи полностью аналогичен $\mathbb{DR}2$ случаю, соотношения выглядят идентично, только кодировки вершин $J$, $F$, $Z$, $Y$ отвечают $\mathbb{DR}3$-цепи. Все рассуждения о вычислении путей полностью аналогичны. Соотношений вводится столько же, сколько для $\mathbb{DR}2$ случая.

\medskip

\subsection{Случай цепи $\mathbb{DL}1$; преобразования 8 и 9}

\begin{figure}[hbtp]
\centering
\includegraphics[width=0.9\textwidth]{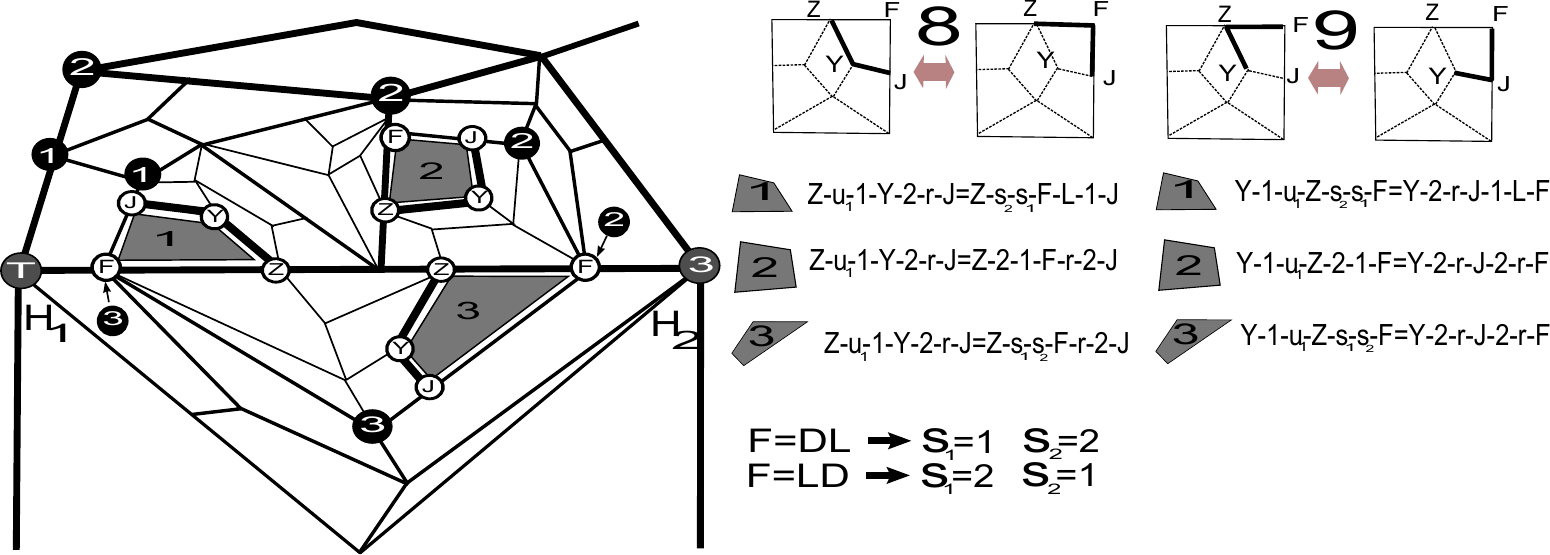}
\caption{Случаи расположения пути вокруг DL1-цепи  и определяющие соотношения для локальных преобразований 8 и 9}
\label{DL1b}
\end{figure}

В правой части рисунка~\ref{DL1b} изображены локальные преобразования $8$ и $9$.

\medskip

Зафиксируем вершину типа $\mathbb{DL}$ третьего уровня и ее начальников.
Черными и серыми кругами отметим вершины, являющиеся начальниками вершин  $Z$, $Y$, $F$, $J$. Два круга попадают в $F$-узлы других расположений, это отмечено стрелками.
Числа в них обозначают, начальниками вершин какого случая является данная вершина. Зная окружение центрального $\mathbb{DL}$-узла, мы можем найти их типы, уровни и окружения (для вершин $H_1$ и $H_2$ -- тип).

Так как начальники вершин  $Z$, $Y$, $F$, $J$ содержатся среди вершин, отмеченных черными кругами, то коды вершин $Z$, $Y$, $F$, $J$, во всех случаях мы можем назвать явно, с точностью до подклееных окружений. Это позволяет ввести определяющие отношения, записанные в правой части рисунка~\ref{DL1b}.

\medskip

{\bf Характеризация.} Составление таблиц полностью аналогично случаю $\mathbb{C}1$-цепи.

\medskip

{\bf Восстановление кода.}
Заметим, что у $Z$ и $F$ в каждом из случаев общие начальники, то есть код каждой из этих вершин может быть вычислен, исходя из кодов остальных трех вершин. $Y$ вычисляется также легко, так как во всех случаях первым начальником является $Z$, а тип второго в каждом случае ясен.

Начальники $J$ вычисляется следующим образом:  случай $1$ -- $\mathbf{Prev}(F)$;  случай $2$ --  первый -- $\mathbf{BottomLeftChain.FBoss}(Z)$ и $\mathbb{A}$-тип второго. Для  случая $3$ заметим, что тип узла, отмеченного символом $T$ может быть установлен, в зависимости от типа ребра, на котором лежат $Z$ и $F$. Так как это ребро, с центром в  $\mathbb{DL}$-вершине, то тип его может быть либо $4$ либо $7$, в соответствии со структурой разбиения макроплитки на подплитки. Если ребро имеет $4$ тип, то вершина $T$ имеет тип $\mathbb{A}$, а первый начальник $J$ -- это вершина из $\mathbb{A}1$ цепи с указателем $2$. Если же ребро имеет тип $7$, то $T$-вершина является вторым начальником $F$ и $Z$, а первый начальник $J$ -- это вершина из $2$-цепи вокруг него с указателем, соответствующим $7$ ребру.

Таким образом, зная код пути $ZYJ$, можно вычислить код пути $ZFJ$, и наоборот. А также зная код пути $FZY$, можно вычислить код пути $FJY$, и наоборот. Таким образом, мы можем осуществить локальное преобразование пути через операцию с его кодом.

\medskip

\subsection{Случай цепи $\mathbb{DL}1$; преобразования 7 и 10}

\begin{figure}[hbtp]
\centering
\includegraphics[width=1\textwidth]{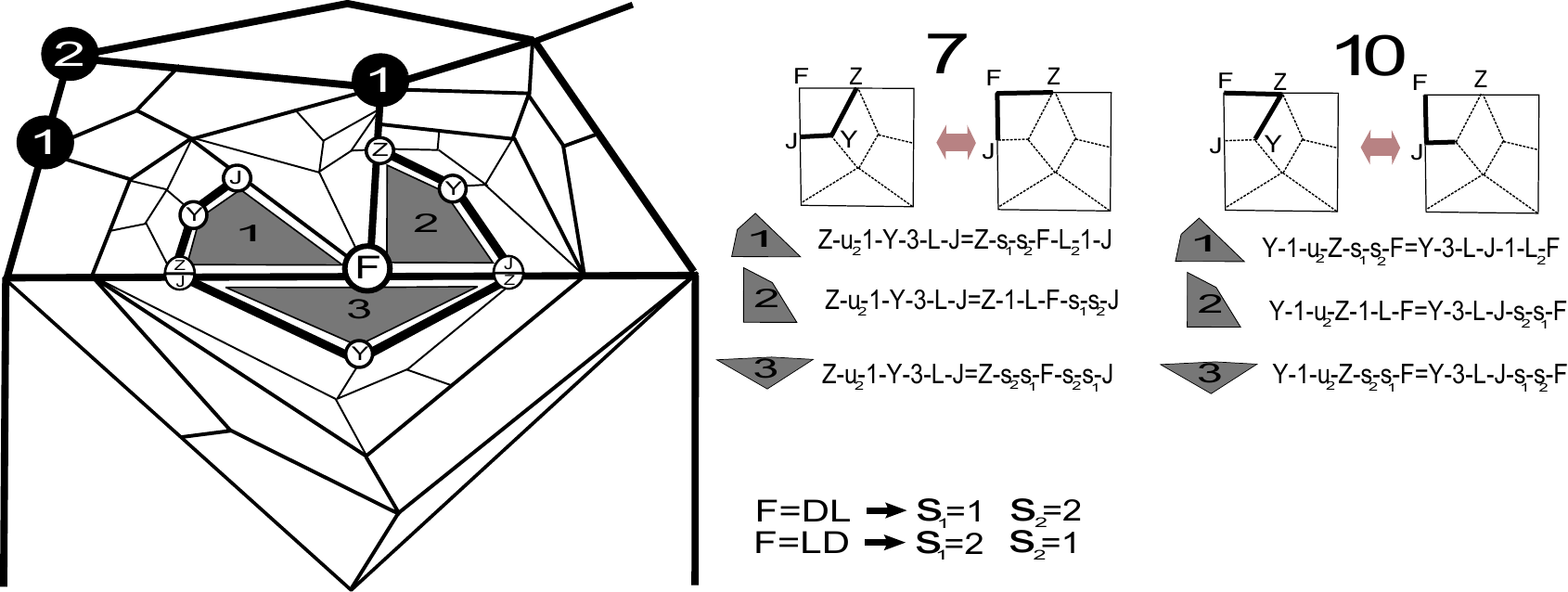}
\caption{Случаи расположения пути вокруг DL1-цепи и соответствующие им определяющие соотношения для локальных преобразований 7 и 10}
\label{DL1a}
\end{figure}

В правой части рисунка~\ref{DL1a} изображены локальные преобразования $7$ и $10$.

\medskip

Фиксируем некоторую вершину типа $\mathbb{DL}$ третьего уровня и ее начальника.
Черными кругами отметим вершины, являющиеся начальниками вершин  $Z$, $Y$, $F$, $J$.
Заметим, что зная окружение центральной вершины, мы можем вычислить коды вершин $Z$, $Y$, $F$, $J$, с точностью до подклееных окружений. Это позволяет ввести определяющие отношения, записанные в правой части рисунка~\ref{DL1a}.

\medskip

{\bf Характеризация.} Составление таблиц полностью аналогично случаю $\mathbb{C}1$-цепи.

\medskip

{\bf  Восстановление кода.}
Вершина $F$ имеет общего начальника с одной из вершин $J$ или $Z$, в каждом из случаев.  То есть, $ZYJ$, можно вычислить код пути $ZFJ$. В обратную сторону, а также для  локального преобразования $10$: очевидно, что зная $F$, можно вычислить коды остальных вершин. Таким образом, мы можем осуществить локальное преобразование пути через операцию с его кодом.

\medskip

\subsection{Случай цепи $\mathbb{DL}2$; преобразования 8 и 9}

\begin{figure}[hbtp]
\centering
\includegraphics[width=0.9\textwidth]{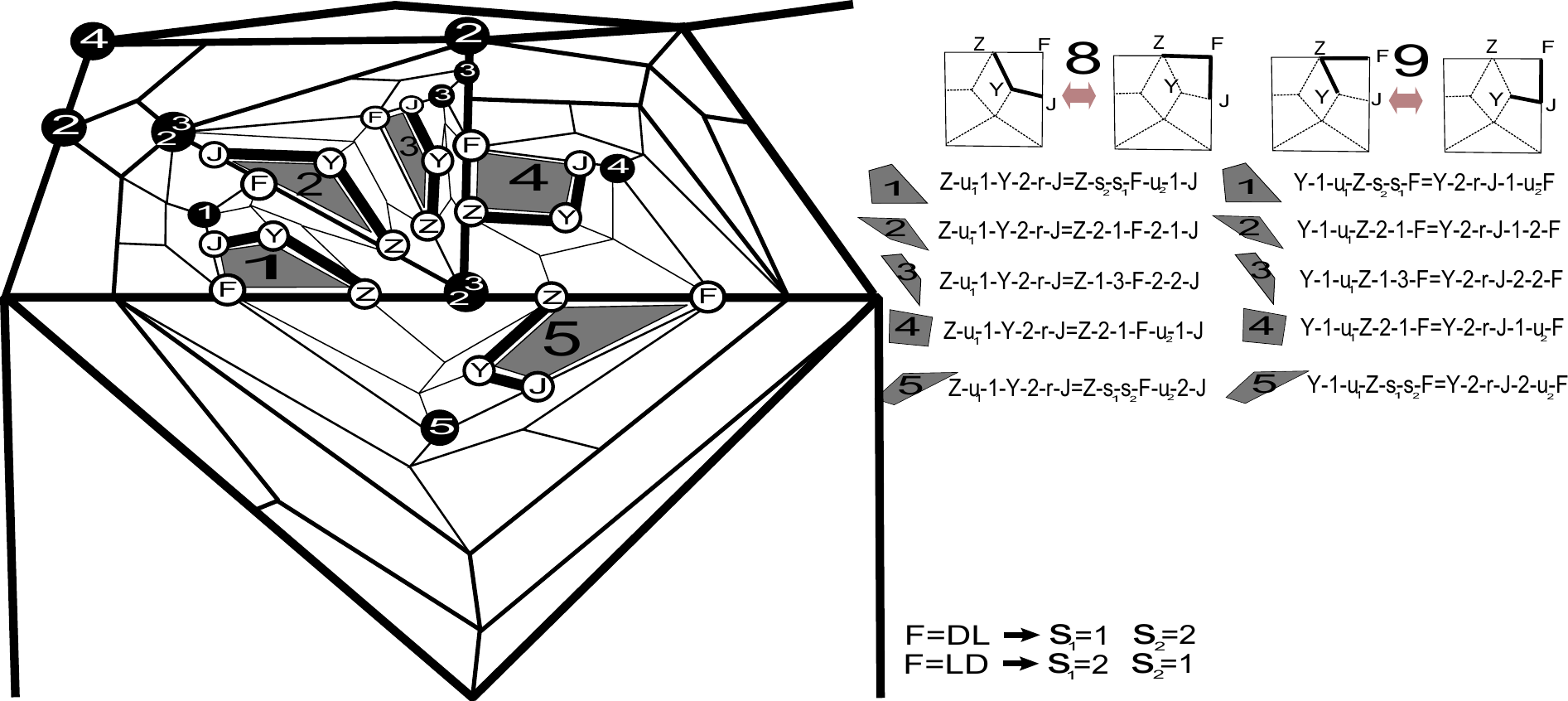}
\caption{Случаи расположения пути вокруг DL2-цепи и определяющие соотношения для локальных преобразований 8 и 9}
\label{DL2b}
\end{figure}

В правой части рисунка~\ref{DL2b} изображены локальные преобразования $8$ и $9$.

\medskip

Зафиксируем вершину типа $\mathbb{DR}$ третьего уровня и ее начальников.
Черными  кругами отметим вершины, являющиеся начальниками вершин  $Z$, $Y$, $F$, $J$.
Заметим, что зная окружение центрального $\mathbb{DL}$-узла, мы можем найти их типы, уровни и окружения. Так как начальники вершин  $Z$, $Y$, $F$, $J$ содержатся среди вершин, отмеченных черными кругами, то коды вершин $Z$, $Y$, $F$, $J$, во всех случаях мы можем назвать явно, с точностью до подклееных окружений. Это позволяет ввести определяющие отношения, записанные в правой части рисунка~\ref{DL2b}.

\medskip

{\bf Характеризация.} Составление таблиц полностью аналогично случаю $\mathbb{C}1$-цепи.

\medskip

{\bf Восстановление кода.}
Заметим, что у $Z$ и $F$ в каждом из случаев общие начальники, то есть код каждой из этих вершин может быть вычислен, исходя из кодов остальных трех вершин. $Y$ вычисляется также легко, так как во всех случаях первым начальником является $Z$, а тип второго в каждом случае ясен.

Начальники $J$ вычисляется следующим образом: в  случаях $1$, $4$, $5$ -- единственный начальник, узел $F$; в  случае $2$ -- начальники те же что и у $F$; в  случае $3$ -- первый и второй начальники как первый и третий у $F$.

Таким образом, зная код пути $ZYJ$, можно вычислить код пути $ZFJ$, и наоборот. А также зная код пути $FZY$, можно вычислить код пути $FJY$, и наоборот. Таким образом, мы можем осуществить локальное преобразование пути через операцию с его кодом.

\medskip

\subsection{Случай цепи $\mathbb{DL}2$; преобразования 7 и 10}

\begin{figure}[hbtp]
\centering
\includegraphics[width=1\textwidth]{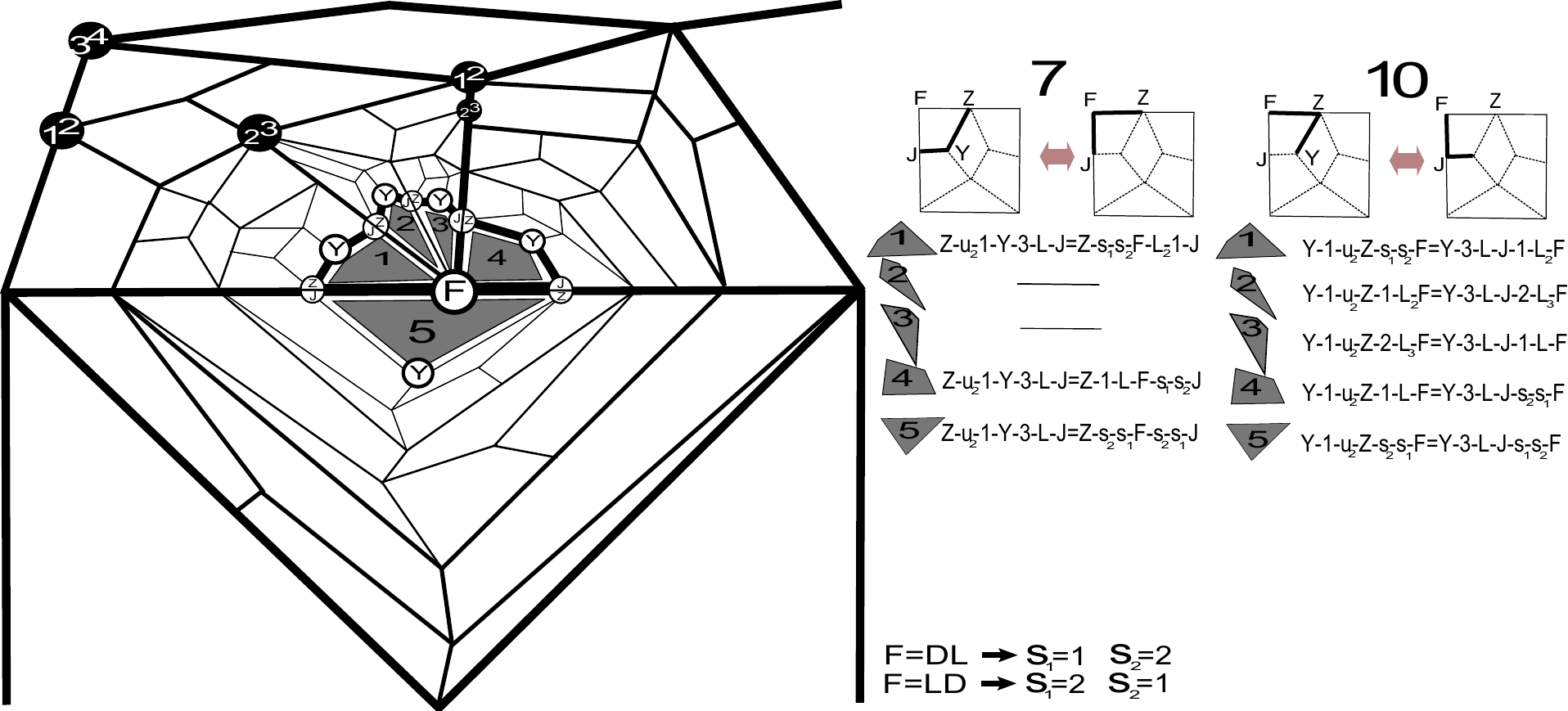}
\caption{Случаи расположения пути вокруг DL2-цепи и определяющие соотношения для локальных преобразований ~7 и 10}
\label{DL2a}
\end{figure}

В правой части рисунка~\ref{DL2a} изображены локальные преобразования $7$ и $10$.

\medskip

Фиксируем некоторую вершину типа $\mathbb{DL}$ третьего уровня и ее начальников.
Черными  кругами отметим вершины, являющиеся начальниками вершин  $Z$, $Y$, $F$, $J$.
Заметим, что зная окружение этой вершины, мы можем вычислить коды вершин $Z$, $Y$, $F$, $J$, с точностью до подклееных окружений. Это позволяет ввести определяющие отношения, записанные в правой части рисунка~\ref{DL2a}.

Для локального преобразования $7$ пути в случаях $2$ и $3$ удовлетворяют условиям мертвого паттерна, и в этих случаях мы соотношения не вводим.

\medskip

{\bf  Характеризация.} Составление таблиц полностью аналогично случаю $\mathbb{C}1$-цепи.

\medskip

{\bf Восстановление кода.}
В случаях $1,4,5$ вершина $F$ имеет общего начальника с одной из вершин $J$ или $Z$. То есть, $ZYJ$, можно вычислить код пути $ZFJ$. В обратную сторону, а также для  локального преобразования $10$: очевидно, что зная $F$, можно вычислить коды остальных вершин. Таким образом, мы можем осуществить локальное преобразование пути через операцию с его кодом.

\medskip

\subsection{Случай цепи $\mathbb{DL}3$; преобразования 7, 8, 9, 10}

Случай $\mathbb{DL}3$-цепи полностью аналогичен $\mathbb{DL}2$ случаю, соотношения выглядят идентично, только кодировки вершин $J$, $F$, $Z$, $Y$ отвечают $\mathbb{DL}3$-цепи. Все рассуждения о вычислении путей полностью аналогичны. Соотношений вводится столько же, сколько для $\mathbb{DL}2$ случая.

\medskip

Далее рассмотрим случаи цепей вблизи края макроплиток, то есть цепи с центрами в вершинах типов $\mathbb{D}$, $\mathbb{U}$, $\mathbb{R}$, $\mathbb{L}$. Рассмотрение этих случаев полностью аналогично случаям обычных боковых узлов, просто используются не все из локальных преобразований. Тем не менее, мы кратко приведем вводимые определяющие соотношения.

\medskip

\subsection{Случай цепи $\mathbb{D}1$; преобразования 8 и 9}

\begin{figure}[hbtp]
\centering
\includegraphics[width=0.9\textwidth]{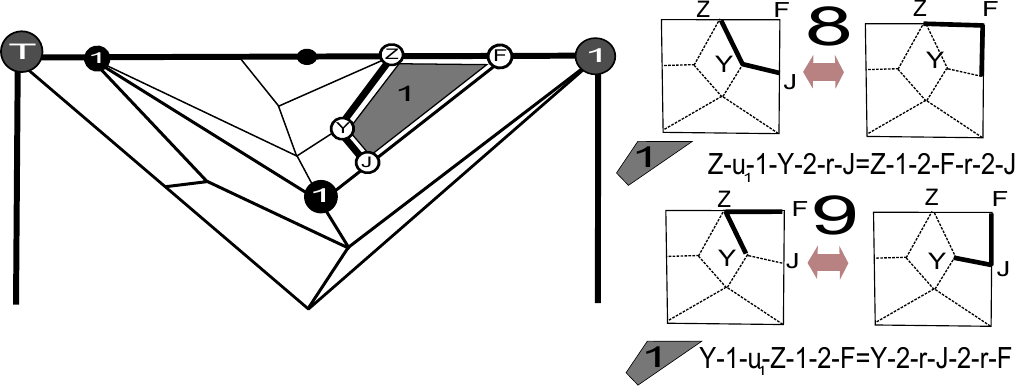}
\caption{Случаи расположения пути вокруг D1-цепи  и соответствующие им определяющие соотношения для локальных преобразований 8 и 9}
\label{D1b}
\end{figure}

В правой части рисунка~\ref{D1b} изображены локальные преобразования $8$ и $9$.

\medskip

Характеризация и восстановление кода полностью аналогично случаю $\mathbb{DL}1$ или  $\mathbb{DR}1$ цепи.

\subsection{Случай цепи $\mathbb{D}1$; преобразования 7 и 10}

\begin{figure}[hbtp]
\centering
\includegraphics[width=1\textwidth]{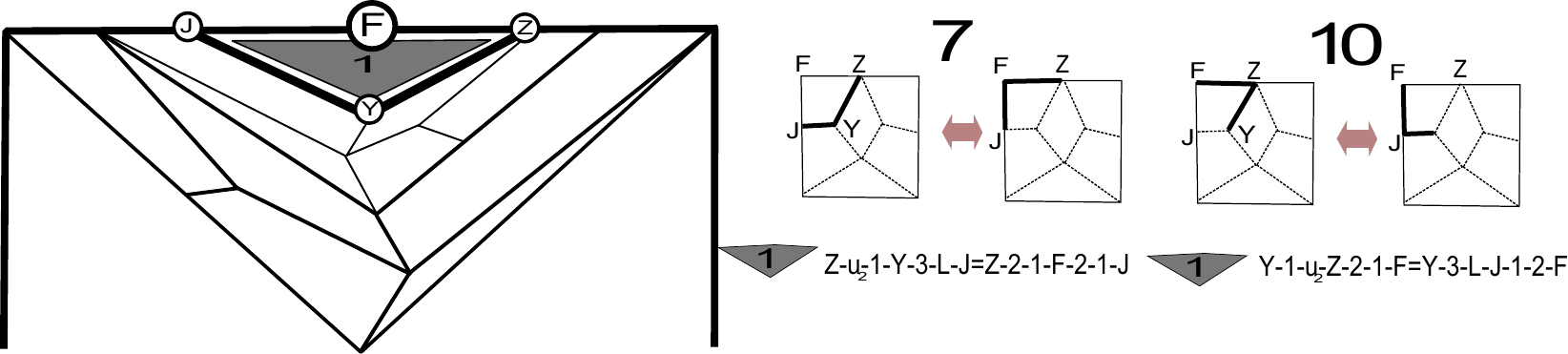}
\caption{Случаи расположения пути вокруг D1-цепи и соответствующие им определяющие соотношения для локальных преобразований 7 и 10}
\label{D1a}
\end{figure}

В правой части рисунка~\ref{D1a} изображены локальные преобразования $7$ и $10$.

\medskip

Характеризация и восстановление кода полностью аналогично случаю $\mathbb{DL}1$ или  $\mathbb{DR}1$ цепи.
\medskip

\subsection{Случай цепи $\mathbb{D}2$; преобразования 8 и 9}

\begin{figure}[hbtp]
\centering
\includegraphics[width=0.9\textwidth]{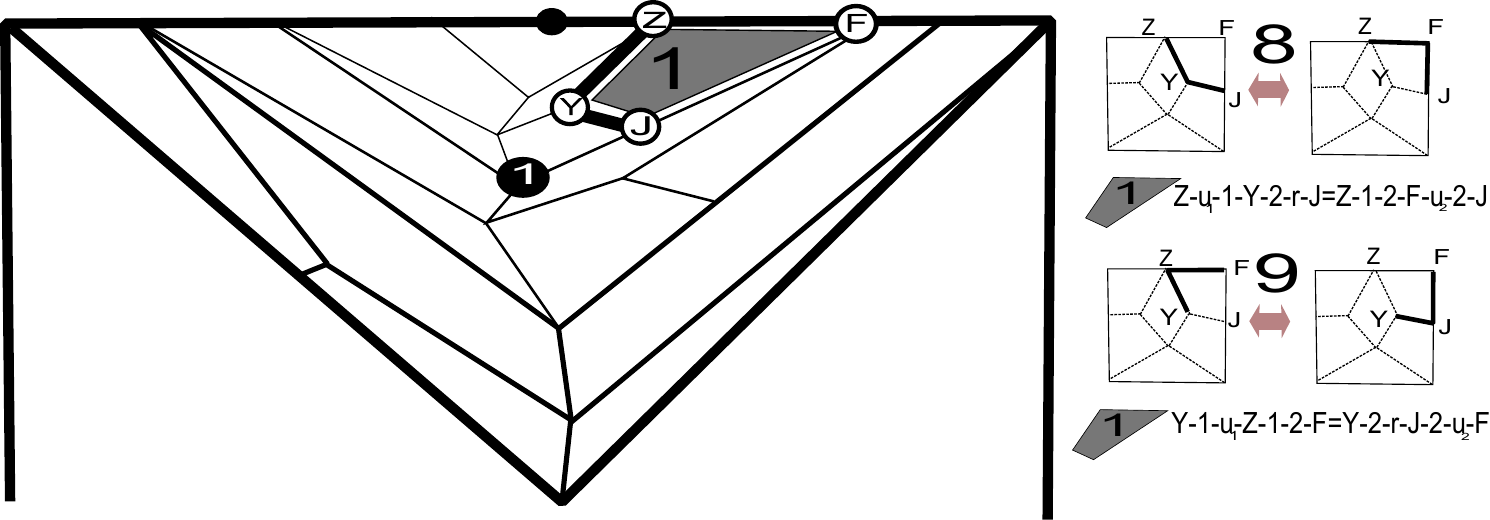}
\caption{Случаи расположения пути вокруг D2-цепи  и соответствующие им определяющие соотношения для локальных преобразований 8 и 9}
\label{D2b}
\end{figure}

В правой части рисунка~\ref{D2b} изображены локальные преобразования $8$ и $9$.

\medskip

Характеризация и восстановление кода полностью аналогично случаю $\mathbb{DL}2$ или  $\mathbb{DR}2$ цепи.

\subsection{Случай цепи $\mathbb{D}2$; преобразования 7 и 10}

\begin{figure}[hbtp]
\centering
\includegraphics[width=1\textwidth]{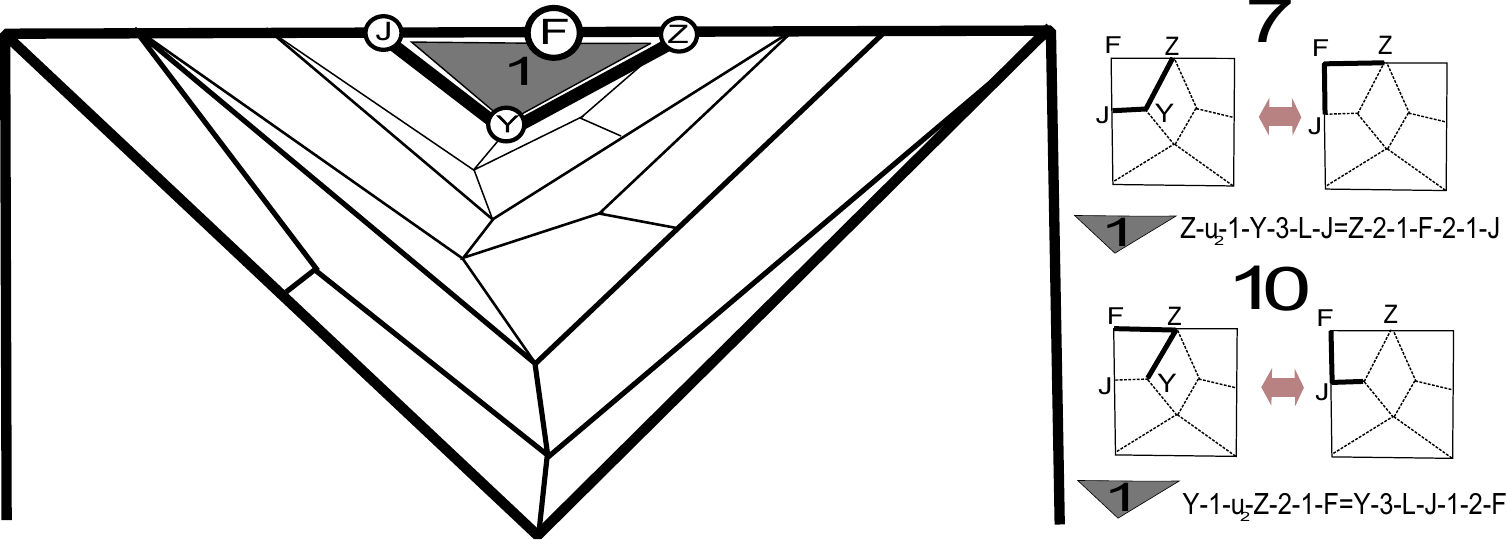}
\caption{Случаи расположения пути вокруг D2-цепи и соответствующие им определяющие соотношения для локальных преобразований 7 и 10}
\label{D2a}
\end{figure}

В правой части рисунка~\ref{D2a} изображены локальные преобразования $7$ и $10$.

\medskip

Характеризация и восстановление кода полностью аналогично случаю $\mathbb{DL}2$ или  $\mathbb{DR}2$ цепи.
\medskip

\subsection{Случай цепи $\mathbb{U}1$; преобразования 8 и 9}

\begin{figure}[hbtp]
\centering
\includegraphics[width=0.9\textwidth]{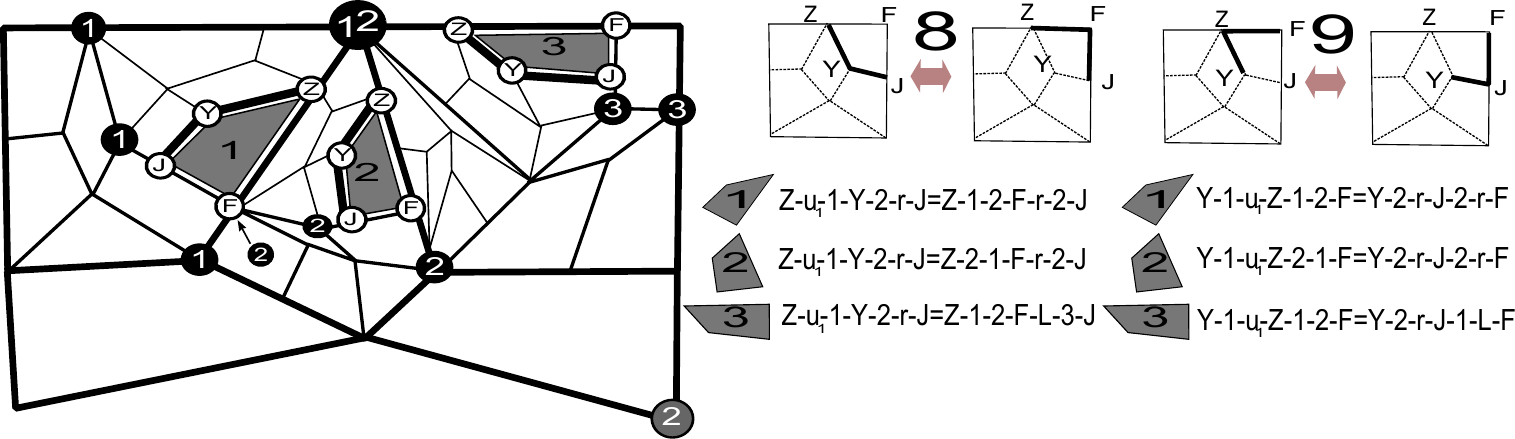}
\caption{Случаи расположения пути вокруг U1-цепи  и соответствующие им определяющие соотношения для 8 и 9 локальных преобразований}
\label{U1b}
\end{figure}

В правой части рисунка~\ref{U1b} изображены локальные преобразования $8$ и $9$.

\medskip

Характеризация и восстановление кода полностью аналогично случаю $\mathbb{UL}1$ или  $\mathbb{UR}1$ цепи.

\medskip

\subsection{Случай цепи $\mathbb{U}1$; преобразования 7 и 10}

\begin{figure}[hbtp]
\centering
\includegraphics[width=1\textwidth]{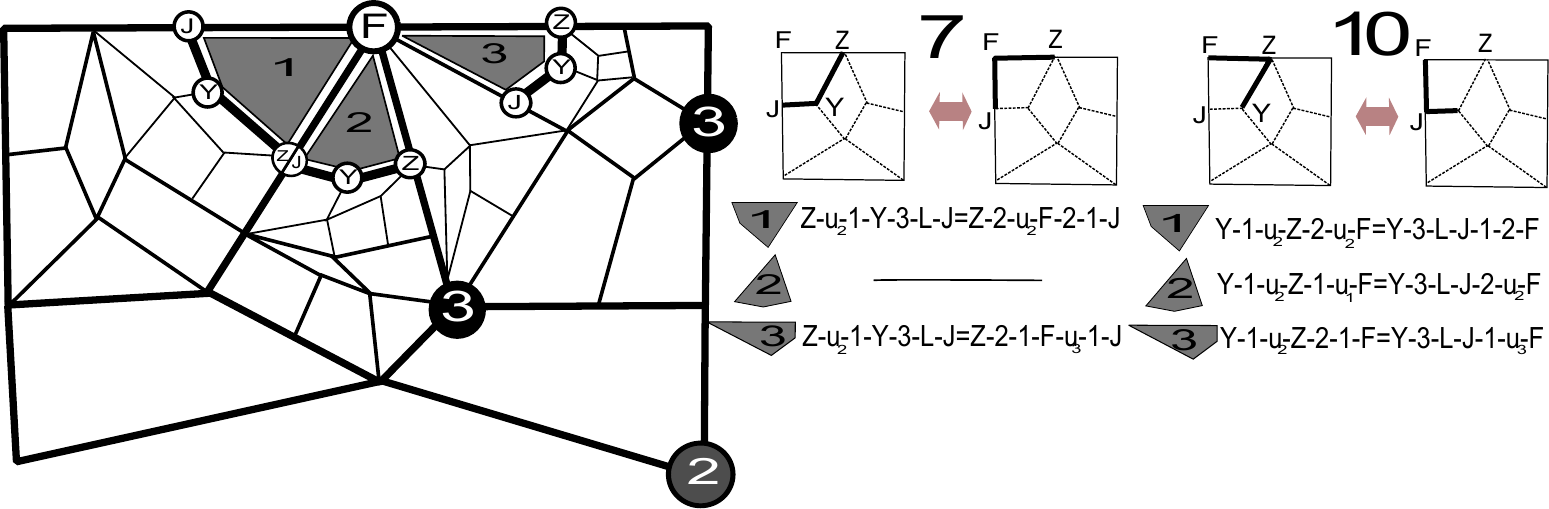}
\caption{Случаи расположения пути вокруг U1-цепи и соответствующие им определяющие соотношения для локальных преобразований 7 и 10}
\label{U1a}
\end{figure}

В правой части рисунка~\ref{U1a} изображены локальные преобразования $7$ и $10$.

\medskip

Характеризация и восстановление кода полностью аналогично случаю $\mathbb{UL}1$ или  $\mathbb{UR}1$ цепи.

\medskip

\subsection{Случай цепи $\mathbb{U}2$; преобразования 8 и 9}

\begin{figure}[hbtp]
\centering
\includegraphics[width=1\textwidth]{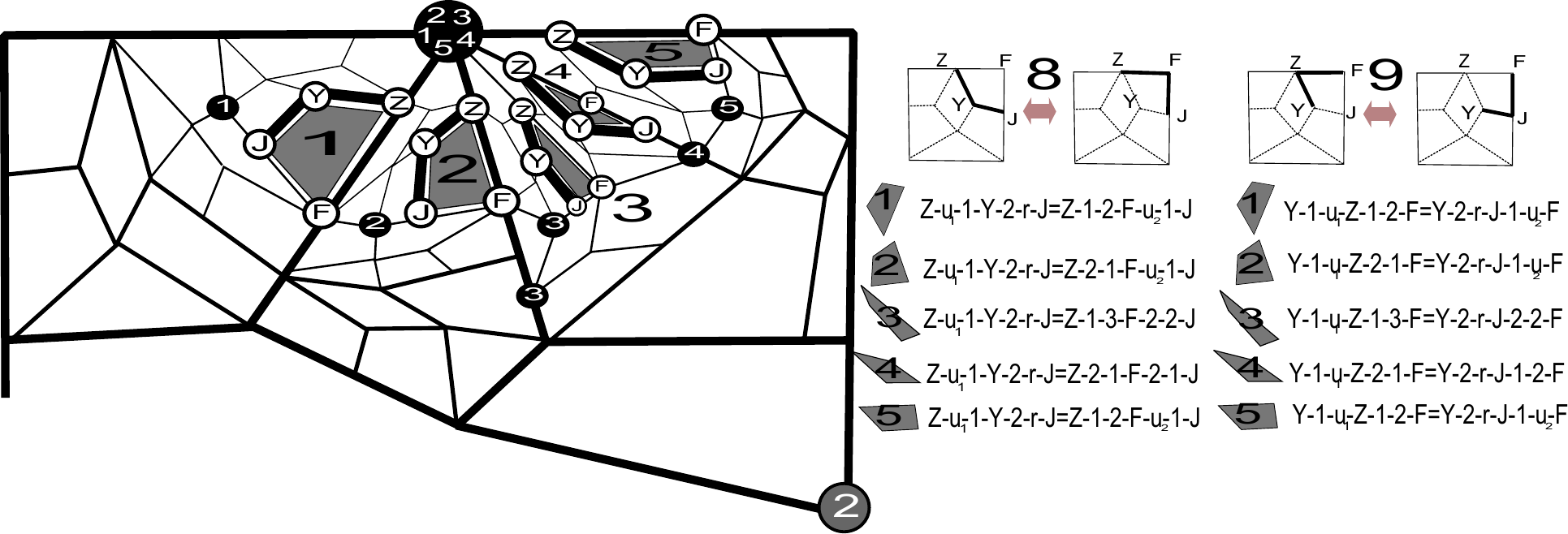}
\caption{Случаи расположения пути вокруг U2-цепи  и соответствующие им определяющие соотношения для локальных преобразований 8 и 9}
\label{U2b}
\end{figure}

В правой части рисунка~\ref{U2b} изображены локальные преобразования $8$ и $9$.

\medskip

Характеризация и восстановление кода полностью аналогично случаю $\mathbb{UL}2$ или  $\mathbb{UR}2$ цепи.

\medskip

\subsection{Случай цепи $\mathbb{U}2$; преобразования 7 и 10}

\begin{figure}[hbtp]
\centering
\includegraphics[width=0.9\textwidth]{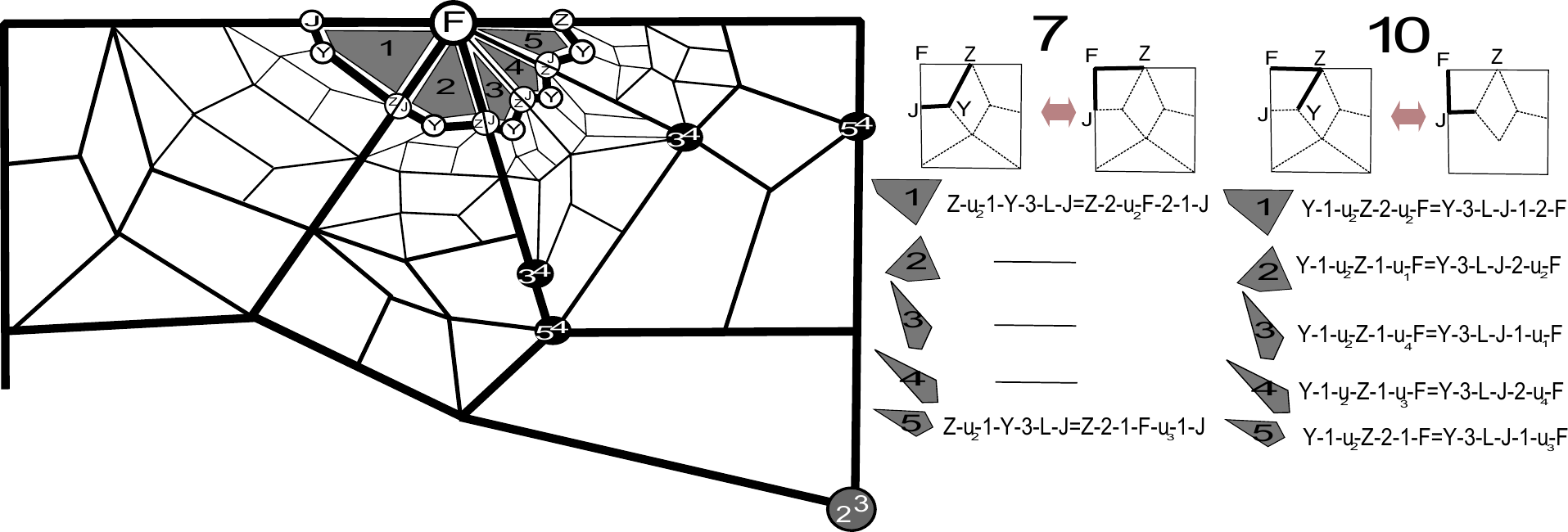}
\caption{Случаи расположения пути вокруг U2-цепи и соответствующие им определяющие соотношения для локальных преобразований 7 и 10}
\label{U2a}
\end{figure}

В правой части рисунка~\ref{U2a} изображены локальные преобразования $7$ и $10$.

\medskip

Характеризация и восстановление кода полностью аналогично случаю $\mathbb{UL}2$ или  $\mathbb{UR}2$ цепи.

\medskip

\subsection{Случай цепи $\mathbb{U}3$; преобразования 7, 8, 9, 10}

Случай $\mathbb{U}3$-цепи полностью аналогичен случаю $\mathbb{U}2$-цепи, соотношения выглядят идентично, только кодировки вершин $J$, $F$, $Z$, $Y$ отвечают $\mathbb{U}3$-цепи. Все рассуждения о вычислении путей полностью аналогичны. Соотношений вводится столько же, сколько для $\mathbb{U}2$ случая.

\medskip

\subsection{Случай цепи $\mathbb{L}1$; преобразования 8 и 9}

\begin{figure}[hbtp]
\centering
\includegraphics[width=0.9\textwidth]{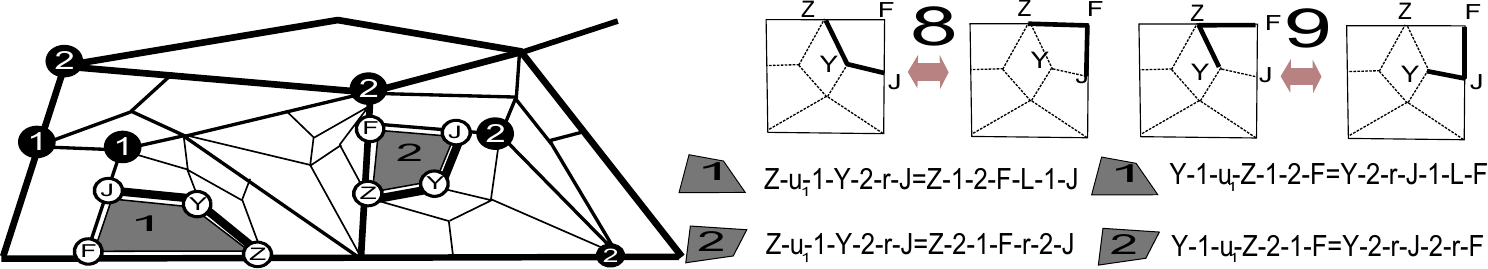}
\caption{Случаи расположения пути вокруг L1-цепи  и соответствующие им определяющие соотношения для локальных преобразований 8 и 9}
\label{L1b}
\end{figure}

В правой части рисунка~\ref{L1b} изображены локальные преобразования $8$ и $9$.

\medskip

Характеризация и восстановление кода полностью аналогично случаю $\mathbb{UL}1$ или  $\mathbb{DL}1$ цепи.

\medskip

\subsection{Случай цепи $\mathbb{L}1$; преобразования 7 и 10}

\begin{figure}[hbtp]
\centering
\includegraphics[width=1\textwidth]{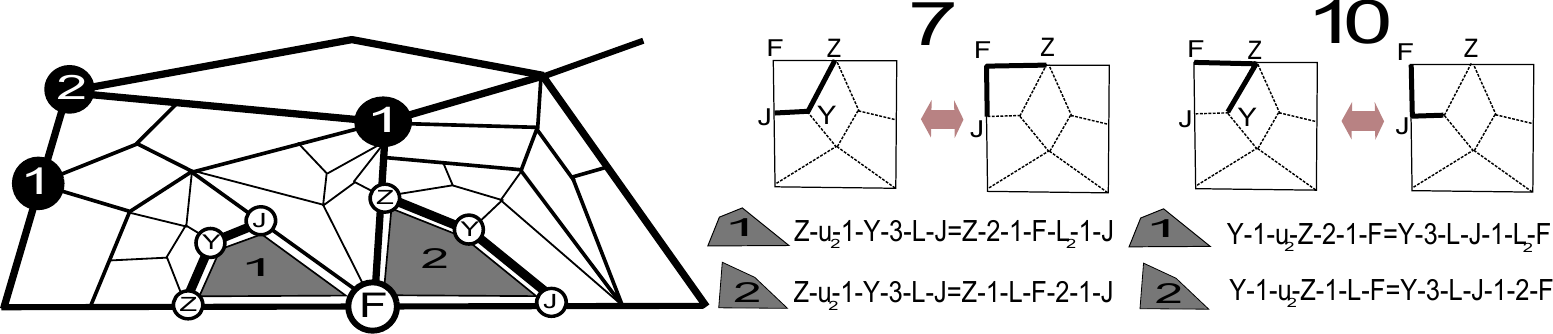}
\caption{Случаи расположения пути вокруг L1-цепи и соответствующие им определяющие соотношения для локальных преобразований 7 и 10}
\label{L1a}
\end{figure}

В правой части рисунка~\ref{L1a} изображены локальные преобразования $7$ и $10$.

\medskip

Характеризация и восстановление кода полностью аналогично случаю $\mathbb{UL}1$ или  $\mathbb{DL}1$ цепи.

\medskip

\subsection{Случай цепи $\mathbb{L}2$; преобразования 8 и 9}

\begin{figure}[hbtp]
\centering
\includegraphics[width=1\textwidth]{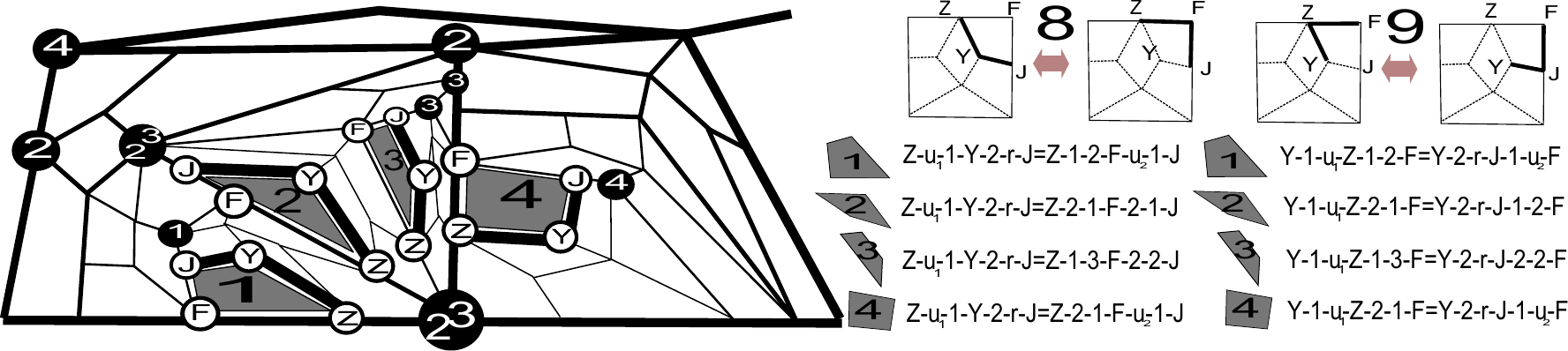}
\caption{Случаи расположения пути вокруг L2-цепи  и соответствующие им определяющие соотношения для локальных преобразований 8 и 9}
\label{L2b}
\end{figure}

В правой части рисунка~\ref{L2b} изображены локальные преобразования $8$ и $9$.

\medskip

Характеризация и восстановление кода полностью аналогично случаю $\mathbb{UL}2$ или  $\mathbb{DL}2$ цепи.

\medskip

\subsection{Случай цепи $\mathbb{L}2$, преобразования 7 и 10}

\begin{figure}[hbtp]
\centering
\includegraphics[width=1\textwidth]{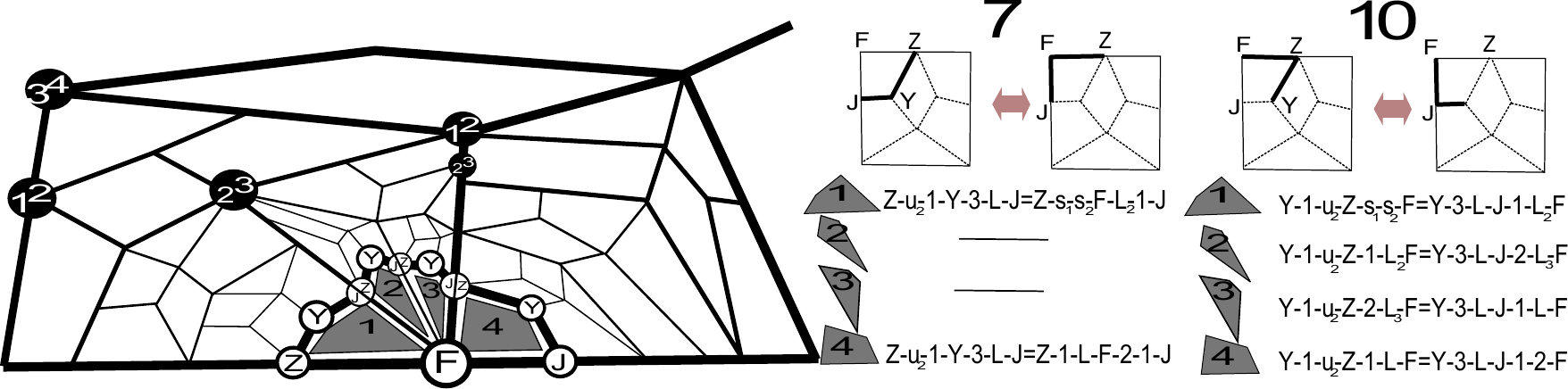}
\caption{Случаи расположения пути вокруг L2-цепи и определяющие соотношения для локальных преобразований 7 и 10}
\label{L2a}
\end{figure}

В правой части рисунка~\ref{L2a} изображены локальные преобразования $7$ и $10$.

\medskip

Характеризация и восстановление кода полностью аналогично случаю $\mathbb{UL}2$ или  $\mathbb{DL}2$ цепи.

\medskip

\subsection{Случай цепи $\mathbb{L}3$; преобразования 7, 8, 9, 10}

Случай $\mathbb{L}3$-цепи полностью аналогичен случаю $\mathbb{L}2$-цепи, соотношения выглядят идентично, только кодировки вершин $J$, $F$, $Z$, $Y$ отвечают $\mathbb{L}3$-цепи. Все рассуждения о вычислении путей полностью аналогичны. Соотношений вводится столько же, сколько для $\mathbb{L}2$ случая.

\medskip

\subsection{Случай цепи $\mathbb{R}1$; преобразования 8 и 9}

\begin{figure}[hbtp]
\centering
\includegraphics[width=1\textwidth]{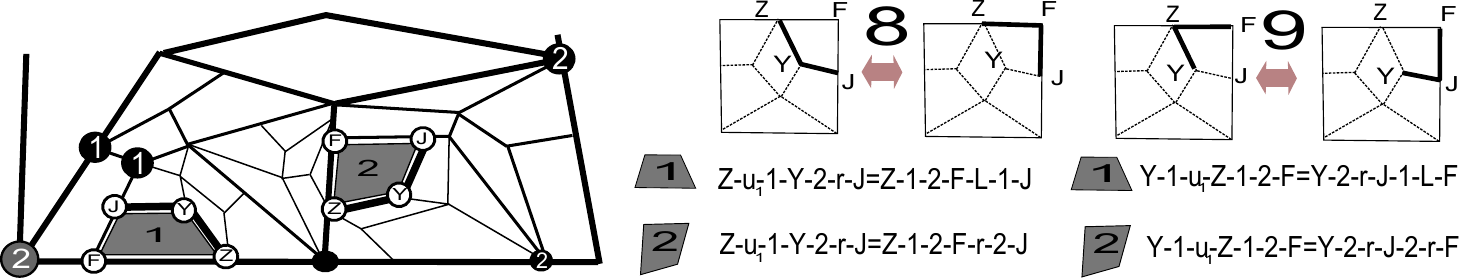}
\caption{Случаи расположения пути вокруг R1-цепи и определяющие соотношения для локальных преобразований 8 и 9}
\label{R1b}
\end{figure}

В правой части рисунка~\ref{R1b} изображены локальные преобразования $8$ и $9$.

\medskip

Характеризация и восстановление кода полностью аналогично случаю $\mathbb{DR}1$ или  $\mathbb{UR}1$ цепи.

\medskip

\subsection{Случай цепи $\mathbb{R}1$; преобразования 7 и 10}

\begin{figure}[hbtp]
\centering
\includegraphics[width=0.9\textwidth]{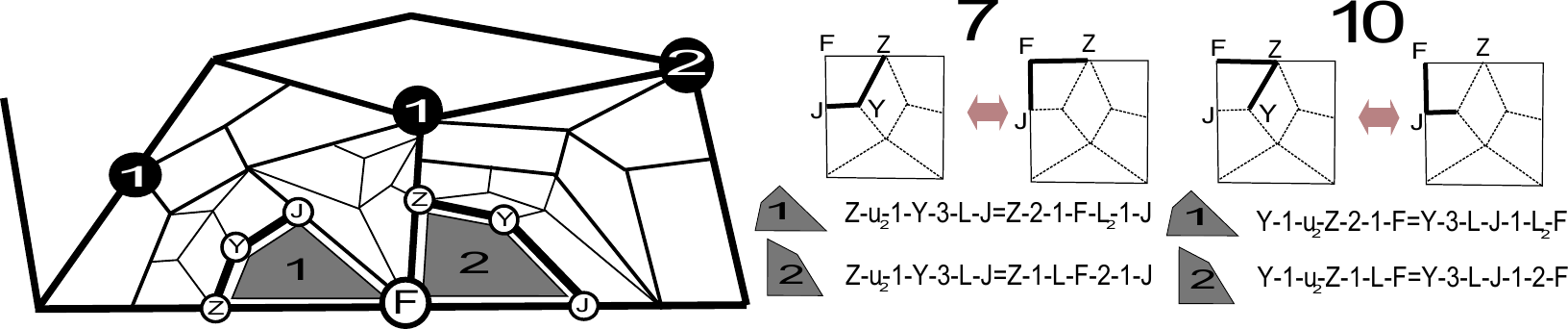}
\caption{Случаи расположения пути вокруг R1-цепи и определяющие соотношения для локальных преобразований 7 и 10}
\label{R1a}
\end{figure}

В правой части рисунка~\ref{R1a} изображены локальные преобразования $7$ и $10$.

\medskip

Характеризация и восстановление кода полностью аналогично случаю $\mathbb{DR}1$ или  $\mathbb{UR}1$ цепи.

\medskip

\subsection{Случай цепи $\mathbb{R}2$; преобразования 8 и 9}

\begin{figure}[hbtp]
\centering
\includegraphics[width=1\textwidth]{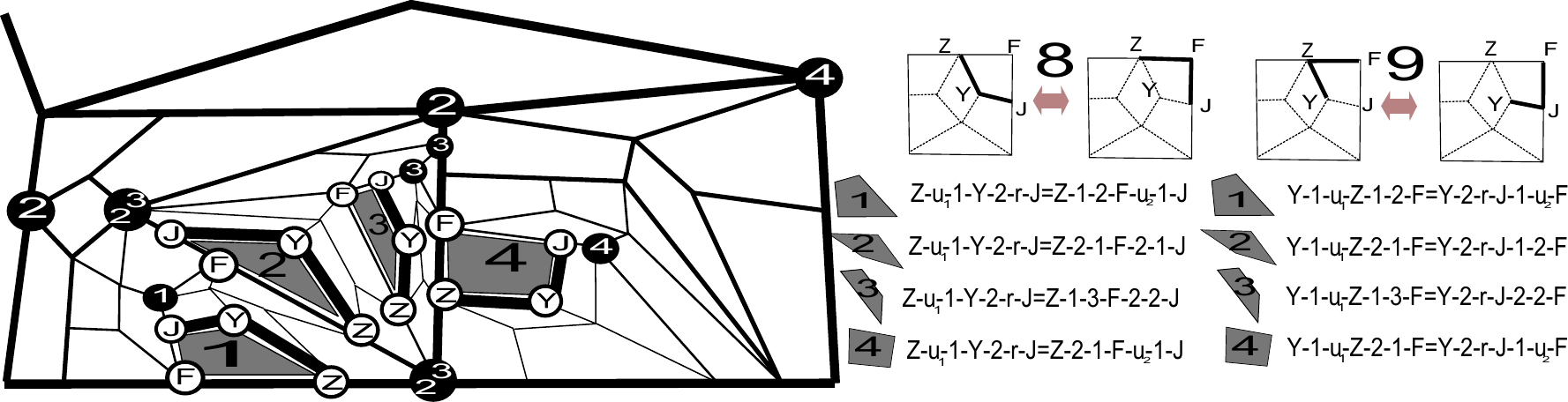}
\caption{Случаи расположения пути вокруг R2-цепи и определяющие соотношения для локальных преобразований 8 и 9}
\label{R2b}
\end{figure}

В правой части рисунка~\ref{R2b} изображены локальные преобразования $8$ и $9$.

\medskip

Характеризация и восстановление кода полностью аналогично случаю $\mathbb{DR}2$ или  $\mathbb{UR}2$ цепи.

\medskip

\subsection{Случай цепи $\mathbb{R}2$; преобразования 7 и 10}

\begin{figure}[hbtp]
\centering
\includegraphics[width=0.9\textwidth]{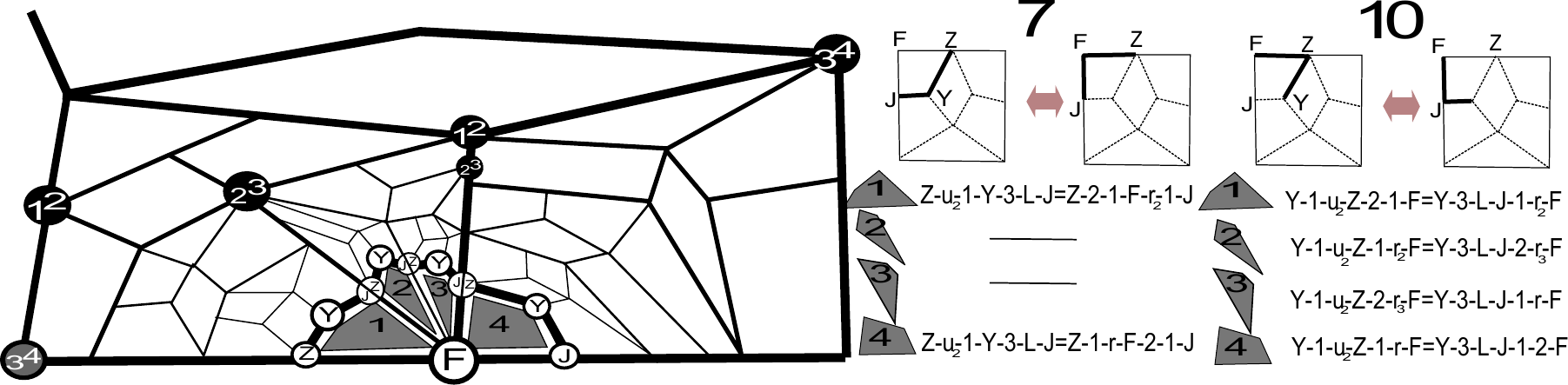}
\caption{Случаи расположения пути вокруг R2-цепи и определяющие соотношения для локальных преобразований 7 и 10}
\label{R2a}
\end{figure}

В правой части рисунка~\ref{R2a} изображены локальные преобразования $7$ и $10$.

\medskip

Характеризация и восстановление кода полностью аналогично случаю $\mathbb{DR}2$ или  $\mathbb{UR}2$ цепи.

\medskip

\subsection{Случай цепи $\mathbb{R}3$; преобразования 7, 8, 9, 10}

Случай $\mathbb{R}3$-цепи полностью аналогичен случаю $\mathbb{R}2$-цепи, соотношения выглядят идентично, только кодировки вершин $J$, $F$, $Z$, $Y$ отвечают $\mathbb{R}3$-цепи. Все рассуждения о вычислении путей полностью аналогичны. Соотношений вводится столько же, сколько для $\mathbb{R}2$ случая.

\medskip

\subsection{Случай цепи $\mathbb{CUL}0$; преобразования 7, 8, 9, 10}

\begin{figure}[hbtp]
\centering
\includegraphics[width=0.9\textwidth]{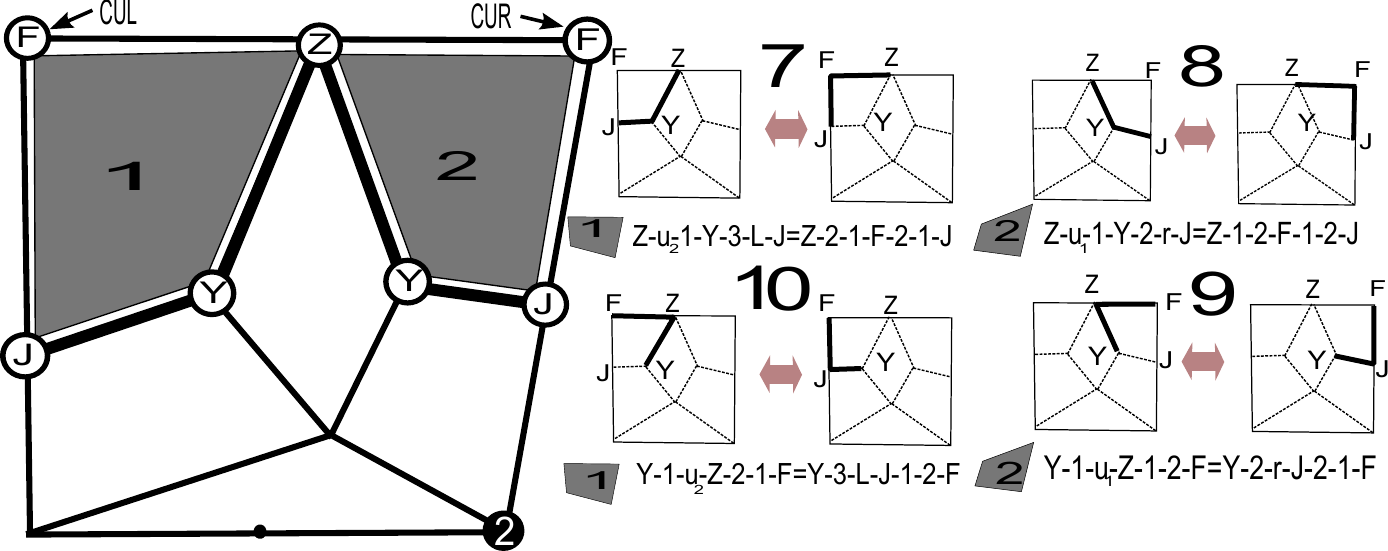}
\caption{Случаи расположения пути вокруг $\mathbb{CUL}0$-цепи и определяющие соотношения для локальных преобразований $7$, $8$, $9$, $10$}
\label{CUL0ab}
\end{figure}

В правой части рисунка~\ref{CUL0ab} изображены локальные преобразования $7$, $8$, $9$ и $10$.

\medskip

Характеризация очевидна.  Восстановление кода четвертой вершины по известным трем другим также ясно во всех случаях.

\medskip

\subsection{Случай цепи $\mathbb{CUL}1$; преобразования  7, 8, 9, 10}

\begin{figure}[hbtp]
\centering
\includegraphics[width=0.9\textwidth]{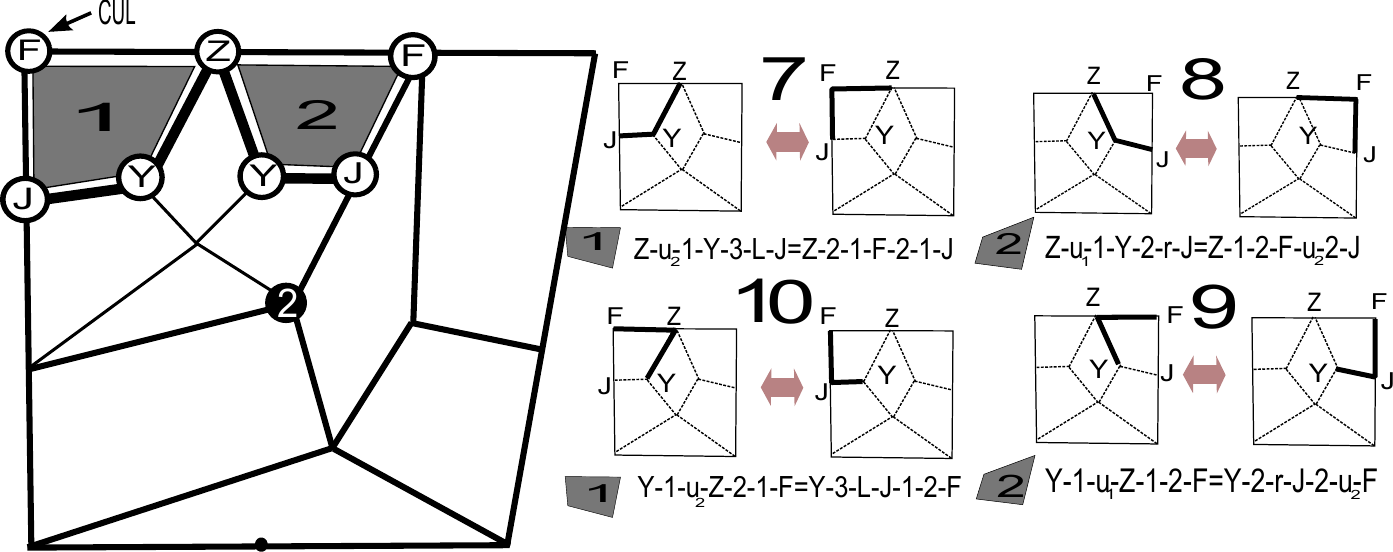}
\caption{Случаи расположения пути вокруг $\mathbb{CUL}1$-цепи и определяющие соотношения для локальных преобразований $7$, $8$, $9$, $10$}
\label{CUL1ab}
\end{figure}

В правой части рисунка~\ref{CUL1ab} изображены локальные преобразования $7$, $8$, $9$ и $10$.

\medskip

Характеризация очевидна.  Восстановление кода четвертой вершины по известным трем другим также ясно во всех случаях.

\medskip

\subsection{Случай цепи $\mathbb{CUR}1$; преобразования  7, 8, 9, 10}

\begin{figure}[hbtp]
\centering
\includegraphics[width=0.9\textwidth]{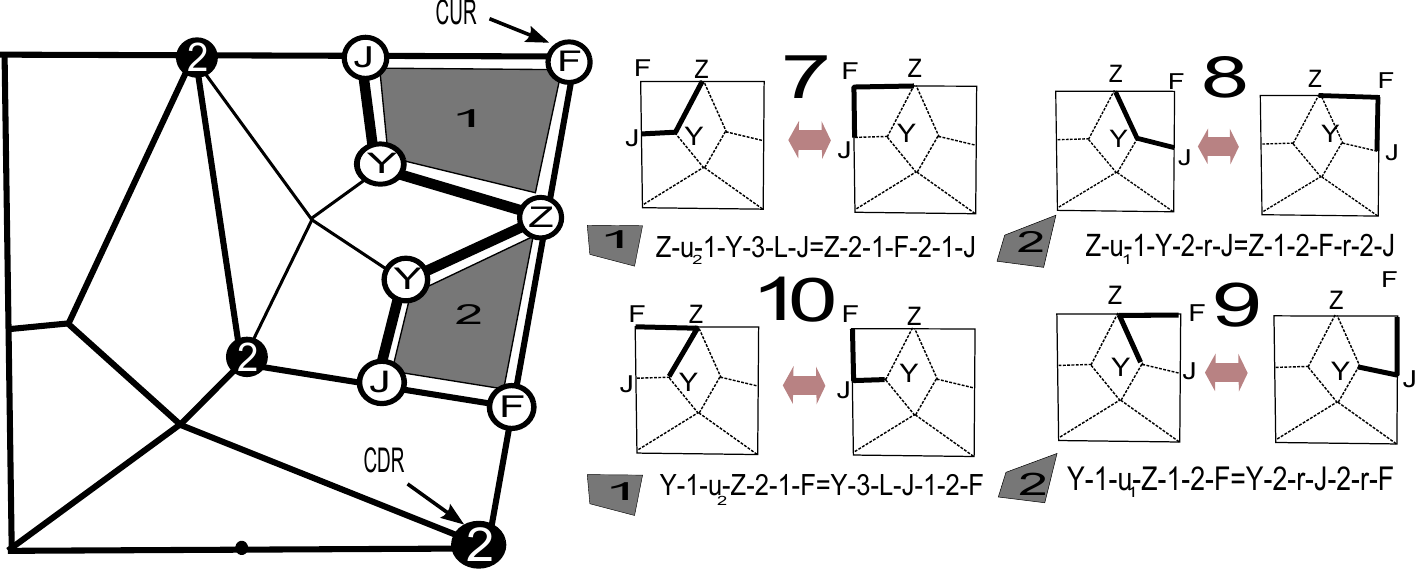}
\caption{Случаи расположения пути вокруг $\mathbb{CUR}1$-цепи и определяющие соотношения для локальных преобразований $7$, $8$, $9$, $10$}
\label{CUR1ab}
\end{figure}

В правой части рисунка~\ref{CUR1ab} изображены локальные преобразования $7$, $8$, $9$ и $10$.

\medskip

Характеризация очевидна.  Восстановление кода четвертой вершины по известным трем другим также ясно во всех случаях.

\medskip

\subsection{Случай цепи $\mathbb{CUR}2$; преобразования 7, 8, 9, 10}

\begin{figure}[hbtp]
\centering
\includegraphics[width=0.9\textwidth]{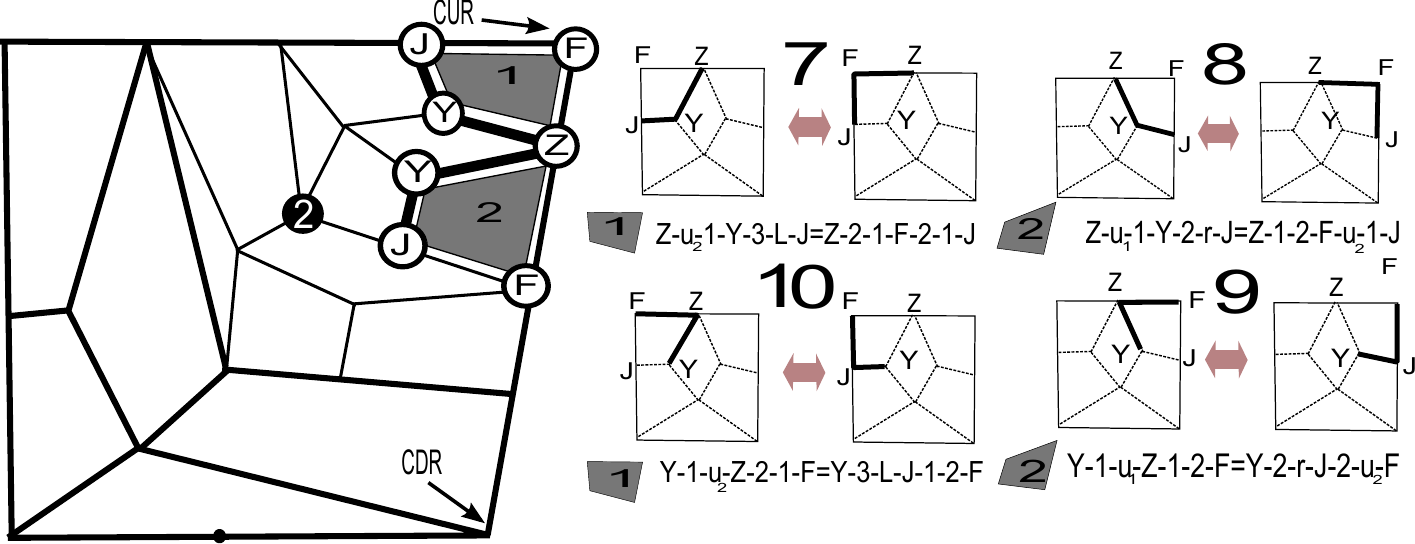}
\caption{Случаи расположения пути вокруг $\mathbb{CUR}2$-цепи  и соответствующие им определяющие соотношения для локальных преобразований $7$, $8$, $9$, $10$}
\label{CUR2ab}
\end{figure}

В правой части рисунка~\ref{CUR2ab} изображены локальные преобразования $8$ и $9$.

\medskip

Характеризация очевидна.  Восстановление кода четвертой вершины по известным трем другим также ясно во всех случаях.

\medskip

\subsection{Случай цепи $\mathbb{CDL}1$; преобразования  7, 8, 9, 10}

\begin{figure}[hbtp]
\centering
\includegraphics[width=0.9\textwidth]{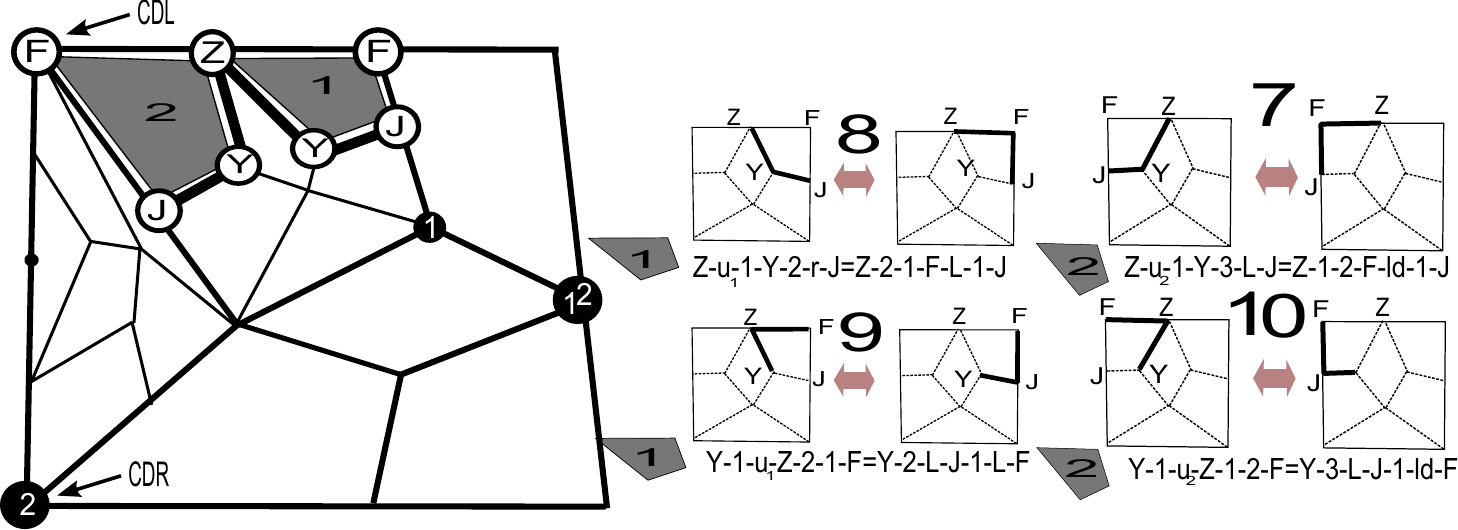}
\caption{Случаи расположения пути вокруг $\mathbb{CDL}1$-цепи и определяющие соотношения для локальных преобразований $7$, $8$, $9$, $10$}
\label{CDL1ab}
\end{figure}

В правой части рисунка~\ref{CDL1ab} изображены локальные преобразования $7$, $8$, $9$ и $10$.

\medskip

Характеризация очевидна. Восстановление кода ясно, так как у $F$ и $J$ -- общее множество начальников, у $Y$ первый начальник -- $Z$, и тип второго во втором случае  -- $\mathbb{A}$, а $Z$ -- краевая вершина без начальников.

\medskip

\subsection{Случай цепи $\mathbb{CDL}2$; преобразования 8 и 9}

\begin{figure}[hbtp]
\centering
\includegraphics[width=0.9\textwidth]{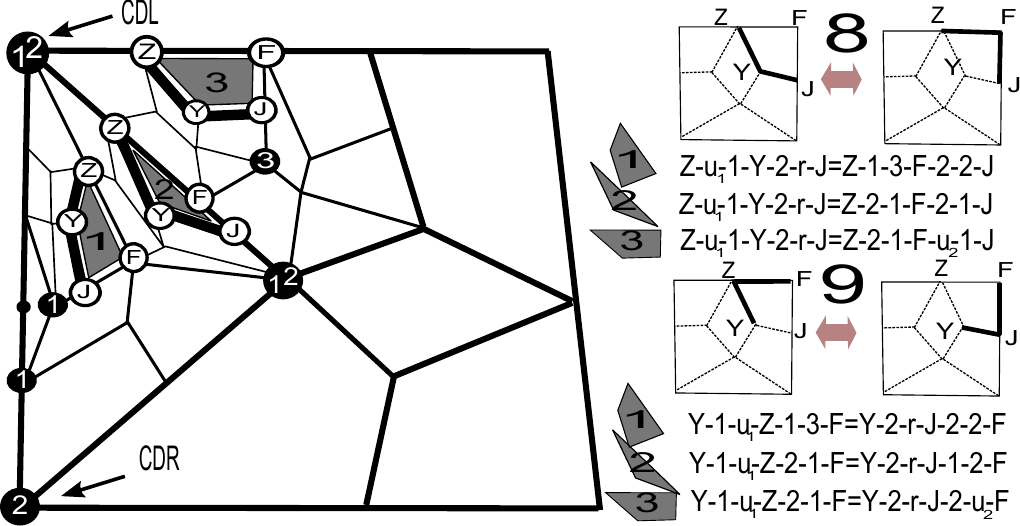}
\caption{Случаи расположения пути вокруг $\mathbb{CDL}2$-цепи и определяющие соотношения для  локальных преобразований 8 и 9}
\label{CDL2b}
\end{figure}

В правой части рисунка~\ref{CDL2b} изображены локальные преобразования $8$ и $9$.

\medskip

Характеризация очевидна. Восстановление кода также ясно, учитывая, что у $F$ и $Z$ всегда общее множество начальников, у $Y$ первый начальник  -- $Z$, а тип второго ясен из расположения. Начальники $J$ во всех случаях либо содержатся среди начальников $F$, либо сам $F$ является начальником.

\medskip

\subsection{Случай цепи $\mathbb{CDL}2$; преобразования 7 и 10}

\begin{figure}[hbtp]
\centering
\includegraphics[width=0.9\textwidth]{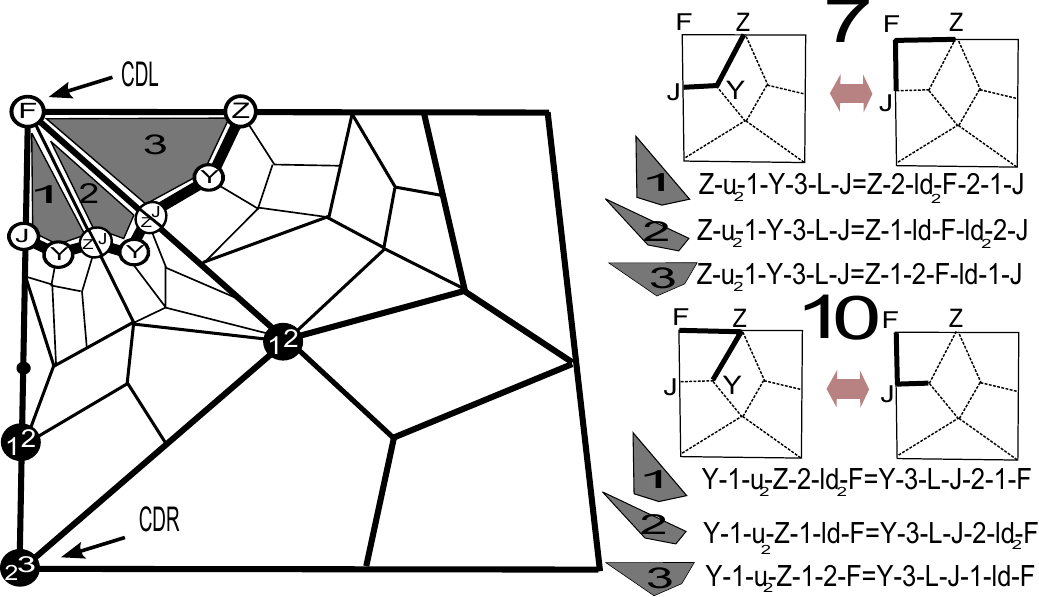}
\caption{Случаи расположения пути вокруг $\mathbb{CDL}2$-цепи и определяющие соотношения для локальных преобразований 7 и 10}
\label{CDL2a}
\end{figure}

В правой части рисунка~\ref{CDL2a} изображены локальные преобразования $7$ и $10$.

\medskip

Характеризация очевидна. Для доказательства восстановления кода достаточно заметить, что окружения четырех вершин в каждом случае очевидны, а окружение начальников (вершины, отмеченных черными кругами), можно явным образом написать.

\medskip

\subsection{Случай цепи $\mathbb{CDL}3$; преобразования 7 и 10}

Случай $\mathbb{CDL}3$-цепи полностью аналогичен случаю $\mathbb{CDL}2$-цепи, соотношения выглядят идентично, только кодировки вершин $J$, $F$, $Z$, $Y$ отвечают $\mathbb{CDL}3$-цепи. Все рассуждения о вычислении путей полностью аналогичны. Соотношений вводится столько же, сколько для $\mathbb{CDL}2$ случая.

\medskip

\subsection{Случай цепи $\mathbb{CDR}1$; преобразования 8 и 9}

\begin{figure}[hbtp]
\centering
\includegraphics[width=0.9\textwidth]{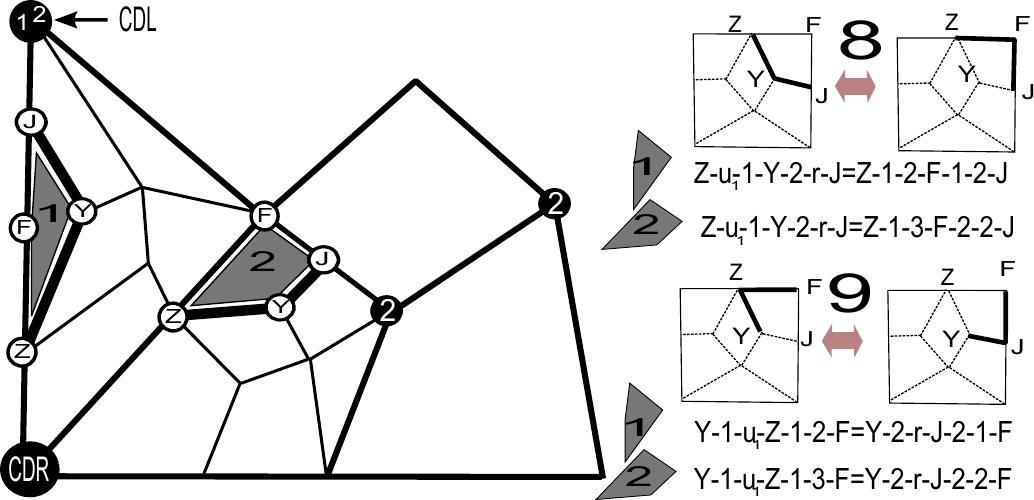}
\caption{Случаи расположения пути вокруг $\mathbb{CDR}1$-цепи и определяющие соотношения для локальных преобразований 8 и 9}
\label{CDR1b}
\end{figure}

В правой части рисунка~\ref{CDR1b} изображены локальные преобразования $8$ и $9$.

\medskip

Характеризация очевидна. Для восстановления кода достаточно заметить, что начальник каждой вершины либо является другой вершиной среди четырех, либо является начальником другой вершины из четырех.

\medskip

\subsection{Случай цепи $\mathbb{CDR}1$; преобразования 7 и 10}

\begin{figure}[hbtp]
\centering
\includegraphics[width=0.9\textwidth]{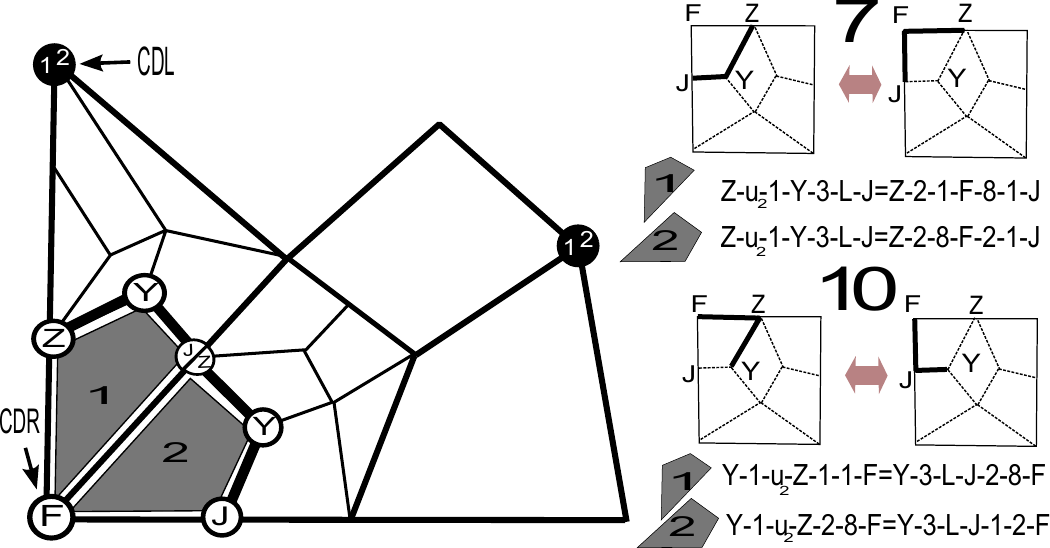}
\caption{Случаи расположения пути вокруг $\mathbb{CDR}1$-цепи и соответствующие им определяющие соотношения для локальных преобразований 7 и 10}
\label{CDR1a}
\end{figure}

В правой части рисунка~\ref{CDR1a} изображены локальные преобразования $7$ и $10$.

\medskip

Характеризация очевидна. Для восстановления кода достаточно заметить, что начальник каждой вершины либо является другой вершиной среди четырех, либо является начальником другой вершины из четырех.

\medskip

\subsection{Случай цепи $\mathbb{CDR}2$; преобразования 8 и 9}

\begin{figure}[hbtp]
\centering
\includegraphics[width=0.9\textwidth]{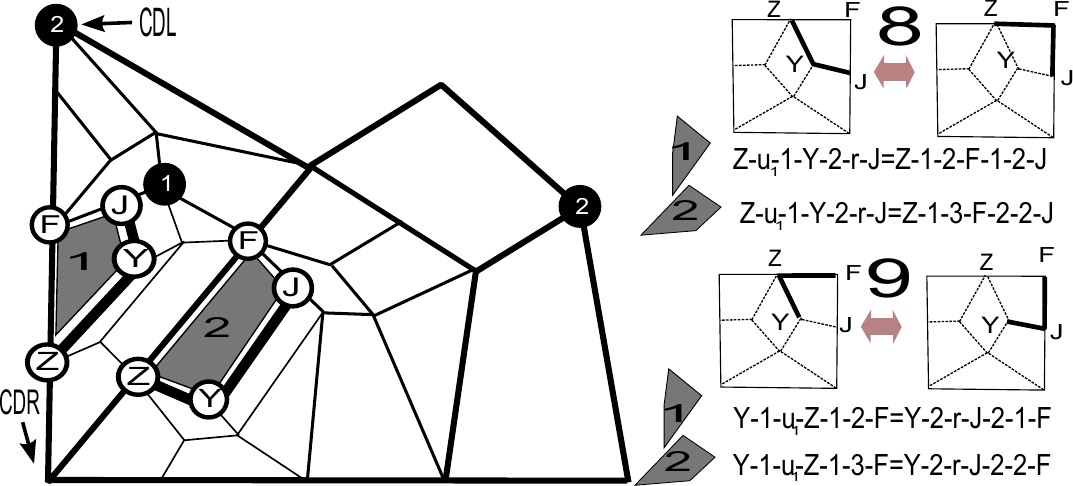}
\caption{Случаи расположения пути вокруг CDR2-цепи  и определяющие соотношения для локальных преобразований 8 и 9}
\label{CDR2b}
\end{figure}

В правой части рисунка~\ref{CDR2b} изображены локальные преобразования $8$ и $9$.

\medskip

Характеризация очевидна. Для восстановления кода достаточно заметить, что начальник каждой вершины либо является другой вершиной среди четырех, либо является начальником другой вершины из четырех.

\medskip

\subsection{Случай цепи $\mathbb{CDR}2$; преобразования 7 и 10}

\begin{figure}[hbtp]
\centering
\includegraphics[width=0.9\textwidth]{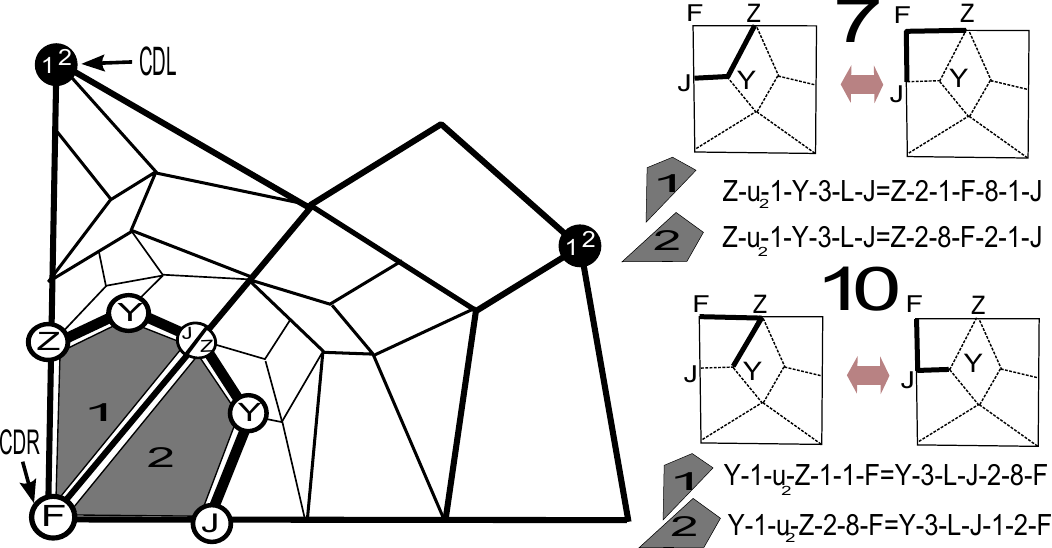}
\caption{Случаи расположения пути вокруг $\mathbb{CDR}2$-цепи и определяющие соотношения для локальных преобразований 7 и 10}
\label{CDR2a}
\end{figure}

В правой части рисунка~\ref{CDR2a} изображены локальные преобразования $7$ и $10$.

\medskip

Характеризация очевидна. Для восстановления кода достаточно заметить, что начальник каждой вершины либо является другой вершиной среди четырех, либо является начальником другой вершины из четырех.

\medskip

\subsection{Подсчет введенных соотношений}
Пусть $F$ -- число различных флагов макроплиток, $P$ -- число различных подклееных окружений, $I$ -- число различных информаций. Пусть также  $\mathbf{Num}(\mathbb{A})$ -- число базовых окружений вершин типа $\mathbb{A}$, аналогично для других типов.
Крайние вершины в пути могут иметь произвольное подклееное окружение, кроме того, соотношения вводятся для каждого заданного флага подклейки.

Для каждого из  локальных преобразований $(1)-(6)$ мы фиксировали один узел с окружением, флагом макроплитки и информацией. Также еще для одного узла могла быть зафиксирована информация. После чего, в каждом из шести случаев расположения мы вводили не более $2$ соотношений для каждого сочетания подклееных окружений крайних узлов пути и флага макроплитки. То есть общее число соотношений для локальных преобразований $(1)-(6)$ не превосходит $72 \cdot I\cdot F \cdot P\cdot \mathbf{Num}$, где $\mathbf{Num}$ -- число всевозможных сочетаний типа, уровня, окружения и информации плоской вершины.
То есть для локальных преобразований $(1)-(6)$ соотношений было введено не более

$$72\cdot 181 \cdot 10^6 \cdot 5 \cdot 10^{18} \cdot 39 \cdot 16 \cdot 10^{12} < 41 \cdot 10^{42}.$$

Теперь посчитаем число введенных соотношений для  локальных преобразований ~$7$, $8$, $9$, $10$.
Напомним, что мы рассматривали вершины из цепей, фиксировали для них информацию, флаг подклейки, а также два подклееных окружения. После этого вводилось не более двух соотношений для каждого локального преобразования.
Таким образом, число введенных соотношений для локальных преобразований ~$7$, $8$, $9$, $10$ не более $8 \cdot F \cdot P^2 \cdot N$, где $N$ -- общее число вершин всех типов, всех базовых окружений, всех возможных информаций.
Учитывая $F \le 16\cdot 10^{18}$ и $P=39$, получаем, что сочетаний базовых окружений и информаций не более $2749 \cdot 181 \cdot 10^6$. Значит число введенных соотношений для ~$7$, $8$, $9$, $10$ локальных преобразований не превосходит $13 \cdot 10^{33}$.

\section{Локальные преобразования при выходе в подклееную макроплитку} \label{pasting_section}

В предыдущей главе мы описали, как вводятся локальные преобразования (аналог определяющих соотношений) для плоских участков путей. В этом главе мы опишем, как это делается для путей, содержащих выход в подклееную макроплитку. То есть мы рассматриваем произвольный путь $X_1e_1e_2X_2e_3e_4X_3$, где $X_1$, $X_2$, $X_3$ -- буквы, отвечающие кодам вершин, а $e_1$, $e_2$, $e_3$, $e_4$ -- буквы, отвечающие ребрам входа и выхода. Хотя бы одно из ребер входа и выхода отвечает ребру, выходящему в подклейку.

\medskip

Для начала докажем следующее утверждение, позволяющее узнать ребра входа и выхода в пути из трех узлов, для случая когда средний узел является ядром некоторой подклейки.

\medskip

{\bf Ядро подклейки.} Пусть $Y$ -- ядро некоторой подклееной макроплитки $T$. Напомним, что в каждом узле $T$ хранится параметр {\it ``флаг подклейки''}, содержащий в себе полный код (тип, уровень, окружение и информацию) соответствующего ядра макроплитки $Y$, а также сочетания типов выходящих из $Y$ ребер, которые лежат на сторонах $T$. например, если $T$ подклеена по ребрам $1$ и $3$ к узлу типа $\mathbb{A}$ c окружением $(\mathbf{7A},\mathbf{7A},\mathbf{3B},\mathbf{4A})$, информацией $\mathbf{Info}$, то параметр {\it флаг подклейки} для узлов $T$ выглядит как $[\mathbb{A},\mathbf{Info},(\mathbf{7A},\mathbf{7A},\mathbf{3B},\mathbf{4A}),1,3]$.
Этот параметр (для узла $X$ внутри $T$) обозначается как $\mathbf{Core}(X)$. Тип выходящего из ядра ребра, соответствующего верхней стороне $T$ обозначается как $\mathbf{T.Core}(X)$. Соответственно, для ребра являющегося левой стороной $T$:  $\mathbf{L.Core}(X)$. Для окружения или типа $Y$ мы используем напрямую обозначение $\mathbf{Core}(X)$.

\medskip

\begin{proposition}[О ребрах около ядра макроплитки.]   \label{tilecore}

Пусть $XYZ$ -- некоторый плоский путь, причем $Y$ -- узел-ядро подклееной макроплитки, на сторонах которой лежат $X$ и $Z$. Тогда ребро выхода из $X$ (это всегда ребро $1$ или ребро $2$) можно установить, зная ребро входа в $Y$. Аналогично, ребро входа в $Z$ устанавливается по известному ребру выхода из $Y$.
\end{proposition}

Напомним, у каждого внутреннего ребра есть А-сторона и B-сторона. При этом первое главное ребро -- то, для которого A-сторона лежит справа.
В соответствии с этим правилом и с определением кодировок ребер входов и выходов можно выписать ребра выхода из $X$ и входа в $Z$ (приведены таблицы~\ref{tableAexits},~\ref{tableBexits},~\ref{tableCexits},~\ref{tableSideexits} для разных типов $Y$):

\begin{table}[hbtp]
\caption{Ребра выхода для $\mathbb{A}$-узла}
\centering
 \begin{tabular}{|c|c|c|c|c|}  \hline
ребро входа в $Y$    & ребро выхода из $X$  \cr
(или выхода из $Y$)  & (или входа в $Z$)  \cr  \hline
$1$ & $1$  \cr \hline
$2$ & $2$  \cr \hline
$3$ & $2$  \cr \hline
$\mathbf{lu}$  & $2$  \cr \hline
$\mathbf{ld}$ & $2$  \cr \hline
  \end{tabular}
\label{tableAexits}
\end{table}

\begin{table}[hbtp]
\caption{Ребра выхода для $\mathbb{B}$-узла}
\centering
 \begin{tabular}{|c|c|c|c|c|}  \hline
ребро входа в $Y$    & ребро выхода из $X$    \cr
(или выхода из $Y$)  & (или входа в $Z$)  \cr  \hline
$1$ & $2$  \cr \hline
$2$ & $1$  \cr \hline
$3$ & $1$  \cr \hline
$\mathbf{ru}$ & $2$  \cr \hline
$\mathbf{mid}$ & $2$  \cr \hline
$\mathbf{rd}$ & $2$  \cr \hline
  \end{tabular}
\label{tableBexits}
\end{table}

\begin{table}[hbtp]
\caption{Ребра выхода для $\mathbb{C}$-узла}
\centering
 \begin{tabular}{|c|c|c|c|c|}  \hline
ребро входа в $Y$    & ребро выхода из $X$    \cr
(или выхода из $Y$)  & (или входа в $Z$)  \cr  \hline
$1$ & $1$  \cr \hline
$2$ & $2$  \cr \hline
$3$ & $1$  \cr \hline
$4$ & $2$  \cr \hline
$\mathbf{ld}_1$ & $1$  \cr \hline
$\mathbf{ld}_2$ & $2$  \cr \hline
$\mathbf{mid}_1$ & $1$  \cr \hline
$\mathbf{mid}_2$ & $2$  \cr \hline
$\mathbf{d}_1$ & $1$  \cr \hline
$\mathbf{d}_2$ & $2$  \cr \hline
  \end{tabular}
\label{tableCexits}
\end{table}

\begin{table}[hbtp]
\caption{Ребра выхода для $\mathbb{UL}$, $\mathbb{UR}$, $\mathbb{DR}$, $\mathbb{DL}$ узлов}
\centering
 \begin{tabular}{|c|c|c|c|c|}  \hline
ребро входа в $Y$    & ребро выхода из $X$    \cr
(или выхода из $Y$)  & (или входа в $Z$)  \cr  \hline
1 & 2  \cr \hline
2 & 1  \cr \hline
$\mathbf{u}_1$ & $1$  \cr \hline
$\mathbf{u}_2$ & $2$  \cr \hline
$\mathbf{u}_3$ & $1$  \cr \hline
$\mathbf{u}_4$ & $2$  \cr \hline
$\mathbf{l}$ & $1$  \cr \hline
$\mathbf{l}_2$ & $1$  \cr \hline
$\mathbf{l}_3$ & $2$  \cr \hline
$\mathbf{r}$ & $2$  \cr \hline
$\mathbf{r}_2$ & $1$  \cr \hline
$\mathbf{r}_3$ & $2$  \cr \hline
  \end{tabular}
\label{tableSideexits}
\end{table}

Таким образом, зная хотя бы один узел подклееной макроплитки $T$ (из тех, которые лежат внутри нее, то есть не на верхней или левой стороне), мы можем узнать ребра входа и выхода для путей вдоль левой или верхней стороны $T$. В дальнейшем, если нам потребуется выписать типы ребер вдоль левой или верхней стороны подклееной макроплитки, мы будем использовать ссылку на предложение~\ref{tilecore}.

\medskip

\begin{proposition}[Об информации узла рядом с ядром макроплитки]   \label{tilecore2}

Пусть $XY$ -- некоторый плоский путь, причем $Y$ -- узел-ядро подклееной макроплитки, на одной из сторон которой лежит $X$. Будем считать, что нам известен код $Y$ (тип, уровень, окружение и информация), а также тип входящего в $Y$ ребра и тип, уровень и окружение $X$. Тогда информация у $X$ может быть восстановлена по этим данным.
\end{proposition}

Доказательство. Допустим сначала, что $X$ имеет один из внутренних типов ($\mathbb{A}$, $\mathbb{B}$, $\mathbb{C}$).

Пусть тип $X$ это $\mathbb{A}$. Рассмотрим ребро входа в $Y$. Имеется пять возможных типов ребер для узла типа  $\mathbb{A}$, это ребра $1$, $2$, $3$, $\mathbf{lu}$ , $\mathbf{ld}$. Если ребро входа в $Y$ это $1$, $2$ или $3$, информация у $X$ будет такой же, как у $Y$. Для ребра $\mathbf{lu}$ , первым начальником $X$ будет $\mathbf{LevelPlus.FBoss}(Y)$, вторым -- $\mathbf{Next.FBoss}(Y)$, третьим -- сам узел $Y$. В случае $\mathbf{ld}$, первым начальником $X$ будет $\mathbf{BottomLeftChain}(Y)$, вторым -- $\mathbb{C}$, с окружением $Y$, третьим -- сам узел $Y$.

 \medskip

Пусть тип $Y$ это $\mathbb{B}$. Имеется шесть возможных ребер из узла типа $\mathbb{B}$, это ребра $1$, $2$, $3$, $\mathbf{ru}$, $\mathbf{rd}$, $\mathbf{mid}$. Если ребро входа в $Y$ имеет тип $1$, $2$ или $3$, то информация у $X$ такая же как у $Y$.
Для ребра $\mathbf{rd}$, первым начальником будет узел в середине ребра 7 для макроплитки $T$, где $Y$ выступает в роли  $\mathbb{B}$-узла. То есть это $1$-цепь вокруг правого нижнего угла $T$, с указателем соответствующим входу по ребру 7. Тип правого нижнего угла $T$ нам известен, так как это второй начальник для $Y$ (это входит в информацию для $Y$). Второй начальник $X$ -- это узел в середине правой стороны $T$. То есть первый начальник $Y$ это $\mathbf{TopFromRight.SBoss}(X)$. Третьим начальником $X$ будет сам узел $Y$.
Для ребра $\mathbf{ru}$, первым начальником будет $1$-цепь вокруг правого верхнего угла $T$, с указателем соответствующим входу по правому ребру. Тип правого верхнего угла $T$ мы можем определить, применив функцию $\mathbf{TopRightType}(Y)$. Второй начальник $X$ совпадает с первым начальником $Y$.  Третьим начальником $X$ будет сам узел $Y$.
Для ребра $\mathbf{mid}$, первым начальником будет $1$-цепь вокруг $\mathbb{A}$-узла $T$ (то есть окружение как у $Y$), с указателем $1$. Второй начальник $X$ это $\mathbb{C}$-узел $T$ (окружение как у $Y$).  Третьим начальником $X$ будет сам узел~$Y$.

 \medskip

Пусть тип $Y$ это $\mathbb{C}$. Имеется десять возможных ребер из узла типа $\mathbb{C}$ ($4$ главных и $6$ неглавных). Можно заметить, что во всех случаях окружения начальников $X$ могут быть восстановлены по информации и окружению $Y$, этот процесс полностью аналогичен описанному выше.

\medskip

Пусть теперь $X$ имеет боковой тип. Если входящее в $Y$ ребро -- главное, то $X$ и $Y$ лежат на одном ребре в некоторой макроплитке, и поэтому у них одинаковая информация.
Пусть входящее в $Y$ ребро -- неглавное. Тогда тип ребра может быть одним из следующих: $\mathbf{u}_1$, $\mathbf{u}_1$, $\mathbf{u}_1$, $\mathbf{u}_1$, $\mathbf{r}$, $\mathbf{r}_2$, $\mathbf{r}_3$, $\mathbf{l}$, $\mathbf{l}_2$, $\mathbf{l}_3$.
Для $\mathbf{u}_1$, $\mathbf{u}_2$ первым начальником $X$ будет $Y$, при этом для $\mathbf{u}_1$ тип второго начальника для $X$ может быть восстановлен по информации и окружению $Y$.
Для $\mathbf{l}$, первым и единственным начальником $X$ будет $\mathbf{Prev}(Y)$, а для $\mathbf{r}$ первый начальник $X$ может быть восстановлен по процедуре  $\mathbf{TopFromRight}(Y)$, тип второго -- по процедуре  $\mathbf{BottomRightTypeFromRight}(Y)$.
Для остальных ребер все начальники $X$ восстанавливаются аналогично случаям внутренних вершин $Y$.

\medskip

\subsection{Обзор случаев}

Пусть есть путь, участвующий в локальном преобразовании, лежащий в макроплитке $T$, причем одно из его ребер -- входящее в подклееную макроплитку. Нужно показать, как происходит преобразование во всех случаях.

\medskip

Заметим, что один конец ребра в подклееную область всегда лежит на верхней или левой стороне макроплитки.

Итак, нужно рассмотреть случаи локальных преобразований, где макроплитка $T$ примыкает к границе левой или верхней стороны подклееной макроплитки. Таким образом, нужно рассмотреть те из случаев $B1-B20$, которые включают граничные стороны $\mathbf{top}$ или $\mathbf{left}$, таблица~\ref{P1P10}.

\medskip
\begin{table}[hbtp]
\caption{Десять окружений макроплиток, содержащих граничные стороны $\mathbf{top}$ или $\mathbf{left}$.}
\centering
 \begin{tabular}{|c|c|c|c|}   \hline
\x{старое \cr обозначение} & \x{окружение \cr макроплитки} & & \x{новое \cr обозначение} \cr \hline
B1.  & $(\mathbf{left},\mathbf{top},\mathbf{right},\mathbf{bottom})$  & & P1 \cr \hline
B2. & $(\mathbf{left},\mathbf{top},\mathbf{1A},\mathbf{3A})$   & & P2 \cr \hline
B3.  & $(\mathbf{7A},x,\mathbf{3B},\mathbf{4A})$,   &$x=\mathbf{left},\mathbf{top}$ & P3    \cr \hline
B4.  & $(\mathbf{7A},x,\mathbf{1A},\mathbf{3A})$,  & $x=\mathbf{left},\mathbf{top}$  & P4 \cr \hline
B5. & $(\mathbf{top},\mathbf{right},\mathbf{6A},\mathbf{2A})$  & & P5 \cr \hline
B6. & $(\mathbf{top},\mathbf{right},\mathbf{1A},\mathbf{3A})$  &  & P6  \cr \hline
B11. & $(x,\mathbf{3B},\mathbf{6A},\mathbf{2A})$,   & $x=\mathbf{left},\mathbf{top}$ & P7 \cr  \hline
B12. & $(x,\mathbf{3B},\mathbf{1A},\mathbf{3A})$,   & $x=\mathbf{left},\mathbf{top}$ & P8 \cr  \hline
B17. & $(x,\mathbf{1A},\mathbf{6A},\mathbf{2A})$,  & $x=\mathbf{left},\mathbf{top}$ & P9  \cr \hline
B18. & $(x,\mathbf{1A},\mathbf{1A},\mathbf{3A})$,  & $x=\mathbf{left},\mathbf{top}$ & P10 \cr \hline
  \end{tabular}
\label{P1P10}
\end{table}

\medskip

Пронумеруем их как $P1-P10$ и рассмотрим каждый отдельно.

\medskip

\begin{figure}[hbtp]
\centering
\includegraphics[width=0.5\textwidth]{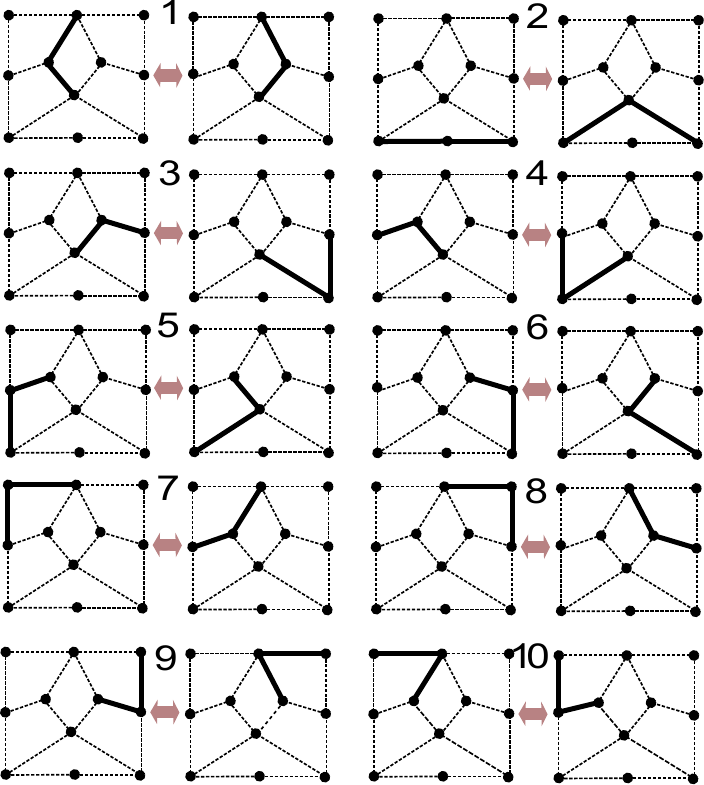}
\caption{Локальные преобразования}
\label{flippaths2}
\end{figure}

\medskip

Для каждого из этих случаев мы рассмотрим те локальные преобразования, где хотя бы один путь имеет узел на верхней или левой стороне. Это преобразования $1$, $2$, $4$, $5$, $7$, $8$, $9$, $10$ (рисунок~\ref{flippaths2}). В остальных случаях ($3$ и $6$) путь не может содержать ребро в подклееную область, так как пути не касаются левой или верхней стороны.

\medskip

{\bf Замечание.} Ниже мы не будем выписывать симметричное к введенному соотношение, отвечающее проходу пути в обратном порядке. Просто будем считать, что соотношений вводится в два раза больше.
Для экономии места в таблицах ниже, ``начальников'' будем называть ``боссами''.

\medskip

\subsection{P1: Макроплитка с окружением $(\mathbf{left},\mathbf{top},\mathbf{right},\mathbf{bottom})$}

В этом параграфе мы рассмотрим случай {\bf P1}, когда макроплитка $T$, внутри которой мы рассматриваем локальное преобразование пути, имеет окружение $(\mathbf{left},\mathbf{top},\mathbf{right},\mathbf{bottom})$. Это происходит, когда макроплитка сама по себе является только что подклееной макроплиткой.

Зафиксируем узел $Q$, являющийся ядром нашей подклееной макроплитки. Его подклееное окружение это $(\mathbf{left},\mathbf{top},\mathbf{right},\mathbf{bottom})$.  У него может быть произвольный тип, кроме углового, произвольное базовое окружение и произвольная информация. Если это боковой или краевой узел, то уровень у него должен быть третий (иначе бы к нему не была подклеена макроплитка).  Сторонами подклейки могут быть любые два ребра, выходящие из нашего узла.

Напомним, что в параметр {\it ядро подклейки} входит тип, уровень, расширенное окружение, информация и типы двух выходящих ребер, являющихся сторонами подклейки. То есть наше ядро подклейки задано.

\medskip

\begin{figure}[hbtp]
\centering
\includegraphics[width=0.9\textwidth]{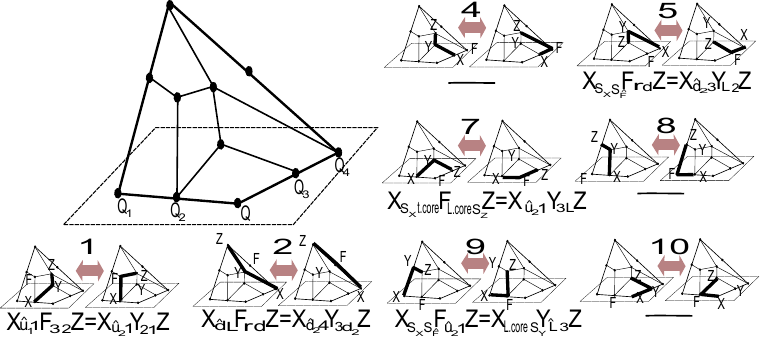}
\caption{Случай P1}
\label{P1}
\end{figure}

\medskip

Определив $Q$, можно также определить и узлы $Q_1$, $Q_2$, $Q_3$, $Q_4$, лежащие на сторонах подклейки. То есть, мы рассматриваем всевозможные сочетания кодов для пятерки вершин, существующих на комплексе. Для каждой такой пятерки, можно определить буквы коды вершин $X$, $Y$, $F$, $Z$, для каждого из указанных на рисунке~\ref{P1} восьми случаев. Также введем обозначения $e_{\mathbf{\mathbf{t.core}}}$ и $e_{\mathbf{\mathbf{l.core}}}$ для типов двух выходящих ребер, являющихся сторонами подклейки.

\medskip

{\bf Характеризация.} По выданному нам коду пути $P_1P_2P_3$, включающему переход по подклееному ребру, мы можем установить, встречается ли этот путь среди перечисленных на рисунке~\ref{P1}. Действительно, первый узел после перехода по подклееному ребру должен содержать параметр {\it ядро подклейки}, соответствующий рассматриваемому нами. Ясно также, что по конфигурации входящих-выходящих ребер и типов вершин можно выявить, какой именно случай из шестнадцати представленных на рисунке~\ref{P1} имеет место.

Итак, мы можем определить по коду пути к какому именно случаю расположения он относится. Теперь покажем, как по этому коду определить код парного к нему пути.

В случаях локальных преобразований $4$, $8$ и $10$ путь удовлетворяет условиям мертвого паттерна и не может быть частью достаточно длинного ненулевого пути, и для этих локальных преобразований мы соотношений не вводим.

\medskip

{\bf Вычисление кода $XFZ$.} Зная код пути $XYZ$ (то есть коды $X$, $Y$, $Z$, а также ребра входа в $Y$, $Z$ и ребра выхода из $X$ и $Y$ вдоль пути $XYZ$), можно выписать код $XFZ$.

Заметим, что коды $X$ и $Z$ мы и так знаем (из пути $XYZ$). Остается выписать ребра выхода из $X$, ребро входа в $F$, код $F$, ребро выхода из $F$ и ребро входа в $Z$ (таблица~\ref{P1a}).

\medskip

\begin{table}[hbtp]
\caption{Случай P1: вычисление кода $XFZ$. }
\centering
 \begin{tabular}{|c|c|c|c|c|c|c|}   \hline
  & ребро из $X$ & ребро в $F$ & \x{тип, уровень \cr и окружение $F$} & \x{информация \cr $F$} & \x{ ребро \cr из $F$ }& \x{ребро \cr в $Z$} \cr \hline
1. & $\widehat{\mathbf{u}_1}$ & $1$ & \x{ $\mathbb{B}$  (окруж. как у $Y$)} & 1 босс -- $X$ & $3$ & $2$  \cr \hline
2. & $\widehat{\mathbf{d}}$ & $\mathbf{l}$ & \x{$\mathbb{D}$, $1$ \cr (окруж. как у $Y$)} & -- & $\mathbf{r}$ & $\mathbf{d}$  \cr \hline
5. & Предл.\ref{tilecore} & Предл.\ref{tilecore} & $\mathbf{Next.FBoss}(Z)$ & как у $X$   & $\mathbf{l}$ & $3$  \cr \hline
7. & Предл.\ref{tilecore} & $\mathbf{T.Core}(Y)$ & $\mathbf{Core}(Y)$ & $\mathbf{Core}(Y)$ & $\mathbf{L.Core}(Y)$ & Предл.\ref{tilecore}   \cr \hline
9. & Предл.\ref{tilecore} & Предл.\ref{tilecore} & $\mathbf{FBoss}(Z)$& как у $X$  & $\widehat{\mathbf{u}_1}$ & $1$  \cr \hline
  \end{tabular}
\label{P1a}
\end{table}

\leftskip=0.0cm

\medskip

{\bf Примечание.} Через $\mathbf{T.Core}(Y)$ и $\mathbf{L.Core}(Y)$ мы обозначаем типы выходящих ребер, соответствующих сторонам подклееной макроплитки. Они содержатся в параметре флаг подклейки у вершины $Y$.

\medskip

{\bf Вычисление кода $XYZ$.} Зная код пути $XFZ$ (то есть коды $X$, $F$, $Z$, а также ребра входа в $F$, $Z$ и ребра выхода из $X$ и $F$ вдоль пути $XFZ$), можно выписать код $XYZ$.

Заметим, что коды $X$ и $Z$ мы и так знаем (из пути $XFZ$). Остается выписать ребра выхода из $X$, ребро входа в $Y$, код $Y$, ребро выхода из $Y$ и ребро входа в $Z$ (таблица~\ref{P1b}).

\medskip
\begin{table}[hbtp]
\caption{Случай P1: вычисление кода $XYZ$. }
\centering
 \begin{tabular}{|c|c|c|c|c|c|c|}   \hline
  & ребро из $X$ & ребро в $Y$ & \x{тип, уровень \cr и окружение $Y$} & информация $Y$ &  \x{ребро \cr из $Y$} & \x{ребро \cr в $Z$} \cr \hline
1. & $\widehat{\mathbf{u}_2}$ & $1$ & $\mathbb{A}$ (окружение как у $F$)& как у $F$ & $2$ & $1$  \cr \hline
2. & $\widehat{\mathbf{d}_2}$ & $4$ & $\mathbb{C}$ (окружение как у $F$)& $1$ босс -- как у $Z$  & $3$ & $\mathbf{d}_2$  \cr
&&&& $2$ босс -$X$, $3$ -$Z$  && \cr  \hline
5. & $\widehat{\mathbf{d}_2}$ & $3$ & $\mathbb{C}$ (окружение как у $Z$) & 1 босс как у $Z$   & $\mathbf{l}$ & $2$  \cr
&&&& $2$ босс -$X$, $3$ -- $\mathbb{CDR}$   &&\cr  \hline
7. & $\widehat{\mathbf{u}_2}$ & $1$ & $\mathbb{A}$, окружение -- & $1$ босс -- $X$ & $3$ & $\mathbf{l}$  \cr
&&& $(\mathbf{left},\mathbf{right},\mathbf{top},\mathbf{bottom})$  &&& \cr  \hline
9. &  $\widehat{\mathbf{r}}$ & $2$ & $\mathbb{R}$,$1$ окружение --  & --  & $\mathbf{r}$ & $2$  \cr
&&& $(\mathbf{left},\mathbf{right},\mathbf{top},\mathbf{bottom})$  &&& \cr  \hline
  \end{tabular}
\label{P1b}
\end{table}

\medskip

Когда все коды определены можно выписать определяющие соотношения, переводящие один путь в другой.

Обозначим результат применения процедуры, описанной в Предложении~\ref{tilecore} к вершине $Q$, связанной с этим ребром как $s_Q$. Как $s_{\widehat{Q}}$ обозначим противоположное главное ребро, то есть, если $s_Q=1$, то $s_{\widehat{Q}}=2$ и $s_Q=2$, то $s_{\widehat{Q}}=1$.  Обозначение $\mathbf{Past}(X)$ используем для обозначения подклееного окружения узла $X$.

\medskip

{\bf Соотношения.} Введем соотношения, реализующие описанные переходы.

\smallskip

$1)$  $X e_{\widehat{u_1}} e_1 F e_3 e_2 Z = X e_{\widehat{u_2}} e_1 Y e_2 e_1 Z   $

\medskip

$2)$  $X e_{\widehat{d}} e_l F e_r e_{d} Z = X e_{\widehat{d_2}} e_4 Y e_3 e_{d_2} Z   $

\medskip

$5)$  $X e_{s_X} e_{s_{\widehat{F}}} F e_l e_{3} Z = X e_{\widehat{d_2}} e_3 Y e_l e_{2} Z   $

\medskip

$7)$  $X e_{s_X} e_{t.core} F e_{l.core} e_{Z} Z = X e_{\widehat{u_2}} e_1 Y e_3 e_{l} Z   $

\medskip

$9)$  $X e_{s_X} e_{s_{\widehat{F}}} F e_{\widehat{u_1}} e_{1} Z = X e_{\widehat{l}} e_u Y e_r e_{2} Z   $

\medskip

\subsection{P2: Макроплитка с окружением $(\mathbf{left},\mathbf{top},\mathbf{1A},\mathbf{3A})$}

В этом параграфе мы рассмотрим случай {\bf P2}, когда макроплитка $T$, внутри которой мы рассматриваем локальное преобразование пути, имеет окружение $(\mathbf{left},\mathbf{top},\mathbf{1A},\mathbf{3A})$. Это происходит, когда макроплитка является прямым потомком подклееной макроплитки.

\begin{figure}[hbtp]
\centering
\includegraphics[width=1\textwidth]{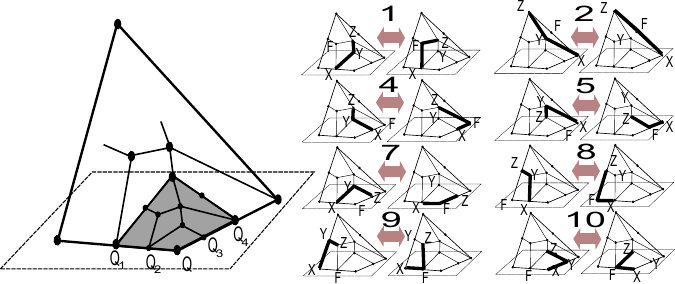}
\caption{Случай P2}
\label{P2}
\end{figure}

Аналогично случаю {\bf P1}, зафиксируем пять вершин $Q$, $Q_1$, $Q_2$, $Q_3$, $Q_4$, занимающих на комплексе положения как на рисунке~\ref{P2}. То есть $Q$ является ядром подклейки, а остальные вершины лежат на сторонах.

Для зафиксированного сочетания кодов этих пяти вершин, можно вычислить коды всех путей, указанных на рисунке~\ref{P2}. Кроме того, по выданному слову (коду пути) можно установить, действительно ли этот путь относится к случаю {\bf P2} и какую конфигурацию из изображенных на рисунке~\ref{P2} он имеет. Действительно по окружению узлов после перехода по подклееному ребру можно установить, что мы имеем дело именно со случаем $P2$, а последовательность входящих-выходящих ребер помогает установить нужную конфигурацию.

\medskip

\medskip

Теперь покажем, как по этому коду определить код парного к нему пути.

{\bf Вычисление кода $XFZ$.} Аналогично случаю {\bf P1}, зная код пути $XYZ$ (то есть коды $X$, $Y$, $Z$, а также ребра входа в $Y$, $Z$ и ребра выхода из $X$ и $Y$ вдоль пути $XYZ$), можно выписать код $XFZ$ (таблица~\ref{P2a}).

\medskip
\begin{table}[hbtp]
\caption{Случай P2: вычисление кода $XFZ$. }
\centering
 \begin{tabular}{|c|c|c|c|c|c|c|}   \hline
  & ребро из $X$ & ребро в $F$ & \x{тип, уровень \cr и окружение $F$} & \x{информация \cr $F$} &  \x{ребро \cr из $F$} & \x{ребро \cr в $Z$}
 \cr \hline
1. & $\widehat{\mathbf{u}_1}$ & $1$ & \x{$\mathbb{B}$  (окружение \cr как у $Y$)} & 1 босс -- $X$ & $3$ & $2$  \cr \hline
2. & $\widehat{\mathbf{l}}$ & $1$ & \x{ $\mathbb{DR}$,$1$, окружение \cr  $3$ ребро }& как у $Z$  & $2$ & $3$  \cr \hline
4. & Предл.\ref{tilecore} & Предл.\ref{tilecore} & $\mathbf{SBoss}(Z)$ & как у $X$  & $\mathbf{l}_2$ & $4$  \cr \hline
5. & Предл.\ref{tilecore} & Предл.\ref{tilecore} & $\mathbf{Next.FBoss}(Z)$ & как у $X$   & $\mathbf{l}$ & $3$  \cr \hline
7. & Предл.\ref{tilecore} & $\mathbf{T.Core}(Y)$ & $\mathbf{Core}(Y)$ & $\mathbf{Core}(Y)$  & $\mathbf{L.Core}(Y)$ & Предл.\ref{tilecore}  \cr \hline
8. & Предл.\ref{tilecore} & Предл.\ref{tilecore} & $\mathbf{FBoss}(Z)$& как у $X$  & $\mathbf{u}_2$ & $2$  \cr \hline
9. & Предл.\ref{tilecore} & Предл.\ref{tilecore} & $\mathbf{FBoss}(Z)$& как у $X$  & $\widehat{\mathbf{u}_1}$ & $1$  \cr \hline
10. & $\mathbf{Top.Core}(Z)$  & Предл.\ref{tilecore} & $\mathbf{FBoss}(Z)$ &  Предл.\ref{tilecore2}  & $\widehat{\mathbf{u}_2}$ & $1$  \cr \hline
  \end{tabular}
\label{P2a}
\end{table}

\medskip

{\bf Вычисление кода $XYZ$.} Зная код пути $XFZ$ (то есть коды $X$, $F$, $Z$, а также ребра входа в $F$, $Z$ и ребра выхода из $X$ и $F$ вдоль пути $XFZ$), можно выписать код $XYZ$ (таблица~\ref{P2b}).

\medskip

\begin{table}[hbtp]
  \caption{Случай P2: вычисление кода $XYZ$. }
\centering
 \begin{tabular}{|c|c|c|c|c|c|c|}   \hline
  & \x{ребро \cr из $X$} & \x{ребро \cr в $Y$} & \x{тип, уровень и \cr окружение $Y$} & информация $Y$ &  \x{ребро \cr из $Y$} & \x{ребро \cr в $Z$} \cr \hline
1. & $\widehat{\mathbf{u}_2}$ & $1$ & $\mathbb{A}$ (окружение как у $F$)& 1 босс -- $X$ & $2$ & $1$  \cr \hline
2. & $\widehat{\mathbf{d}_2}$ & $4$ & $\mathbb{C}$ (окружение как у $F$)& \x{  $1$ босс -- как у $Z$\cr $2$ босс -$X$, \cr $3$ босс -- $Z$  }  & $3$ & $\mathbf{d}_2$  \cr \hline
4. & $\widehat{\mathbf{l}}$ & $3$ & $\mathbb{A}$ (окружение как у $Z$)& 1 босс как у $Z$  & $2$ & $1$  \cr \hline
5. & $\widehat{\mathbf{d}_2}$ & $3$ & $\mathbb{C}$ (окружение как у $Z$) & \x{  1 босс как у $Z$ \cr  $2$ босс -$X$, \cr  $3$ босс -- $\mathbb{A}$, \cr окружение -- \cr $\mathbf{Past}(X)$ } & $\mathbf{l}$ & $2$  \cr \hline
7. & $\widehat{\mathbf{u}_2}$ & $1$ & \x{ $\mathbb{A}$, окружение -- \cr $(\mathbf{left},\mathbf{right},\mathbf{\mathbf{1A}},\mathbf{3A})$ }& $1$ босс -- $X$ & $3$ & $\mathbf{l}$  \cr \hline
8. & $\widehat{\mathbf{u}_1}$ & $1$ & \x{$\mathbb{B}$ окружение -- \cr $(\mathbf{left},\mathbf{right},\mathbf{1A},\mathbf{3A})$} & \x{ $1$ босс -- $X$ \cr тип $2$ -- $\mathbb{A}$ }  & $2$ & $\mathbf{r}$  \cr \hline
9. &  $\widehat{\mathbf{u}_2}$ & $2$ & \x{ $\mathbb{RU}$, окружение -- \cr $0$-цепь вокруг  $\mathbb{A}$ (ук $1$) \cr c окружением $\mathbf{Past}(X)$ \cr уровень 1 }& $1$ босс -- $X$  & $\mathbf{r}$ & $2$  \cr \hline
10. & $\mathbf{L.Core}(Z)$ & Предл.\ref{tilecore} & $\mathbf{Next.Past.FBoss}(Z)$& Предл.\ref{tilecore2}  & $\widehat{\mathbf{l}}$ & $3$  \cr \hline
  \end{tabular}
\label{P2b}
\end{table}

{\bf Примечание.}
Запись $\mathbf{Next.Past.FBoss}(Z)$ означает, что берется следующий узел в цепи для первого начальника $Z$, причем рассматривается его подклееное окружение (а не основное).

\medskip

{\bf Соотношения.} Вводятся аналогично случаю $P1$, коды вершин и ребер указаны в таблицах выше.

\subsection{P3: Макроплитка с окружением $(\mathbf{7A},x,\mathbf{3B},\mathbf{4A})$; $x=\mathbf{left},\mathbf{top}$}

В этом параграфе мы рассмотрим случай {\bf P3}, когда макроплитка $T$, внутри которой мы рассматриваем локальное преобразование пути, имеет окружение $(\mathbf{7A},x,\mathbf{1A},\mathbf{3A})$; $x=\mathbf{left},\mathbf{top}$. Это происходит, когда макроплитка является левой нижней подплиткой и верхней стороной выходит на левую или верхнюю сторону некоторой подклееной макроплитки.

\begin{figure}[hbtp]
\centering
\includegraphics[width=1\textwidth]{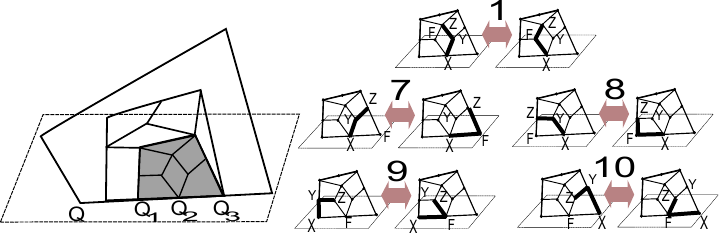}
\caption{Случай P3}
\label{P3}
\end{figure}

Зафиксируем пять вершин $Q$, $Q_1$, $Q_2$, $Q_3$, занимающих на комплексе положения как на рисунке~\ref{P3}. То есть $Q$ является ядром подклейки, а остальные вершины лежат на сторонах.

Для зафиксированного сочетания кодов этих пяти вершин, можно вычислить коды всех путей, указанных на рисунке~\ref{P3}. Кроме того, по выданному слову (коду пути) можно установить, действительно ли этот путь относится к случаю $P2$ и какую конфигурацию из изображенных на рисунке~\ref{P3} он имеет.

\medskip

\medskip

{\bf Вычисление кода $XFZ$.} Зная код пути $XYZ$ (то есть коды $X$, $Y$, $Z$, а также ребра входа в $Y$, $Z$ и ребра выхода из $X$ и $Y$ вдоль пути $XYZ$), можно выписать код $XFZ$. Обозначение $\mathbf{Past}(X)$ используем для обозначения подклееного окружения узла $X$ (таблица~\ref{P3a}).

\medskip
\begin{table}[hbtp]
  \caption{Случай P3: вычисление кода $XFZ$.}
\centering
 \begin{tabular}{|c|c|c|c|c|c|c|}   \hline
  & ребро из $X$ & ребро в $F$ & \x{тип, уровень и \cr окружение $F$} & информация $F$ &  ребро из $F$ & \x{ребро \cr в $Z$} \cr \hline
1. & $\widehat{\mathbf{u}_1}$ & $1$ & \x{$\mathbb{B}$ \cr (окруж как у $Y$) }& \x{ 1 босс -- $X$ \cr тип 2 -- $\mathbb{A}$ } & $3$ & $2$  \cr \hline
7. &Предл.\ref{tilecore} & Предл.\ref{tilecore} & $\mathbf{SBoss}(Z)$ & как у $X$  & \x{ $\widehat{\mathbf{u}_3}$ если \cr $\mathbf{Past}(F)=\mathbb{U}$ \cr  $\widehat{\mathbf{l}_2}$ если \cr $\mathbf{Past}(F)=\mathbb{L}$ } & $\mathbf{l}$   \cr \hline
8. & Предл.\ref{tilecore} & Предл.\ref{tilecore} & $\mathbf{Next.FBoss}(Z)$ & как у $X$  & $\mathbf{l}$ & $1$  \cr \hline
9. & Предл.\ref{tilecore} & Предл.\ref{tilecore} & $\mathbf{FBoss}(Z)$ & как у $X$  & $\widehat{\mathbf{u}_1}$ & $1$  \cr \hline
10. & Предл.\ref{tilecore} & Предл.\ref{tilecore} & $\mathbf{FBoss}(Z)$& как у $X$  & $\widehat{\mathbf{u}_2}$ & $1$  \cr \hline
  \end{tabular}
\label{P3a}
\end{table}

\medskip

{\bf Вычисление кода $XYZ$.} Зная код пути $XFZ$ (то есть коды $X$, $F$, $Z$, а также ребра входа в $F$, $Z$ и ребра выхода из $X$ и $F$ вдоль пути $XFZ$), можно выписать код $XYZ$ (таблица~\ref{P3b}).

\medskip
\begin{table}[hbtp]
 \caption{Случай P3: вычисление кода $XYZ$. }
\centering
 \begin{tabular}{|c|c|c|c|c|c|c|}   \hline
& \x{ребро \cr из $X$} & \x{ребро \cr в $Y$} & \x{тип, уровень \cr и окружение $Y$} & информация $Y$ & \x{ ребро \cr из $Y$} & \x{ребро \cr в $Z$ }\cr \hline
1. & $\widehat{\mathbf{u}_2}$ & $1$ & $\mathbb{A}$ (окруж. как у $F$)& 1 босс -- $X$ & $2$ & $1$  \cr \hline
7. & $\widehat{\mathbf{u}_2}$ & $1$ & \x{ $\mathbb{A}$, окружение -\cr  $(\mathbf{7A},x,\mathbf{3B},\mathbf{4A})$, \cr $x=\mathbf{left},\mathbf{top}$  }& $1$ босс -- $X$ & $3$ & $\mathbf{l}$  \cr  \hline
8. & $\widehat{\mathbf{u}_1}$ & $1$ &  \x{  $\mathbb{B}$, окружение -- \cr $(\mathbf{7A},x,\mathbf{3B},\mathbf{4A})$, \cr  $x=\mathbf{left},\mathbf{top}$  }& $1$ босс -- $X$  & $2$ & $\mathbf{r}$  \cr \hline
9. &  $\widehat{\mathbf{l}}$ & $2$ & \x{  $\mathbb{DR}$,1 окружение -- $3$ } & $1$ босс -- $Y$  & $\mathbf{r}$ & $2$  \cr \hline
10. & $\widehat{\mathbf{l}_2}$ или $\widehat{\mathbf{u}_1}$  & $1$ & $\mathbf{Next.Past.FBoss}(Z)$ & \x{ 1 босс -- \cr $\mathbf{TopFromCorner.Past}(X)$ \cr  $2$ босс -- $X$, $3$ босс -- \cr $\mathbf{RightCorner.Past}(X)$} & $\mathbf{l}$ &  $3$   \cr  \hline
  \end{tabular}
\label{P3b}
\end{table}

\medskip

{\bf Примечание}. Ребро из $X$ в случае $10$ имеет тип $\widehat{\mathbf{l}_2}$ если $\mathbf{Past}(X)=\mathbb{L}$ и $\widehat{\mathbf{u}_1}$ если $\mathbf{Past}(X)=\mathbb{U}$.

Запись $\mathbf{Next.Past.FBoss}(Z)$ означает, что берется следующий узел в цепи для первого начальника $Z$, причем рассматривается его подклееное окружение (а не основное).

\medskip

\subsection{P4: Макроплитка с окружением $(\mathbf{7A},x,\mathbf{1A},\mathbf{3A})$; $x=\mathbf{left},\mathbf{top}$}

В этом параграфе мы рассмотрим случай {\bf P4}, когда макроплитка $T$, внутри которой мы рассматриваем локальное преобразование пути, имеет окружение $(\mathbf{7A},x,\mathbf{1A},\mathbf{3A})$; $x=\mathbf{left},\mathbf{top}$. Это происходит, когда макроплитка является прямым потомком некоторой макроплитки, являющейся левой нижней подплиткой.  Верхней стороной $T$ выходит на левую или верхнюю сторону некоторой подклееной макроплитки.

Определение кодов вершин можно провести аналогично предыдущим случаям.

\medskip

\begin{figure}[hbtp]
\centering
\includegraphics[width=1\textwidth]{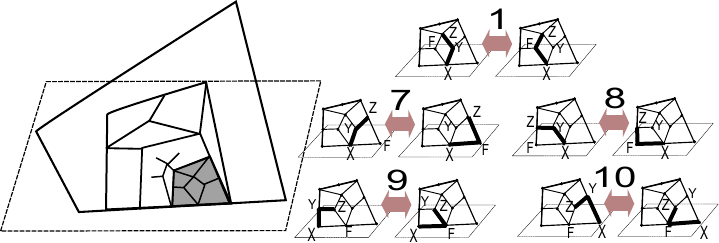}
\caption{Случай P4}
\label{P4}
\end{figure}

\medskip

{\bf Вычисление кода $XFZ$.} Зная код пути $XYZ$ (то есть коды $X$, $Y$, $Z$, а также ребра входа в $Y$, $Z$ и ребра выхода из $X$ и $Y$ вдоль пути $XYZ$), можно выписать код $XFZ$ (таблица~\ref{P4a}).

\medskip

\begin{table}[hbtp]
  \caption{Случай P4: вычисление кода $XFZ$. }
\centering
 \begin{tabular}{|c|c|c|c|c|c|c|}   \hline
  & ребро из $X$ & ребро в $F$ & \x{тип, уровень и \cr окружение $F$} & \x{информация \cr $F$} &  ребро из $F$ & \x{ребро \cr в $Z$}
 \cr \hline
1. & $\widehat{\mathbf{u}_1}$ & $1$ & \x{  $\mathbb{B}$ \cr  (окруж как у $Y$) }& \x{ 1 босс -- $X$ \cr тип 2 -- $\mathbb{A}$ }& $3$ & $2$  \cr \hline
7. & Предл.\ref{tilecore} & Предл.\ref{tilecore} & $\mathbf{SBoss}(Z)$ & как у $X$  & \x{ $\widehat{\mathbf{u}_3}$ если \cr $\mathbf{Past}(F)=\mathbb{U}$ \cr $\widehat{\mathbf{l}_2}$ если \cr $\mathbf{Past}(F)=\mathbb{L}$ } & $\mathbf{l}$   \cr \hline
8. & Предл.\ref{tilecore} & Предл.\ref{tilecore} & $\mathbf{FBoss}(Z)$ & как у $X$  & $\widehat{\mathbf{u}_2}$ & $2$  \cr \hline
9. & Предл.\ref{tilecore} & Предл.\ref{tilecore} & $\mathbf{FBoss}(Z)$ & как у $X$  & $\widehat{\mathbf{u}_1}$ & $1$  \cr \hline
10. & Предл.\ref{tilecore} & Предл.\ref{tilecore} & $\mathbf{FBoss}(Z)$& как у $X$  & $\widehat{\mathbf{u}_2}$ & $1$  \cr \hline
  \end{tabular}
\label{P4a}
\end{table}

\medskip

{\bf Вычисление кода $XYZ$.} Зная код пути $XFZ$ (то есть коды $X$, $F$, $Z$, а также ребра входа в $F$, $Z$ и ребра выхода из $X$ и $F$ вдоль пути $XFZ$), можно выписать код $XYZ$ (таблица~\ref{P4b}).

\medskip

\begin{table}[hbtp]
  \caption{Случай P4: вычисление кода $XYZ$. }
\centering
 \begin{tabular}{|c|c|c|c|c|c|c|}   \hline
& \x{ребро \cr из $X$} & \x{ребро \cr в $Y$} & \x{тип, уровень и \cr окружение $Y$} & информация $Y$ &  \x{ ребро \cr из $Y$ }& \x{ребро \cr в $Z$} \cr \hline
1. & $\widehat{\mathbf{u}_2}$ & $1$ & \x{$\mathbb{A}$, \cr (окружение как у $F$)} & 1 босс -- $X$ & $2$ & $1$  \cr \hline
7. & $\widehat{\mathbf{u}_2}$ & $1$ & \x{  $\mathbb{A}$, окружение -- \cr  $(\mathbf{7A},x,\mathbf{1A},\mathbf{3A})$, \cr $x=\mathbf{left},\mathbf{top}$ } & $1$ босс -- $X$ & $3$ & $\mathbf{l}$  \cr \hline
8. & $\widehat{\mathbf{u}_1}$ & $1$ & \x{ $\mathbb{B}$, окружение -- \cr  $(\mathbf{7A},x,\mathbf{1A},\mathbf{3A})$, \cr $x=\mathbf{left},\mathbf{top}$ }& $1$ босс -- $X$  & $2$ & $\mathbf{r}$  \cr  \hline
9. &  $\widehat{\mathbf{u}_2}$ & $2$ & \x{  $\mathbb{RU}$,1, окружение -- \cr $\mathbb{A}0$-цепь  (указ $1$) \cr c окруж $\mathbf{Past}(X)$ } & $1$ босс -- $X$  & $\mathbf{r}$ & $2$  \cr \hline
10. & $\widehat{\mathbf{l}_2}$ или $\widehat{\mathbf{u}_1}$  & $1$ & $\mathbf{Next.Past.FBoss}(Z)$ & \x{  1 босс -- \cr $\mathbf{TopFromCorner.Past}(X)$ \cr 2 босс -- $X$,  $3$ босс -- \cr $\mathbf{RightCorner.Past}(X)$ }& $\mathbf{l}$ & $3$  \cr  \hline
  \end{tabular}
\label{P4b}
\end{table}

\medskip

{\bf Примечание.} Ребро из $X$ в случае $10$ имеет тип $\widehat{\mathbf{l}_2}$ если $\mathbf{Past}(X)=\mathbb{L}$, и $\widehat{\mathbf{u}_1}$ если $\mathbf{Past}(X)=\mathbb{U}$.

Запись $\mathbf{Next.Past.FBoss}(Z)$ означает, что берется следующий узел в цепи для первого начальника $Z$, причем рассматривается его подклееное окружение (а не основное).

\medskip

\subsection{P5: Макроплитка с окружением $(\mathbf{top},\mathbf{right},\mathbf{6A},\mathbf{2A})$}

В этом параграфе мы рассмотрим случай {\bf P5}, когда макроплитка $T$, внутри которой мы рассматриваем локальное преобразование пути, имеет окружение $(\mathbf{top},\mathbf{right},\mathbf{6A},\mathbf{2A})$. Это происходит, когда макроплитка является правой верхней подплиткой некоторой подклееной макроплитки.

Определение кодов вершин можно провести аналогично предыдущим случаям.

\medskip

\begin{figure}[hbtp]
\centering
\includegraphics[width=1\textwidth]{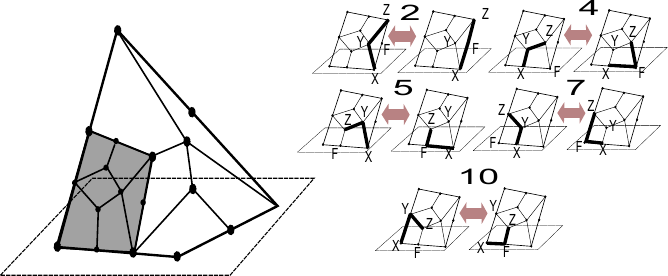}
\caption{Случай P5}
\label{P5}
\end{figure}

\medskip

{\bf Вычисление кода $XFZ$.} Зная код пути $XYZ$ (то есть коды $X$, $Y$, $Z$, а также ребра входа в $Y$, $Z$ и ребра выхода из $X$ и $Y$ вдоль пути $XYZ$), можно выписать код $XFZ$ (таблица~\ref{P5a}).

Для преобразования $7$ пути удовлетворяют условиям мертвого паттерна. Для них мы соотношения не вводим.

\medskip
\begin{table}[hbtp]
\caption{Случай P5: вычисление кода $XFZ$. }
\centering
 \begin{tabular}{|c|c|c|c|c|c|c|}   \hline
  & ребро из $X$ & ребро в $F$ & \x{тип, уровень и \cr окружение $F$} & информация $F$ &  ребро из $F$ & ребро в $Z$
 \cr  \hline
2. & $\widehat{\mathbf{u}_1}$ & $1$ & \x{ $\mathbb{DR}$,1, \cr окружение -- $3$ } & \x{ 1 босс -- $X$ \cr тип 2 -- $\mathbb{CDR}$  } & $1$ & $2$  \cr \hline
4. & Предл.\ref{tilecore} &  Предл.\ref{tilecore} & $\mathbf{SBoss}(Z)$ & как у $X$  & $\widehat{\mathbf{u}_3}$  & $4$   \cr \hline
5. &  Предл.\ref{tilecore} &  Предл.\ref{tilecore} & $\mathbf{Next.FBoss}(Z)$ & как у $X$  & $\widehat{\mathbf{l}}$ & $3$  \cr \hline
10. &  Предл.\ref{tilecore} &  Предл.\ref{tilecore} & $\mathbf{Next.FBoss}(Z)$ & как у $X$  & $\widehat{\mathbf{l}}$ & $3$  \cr \hline
  \end{tabular}
\label{P5a}
\end{table}

\medskip

{\bf Вычисление кода $XYZ$.} Зная код пути $XFZ$ (то есть коды $X$, $F$, $Z$, а также ребра входа в $F$, $Z$ и ребра выхода из $X$ и $F$ вдоль пути $XFZ$), можно выписать код $XYZ$ (таблица~\ref{P5b}).

\medskip

\begin{table}[hbtp]
\caption{Случай P5: вычисление кода $XYZ$. }
\centering
 \begin{tabular}{|c|c|c|c|c|c|c|}   \hline
& \x{ребро \cr из $X$} & \x{ребро \cr в $Y$} & \x{тип, уровень и \cr окружение $Y$} & информация $Y$ &  \x{ребро \cr из $Y$} & \x{ребро \cr в $Z$} \cr \hline
2. & $\widehat{\mathbf{u}_3}$ & $4$ & \x{ $\mathbb{C}$, окружение как у $T$} & \x{1 босс -- $\mathbb{U}$, \cr окружение как у $T$, \cr $(\mathbf{top},\mathbf{right},\mathbf{6A},\mathbf{2A})$, \cr  2 босс -- $X$, \cr 3 босс -- $\mathbb{B}$, окружение \cr $(\mathbf{left},\mathbf{top},\mathbf{right},\mathbf{bottom})$}  & $2$ & $1$  \cr \hline
4. & $\widehat{\mathbf{l}}$ & $3$ & $\mathbb{A}$, окружение как у $T$ & $1$ босс как у $Z$ & $2$ & $1$  \cr \hline
5. & $\widehat{\mathbf{u}_3}$ & $4$ & $\mathbb{C}$, окружение как у $T$ & \x{ $1$ босс -- как у $Z$ \cr 2 босс -- $X$,\cr 3 босс --  $\mathbb{B}$, окружение \cr $(\mathbf{left},\mathbf{top},\mathbf{right},\mathbf{bottom})$ }  & $1$ & $2$  \cr  \hline
10. & $\widehat{\mathbf{r}}$ & $2$ & $\mathbf{FBoss}(Z)$ & --   & $\mathbf{u}_2$ & $1$  \cr \hline
  \end{tabular}
\label{P5b}
\end{table}

\medskip

\subsection{P6: Макроплитка с окружением $(\mathbf{top},\mathbf{right},\mathbf{1A},\mathbf{3A})$}

В этом параграфе мы рассмотрим случай {\bf P6}, когда макроплитка $T$, внутри которой мы рассматриваем локальное преобразование пути, имеет окружение $(\mathbf{top},\mathbf{right},\mathbf{1A},\mathbf{3A})$. Это происходит, когда макроплитка является прямым потомком некоторой макроплитки, являющейся правой верхней подплиткой некоторой подклееной макроплитки.

Определение кодов вершин можно провести аналогично предыдущим случаям.

\medskip

\begin{figure}[hbtp]
\centering
\includegraphics[width=1\textwidth]{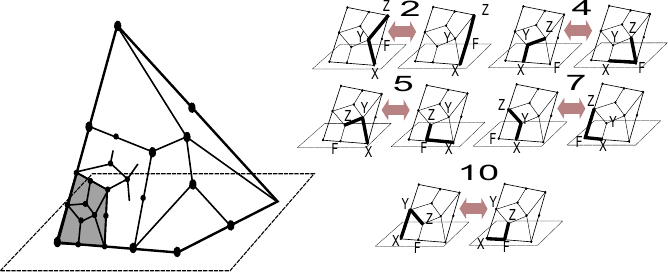}
\caption{Случай P6}
\label{P6}
\end{figure}

\medskip

{\bf Вычисление кода $XFZ$.} Зная код пути $XYZ$ (то есть коды $X$, $Y$, $Z$, а также ребра входа в $Y$, $Z$ и ребра выхода из $X$ и $Y$ вдоль пути $XYZ$), можно выписать код $XFZ$ (таблица~\ref{P6a}).

Для преобразования $7$ пути удовлетворяют условиям мертвого паттерна. Для них мы соотношения не вводим.

\medskip

\begin{table}[hbtp]
  \caption{Случай P6: вычисление кода $XFZ$. }
\centering
\begin{tabular}{|c|c|c|c|c|c|c|}   \hline
  & ребро из $X$ & ребро в $F$ & \x{тип, уровень и \cr окружение $F$} & информация $F$ &  ребро из $F$ & ребро в $Z$ \cr \hline
2. & $\widehat{\mathbf{l}}$ & $1$ & \x{ $\mathbb{DR}$,1, \cr окружение -- $3$ } & 1 босс как у $Z$ & $2$ & $3$  \cr \hline
4. &  Предл.\ref{tilecore} &  Предл.\ref{tilecore} & $\mathbf{SBoss}(Z)$ & как у $X$  & $\widehat{\mathbf{l}_3}$  & $4$   \cr \hline
5. &  Предл.\ref{tilecore} &  Предл.\ref{tilecore} & $\mathbf{Next.FBoss}(Z)$ & как у $X$  & $\widehat{\mathbf{l}}$ & $3$  \cr \hline
10. &  Предл.\ref{tilecore} &  Предл.\ref{tilecore} & $\mathbf{Next.FBoss}(Z)$ & как у $X$  & $\widehat{\mathbf{l}}$ & $3$  \cr \hline
  \end{tabular}
\label{P6a}
\end{table}

\medskip

{\bf Вычисление кода $XYZ$.} Зная код пути $XFZ$ (то есть коды $X$, $F$, $Z$, а также ребра входа в $F$, $Z$ и ребра выхода из $X$ и $F$ вдоль пути $XFZ$), можно выписать код $XYZ$ (таблица~\ref{P6b}).

\medskip

\begin{table}[hbtp]
\caption{Случай P6: вычисление кода $XYZ$. }
\centering
 \begin{tabular}{|c|c|c|c|c|c|c|}   \hline
& \x{ребро \cr из $X$} & \x{ребро \cr в $Y$} & \x{тип, уровень и \cr окружение $Y$} & информация $Y$ &  \x{ребро \cr из $Y$} & \x{ребро \cr в $Z$} \cr \hline
2. & $\widehat{\mathbf{l}_2}$ & $4$ & \x{$\mathbb{C}$, \cr окружение -- как у $T$} & \x{ 1 босс -- $\mathbb{U}$, \cr окружение как у $T$\cr 2 босс -- $X$, \cr 3 босс -- $\mathbb{A}$, окружение \cr как у $\mathbf{FBoss}(F)$  } & $2$ & $1$  \cr \hline
4. & $\widehat{\mathbf{l}}$ & $3$ & $\mathbb{A}$, окружение как у $T$ & $1$ босс как у $Z$ & $2$ & $1$  \cr \hline
5. & $\widehat{\mathbf{u}_3}$ & $4$ & $\mathbb{C}$, окружение как у $T$ & \x{ $1$ босс -- как у $Z$\cr 2 босс -- $X$, \cr 3 босс --  $\mathbb{A}$,окружение \cr $\mathbf{RightCorner.Past}(X)$  }  & $1$ & $2$  \cr  \hline
10. & $\widehat{\mathbf{r}}$  & $2$ & $\mathbf{FBoss}(Z)$ & 1 босс -- $X$   & $\mathbf{u}_2$ & $1$  \cr \hline
  \end{tabular}
\label{P6b}
\end{table}

\medskip

\subsection{P7: Макроплитка с окружением $(x,\mathbf{3B},\mathbf{6A},\mathbf{2A})$; $x=\mathbf{left},\mathbf{top}$}

В этом параграфе мы рассмотрим случай {\bf P7}, когда макроплитка $T$, внутри которой мы рассматриваем локальное преобразование пути, имеет окружение $(x,\mathbf{3B},\mathbf{6A},\mathbf{2A})$; $x=\mathbf{left},\mathbf{top}$. Это происходит, когда макроплитка является правой верхней подплиткой некоторой макроплитки, являющейся левой нижней подплиткой подклееной макроплитки.

Определение кодов вершин можно провести аналогично предыдущим случаям.

\medskip

\begin{figure}[hbtp]
\centering
\includegraphics[width=1\textwidth]{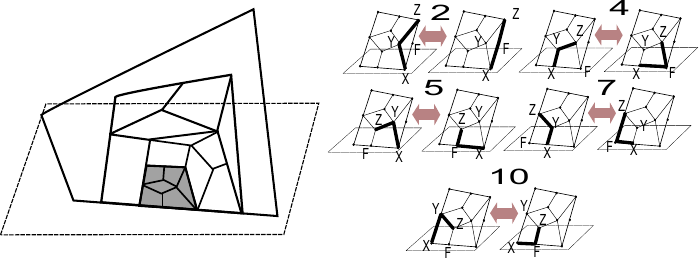}
\caption{Случай P7}
\label{P7}
\end{figure}

\medskip

{\bf Вычисление кода $XFZ$.} Зная код пути $XYZ$ (то есть коды $X$, $Y$, $Z$, а также ребра входа в $Y$, $Z$ и ребра выхода из $X$ и $Y$ вдоль пути $XYZ$), можно выписать код $XFZ$ (таблица~\ref{P7a}).

\medskip

\begin{table}[hbtp]
  \caption{Случай P7: вычисление кода $XFZ$. }
\centering
 \begin{tabular}{|c|c|c|c|c|c|c|}   \hline
  & ребро из $X$ & ребро в $F$ & \x{тип, уровень и \cr окружение $F$} & информация $F$ &  ребро из $F$ & ребро в $Z$ \cr \hline
2. & $\widehat{\mathbf{u}_1}$ & $1$ &  \x{ $\mathbb{DR}$,1, \cr окружение -- $2$  } & \x{ 1 босс -- $X$ \cr  тип 2 -- $\mathbb{A}$} & $2$ & $1$  \cr \hline
4. & Предл.\ref{tilecore} & Предл.\ref{tilecore}& $\mathbf{SBoss}(Z)$ & как у $X$  & $\widehat{\mathbf{u}_3}$  & $4$   \cr \hline
5. & Предл.\ref{tilecore} & Предл.\ref{tilecore} & $\mathbf{Next.FBoss}(Z)$ & как у $X$  & $\widehat{\mathbf{l}}$ & $3$  \cr \hline
7. & Предл.\ref{tilecore} & Предл.\ref{tilecore} & $\mathbf{Next.FBoss}(Z)$ & как у $X$  & $\widehat{\mathbf{l}}$ & $1$  \cr \hline
10. & Предл.\ref{tilecore} & Предл.\ref{tilecore} & $\mathbf{Next.FBoss}(Z)$ & как у $X$  & $\widehat{\mathbf{l}}$ & $3$  \cr \hline
  \end{tabular}
\label{P7a}
\end{table}

\medskip

{\bf Вычисление кода $XYZ$.} Зная код пути $XFZ$ (то есть коды $X$, $F$, $Z$, а также ребра входа в $F$, $Z$ и ребра выхода из $X$ и $F$ вдоль пути $XFZ$), можно выписать код $XYZ$ (таблица~\ref{P7b}).

\medskip

\begin{table}[hbtp]
\caption{Случай P7: вычисление кода $XYZ$. }
\centering
\begin{tabular}{|c|c|c|c|c|c|c|}   \hline
& ребро  & ребро & \x{тип, уровень и \cr окружение $Y$} & информация $Y$ &  ребро & ребро  \cr
&из $X$ &в $Y$&&& из $Y$ & в $Z$ \cr \hline
2. & $\widehat{\mathbf{u}_3}$ & $4$ & $\mathbb{C}$, окружение -- как у $T$ & \x{ 1 босс -- \cr $\mathbf{TopFromCorner.Past}(X)$,\cr 2 босс -- $X$, \cr 3 босс -- $Z$  } & $3$ & $\mathbf{ru}$  \cr \hline
4. & $\widehat{\mathbf{l}}$ & $3$ & $\mathbb{A}$ окружение как у $T$ & $1$ босс как у $Z$ & $2$ & $1$  \cr \hline
5. & $\widehat{\mathbf{u}_3}$ & $4$ & $\mathbb{C}$ окружение как у $T$ &  \x{$1$ босс -- как у $Z$\cr 2 босс -- $X$, \cr 3 босс --  $\mathbb{B}$,окружение \cr  $(\mathbf{7A},x,\mathbf{3B},\mathbf{4A})$  } & $2$ & $1$  \cr \hline
7. &  $\widehat{\mathbf{l}}$ & $3$ & $\mathbb{A}$ окружение как у $T$  & $1$ босс -- $Z$  & $1$ & $\mathbf{u}_2$  \cr\hline
10. & $\widehat{\mathbf{l}}$  & $1$ & $\mathbf{FBoss}(Z)$ & 1 босс -- $\mathbf{Prev}(X)$  & $\mathbf{u}_2$ & $1$  \cr \hline
  \end{tabular}
\label{P7b}
\end{table}

\medskip

\subsection{P8: Макроплитка с окружением $(x,\mathbf{3B},\mathbf{1A},\mathbf{3A})$; $x=\mathbf{left},\mathbf{top}$}

В этом параграфе мы рассмотрим случай {\bf P8}, когда макроплитка $T$, внутри которой мы рассматриваем локальное преобразование пути, имеет окружение $(x,\mathbf{3B},\mathbf{1A},\mathbf{3A})$; $x=\mathbf{left},\mathbf{top}$. Это происходит, когда макроплитка является прямым потомком правой верхней подплитки некоторой макроплитки, являющейся левой нижней подплиткой подклееной макроплитки.

Определение кодов вершин можно провести аналогично предыдущим случаям.

\medskip

\begin{figure}[hbtp]
\centering
\includegraphics[width=1\textwidth]{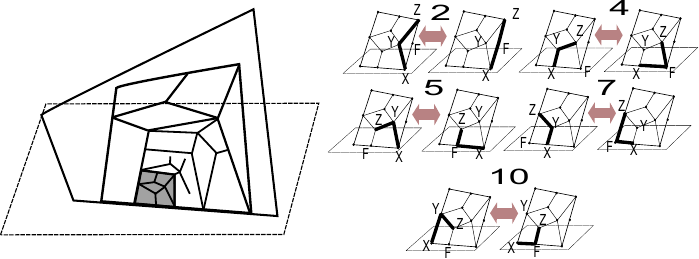}
\caption{Случай P8}
\label{P8}
\end{figure}

\medskip

{\bf Вычисление кода $XFZ$.} Зная код пути $XYZ$ (то есть коды $X$, $Y$, $Z$, а также ребра входа в $Y$, $Z$ и ребра выхода из $X$ и $Y$ вдоль пути $XYZ$), можно выписать код $XFZ$ (таблица~\ref{P8a}).

\medskip

\begin{table}[hbtp]
\caption{Случай P8: вычисление кода $XFZ$. }
\centering
 \begin{tabular}{|c|c|c|c|c|c|c|}   \hline
  & ребро из $X$ & ребро в $F$ & \x{тип, уровень и \cr окружение $F$} & информация $F$ &  ребро из $F$ & ребро в $Z$ \cr \hline
2. & $\widehat{\mathbf{l}}$ & $1$ & \x{  $\mathbb{DR}$,1, \cr окружение -- $3$  }& 1 босс -- как у $Z$ & $2$ & $3$  \cr \hline
4. & Предл.\ref{tilecore} & Предл.\ref{tilecore} & $\mathbf{SBoss}(Z)$ & как у $X$  & $\widehat{\mathbf{l}_2}$  & $4$   \cr \hline
5. & Предл.\ref{tilecore} & Предл.\ref{tilecore} & $\mathbf{Next.FBoss}(Z)$ & как у $X$  & $\widehat{\mathbf{l}}$ & $3$  \cr \hline
7. & Предл.\ref{tilecore} & Предл.\ref{tilecore} & $\mathbf{Next.FBoss}(Z)$ & как у $X$  & $\widehat{\mathbf{l}}$ & $1$  \cr \hline
10. & Предл.\ref{tilecore} & Предл.\ref{tilecore} & $\mathbf{Next.FBoss}(Z)$ & как у $X$  & $\widehat{\mathbf{l}}$ & $3$  \cr \hline
  \end{tabular}
\label{P8a}
\end{table}

\medskip

{\bf Вычисление кода $XYZ$.} Зная код пути $XFZ$ (то есть коды $X$, $F$, $Z$, а также ребра входа в $F$, $Z$ и ребра выхода из $X$ и $F$ вдоль пути $XFZ$), можно выписать код $XYZ$ (таблица~\ref{P8b}).

\medskip

\begin{table}[hbtp]
 \caption{Случай P8: вычисление кода $XYZ$. }
\centering
\begin{tabular}{|c|c|c|c|c|c|c|}   \hline
& \x{ребро \cr из $X$} & \x{ребро \cr в $Y$} & \x{тип, уровень и \cr окружение $Y$} & информация $Y$ &  \x{ребро \cr из $Y$} & \x{ребро \cr в $Z$} \cr \hline
2. & $\widehat{\mathbf{l}_2}$ & $4$ & $\mathbb{C}$, окружение -- как у $T$ & \x{1 босс -- \cr  $\mathbf{TopFromCorner.Past}(X)$, \cr 2 босс -- $X$ \cr 3 босс -- $Z$   } & $3$ & $\mathbf{ru}$  \cr \hline
4. & $\widehat{\mathbf{l}}$ & $3$ & $\mathbb{A}$, окружение как у $T$ & $1$ босс как у $Z$ & $2$ & $1$  \cr \hline
5. & $\widehat{\mathbf{l}_2}$ & $4$ & $\mathbb{C}$, окружение как у $T$ & \x{ $1$ босс -- как у $Z$\cr 2 босс -- $X$, \cr 3 босс --  $\mathbb{A}$, окружение \cr $\mathbf{RightCorner.Past}(X)$ } & $1$ & $2$  \cr \hline
7. &  $\widehat{\mathbf{l}}$ & $3$ & $\mathbb{A}$, окружение как у $T$  & $1$ босс -- $Z$  & $1$ & $\mathbf{u}_2$  \cr\hline
10. & $\widehat{\mathbf{l}}$  & $2$ & $\mathbf{FBoss}(Z)$ & 1 босс -- $\mathbf{Prev}(X)$  & $\mathbf{u}_2$ & $1$  \cr \hline
  \end{tabular}
\label{P8b}
\end{table}

\medskip

\subsection{P9: Макроплитка с окружением $(x,\mathbf{1A},\mathbf{6A},\mathbf{2A})$; $x=\mathbf{left},\mathbf{top}$}

В этом параграфе мы рассмотрим случай {\bf P9}, когда макроплитка $T$, внутри которой мы рассматриваем локальное преобразование пути, имеет окружение $(x,\mathbf{3B},\mathbf{6A},\mathbf{2A})$; $x=\mathbf{left},\mathbf{top}$. Это происходит, когда макроплитка является правой верхней подплиткой некоторой макроплитки, являющейся левой верхней подплиткой подклееной макроплитки.

Определение кодов вершин можно провести аналогично предыдущим случаям.

\medskip

\begin{figure}[hbtp]
\centering
\includegraphics[width=1\textwidth]{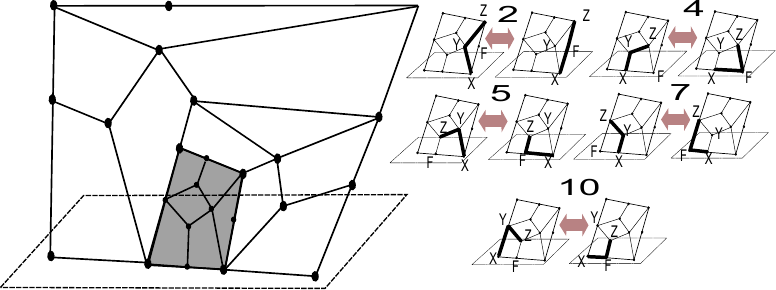}
\caption{Случай P9}
\label{P9}
\end{figure}

\medskip

{\bf Вычисление кода $XFZ$.} Зная код пути $XYZ$ (то есть коды $X$, $Y$, $Z$, а также ребра входа в $Y$, $Z$ и ребра выхода из $X$ и $Y$ вдоль пути $XYZ$), можно выписать код $XFZ$ (таблица~\ref{P9a}).

\medskip
\begin{table}[hbtp]
  \caption{Случай P9: вычисление кода $XFZ$. }
\centering
 \begin{tabular}{|c|c|c|c|c|c|c|}   \hline
  & ребро из $X$ & ребро в $F$ & \x{тип, уровень и \cr окружение $F$} & информация $F$ &  ребро из $F$ & ребро в $Z$ \cr \hline
2. & $\widehat{\mathbf{u}_1}$ & $1$ & \x{  $\mathbb{DR}$,1, \cr окружение -- $2$ }& \x{ 1 босс -- $X$ \cr  тип 2 -- $\mathbb{A}$} & $2$ & $1$  \cr \hline
4. & Предл.\ref{tilecore} & Предл.\ref{tilecore} & $\mathbf{SBoss}(Z)$ & как у $X$  & $\widehat{\mathbf{u}_3}$  & $4$   \cr \hline
5. & Предл.\ref{tilecore} & Предл.\ref{tilecore} & $\mathbf{Next.FBoss}(Z)$ & как у $X$  & $\widehat{\mathbf{l}}$ & $3$  \cr \hline
7. & Предл.\ref{tilecore} & Предл.\ref{tilecore} & $\mathbf{Next.FBoss}(Z)$ & как у $X$  & $\widehat{\mathbf{u}_2}$ & $1$  \cr \hline
10. & Предл.\ref{tilecore} & Предл.\ref{tilecore} & $\mathbf{Next.FBoss}(Z)$ & как у $X$  & $\widehat{\mathbf{l}}$ & $3$  \cr \hline
  \end{tabular}
\label{P9a}
\end{table}

\medskip

{\bf Вычисление кода $XYZ$.} Зная код пути $XFZ$ (то есть коды $X$, $F$, $Z$, а также ребра входа в $F$, $Z$ и ребра выхода из $X$ и $F$ вдоль пути $XFZ$), можно выписать код $XYZ$ (таблица~\ref{P9b}).

\medskip

\begin{table}[hbtp]
\caption{Случай P9: вычисление кода $XYZ$.}
\centering
\begin{tabular}{|c|c|c|c|c|c|c|}   \hline
& ребро  & ребро  & \x{тип, уровень и \cr окружение $Y$} & информация $Y$ &  ребро  & ребро  \cr
&из $X$ &в $Y$&&& из $Y$ & в $Z$ \cr \hline
2. & $\widehat{\mathbf{u}_3}$ & $4$ & $\mathbb{C}$, окружение -- как у $T$ & \x{ 1 босс -- \cr  $\mathbf{TopFromCorner.Past}$(X), \cr 2 босс -- $X$, \cr 3 босс -- $Z$ } & $3$ & $\mathbf{ru}$  \cr \hline
4. & $\widehat{\mathbf{l}}$ & $3$ & $\mathbb{A}$, окружение как у $T$ & $1$ босс как у $Z$ & $2$ & $1$  \cr \hline
5. & $\widehat{\mathbf{u}_3}$ & $4$ & $\mathbb{C}$, окружение как у $T$ & \x{  $1$ босс -- как у $Z$\cr 2 босс -- $X$, \cr 3 босс --  $\mathbb{B}$, окружение \cr $(\mathbf{left},\mathbf{top},\mathbf{1A},\mathbf{3A})$ } & $1$ & $2$  \cr \hline
7. &  $\widehat{\mathbf{l}}$ & $3$ & $\mathbb{A}$, окружение как у $T$  & $1$ босс -- $Z$  & $1$ & $\mathbf{u}_2$  \cr\hline
10. & $\widehat{\mathbf{l}}$  & $1$ & $\mathbf{FBoss}(Z)$ & 1 босс -- $X$  & $\mathbf{u}_2$ & $1$  \cr \hline
  \end{tabular}
\label{P9b}
\end{table}

\medskip

\subsection{P10: Макроплитка с окружением $(x,\mathbf{1A},\mathbf{1A},\mathbf{3A})$; $x=\mathbf{left},\mathbf{top}$}

В этом параграфе мы рассмотрим случай {\bf P10}, когда макроплитка $T$, внутри которой мы рассматриваем локальное преобразование пути, имеет окружение $(x,\mathbf{3B},\mathbf{1A},\mathbf{3A})$; $x=\mathbf{left},\mathbf{top}$. Это происходит, когда макроплитка является прямым потомком правой верхней подплитки некоторой макроплитки, являющейся левой верхней подплиткой подклееной макроплитки.

Определение кодов вершин можно провести аналогично предыдущим случаям.

\medskip

\begin{figure}[hbtp]
\centering
\includegraphics[width=1\textwidth]{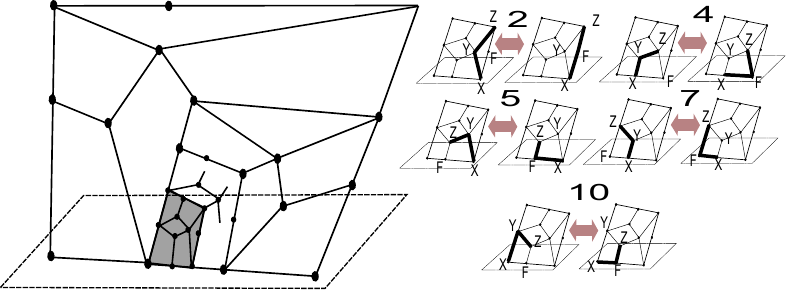}
\caption{Случай P10}
\label{P10}
\end{figure}

\medskip

{\bf Вычисление кода $XFZ$.} Зная код пути $XYZ$ (то есть коды $X$, $Y$, $Z$, а также ребра входа в $Y$, $Z$ и ребра выхода из $X$ и $Y$ вдоль пути $XYZ$), можно выписать код $XFZ$ (таблица~\ref{P10a}).

\medskip

\begin{table}[hbtp]
\caption{Случай P10: вычисление кода $XFZ$. }
\centering
 \begin{tabular}{|c|c|c|c|c|c|c|}   \hline
  & ребро из $X$ & ребро в $F$ & \x{тип, уровень и \cr окружение $F$} & информация $F$ &  ребро из $F$ & ребро в $Z$ \cr \hline
2. & $\widehat{\mathbf{l}}$ & $1$ & \x{  $\mathbb{DR}$ окружение -- $3$ }& 1 босс как у $Z$ & $2$ & $3$  \cr \hline
4. & Предл.\ref{tilecore} & Предл.\ref{tilecore} & $\mathbf{SBoss}(Z)$ & как у $X$  & $\widehat{\mathbf{l}_2}$  & $4$   \cr \hline
5. & Предл.\ref{tilecore} & Предл.\ref{tilecore} & $\mathbf{Next.FBoss}(Z)$ & как у $X$  & $\widehat{\mathbf{l}}$ & $3$  \cr \hline
7. & Предл.\ref{tilecore} & Предл.\ref{tilecore} & $\mathbf{FBoss}(Z)$ & как у $X$  & $\widehat{\mathbf{u}_2}$ & $2$  \cr \hline
10. & Предл.\ref{tilecore} & Предл.\ref{tilecore} & $\mathbf{Next.FBoss}(Z)$ & как у $X$  & $\widehat{\mathbf{l}}$ & $3$  \cr \hline
  \end{tabular}
\label{P10a}
\end{table}

\medskip

{\bf Вычисление кода $XYZ$.} Зная код пути $XFZ$ (то есть коды $X$, $F$, $Z$, а также ребра входа в $F$, $Z$ и ребра выхода из $X$ и $F$ вдоль пути $XFZ$), можно выписать код $XYZ$ (таблица~\ref{P10b}).

\medskip
\begin{table}[hbtp]
  \caption{Случай P10: вычисление кода $XYZ$. }
\centering
\begin{tabular}{|c|c|c|c|c|c|c|}   \hline
& \x{ребро \cr из $X$} & \x{ребро \cr в $Y$} & \x{тип, уровень и \cr окружение $Y$} & информация $Y$ &  \x{ребро \cr из $Y$ }& \x{ребро \cr в $Z$} \cr \hline
2. & $\widehat{\mathbf{l}_2}$ & $4$ & $\mathbb{C}$, окружение -- как у $T$ & \x{ 1 босс -- \cr $\mathbf{LevelPlus.FBoss}(Z)$, \cr 2 босс -- $X$, \cr 3 босс -- $Z$} & $3$ & $\mathbf{ru}$  \cr \hline
4. & $\widehat{\mathbf{l}}$ & $3$ & $\mathbb{A}$, окружение как у $T$ & $1$ босс как у $Z$ & $2$ & $1$  \cr \hline
5. & $\widehat{\mathbf{l}_2}$ & $4$ & $\mathbb{C}$, окружение как у $T$ & \x{ $1$ босс -- как у $Z$\cr 2 босс -- $X$, \cr 3 босс --  $\mathbb{A}$, окружение \cr $\mathbf{RightCorner.Past}(X)$ } & $1$ & $2$  \cr  \hline
7. &  $\widehat{\mathbf{l}}$ & $3$ & $\mathbb{A}$, окружение как у $T$  & $1$ босс -- $Z$  & $1$ & $\mathbf{u}_2$  \cr\hline
10. & $\widehat{\mathbf{u}_2}$  & $2$ & $\mathbf{FBoss}(Z)$ & 1 босс -- $X$  & $\mathbf{u}_2$ & $1$  \cr \hline
  \end{tabular}
\label{P10b}
\end{table}

\medskip

Таким образом, во всех случаях можно провести локальное преобразование, как и для плоских путей. Это дает нам возможности преобразовывать пути, меняя локальные участки на им эквивалентные и приводить слово к канонической форме.

\medskip

\subsection{Оценка числа введенных соотношений}

В каждом случае мы выбирали ядро подклейки, после чего выбирали не более четырех вершин ($Q_1$, $Q_2$, $Q_3$, $Q_4$) на сторонах подклееной макроплитки, для которых фиксировали их коды. После чего записывали не более $16$ соотношений, для каждого выбранного случая. Учитывая конструкцию подклейки, можно заметить, что информация вершин $Q_i$ может быть вычислена по информации ядра подклейки и известным ребрам-сторонам. То есть соотношений будет не более чем $16 \cdot F \cdot H^4$, где $F$ -- количество флагов подклейки, $H$ -- число расширенных окружений. То есть не более чем $16 \cdot 82485^4 \cdot 5 \cdot 10^{18}<4\cdot 10^{39}$ соотношений.

\medskip

\section{Приведение к канонической форме} \label{canonic}

В этой главе мы покажем, как проводить преобразования для приведения слова к канонической форме. Сначала мы рассмотрим пути, лежащие на одном ребре некоторой макроплитки, а затем разберем общий случай.

Пусть $X_i$ -- буквы, кодирующие входящие ребра, $Y_i$ -- буквы, кодирующие узлы (их типы, окружения и информации), $Z_i$ -- буквы, кодирующие выходящие ребра. Фактически, $X_i$ и $Z_i$ выбираются из одного алфавита ребер входа-выхода.
Напомним, что слово $W$ имеет {\it $\mathbf{CODE}$-форму}, если в нем справа от любой (не последней в слове) буквы семейства $X$ обязательно стоит буква семейства $Y$, справа от любой (не последней в слове)  буквы семейства $Y$ стоит буква семейства $Z$, а справа от любой (не последней в слове)  буквы семейства $Z$ стоит буква семейства $X$.
Рассмотрим полугруппу $S$ с нулем $0$ и порождающими $\{X_i,Y_i,Z_i \}$.
Будем считать, что в полугруппе $S$ введены определяющие соотношения следующих категорий:

\medskip

{\bf Категория 1}: соотношения, обеспечивающие $\mathbf{CODE}$-форму: $X_iZ_j=0$, $X_iX_j=0$, $Y_iX_j=0$, $Z_iY_j=0$, $Y_iY_j=0$, $Z_iZ_j=0$, для всевозможных пар $(i,j)$.

\medskip

{\bf Категория 2}: $U_i=0$, где $U_i$ пробегает всевозможные слова длины не более $4$, не являющиеся кодировкой никакого пути на комплексе.

\medskip

{\bf Категория 3}: $U_i=0$, где $U_i$ пробегает всевозможные слова, являющиеся кодировкой какого-либо пути, представляющего нулевую форму (путь туда и обратно по некоторому ребру).

\medskip

{\bf Категория 4}: $U_i=0$, где $U_i$ пробегает всевозможные слова, длины не более $6$, являющиеся кодировкой какого-либо пути, представляющего путь с мертвым паттерном (одним из мертвых паттернов, перечисленных в предложении~\ref{DeadPaterns}).

\medskip

{\bf Категория $5$}: $U_i=V_j$, где $U_i$ и $V_j$ -- кодировки путей, участвующих в некотором локальном преобразовании. Конечное число таких пар приведено в главах~\ref{flip_section} и~\ref{pasting_section}.

\medskip

{\bf Оценка числа соотношений.} Букв в алфавите не более $7 \cdot 10^{36}$. Всевозможных слов длины не более $4$ существует не более $2{,}401 \cdot 10^{147}$. Учитывая, что в других категориях соотношений значительно меньше, можно считать, что общее число соотношений, введенных в полугруппе, не превосходит $3 \cdot 10^{147}$.

\smallskip

{\bf Замечание.} Введение соотношений можно значительно оптимизировать, сильно уменьшив их число. Но это выходит за рамки этой работы.

\medskip

Далее мы рассмотрим следствия из этих соотношений.

\subsection{Свойства путей, лежащих на одном ребре}

Пусть $W$ -- последовательность боковых узлов, принадлежащих одному внутреннему ребру некоторой макроплитки.
Поскольку мы ввели мономиальные соотношения для всех слов, не длиннее $10$, не являющихся кодировками путей на комплексе, мы можем установить, какие слова в принципе могут быть разрешены в полугруппе $S$.

\begin{definition}
Определим {\it последовательность уровней на ребре} $\mathbf{EdgeLevels}(k)$ в зависимости от натурального параметра $k$ следующим образом:

\smallskip

$\mathbf{EdgeLevels}(1)$=$1$;

$\mathbf{EdgeLevels}(2)$=$121$;

$\mathbf{EdgeLevels}(3)$=$1213121$;

Для $k>3$, $\mathbf{EdgeLevels}(k)$=$\mathbf{EdgeLevels}(k-1)$ $3$ $\mathbf{EdgeLevels}(k-1)$.

\end{definition}

Рассмотрим последовательность $\mathbf{EdgeLevels}(k)$.  Заменим $1$ на $2$, $2$ на $3$, $3$ оставляем как есть. После этого в начале и конце последовательности, а также между любыми двумя членами добавим $1$. Можно показать по индукции, что при этом получится $\mathbf{EdgeLevels}(k+1)$. также можно заметить, что последовательность периодическая, с периодом $3121$ и предпериодом $121$.
Ясно, что для любого бокового ребра существует такое $k$, что последовательность $\mathbf{EdgeLevels}(k)$ отвечает значениям уровней боковых узлов на этом ребре.

\medskip

\begin{proposition}   \label{edgepaths}
Пусть $W$ -- слово в $S$, в котором все буквы, соответствующие узлам, кодируют боковые узлы, лежащие на одном и том же ребре, одинаковым для всех параметром информации. Кроме того, пусть все реберные буквы одинаковы и равны $1$ или $2$ (то есть все переходы по главным ребрам).  Тогда:

$1$. Уровни узлов на ребре идут в последовательности $\mathbf{EdgeLevels}(k)$ для некоторого натурального $k$.

$2$. Для каждого узла $X$, уровень которого выше $1$, существует подклееная макроплитка $T$, такая что тип $X$ относительно $T$ (тип в подклейке) равен либо $\mathbb{CUR}$ либо $\mathbb{CDL}$, причем оба таких типа (для разных макроплиток) могут быть только у одного узла в $W$, находящегося в середине внутреннего ребра макроплитки.

\end{proposition}

{\bf Доказательство}. Первое утверждение следует из того, что в любом комплексе уровни идут именно в такой последовательности. То есть, все слова длины не более $20$, не являющиеся подсловами такой периодической последовательности, присутствуют среди мономиальных соотношений. Из этого следует, что если последовательность узлов не содержит запрещенных подслов, последовательность уровней узлов соответствует указанной.

Для доказательства второго утверждения рассмотрим, как устроены последовательности узлов на внутренних ребрах построенного комплекса. После операции разбиения, узлы первого уровня становятся узлами второго. При операции подклейки, каждый узел второго уровня становится либо правым верхним, либо левым нижним углом новой подклееной макроплитки, причем в дальнейшем к этим узлам макроплитки углами не подклеиваются.
Единственным узлом, в подклееном типе которого может присутствовать и $\mathbb{CUR}$ и $\mathbb{CDL}$, является узел в середине внутреннего ребра, так как он может участвовать в нескольких подклееных макроплитках.

\medskip

{\bf Ранги боковых узлов}.
Пусть $W$ -- последовательность боковых узлов, принадлежащих одному внутреннему ребру некоторой макроплитки.
Будем говорить, что {\it боковой узел имеет ранг $1$} (или, соответственно, ранг $2$), если этот боковой узел -- первого (или, соответственно, второго) уровня. {\it Регулярным словом} ранга $1$ будем считать слово из одного узла ранга $1$.

Теперь индуктивно определим остальные ранги узлов. {\it Регулярным словом} ранга $n$ будем называть последовательность узлов $XYX$, где $X$ -- регулярное слово ранга $n-1$, а $Y$ -- узел ранга $n$. То есть, например, если ранги узлов идут в последовательности $121$, то это регулярное слово ранга $2$.
Будем говорить, что {\it боковой узел $X$ имеет ранг $n$}, если выполнены следующие условия:

$1.$ $X$ находится в середине подслова $\widehat{W} \subset W$, состоящего из $2^{n-1}+1$ букв;

$2.$ Слово $\widehat{W}$ представляется в форме $YUXVZ$, где $Y$, $Z$ -- боковые узлы ранга $n-1$, $U$, $V$ -- регулярные слова ранга $n-2$.

\medskip

Слово $\widehat{W}$ будем называть {\it словом-представителем} для буквы $X$ ранга $n$.

\medskip

{\it Диаграммой рангов} слова $U$ будем называть слово $\widehat{U}$ в счетном алфавите $\{ 1,2,\dots \}$, получаемое из слова $U$ заменой каждой буквы на ее ранг.
Например, узел $X$ имеет ранг $3$, если в $W$ есть подслово $Z_1Y_1XY_2Z_2$, где у $Z_1$ и $Z_2$ ранг $2$, а у $Y_1$ и $Y_2$ ранг $1$. Тогда диаграмма рангов выглядит так: $21312$.

%Будем называть последовательность узлов $W$ {\it рангово-регулярной}, если она является подсловом некоторого регулярного слова ранга $n$.

\medskip

{\bf Замечание.}
Смысл понятия {\it слово-представитель $X$} в том что слово-представитель является половиной периметра (верхняя плюс левая сторона) подклееной макроплитки, ядром которой является $X$.

\medskip

\begin{figure}[hbtp]
\centering
\includegraphics[width=1\textwidth]{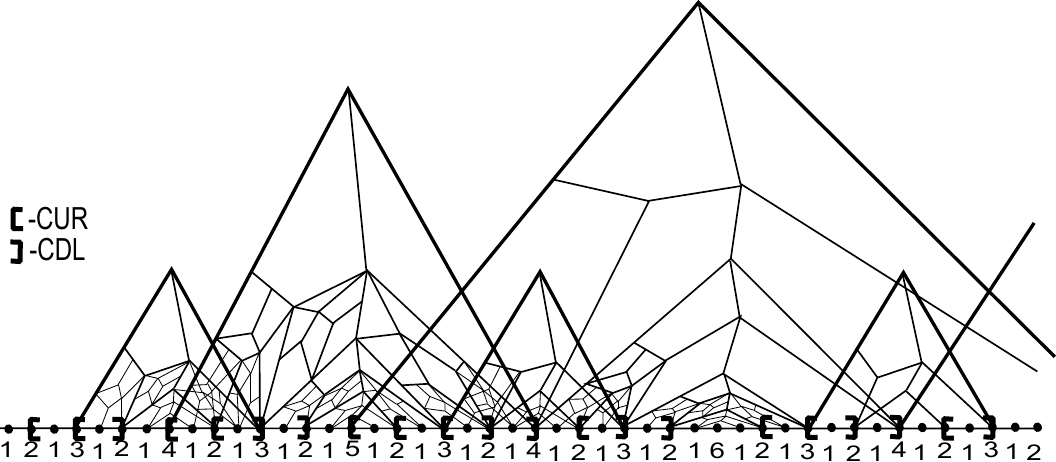}
\caption{Устройство подклеек на ребре}
\label{rankword}
\end{figure}

\medskip

\begin{proposition}   \label{skobka}
Пусть $W$ -- последовательность узлов на одном внутреннем ребре некоторой макроплитки и узел $X$ имеет ранг $n$. Рассмотрим слова длиной $2^k+1$,  с центром в $X$.  Его слово-представитель, имеющее длину $2^{n-1}+1$, является наименьшим среди таких слов, левый узел которых имеет в подклееном типе $\mathbb{CUR}$, а правый --  $\mathbb{CDL}$.

\end{proposition}

Это утверждение доказывается индукцией по $n$.

\medskip

\begin{proposition}   \label{mainwords}
Пусть $W$ -- слово в $S$, в котором все буквы, соответствующие узлам, кодируют боковые узлы, причем все они лежат на ребре одного и того же типа, и параметр информации у всех одинаковый.

Пусть $X$ -- буква в слове $W$, уровень которой равен $3$, причем между краями $W$ и $X$ находятся, как минимум, $4$ буквы, отвечающие кодам узлов. Тогда либо некоторое подслово $\widehat{W}\in W$, содержащее $X$, приводится к нулю с помощью определяющих соотношений в $S$, либо существует подслово в $W$, являющееся словом-представителем $X$ и оно может быть реализовано в качестве кодировки пути на некотором ребре комплекса.

\end{proposition}

{\bf Доказательство}. Будем считать, что мы читаем слова в положительном (относительно несущего ребра) направлении.  Пусть уровень $X$  (то есть и ранг тоже) равен $3$. Рассмотрим подслово $W$, содержащее $5$ узловых букв, с центром в $X$.

Допустим, первая и последняя буквы содержат в подклееном типе, соответственно, $\mathbb{CDL}$ и $\mathbb{CUR}$.
Если $W$ не встречается на комплексе, оно должно присутствовать среди мономиальных определяющих соотношений. Иначе оно автоматически реализуется в качестве кодировки пути.
 В этом случае $W$ -- кодировка пути по верхней и левой стороне подклееной макроплитки уровня $2$, и $W$ -- слово-представитель $X$.

Пусть, теперь, первая и последняя буквы содержат в подклееном типе, соответственно, $\mathbb{CUR}$ и $\mathbb{CDL}$.  (Иные варианты не встречаются на комплексе в качестве кодировок и потому приводятся к нулю).

\smallskip

\begin{figure}[hbtp]
\centering
\includegraphics[width=1\textwidth]{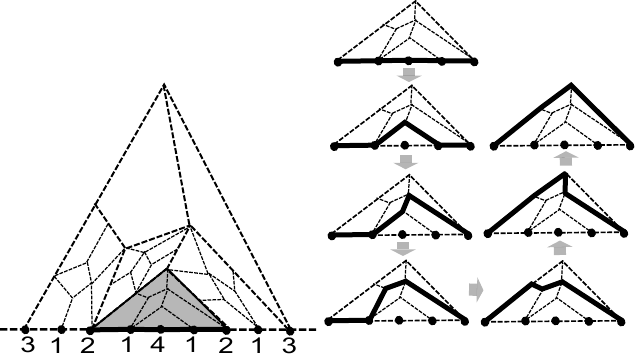}
\caption{Преобразование слова-представителя}
\label{enlarge}
\end{figure}

В этом случае ранг центральной буквы -- $4$ или выше. Рассмотрим слово $V_1$ из девяти букв с центром в $X$. Диаграмма рангов $V_1$ имеет вид $3121x1213$. Слово $V_1$ встречается на комплексе, иначе все можно сразу привести к нулю. Рассмотрим типы (в подклееной части) для четвертой и шестой букв этого слова. Это должны быть разные типы, а учитывая положительное прочтение слова, четвертая должна иметь тип $\mathbb{U}$, а шестая $\mathbb{L}$ (речь о типах в подклееной области). То есть, к подслову из четвертой, пятой и шестой букв может быть применено  соотношение для локального преобразования из числа введенных в случае $P2$. Затем можно применить другие локальные преобразования $7$ из числа введенных в случае $P2$, для соответствующих трехбуквенных слов. В результате мы преобразуем слово $V_1$ в кодировку пути, изображенного на рисунке~\ref{enlarge}.

\medskip

\begin{figure}[hbtp]
\centering
\includegraphics[width=0.9\textwidth]{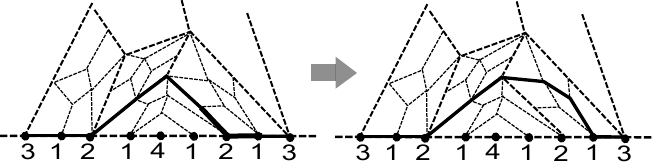}
\caption{Преобразование слова-представителя}
\label{enlarge2}
\end{figure}

\medskip

Рассмотрим трехбуквенное подслово, отмеченное в левой части на рисунке~\ref{enlarge2}. Заметим, что либо путь с таким кодом не может встретиться на комплексе (и тогда оно равно нулю), либо к нему можно применить соотношение $P3-8$. Действительно, узел $X$ обязан иметь уровень $1$, значит, базовый тип его $\mathbb{UL}$ или $\mathbb{LU}$, а тип в подклейке -- $\mathbb{U}$ (если будет $\mathbb{L}$, то такого участка из трех узлов не может встретиться на комплексе). Затем, аналогично применяем соотношения $B3-6$, $P3-1$, и получаем кодировку пути, изображенного в правой части рисунка~\ref{enlarge2}.

Далее возможно два варианта. Последняя буква слова $V_1$ может содержать в подклееном типе либо $\mathbb{CDL}$ либо $\mathbb{CUR}$.

\begin{figure}[hbtp]
\centering
\includegraphics[width=0.9\textwidth]{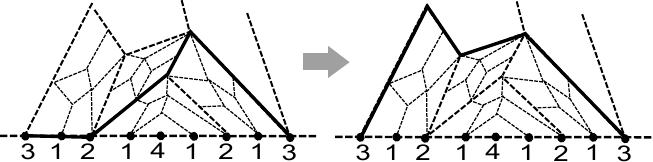}
\caption{Преобразование слова-представителя}
\label{enlarge3}
\end{figure}

В первом случае мы можем применить соотношение $P3-10$, далее $B3-5$, $B3-2$. Применяя далее подходящие соотношения, можно получить из слова кодировку пути, изображенного в левой части на рисунке~\ref{enlarge3}. Согласно предложению~\ref{skobka}, в этом случае рассмотренное девятибуквенное слово и есть слово-представитель для $X$. Применив далее несколько соотношений, мы можем преобразовать это слово в кодировку пути, указанного в правой части рисунка~\ref{enlarge3}.

\smallskip

\begin{figure}[hbtp]
\centering
\includegraphics[width=0.7\textwidth]{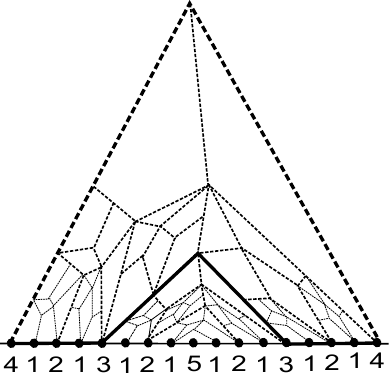}
\caption{Преобразование слова-представителя}
\label{enlarge4}
\end{figure}

Во втором случае мы продолжаем применять соотношения, пока не получим из нашего девятибуквенного слова кодировку указанного на рисунке~\ref{enlarge4} пути.

Теперь рассмотрим $17$-буквенное слово $V_2$ с центром в $X$ (рисунок~\ref{enlarge4}). Мы опять можем брать трехбуквенные слова на границе $17$-буквенного слова так, что либо эти трехбуквенные слова будут встречаться в мономиальных соотношениях, либо в соотношениях описанных в главах про локальные преобразования путей. После проведения необходимых преобразований, мы опять получим альтернативу: последняя буква слова $V_2$ может содержать в подклееном типе либо $\mathbb{CDL}$ либо $\mathbb{CUR}$. Аналогично $9$-буквенному слову, мы либо получаем кодировку пути, полученного из слова-представителя, либо возможность продолжить процесс и перейти к рассмотрению $33$ буквенного слова.

Продолжая процесс преобразования, мы дойдем до того момента, когда сможем окончательно выписать кодировку пути, полученного из слово-представителя, либо получим в какой-то момент нулевое слово.

\medskip

\begin{proposition}   \label{mainwords2}

Пусть $W$ -- слово в $S$, удовлетворяющее следующим условиям:

$1$. $W$ состоит из трех частей $W=V_1XV_2$, где

$X$- узловая буква, соответствующая внутреннему узлу ($\mathbb{A}$, $\mathbb{B}$ или $\mathbb{C}$), либо $X$ -- любая узловая буква, но ребро входа в $X$ или ребро выхода из $X$ -- неглавное;

$V_1$,$V_2$ -- слова, все узловые буквы которых обладают одной и той же информацией (своей для $V_1$ и для $V_2$), а все реберные буквы (кроме, возможно, входящего в $X$ и выходящего из $X$)  соответствуют главным ребрам, причем тип ребра постоянен для $V_1$ и для $V_2$.

Тогда либо некоторое подслово $\widehat{W}\in W$, содержащее $X$, приводится к нулю с помощью определяющих соотношений в $S$, либо существует подслово в $W$ являющееся словом-представителем $X$ и оно может быть реализовано в качестве кодировки пути на некотором ребре комплекса.

\end{proposition}

Доказательство полностью аналогично предложению~\ref{mainwords}. В данном случае в процессе доказательства используются соотношения, связанные с подклееной макроплиткой с ядром в $X$.

\medskip

\begin{proposition}   \label{checkedges}

Пусть $W$ -- слово в $S$ и каждая узловая буква $W$ кодирует боковой узел, причем информация всех таких узлов совпадает, а каждая реберная буква соответствует главному ребру (то есть для всех $1$ или для всех $2$). Тогда либо $W$ можно привести к нулю с помощью определяющих соотношений в $S$, либо $W$ представляется в виде кодировки некоторого пути на комплексе.

\end{proposition}

{\bf Доказательство}. Заметим, что если уровень всех узловых букв не выше $2$, то $W$ не может быть содержать более $3$ узловых букв, в этом случае $W$ имеет малую длину и если в нем не содержится запрещенных слов, значит оно реализуется в качестве кодировки пути.

Выберем некоторую букву $X$, уровня $3$. Пусть ее ранг равен $k$. Применим к $X$ предложение~\ref{mainwords}. Получим $V(X)$ -- слово-представитель $X$. Обозначим крайнюю справа букву $V(X)$ как $X_1$. Ее ранг будет меньше на $1$, чем ранг $X$. Применим к $X_1$ предложение~\ref{mainwords}, опять возьмем в $V(X_1)$ крайнюю справа букву. Продолжим так применять предложение~\ref{mainwords}, пока ранг и уровень всех буквы $X_i$,  которой равны $3$. Если в какой-то момент слово $W$ закончилось, мы можем рассмотреть аналогичный процесс, взяв на каждом шаге крайние слева буквы слов-представителей.

%Пусть процесс закончился с обеих сторон. Вернемся к началу процесса. Проделаем с $V(X)$ локальные преобразования, показанные на рисунке~\ref{enlarge}. Теперь путь приходит в вершину $X_1$ из подклееной макроплитки, соответствующей слову представителю $V(X)$ и уходит далее по главному ребру. Рассмотрим вершину, соответствующую букве $X_1$, как ядро подклееной макроплитки, со сторонами на ребрах, по которым приходит и уходит наш путь. В этой макроплитке также можно применить локальные преобразования, показанные на рисунке~\ref{enlarge}. В целом эти локальные преобразования аналогичны тем, что применяются в предложении~\ref{mainwords}, просто используется другая подклееная макроплитка. Путь можно ``развернуть'' до вершины $X_2$. Далее, можно аналогичную операцию для $X_2$ и для остальных вершин $X_i$. На каждом шаге процесса мы увеличиваем ту часть пути, которая проходит не по изначальному главному пути. Учитывая, что все вершины $X_i$ имеют уровень $3$, вплоть до края слова, описанный процесс позволяет нам получить путь, эквивалентный исходному, но такой, что только его концы лежат на главном ребре, причем в крайних узловых точках слова. То есть, мы получаем реализацию слова в качестве пути на комплексе.

\medskip

Будем пока считать, что слово не заканчивается. Тогда мы в итоге получим букву $X_i$, ранг и уровень которой равны $2$. Следующая за $X_i$ узловая буква должна иметь уровень $1$, иначе возникнет запрещенное слово. Рассмотрим еще одну следующую узловую букву $Y$. Ее уровень равен $3$, применим к ней предложение~\ref{mainwords}. Можно заметить, что левый конец полученного слова-представителя $V(Y)$ -- это буква $X$. То есть, ранг $Y$ на один выше, чем ранг $X$.

Таким образом, с помощью указанного процесса, мы можем по заданной узловой букве либо найти букву с большим рангом, либо дойти до конца слова (в обе стороны). Во втором случае указанная буква будет обладать наивысшим рангом из всех узловых букв.

Пусть буква $X$ -- обладает максимальным рангом из всех узловых букв $W$. Теперь применим предложение~\ref{mainwords} к $X$ и рассмотрим крайние буквы $X_L$ и $X_R$ полученного слова-представителя. Заметим, что их ранг на $1$ меньше, чем ранг $X$.  Для каждой из этих букв зафиксируем пару соседних с ней реберных букв так, что одно ребро уходит в подклейку (соответствующую слову-представителю $X$), а другая -- соответствует главному ребру. Применим к этим буквам с такими парами ребер предложение~\ref{mainwords2}, получая расширение подслова, представимого в качестве кодировки некоторого пути. Продолжим далее аналогичным образом, пока крайние буквы не будут иметь ранг $2$. Следующие две буквы с обоих концов должны иметь ранг $1$, и полученное таким образом слово представляется в виде кодировки пути на комплексе. Заметим, что если после этого слева или справа имеется какая-либо буква, она должна иметь больший ранг, чем $X$ (в этом можно убедиться, применив для этой буквы предложение~\ref{mainwords}). Таким образом, все слово может быть представлено в виде кодировки пути на комплексе.

Пусть у нас получился путь $P$ на комплексе. Теперь можно проделать над ним все локальные преобразования в обратном порядке, возвращая слову форму, где все узловые буквы принадлежат одному ребру и все переходы -- главные. Путь $P$ преобразуется в эквивалентный путь, идущий по одному внутреннему ребру некоторой макроплитки.

\medskip

\begin{proposition}   \label{checkedges2}

Пусть $W$ -- слово в $S$, удовлетворяющее следующим условиям:

$W$ состоит из трех частей $W=V_1XV_2$, где

$X$ -- узловая буква, соответствующая внутреннему узлу ($\mathbb{A}$, $\mathbb{B}$ или $\mathbb{C}$), либо $X$ -- любая узловая буква, но ребро входа в $X$ или ребро выхода из $X$ -- неглавное;

$V_1$, $V_2$ -- слова, все узловые буквы которых обладают одной и той же информацией (своей для $V_1$ и для $V_2$), а все реберные буквы (кроме, возможно, входящего в $X$ и выходящего из $X$) соответствуют главным ребрам, причем тип ребра постоянен для $V_1$ и для $V_2$.

Тогда либо $W$ можно привести к нулю с помощью определяющих соотношений в $S$, либо $W$ представляется в виде кодировки некоторого пути на комплексе.

\end{proposition}

Доказательство аналогично предложению~\ref{checkedges}.

\medskip

\begin{proposition}[О пути внутри макроплитки]   \label{innerpath}

Пусть $W$-- слово в $S$, отвечающее следующим условиям:

{\bf $1$.} В слове  $W$ отсутствуют буквы, кодирующие входы-выходы в подклейку;

{\bf $2$.} Первая и последняя буквы слова $W$ -- узловые  узлы;

{\bf $3$.} Вторая буква слова $W$ кодирует неглавное выходящее ребро;

{\bf $4$.} Предпоследняя буква слова $W$ кодирует неглавное входящее ребро;

{\bf $5$.} Кроме второй и предпоследней букв в слове $W$ ни одна буква не кодирует неглавных входящих или выходящих ребер (то есть все остальные буквы, кодирующие ребра -- отвечают главным ребрам).

Тогда $W$ либо может быть приведено к нулю, либо $W$ является кодировкой некоторого пути, соединяющего два узла на периметре некоторой макроплитки.

\end{proposition}

{\bf Доказательство}.  Можно считать, что слово $W$ имеет $\mathbf{CODE}$-форму, то есть буквы, кодирующие входящие ребра, узлы и выходящие ребра, встречаются в слове в правильном порядке. В противном случае слово сразу можно привести к нулю, c помощью соотношений категории $1$. Также, будем считать что всякое подслово длины $4$ в нашем слове является регулярным, иначе,  можно применить мономиальное соотношение из категории~$2$.

Согласно $\mathbf{CODE}$-форме, буквы $1$, $4$, $7$, $\dots$, $3k+1$ слова $W$ кодируют узлы. Буквы $2$, $5$, $8$, $\dots$, $3k+2$ -- кодируют выходящие ребра, а буквы $3$, $6$, $9$, $\dots$, $3k$ -- кодируют входящие ребра.

Неглавное ребро входа означает, что путь входит внутрь некоторой макроплитки.
Заметим, всего существует шесть возможных ``входов'' внутрь макроплитки: по ребрам $1$ и $2$ (из середины верхней стороны), по ребру $3$ из середины левой стороны, по ребру $6$ из середины правой стороны, по  ребру $7$ из левого нижнего угла и по ребру $8$ из правого нижнего угла.

\medskip

По второй букве $W$, кодирующей вход в макроплитку, можно установить, с каким из указанных шести входов мы имеем дело. Допустим, сначала, что это вход по ребру $3$ из середины левой стороны, а выход (предпоследняя буква) -- по ребру $1$ в середину верхней стороны.

Допустим, четвертая буква $W$ соответствует черному внутреннему узлу ($\mathbb{A}$, $\mathbb{B}$ или $\mathbb{C}$). Рассмотрим седьмую букву, соответствующую третьему по очередности узлу вдоль нашего пути. Если эта буква является последней в $W$, то $W$ представляет собой путь из середины левой стороны через узел типа $\mathbb{A}$ в середину верхней стороны. Все слова не соответствующие этому, не кодируют никакого пути на комплексе и поэтому содержатся среди мономиальных соотношений категории $2$. Значит, седьмая буква не является последней. Заметим, что если она соответствует боковому узлу, то подслово $W$ из первых семи букв будет кодировать несуществующий путь. В случае, когда седьмая буква соответствует внутреннему узлу, можно рассмотреть возможные типы десятой и тринадцатой букв. Можно заметить, что если путь с такой кодировкой существует, в нем содержится участок, проходимый туда и обратно (нулевая форма) и в этом случае, в нашем слове будет подслово из категории $3$ мономиальных соотношений.

\medskip

Теперь пусть четвертая буква $W$ кодирует боковой узел. Пусть в слове $W$ все буквы кодируют боковые узлы. Участок пути из одного бокового узла в другой по главному ребру не меняет тип этого ребра. В нашем случае это ребро $3$, причем с левой стороны -- сторона $\mathbf{3A}$, а с правой -- $\mathbf{3B}$. Но при выходе (когда мы приходим в середину верхней стороны) ребро должно быть типа $1$, то есть найдется локальный участок, где тип ребра меняется, значит можно применить соотношение из категории $2$.

Итак, после цепочки боковых узлов, нам должен встретиться внутренний узел. Поскольку у нас ребро $3$, это может быть только узел типа $\mathbb{A}$, для других не будет реализующих путей на комплексе. После узла $\mathbb{A}$, мы можем продолжать путь по ребру $1$, либо по ребру~$4$. (Если путь далее продолжается по ребру $3$, то образуется участок с нулевой формой и опять же можно применить мономиальное соотношение.)

В случае ребра $4$, мы продолжим изучение последующих букв, кодирующих узлы. Аналогично можно установить, что после нескольких боковых узлов, мы должны встретить узел типа $\mathbb{C}$, иначе можно будет применить мономиальное соотношение. После узла $\mathbb{C}$, дальнейший путь может продолжиться по одному из ребер $7$, $8$ или $5$. В первых двух случаях, если после цепочки боковых, мы встретим внутренний узел, можно будет применить мономиальное соотношение, так как такого участка пути не может существовать на комплексе.

Если же дальнейший путь идет по ребру $5$, опять рассмотрим цепочку боковых узлов до первого внутреннего узла (опять же, боковые узлы не могут поменять тип ребра, а в конце ребро должно быть первого типа). Это может быть только узел типа $\mathbb{B}$ и дальнейший путь может идти только по ребрам $2$ или $6$. В обоих случаях, для любого встреченного узла внутреннего типа можно будет применить мономиальное соотношение.

\medskip

Таким образом, наше слово должно содержать цепочку боковых узлов на ребре $1$, потом $\mathbb{A}$-узел, потом цепочка боковых узлов типа $3$. Если указанные цепочки содержат не более $10$ букв, кодирующих узлы, то можно напрямую либо найти все наше слово среди запрещенных, а если его там нет, значит слово реализуется в качестве кодировки. Пусть указанные цепочки достаточно длинные. Временно выбросим из слова три первые и три последние буквы. Останется слово, к которому мы можем применить предложение~\ref{checkedges2}, выбрав в качестве $X$ наш $\mathbb{A}$-узел.

То есть, урезанное слово без первых трех и последних трех букв может быть представлено в качестве кодировки пути на комплексе. Пусть это будет путь внутри некоторой макроплитки $T$, с заданными окружениями тех узлов на границе $T$, которые являются начальниками узлов из $W$.

При этом первые три и последние три буквы кодируют вход и выход в макроплитку, плюс информацию двух узлов на границе. Если вход или выход не соответствует ребрам, по которым далее должен идти путь, мы сразу можем применить мономиальное соотношение. Значит, окружение макроплитки и окружение и информация двух входных узлов могут сочетаться в одном объекте на комплексе.
Тогда получаем, что наше слово  реализуется в качестве кодировки некоторого пути на комплексе.

\medskip

{\bf Замечание $1$}. Если вход и выход осуществляются не по ребрам $3$ и $1$, можно провести полностью аналогичный побуквенный анализ слова, и найти подслово категорий $2$ или $3$, либо установить, что слово есть кодировка некоторого пути в макроплитке. Таким образом, предложение можно доказать для входов внутрь макроплитки по ребрам любых внутренних типов.

\medskip
{\bf Замечание $2$}. Неглавные ребра в начале и конце слова могут также быть ребрами в подклейку, рассуждения при этом не меняются.

\begin{proposition}[О сокращении подпути по подклееной макроплитке]   \label{deletepath}

Пусть $W$ -- слово в $S$, причем первая и последняя буквы кодируют некоторые узлы (то есть это некоторые из букв $Y_i$), а вторая и предпоследняя -- кодируют, соответственно входящее ребро в подклейку и выходящее ребро из подклейки, причем кроме этих двух, больше входов и выходов в подклееные области в слове $W$ не встречаются. Тогда $W$ либо может быть приведено к нулю, либо к форме, содержащей менее двух входов или выходов в подклееные области.

\end{proposition}

{\bf Доказательство}.  Можно считать, что $W$ имеет $\mathbf{CODE}$-форму, то есть буквы, кодирующие входящие ребра, узлы и выходящие ребра, встречаются в слове в правильном порядке.

{\it Cкобочной структурой} слова будем называть слово в алфавите из двух букв:$\{ \mathbf{[}$, $\mathbf{]}\}$  открывающей и закрывающей скобок, составленное следующим образом: -- каждая буква, кодирующая выходящее неглавное ребро, или ребро, ведущее в подклееную макроплитку, заменяется на ``$[$'', а каждая буква, кодирующая входящее неглавное ребро, или ребро, ведущее из подклееной макроплитки, заменяется на ``$]$''. Остальные буквы заменяются пустыми буквами.

\medskip

Рассмотрим скобочную структуру $W$. Она начинается открывающей скобкой и заканчивается закрывающей. Очевидно, что существует подслово в $W$, которому в скобочной структуре соответствует подслово ``$[$ $]$'' (из открывающей и закрывающей скобок). Это подслово $\widehat{W}$ удовлетворяет условиям предложения~\ref{innerpath}, то есть, в некоторой макроплитке есть путь, соединяющий точки периметра, с таким кодом. Но тогда есть цепочка локальных преобразований, приводящая этот путь к эквивалентному, но идущему по периметру той же макроплитки. Поскольку локальным преобразованиям отвечают соотношения из категории $5$, подслово $\widehat{W}$ можно привести к форме, скобочная структура которой содержит не более одной скобки (а не две). Отметим также, что начальная и конечная скобки не исчезают при такой операции.

Будем проводить такие операции, пока во всем слове $W$ не останется только открывающей и закрывающей скобок. Теперь наше слово соответствует пути, все ребра которого -- главные, кроме первого и последнего. Данная ситуация аналогична условиям предложения~\ref{innerpath}. Таким образом, если наше слово ненулевое, оно представляет некоторый путь и мы можем провести цепочку локальных преобразований, чтобы этот путь шел по периметру макроплитки. Такой путь не может содержать более одного ребра в подклейки. Таким образом, в результате преобразований, в слове понизилось количество ребер в подклейки.

\medskip

\begin{proposition}[О двух выходах в подклейку подряд]   \label{twoexits}

Пусть $W$-- слово в полугруппе $S$, представляющееся в виде $W=xV_1yV_2$, где $x,y,z$ -- буквы, кодирующие выходящие в подклейку ребра, а $V_1$, $V_2$ -- слова,  не содержащие букв, кодирующих возвращающиеся из подклеек ребра.

Тогда если количество узловых букв в $V_2$ более чем в $5$ раз превышает количество узловых букв в $V_1$, то слово может быть приведено к нулю в $S$.
\end{proposition}

{\bf Доказательство.}  Рассмотрим подслово $V_2$.  Допустим, все его реберные буквы кодируют главные ребра. Если среди узловых букв $V_2$  нет букв, соответствующим узлам $\mathbb{A}$, $\mathbb{B}$, $\mathbb{C}$, $\mathbb{R}$, $\mathbb{D}$, $\mathbb{CDR}$, то все узлы $V_2$ боковые и можно применить предложение~\ref{checkedges2} к слову $W=xV_1yV_2$, где в качестве узла $X$ из формулировки предложения~\ref{checkedges2} используется последний узел в $V_1$ (узел, где мы выходим в подклейку второй раз).  Так как одно ребро длиннее другого, процесс преобразований приведет к ситуации, когда подслово из трех узлов, к которому мы хотим применить локальное преобразование, представляет собой запрещенный участок пути, то есть можно применить мономиальное соотношение.

Если среди узловых букв $V_2$  есть буквы, соответствующие узлам $\mathbb{A}$, $\mathbb{B}$, $\mathbb{C}$, $\mathbb{R}$, $\mathbb{D}$, $\mathbb{CDR}$,  то можно заметить, что этих узлов не более $4$. Иначе можно проанализировать, как путь ходит по макроплитке и привести кусок пути к нулевой форме (рассуждения, аналогичные используемым в предложении~\ref{innerpath}). Заметим также, что расстояние между двумя такими узлами (если считать его в количестве узловых букв) будет одинаковым. Иначе можно применить предложение~\ref{checkedges2} взяв в качестве $X$ букву на границе между неравными участками.  Таким образом, общая длина $V_2$ не более чем в пять раз превышает длину начального участка (от начала $V_2$ до первого узла из числа указанных). Длина начального участка не превышает длины $V_1$ (по тем же причинам, как если в $V_2$ есть только боковые узлы).

\medskip

Пусть среди реберных букв $V_2$ есть переходы по неглавному ребру (то есть входы в подплитку). Будем считать, что есть только входы по неглавному ребру, так как если есть выходы, то можно применить сокращение пары скобок (предложение~\ref{innerpath}).
Отметим букву с каждым таким входом угловой скобкой ``<''. Согласно предложению~\ref{longpath}, путь, проходящий по макроплитке уровня $n$, имеет длину менее $5\cdot 2^{n-2}$. То есть, каждый кусок пути между узлами, отмеченными угловой скобкой, будет менее $5\cdot 2^{n-2}$, где $n$ -- уровень макроплитки, где проходит этот участок. С каждым прохождением угловой скобки мы попадаем в макроплитку, уровень которой как минимум, на $1$ меньше чем был до этого. Значит, общая длина пути не может превышать $5\cdot 2^{k-2}$, где $k$ -- уровень макроплитки, в которую мы переходим в узле $V_2$.

При этом, учитывая предложение~\ref{pasting_distance}, количество узловых букв (или длина по ребрам) в подслове $V_1$ не менее $2^{k-1}$. Отсюда следует наше утверждение.

\medskip

\begin{proposition}[О безквадратности пути по одному ребру]   \label{nosquares}

Пусть $W$-- путь на комплексе, целиком проходящий внутри некоторого ребра (внутреннего или граничного) и не содержащий узлов с типами  $\mathbb{DL}$,$\mathbb{LD}$, $\mathbb{DR}$, $\mathbb{RD}$, $\mathbb{RU}$, $\mathbb{UR}$, $\mathbb{R}$, $\mathbb{D}$.
Тогда соответствующее слово в полугруппе $S$ не содержит двух одинаковых подслов, идущих подряд.
\end{proposition}

{\bf Доказательство}. Рассмотрим подстановочную систему (или DOLL-систему) в алфавите $\{ U_1,L_1, U, L \}$, определяемую следующим образом: $U_1\rightarrow U_1UL_1$, $L_1\rightarrow U_1LL_1$, $U\rightarrow U$, $L\rightarrow L$.

Такая система, в частности, генерирует следующую последовательность слов:

$$U_1\rightarrow U_1UL_1 \rightarrow  U_1UL_1 U U_1LL_1 \rightarrow U_1UL_1 U U_1LL_1 U U_1UL_1 L U_1LL_1 \rightarrow \dots $$

Пусть $X$ -- некоторое слово. Обозначим описанную подстановочную замену как $f(X)$.

\medskip

Пусть сначала путь отвечающий $W$ проходит по внутреннему ребру. Рассмотрим кодировку нашего слова $W$ в полугруппе $S$. Проведем над ним следующие операции:

\medskip

{\it каждую букву, кодирующую $\mathbb{UL}$-узел, уровня выше первого, заменим на букву $U$;}

{\it каждую букву, кодирующую $\mathbb{LU}$-узел, уровня выше первого, заменим на букву $L$;}

{\it каждую букву, кодирующую $\mathbb{UL}$-узел первого уровня, заменим на букву $U_1$;}

{\it каждую букву, кодирующую $\mathbb{LU}$-узел первого уровня, заменим на букву $L_1$;}

{\it все буквы, кодирующие не узловые буквы, удалим из слова.}

\medskip

Получившееся слово в алфавите $\{ U_1,L_1, U, L \}$ будем называть $W'$. Заметим, что если $W$ содержит два одинаковых подслова, идущих подряд, то таким же свойством будет обладать и $W'$.

Покажем, что $W'$ является подсловом $f^n(U_1)$ или $f^n(L_1)$ для некоторого $n$. Действительно, выберем в нашем слове $W$ букву, $X$ кодирующую узел с самым большим рангом $n$ (узел, появившийся ранее всех при разбиениях). Такая буква только одна, так как если бы их было хотя бы две, то согласно правилам разбиения, между ними нашелся бы узел с более высоким рангом.

Допустим, ради определенности, что узел $X$ имеет тип $\mathbb{UL}$. Можно заметить, что последовательность узловых букв в слове $W'$ является подсловом в слове $f^{n-1}(U_1)$, это легко можно установить индукцией по $n$, подстановочная последовательность $f$ полностью соответствует образованию узловых вершин при разбиениях.

\medskip

Таким образом, для доказательства предложения, достаточно показать, что $f^{n}(U_1)$ не содержит двух одинаковых подслов подряд для любого $n$.

Допустим, это не так. Возьмем наименьшее $n$, при котором такие подслова появляются.
Пусть $f^{n}(U_1)$ содержит подслово $QQ$. Разметим наше слово следующим образом. Выделим в $f^{n}(U_1)$ обычными (не квадратными) скобками участки, соответствующие подсловам $U_1UL_1$, $U_1LL_1$, $U$, $V$, образовавшимся при последнем переходе
$f^{n-1}(U_1) \rightarrow f^{n}(U_1)$. Сами скобки не принадлежат алфавиту, а являются просто инструментом разметки. Если при этом в $Q$ вошло целое количество участков, начинающихся и заканчивающихся скобкой,
то в слове $f^{n-1}(U_1)$ также будет два одинаковых подслова $f^{-1}(Q)$ идущих подряд, что дает противоречие с выбором $n$.

\medskip

Значит есть подслово $(U_1UL_1)$ или $(U_1LL_1)$, часть которого лежит вне $Q$, а часть внутри. То есть, $Q$ начинается с $UL_1)$ или $LL_1)$,  либо с $U_1)$ или $L_1)$.

Рассмотрим первый случай. $Q$ должно заканчиваться на $(U_1$, то есть $Q$ представляется в виде $UL_1)Y(U_1$, где $Y$ -- слово содержащее целое количество участков в скобках и к нему можно применить $f^{-1}$. Тогда $f^{n}(U_1)$ представляется в виде
$$W_1 (U_1 \bf{UL_1)Y(U_1UL_1)Y(U_1} UL_1) W_2,$$

где $W_1$, $W_2$ -- некоторые подслова, содержащие целое число участков в скобках.

Можно заметить, что при применении $f^{-1}$ к подслову $(U_1 \bf{UL_1)Y(U_1UL_1)Y(U_1} UL_1)$ там также образуется два одинаковых подслова идущих подряд. То есть такие слова встретятся в $f^{n-1}(U_1)$, что невозможно.

Остальные случаи: $Q$ начинается с $LL_1)$; $U_1)$; $L_1)$ рассматриваются полностью аналогично.

\medskip

\begin{proposition}   \label{powerpath}
Пусть $W$ -- произвольное слово в $S$. Тогда слово $W^9$ приводится к нулю.

\end{proposition}

{\bf Доказательство.} Отметим в слове $W$ открывающими скобками -- ребра, кодирующие входы в подклейки, а закрывающими -- выходы. Если подслово начинается с открывающей скобки, заканчивается закрывающей, и не содержит внутри скобок, будем применять предложение~\ref{deletepath}. Таким образом, в какой-то момент мы преобразуем наше слово $W$ в форму $\widehat{W}=W_1W_2$, где $W_1$ содержит только закрывающие скобки, а $W_2$ только открывающие.

Заметим, что $W^{n+1}=W_1 (W_2W_1)^n W_2$.
Рассмотрим слово $W_2W_1$. Применяя к нему такой же процесс (предложение~\ref{deletepath} к каждой паре скобок), мы в итоге получим слово $V$, либо не содержащее скобок вообще, либо содержащее только один их вид.

Пусть, например, $V$ содержит только открывающие скобки, то для $n\le 7$, $V^n$ содержит подслово $xV_1yV_2$, где $x$, $y$ -- открывающие скобки (буквы, кодирующие ребра в подклейки), а количество узловых букв в $V_2$ превышает число узловых букв в $V_1$ в более чем в пять раз. Таким образом, учитывая предложение~\ref{twoexits}, в этом случае $V^n$ можно привести к нулю, а значит и $W^9$ тоже.

Пусть $V$ не содержит скобок вообще. В этом случае мы имеем дело с плоской ситуацией, без выходов в подклейки. Расставим скобки новым образом: открывающие для неглавных выходящих из узла ребер, а закрывающие -- для неглавных входящих в узел ребер. Фактически, теперь открытие и закрытие скобок символизирует вход и выход в макроплитку. Теперь будем применять предложение~\ref{innerpath} для участков внутри пары правильно расположенных скобок. Таким образом, мы приведем $V$ к форме $V_1V_2$, где в $V_1$ нет открывающих скобок, в $V_2$ нет закрывающих.

Заметим, что $V^n=V_1 (V_2V_1)^{n-1} V_2$.
Рассмотрим слово $V_2V_1$. Применяя к нему такой же процесс (предложение~\ref{innerpath} к каждой паре скобок), мы в итоге получим слово $Q$, либо не содержащее скобок вообще, либо содержащее только один их вид.

Пусть, например, $Q$ содержит только открывающие скобки, то $Q^7$ содержит подслово $xQ_1yQ_2$, где $x$, $y$ -- открывающие скобки (буквы, кодирующие ребра в подклейки), а количество узловых букв в $Q_2$ превышает число узловых букв в $Q_1$ более чем в пять раз. Таким образом, учитывая предложение~\ref{twoexits}, в этом случае $Q^7$ можно привести к нулю, а значит и $V^8$ тоже.

Пусть $Q$ не содержит скобок вообще.  В этом случае, все ребра главные. Будем следить, по каким ребрам проходит предполагаемый путь. Допустим, $Q$ содержит узловую букву кодирующую не боковой узел (то есть, узел типа $\mathbb{A}$, $\mathbb{B}$, $\mathbb{C}$, или угловой), тогда путь должен проходить по нескольким ребрам одной макроплитки. Заметим, что нельзя посетить более четырех внутренних ребер макроплитки не повернув ни в одном месте назад и не выйдя на периметр. Но $Q^7$ должен содержать больше четырех таких участков, значит, в каком-то из них можно будет применить мономиальное соотношение.

Пусть теперь все узловые буквы $Q$ кодируют боковые узлы. То есть предполагаемый путь проходит по одному и тому же ребру некоторой макроплитки. Для такого случая можно применить предложение~\ref{checkedges}. То есть мы применяем несколько преобразований, и если в процессе нельзя было применить мономиальное соотношение, в конце получится представление нашего слова в виде кодировки некоторого пути на комплексе. То есть, наш путь лежит полностью внутри некоторого ребра на комплексе.

\medskip

Внутри каждого ребра на комплексе нет более трех узлов с типами, отличными от $\mathbb{UL}$ и $\mathbb{LU}$. Так как наш путь содержит девять одинаковых участков подряд, в них не может встретиться вершин с типами, отличными от $\mathbb{UL}$ и $\mathbb{LU}$.
Но тогда в слове не может быть даже двух одинаковых участков подряд (предложение~\ref{nosquares}).
Таким образом, $W^9$ всегда приводится к нулю.

\medskip

Теперь мы можем завершить доказательство нашей основной теоремы:

\begin{theorem}
В полугруппе $S$, заданной конечным числом определяющих соотношений, существует бесконечное число различных слов, не равных нулю. При этом для любого слова $W$ его девятая степень приводится к нулю: $W^9=0$.

\end{theorem}

{\bf Доказательство}. Каждому регулярному пути на комплексе можно привести в соответствие его кодировку -- слово в полугруппе. При этом, при преобразовании пути его длина не меняется. Таким образом, если некоторый путь приводится к нулю, это значит, что либо в нем нашлось запрещенное подслово, либо локально некратчайший кусок пути, либо мертвый паттерн. Во всех этих случаях путь не является регулярным. То есть, все регулярные пути к нулю привести нельзя.

Вторая часть утверждения следует из предложения~\ref{powerpath}. Таким образом, построенная конечно определенная полугруппа $S$ содержит бесконечное множество различных слов, и при этом является нильполугруппой, где каждое слово в девятой степени приводится к нулю.

\medskip

\medskip

\section{ Приложение. Подсчет числа окружений} \label{Appendix}

В этой главе мы посчитаем число различных четверок-комбинаций типов ребер.

\medskip

\begin{proposition}[О нижней стороне макроплитки]   \label{bottomside}

Пусть $X$ и $Y$ -- вершины, соответствующие левому нижнему и правому нижнему углу макроплитки $T$. Кроме того, пусть $Z$ -- боковая вершина, лежащая на нижней стороне $T$. Тогда начальники у $Z$ -- это пара $(X,Y)$.
\end{proposition}

{\bf Доказательство.} Рассмотрим положение $T$ в ее родительской макроплитке $T'$. Заметим, что нижняя сторона $T$ в любом случае не попадает на границу $T'$, то есть обязательно является одним из внутренних ребер в $T'$. Тогда его концы $X$, $Y$ являются начальниками для всех боковых вершин лежащих внутри ребра.

\medskip

{\bf Определение.} Код пути $P$ будем называть {\it мертвым}, если существует натуральное $N$, такое что для любых путей $W_1, W_2$ длины более $N$ путь $W_1PW_2$ может быть приведен к нулевой форме.

\subsection{ Варианты расположения макроплиток}

В этой главе мы разберем различные случаи расположения макроплитки на мозаике и определим, какие возможные комбинации плиточной информации (то есть сочетания типов ее граничных ребер) она может иметь.

{\bf Определение.}  Назовем макроплитку $X$ {\it прямым потомком} для макроплитки $Y$, если $X$ является левой верхней подплиткой $Y$, либо является левой верхней подплиткой некоторой макроплитки $Z$, которая уже является прямым потомком для $Y$.

\medskip

\begin{proposition}[О комбинациях ребер при примыкании к краю]   \label{combos_edges}

Рассмотрим некоторую макроплитку $T$ с типами граничных ребер $t=\mathbf{Top}(T)$,  $l=\mathbf{Left}(T)$, $r=\mathbf{Right}(T)$ и $b=\mathbf{Bottom}(T)$. Пусть также среди типов ребер $t$  и $l$ хотя бы один обозначает краевое ребро ($\mathbf{top}$, $\mathbf{right}$, $\mathbf{bottom}$, $\mathbf{left}$). Тогда всевозможные варианты значений $t$, $l$, $r$, $b$ исчерпываются описанными ниже пятью случаями:

1. $T$ -- подклееная макроплитка. $l=\mathbf{left}$, $t=\mathbf{top}$, $r=\mathbf{right}$, $b=\mathbf{bottom}$. То есть, вся четверка
($\mathbf{left},\mathbf{top},\mathbf{right},\mathbf{bottom}$).

2. $T$ -- непосредственная подплитка подклееной плитки. В зависимости от ее типа четверка определяется следующим образом:

\medskip

\begin{enumerate}

\item {\it левое нижнее положение}  -- $(\mathbf{7A}, \mathbf{left}, \mathbf{3B},\mathbf{4A})$
\item {\it левое верхнее положение} -- $(\mathbf{left},\mathbf{top},\mathbf{1A},\mathbf{3A})$

\item {\it правое верхнее положение} -- $(\mathbf{top},\mathbf{right},\mathbf{6A},\mathbf{2A})$
\item {\it правое нижнее положение}  --  $(\mathbf{right},\mathbf{8A},\mathbf{5A},\mathbf{6B})$
\item  {\it нижнее положение}   -- $(\mathbf{8B},\mathbf{bottom},\mathbf{bottom},\mathbf{7B})$
\end{enumerate}

  \medskip

3. $T$ -- прямой потомок некоторой макроплитки $T'$, которая является непосредственной подплиткой подклееной макроплитки, либо $T$ -- правая нижняя подплитка макроплитки, являющейся нижней подплиткой в подклееной макроплитке, либо прямой потомок такой макроплитки.

$(\mathbf{7A}, \mathbf{left}, \mathbf{1A},\mathbf{3A})$

$(\mathbf{top},\mathbf{right},\mathbf{1A},\mathbf{3A})$

$(\mathbf{right}(\mathbf{bottom}),\mathbf{8A},\mathbf{1A},\mathbf{3A})$

$(\mathbf{bottom},\mathbf{8A},\mathbf{5A},\mathbf{6B})$

$(\mathbf{8B},\mathbf{bottom},\mathbf{1A},\mathbf{3A})$

 \medskip

4.  $T$ примыкает левой стороной к краю подклееной макроплитки, то есть $l=\mathbf{left},\mathbf{top},\mathbf{right},\mathbf{bottom}$. Тогда четверка $(l,t,r,b)$ может быть следующей:

\medskip

$(l,\mathbf{3B}, \mathbf{6A}, \mathbf{2A})$,  $(\mathbf{right}, \mathbf{6A}, \mathbf{6A}, \mathbf{2A})$, $(\mathbf{bottom}, \mathbf{6A}, \mathbf{6A}, \mathbf{2A})$,

$(l,\mathbf{1A}, \mathbf{6A}, \mathbf{2A})$,  $(\mathbf{right}, \mathbf{6A}, \mathbf{1A}, \mathbf{3A})$, $(\mathbf{bottom}, \mathbf{6A}, \mathbf{1A}, \mathbf{3A})$,

$(l,\mathbf{3B}, \mathbf{1A}, \mathbf{3A})$,  $(\mathbf{bottom}, \mathbf{bottom}, \mathbf{6A}, \mathbf{2A})$,

$(l,\mathbf{1A}, \mathbf{1A}, \mathbf{3A})$,  $(\mathbf{bottom}, \mathbf{bottom}, \mathbf{1A}, \mathbf{3A})$

где $l=\mathbf{left},\mathbf{top},\mathbf{right},\mathbf{bottom}$

\medskip

5. $T$ примыкает верхней стороной к краю подклееной макроплитки, то есть $t=\mathbf{left},\mathbf{top},\mathbf{right},\mathbf{bottom}$. Тогда четверка $(l,t,r,b)$ может быть следующей:

$(\mathbf{7A}, t, \mathbf{3B}, \mathbf{4A})$,  $(\mathbf{bottom}, \mathbf{bottom}, \mathbf{6A}, \mathbf{2A})$,

$(\mathbf{7A}, t, \mathbf{1A}, \mathbf{3A})$,  $(\mathbf{bottom}, \mathbf{bottom}, \mathbf{1A}, \mathbf{3A})$,

где $t=\mathbf{left},\mathbf{top},\mathbf{right},\mathbf{bottom}$

\medskip

Общий список возможных сочетаний перечислен в таблице~\ref{Bcases}.

\medskip

\begin{table}[hbtp]
\caption{Все варианты окружений, включающие примыкание к краю}
\centering
 \begin{tabular}{|c|c|c|}   \hline
B1.  & $(\mathbf{left},\mathbf{top},\mathbf{right},\mathbf{bottom})$  &  \cr \hline
B2. & $(\mathbf{left},\mathbf{top},\mathbf{1A},\mathbf{3A})$   & \cr \hline
B3.  & $(\mathbf{7A},x,\mathbf{3B},\mathbf{4A})$,   &$x=\mathbf{left},\mathbf{top},\mathbf{right},\mathbf{bottom}$    \cr \hline
B4.  & $(\mathbf{7A},x,\mathbf{1A},\mathbf{3A})$,  & $x=\mathbf{left},\mathbf{top},\mathbf{right},\mathbf{bottom}$   \cr \hline
B5. & $(\mathbf{top},\mathbf{right},\mathbf{6A},\mathbf{2A})$  &  \cr \hline
B6. & $(\mathbf{top},\mathbf{right},\mathbf{1A},\mathbf{3A})$  &    \cr \hline
B7. & $(x,\mathbf{8A},\mathbf{5A},\mathbf{6B})$  & $x=\mathbf{right},\mathbf{bottom}$    \cr \hline
B8. & $(x,\mathbf{8A},\mathbf{1A},\mathbf{3A})$  & $x=\mathbf{right},\mathbf{bottom}$    \cr \hline
B9. & $(\mathbf{8B},\mathbf{bottom},\mathbf{bottom},\mathbf{7B})$    &\cr \hline
B10. & $(\mathbf{8B},\mathbf{bottom},\mathbf{1A},\mathbf{3A})$    &\cr \hline
B11. & $(x,\mathbf{3B},\mathbf{6A},\mathbf{2A})$,   & $x=\mathbf{left},\mathbf{top},\mathbf{right},\mathbf{bottom}$ \cr  \hline
B12. & $(x,\mathbf{3B},\mathbf{1A},\mathbf{3A})$,   & $x=\mathbf{left},\mathbf{top},\mathbf{right},\mathbf{bottom}$\cr  \hline
B13. & $(\mathbf{right},\mathbf{6A},\mathbf{6A},\mathbf{2A})$   &  \cr \hline
B14. & $(\mathbf{right},\mathbf{6A},\mathbf{1A},\mathbf{3A})$    &   \cr \hline
B15. & $(\mathbf{bottom},\mathbf{6A},\mathbf{6A},\mathbf{2A})$   &\cr \hline
B16. & $(\mathbf{bottom},\mathbf{6A},\mathbf{1A},\mathbf{3A})$   &\cr \hline
B17. & $(x,\mathbf{1A},\mathbf{6A},\mathbf{2A})$,  & $x=\mathbf{left},\mathbf{top},\mathbf{right},\mathbf{bottom}$   \cr \hline
B18. & $(x,\mathbf{1A},\mathbf{1A},\mathbf{3A})$,  & $x=\mathbf{left},\mathbf{top},\mathbf{right},\mathbf{bottom}$  \cr \hline
B19 & $(\mathbf{bottom}, \mathbf{bottom},\mathbf{6A},\mathbf{2A})$   &  \cr \hline
B20. & $(\mathbf{bottom}, \mathbf{bottom},\mathbf{1A},\mathbf{3A})$  &  \cr \hline
  \end{tabular}
\label{Bcases}
\end{table}

\end{proposition}

\medskip

{\bf Доказательство.} Очевидно, что если $T$ сама является подклееной плиткой, то типы ребер будут как в случае 1. Если $T$ является непосредственной подплиткой подклееной макроплитки, то возможны
пять случаев (кроме средней подплитки, так как ее границы не выходят на край) и типы ребер будут как в случае 2.

\medskip

Теперь рассмотрим вариант, когда $T$ не является ни подклееной плиткой, ни непосредственной подплиткой подклееной плитки. Пусть левый верхний угол $T$ является углом подклееной макроплитки (но сама $T$ подклееной не является). Тогда возможны два случая. Первый -- $T$ является правой нижней подплиткой некоторой плитки, которая является нижней подплиткой в подклееной макроплитке. Тогда четверка $(l,t,r,b)$ принимает значения $\mathbf{bottom},\mathbf{8A},\mathbf{5A},\mathbf{6B}$. Второй -- $T$ есть прямой потомок некоторой плитки $X$, которая является одной из четырех непосредственных подплиток подклееной макроплитки ($X$ не может быть средней подплиткой, так как она не примыкает к краевым сторонам, и левой верхней подплиткой, так как тогда $T$ была бы прямым потомком подклееной плитки). В этом случае пара $(l, t)$ принимает следующие значения, соответствующие четырем положениям $X$: $(\mathbf{7A},\mathbf{left})$, $(\mathbf{top},\mathbf{right})$, $(\mathbf{right},\mathbf{8A})$, $(\mathbf{8B},\mathbf{bottom})$. Так как $T$ является прямым потомком $X$, то правое и нижнее ребро будут иметь типы $(\mathbf{1A},\mathbf{3A})$. Получаем перечисленные в случае $3$ варианты.

\medskip

\begin{table}[hbtp]
\caption{Комбинации ребер, возможные при внутреннем расположении макроплитки}
\centering
 \begin{tabular}{|c|c|}   \hline
для любого внутреннего типа $t$ & $(\mathbf{7A},t,\mathbf{3B},\mathbf{4A})$ \cr \hline
если $t \in \{ \mathbf{2A},\mathbf{3A},\mathbf{4A},\mathbf{5B},\mathbf{6B},\mathbf{7B} \}$  & $(\mathbf{8B},t,t,\mathbf{7B})$, $(t,t,\mathbf{6A},\mathbf{2A})$,   \cr
(примыкание нижней стороной) &   \cr \hline
если $t=\mathbf{1A}$ & $(l,\mathbf{1A},\mathbf{6A},\mathbf{2A})$, \cr
& ($\mathbf{l}$-любой внутренний тип)   \cr \hline
если $t=\mathbf{1B}$ & $(\mathbf{4B},\mathbf{1B},\mathbf{2B},\mathbf{5B})$,  \cr \hline
если $t=\mathbf{2B}$ & $(\mathbf{1B},\mathbf{2B},\mathbf{6A},\mathbf{2A})$,  \cr  \hline
если $t=\mathbf{3B}$ & $(l,\mathbf{3B},\mathbf{6A},\mathbf{2A})$,    \cr
& ($\mathbf{l}$-любой внутренний тип)   \cr \hline
если $t=\mathbf{5A}$ & $(\mathbf{8A},\mathbf{5A},\mathbf{6A},\mathbf{2A})$,  \cr \hline
если $t=\mathbf{6A}$ & $(l,\mathbf{6A},\mathbf{6A},\mathbf{2A})$, \cr
& ($\mathbf{l}$-любой внутренний, кроме $\mathbf{1B}, \mathbf{4B}, \mathbf{7A}, \mathbf{8A}, \mathbf{8B}$.)   \cr \hline
если $t=\mathbf{8A}$ & $(l,\mathbf{8A},\mathbf{5A},\mathbf{6B})$,  \cr
& ($\mathbf{l}$-любой внутренний, кроме $\mathbf{1B}, \mathbf{4B}, \mathbf{7A}, \mathbf{8A}, \mathbf{8B}$.)   \cr \hline
\x{если $(l,t,r,d)$- \cr одна из четверок выше }& $(l,r,\mathbf{1A},\mathbf{3A})$  \cr \hline
  \end{tabular}
\label{Acombos}
\end{table}

Итак, теперь пусть теперь левый верхний угол $T$ не является углом подклееной плитки. Тогда одновременно левая и верхняя стороны $T$ примыкать к краям не могут. Пусть сначала
$T$ примыкает левой стороной к краю подклееной макроплитки (но не в условиях случаев $1$-$3$).
Пусть $X$ это максимальная плитка, прямым потомком которой является наша макроплитка $T$ (либо, если $T$ -- не левая верхняя подплитка, то $X=T$).
Заметим, что $X$ не может быть подклееной плиткой или одной из шести подплиток подклееной плитки, так как тогда реализовался бы один из предыдущих случаев. Заметим, что $X$ может занимать только правое верхнее положение в своей родительской макроплитке $X'$, так как в остальных положениях левая сторона не будет примыкать к краю, либо можно было бы выбрать более крупную $X$, либо левый верхний угол $T$ попадет в угол подклееной плитки. То есть, $X'$ примыкает к краю своей верхней стороной, а верхняя сторона нашей плитки $T$ лежит на правой стороне $X'$. В свою очередь, $X'$ может занимать левое нижнее, левое верхнее, правое верхнее либо нижнее положение в своей родительской макроплитке. В зависимости от этого положения $t$ принимает значения $\mathbf{3B}$, $\mathbf{1A}$, $\mathbf{6A}$, $\mathbf{bottom}$. Заметим, что если $X'$ занимает нижнее положение в своей родительской макроплитке, (последний случай) то эта родительская макроплитка должна быть подклееной, так как иначе нижняя сторона не будет выходить на край. Аналогично, если  $X'$ занимает правое верхнее положение, то либо родительская плитка -- подклееная, и тогда край -- правый, либо родительская плитка -- нижняя подплитка подклееной и тогда край -- нижний.  Значения $(r,b)$ определяются как $(\mathbf{1A},\mathbf{3A})$, если $X\ne T$, и $(\mathbf{6A},\mathbf{2A})$, если $T=X$. Получаем весь список значений из случая~4.

\medskip
\begin{table}[hbtp]
\caption{Все варианты окружений для внутренних расположений}
\centering
 \begin{tabular}{|c|c|c|}   \hline
A1.  & $(\mathbf{7A},t,\mathbf{3B},\mathbf{4A})$  & для любого внутреннего типа $t$ \cr \hline
A2. & $(\mathbf{7A},t,\mathbf{1A},\mathbf{3A})$   &для любого внутреннего типа $t$ \cr \hline
A3.  & $(\mathbf{8B},t,t,\mathbf{7B})$,  &  где $t \in \{ \mathbf{2A},\mathbf{3A},\mathbf{4A},\mathbf{5B},\mathbf{6B},\mathbf{7B} \}$  \cr \hline
A4.  & $(\mathbf{8B},t,\mathbf{1A},\mathbf{3A})$,  &  где $t \in \{ \mathbf{2A},\mathbf{3A},\mathbf{4A},\mathbf{5B},\mathbf{6B},\mathbf{7B} \}$  \cr \hline
A5. & $(t,t,\mathbf{6A},\mathbf{2A})$,  &где $t \in \{ \mathbf{2A},\mathbf{3A},\mathbf{4A},\mathbf{5B},\mathbf{6B},\mathbf{7B} \}$    \cr \hline
A6. & $(t,t,\mathbf{1A},\mathbf{3A})$,  &где $t \in \{ \mathbf{2A},\mathbf{3A},\mathbf{4A},\mathbf{5B},\mathbf{6B},\mathbf{7B} \}$    \cr \hline
A7. & $(l,\mathbf{1A},\mathbf{6A},\mathbf{2A})$,   &($l$-любой внутренний тип)   \cr \hline
A8. & $(l,\mathbf{1A},\mathbf{1A},\mathbf{3A})$,   &($l$-любой внутренний тип)   \cr \hline
A9. & $(\mathbf{4B},\mathbf{1B},\mathbf{2B},\mathbf{5B})$,   &\cr \hline
A10. & $(\mathbf{4B},\mathbf{1B},\mathbf{1A},\mathbf{3A})$,   &\cr \hline
A11. & $(\mathbf{1B},\mathbf{2B},\mathbf{6A},\mathbf{2A})$,   &\cr  \hline
A12. & $(\mathbf{1B},\mathbf{2B},\mathbf{1A},\mathbf{3A})$,   &\cr  \hline
A13. & $(l,\mathbf{3B},\mathbf{6A},\mathbf{2A})$,    &($l$-любой внутренний тип)   \cr \hline
A14. & $(l,\mathbf{3B},\mathbf{1A},\mathbf{3A})$,    &($l$-любой внутренний тип)   \cr \hline
A15. & $(\mathbf{8A},\mathbf{5A},\mathbf{6A},\mathbf{2A})$,   &\cr \hline
A16. & $(\mathbf{8A},\mathbf{5A},\mathbf{1A},\mathbf{3A})$,   &\cr \hline
A17. & $(l,\mathbf{6A},\mathbf{6A},\mathbf{2A})$,  &($l$-любой внутренний, кроме $\mathbf{1B}, \mathbf{4B}, \mathbf{7A}, \mathbf{8A}, \mathbf{8B}$.)   \cr \hline
A18. & $(l,\mathbf{6A},\mathbf{1A},\mathbf{3A})$,  &($l$-любой внутренний, кроме $\mathbf{1B}, \mathbf{4B}, \mathbf{7A}, \mathbf{8A}, \mathbf{8B}$.)   \cr \hline
A19 & $(l,\mathbf{8A},\mathbf{5A},\mathbf{6B})$,   &($l$-любой внутренний, кроме $\mathbf{1B}, \mathbf{4B}, \mathbf{7A}, \mathbf{8A}, \mathbf{8B}$.)   \cr \hline
A20. & $(l,\mathbf{8A},\mathbf{1A},\mathbf{3A})$,  & ($l$-любой внутренний, кроме $\mathbf{1B}, \mathbf{4B}, \mathbf{7A}, \mathbf{8A}, \mathbf{8B}$.)   \cr \hline
  \end{tabular}
\label{Acases}
\end{table}

Итак {\bf остался случай}, когда $T$ примыкает верхней стороной к краю.
Пусть $X$ это максимальная плитка, прямым потомком которой является наша макроплитка $T$ (либо, если $T$ -- не левая верхняя подплитка, то $X=T$). Заметим, что $X$ может занимать левое нижнее либо правое верхнее положение в своей родительской макроплитке $X'$, так как в остальных положениях левая сторона не будет примыкать к краю, либо можно было бы выбрать более крупную $X$, либо левый верхний угол $T$ попадет в угол подклееной плитки. В случае правого верхнего положения, $X'$ должна быть нижней подплиткой подклееной плитки, иначе левый верхний угол $T$ попадет в угол подклееной плитки. Но тогда левая сторона $T$ тоже попадет на край, и это будет один из уже описанных случаев: $(\mathbf{bottom}, \mathbf{bottom}, \mathbf{1A}, \mathbf{3A})$ и $(\mathbf{bottom}, \mathbf{bottom}, \mathbf{1A}, \mathbf{3A})$. Итак, остается вариант, когда $X$ является левой нижней подплиткой $X'$. Тогда $l=\mathbf{7A}$, а пара $r,b$ определяется опять в зависимости от того, совпадает ли $T$ c $X$: либо ($\mathbf{3B},\mathbf{4A}$) если совпадает, либо $(\mathbf{1A},\mathbf{3A})$ если нет.
Получаем перечисленные варианты в случае $5$.

\medskip

\begin{proposition}   \label{combos_inner}
{\bf О комбинациях ребер при внутреннем расположении.}
Рассмотрим некоторую макроплитку $T$ с типами граничных ребер $t=\mathbf{Top}(T)$,  $l=\mathbf{Left}(T)$, $r=\mathbf{Right}(T)$ и $b=\mathbf{Bottom}(T)$, причем типы $l$ и $t$ ребер -- внутренние. Тогда возможны комбинации типов ребер перечисленные в таблице~\ref{Acombos}.

Общий список возможных сочетаний перечислен в таблице~\ref{Acases}.

\end{proposition}

\medskip

{\bf Доказательство.} Плитка $T$ может занимать одно из $6$ возможных мест в родительской макроплитке. В четырех из них можно сразу определить сочетание ребер и случай, к которому относится данная ситуация (таблица~\ref{innercombos2}).

\begin{table}[hbtp]
\caption{все варианты окружений, включающие примыкание к краю}
\centering
 \begin{tabular}{|c|c|c|c|}   \hline
 положение $T$ & комбинация ребер&& \x{итоговый \cr случай} \cr \hline
левое нижнее  & $(\mathbf{7A},x,\mathbf{3B},\mathbf{4A})$  & для любого внутреннего типа $x$ & A1 \cr \hline
правое нижнее  & $(x,\mathbf{8A},\mathbf{5A},\mathbf{6B})$  & \x{ $x$-правая сторона \cr  (любая, кроме \cr $\mathbf{1B}, \mathbf{4B}, \mathbf{7A}, \mathbf{8A}, \mathbf{8B}$) }  & A19 \cr \hline
среднее  & $(\mathbf{4B},\mathbf{1B},\mathbf{2B},\mathbf{5B})$  &  & A9 \cr \hline
нижнее & $(\mathbf{8B},x,x,\mathbf{7B})$& \x{ $x$-нижняя сторона, \cr (то есть \mathbf{2A},\mathbf{3A},\mathbf{4A},\mathbf{5B},\mathbf{6B},\mathbf{7B}) } & A3 \cr \hline
  \end{tabular}
\label{innercombos2}
\end{table}

\begin{table}[hbtp]
\caption{все варианты окружений, включающие примыкание к краю}
\centering
 \begin{tabular}{|c|c|c|c|}   \hline
 \x{положение \cr левого верхнего угла  }  & \x{ комбинация \cr ребер} &&  \cr \hline
середина верхней стороны  & \x{$(x,\mathbf{1A},\mathbf{6A},\mathbf{2A})$ \cr $(\mathbf{1B},\mathbf{2B},\mathbf{6A},\mathbf{2A})$} & \x{ $x$ -- верхняя сторона, \cr то есть любой тип} & \x{A7 \cr A11}  \cr\hline
середина левой стороны  & $(x,\mathbf{3B},\mathbf{6A},\mathbf{2A})$  & \x{ $x$ -- левая сторона, \cr то есть любой тип}  & A13 \cr\hline
середина правой стороны & $(x,\mathbf{6A},\mathbf{6A},\mathbf{2A})$  & \x{ $x$ -- правая сторона, \cr (любой, кроме \cr $\mathbf{1B},\mathbf{4B},\mathbf{7A},\mathbf{8A},\mathbf{8B}$)}  & A17 \cr\hline
внутренний $\mathbb{C}$-узел & $(\mathbf{8A},\mathbf{5A},\mathbf{6A},\mathbf{2A})$&  & A15 \cr \hline
середина нижней стороны & $(x,x,\mathbf{6A},\mathbf{2A})$& \x{ $x$ -- нижняя сторона,  \cr то есть \mathbf{2A},\mathbf{3A},\mathbf{4A},\mathbf{5B},\mathbf{6B},\mathbf{7B} }   & A5 \cr \hline
  \end{tabular}
\label{edgecombos}
\end{table}

\medskip

Пусть $T$ занимает правое верхнее положение. Тогда возможны ситуации, указанные в таблице~\ref{edgecombos}.

\medskip

Если $T$ занимает левое верхнее положение, рассмотрим максимального макроплитку $T'$, являющуюся прямым предком $T$ (то есть такую, что из $T'$ переходами к левой верхней подплитке можно за несколько операций получить $T$). Заметим, что в $T'$ левая и верхняя стороны такие же как в $T$, и при этом у $T'$ есть родительская макроплитка, в которой она занимает одно из пяти положений (кроме левого верхнего).

Таким образом, тип левого и верхнего ребер $T'$ будет как в одном из случаев рассмотренных выше, и общая четверка для $T$ будет состоять из этих двух ребер плюс \mathbf{1A} и \mathbf{3A} для правого и нижнего ребер. То есть получим один из четных случаев выписанных выше.

%\vfil \FloatBarrier

\section{Дальнейшие вопросы и соображения} \label{future}

Методы данной работы, возможно, помогут построить и иные объекты. В частности, представляет интерес следующий

\medskip {\bf Вопрос.} {\it Существует ли  конечно определенная полугруппа с размерностью Гель\-фан\-да--Ки\-рил\-ло\-ва $2.5$.}

\smallskip

Л.~Н.~Шеврин поставил следующие вопросы:

\medskip
{\bf Вопрос 1.}
{\it Каково наименьшее возможное число элементов базиса бесконечной конечно определенной нильполугруппы конечного индекса?}

\medskip
{\bf Вопрос 2.} {\it Каков наименьший возможный индекс бесконечной конечно определенной нильполугруппы конечного индекса?}
\medskip

 Можно ли опустить в формулировке обоих вопросов слова «конечного индекса»? А именно, эти слова нужно будет опустить, если будет дан положительный ответ на следующий вопрос:

\medskip
{\bf Вопрос 3.} {\it Любая ли бесконечная конечно определенная полугруппа имеет конечный индекс?}

Нам представляются чрезвычайно важным для дальнейшего обобщение  результата Х.~Гудмана-Штраусса. В частности, представляют интерес два вопроса:

\medskip
{\bf А.} Топологическое обобщение результата Х.~Гудмана-Штраусса. Мозаика рассматривается как плоский граф, где плитками являются простые циклы. Нужно доказать, что заданием конечного числа локальных правил можно добиться заданного иерархического правила построения.

\medskip
{\bf В.} Обобщение иерархического правила. Нужно доказать, что заданием конечного числа локальных правил можно добиться заданного рекурсивного правила построения мозаики.

\medskip

Последнее нам представляется важным как в связи с третьим вопросом Л.~Н.~Шеврина (мы предполагаем отрицательный ответ), так и для построениия иных рекурсивных объектов. По всей видимости, эти вопросы должны предшествовать теоретико-кольцевым исследованиям.

И, наконец, возникает вопрос, на который обратил наше внимание Л.~Н.~Шеврин: {\it какие максимальные ненулевые степени элементов возможны в построенной нами полугруппе~$S$?}

Все перечисленные вопросы представляются полезными для получения идей, позволяющих работать в кольцевой тематике.

%\section{Библиография}

%\addcontentsline{toc}{chapter}{\Numline {}Библиография}

%\markboth{}{Библиография}

%\begin{enumerate}

% Доклады, сделанные по материалам работы, приведены в конце списка литературы.

\end{document}